\documentclass[12pt,a4paper,titlepage]{report}
\usepackage[mathscr,mathcal]{eucal}
\usepackage{amsfonts}
\usepackage{amssymb}
\usepackage{xypic}
\newtheorem{theorem}{Theorem}[section]
\newtheorem{definition}[theorem]{Definition}
\newtheorem{ypoth}[theorem]{Assumptions}
\newtheorem{ypothesh}[theorem]{Assumption}
\newtheorem{lemma}[theorem]{Lemma}
\newtheorem{limma}{Lemma}
\newtheorem{corollary}[theorem]{Corollary}
\newtheorem{proposition}[theorem]{Proposition}
\newtheorem{thedef}[theorem]{Theorem-Definition}

\newtheorem{px}[theorem]{Example}

\newcommand{\ie}{\emph{i.e., }}

\title{On the cyclotomic Hecke algebras of complex reflection
groups}

\begin{document}

 $$\textrm{\textbf{UNIVERSIT{\'E} PARIS 7 - DENIS DIDEROT}}$$
                  $$\textrm{UFR Math{\'e}matiques}$$
                                 \\
               $$\textrm{\textbf{TH{\`E}SE DE DOCTORAT}}$$
              $$\textrm{Discipline: Math{\'e}matiques}$$
                                 \\
              $$\textrm{Pr{\'e}sent{\'e}e par}$$
         $$\textrm{\textbf{MARIA CHLOUVERAKI}}$$
  $$\textrm{pour obtenir le grade de Docteur de l'Universit{\'e} Paris 7}$$
                                 \\ \\
$$\textrm{\large SUR LES ALG{\`E}BRES DE HECKE CYCLOTOMIQUES}$$
  $$\textrm{\large DES GROUPES DE R{\'E}FLEXIONS COMPLEXES}$$

$$\textrm{\Huge $\widetilde{}$}$$
 $$\textrm{\large ON THE CYCLOTOMIC HECKE ALGEBRAS}$$
  $$\textrm{\large OF COMPLEX REFLECTION GROUPS}$$
                                  \\
   $$\textrm{\textbf{Th{\`e}se dirig{\'e}e par M. Michel BROU{\'E}}}$$

$$\textrm{Soutenue le 21 septembre 2007 devant le jury compos{\'e}
de}$$

$$\begin{array}{lccl}
  \textrm{\textbf{M. Michel BROU{\'E}}} & & & \textrm{directeur de th{\`e}se} \\
  \textrm{\textbf{M. Meinolf GECK}} & & & \textrm{rapporteur} \\
  \textrm{\textbf{M. Iain GORDON}} & & &  \textrm{examinateur}\\
  \textrm{\textbf{M. Gunter MALLE}} & & & \textrm{rapporteur} \\
  \textrm{\textbf{M. Jean MICHEL}} & & &  \textrm{examinateur}\\
  \textrm{\textbf{M. Rapha\"{e}l ROUQUIER}} & & &  \textrm{examinateur}\\
\end{array}$$

\newpage
 ~
\newpage

\chapter*{Introduction}

Les travaux de G.Lusztig sur les caract{\`e}res irr{\'e}ductibles
des groupes r{\'e}ductifs sur les corps finis ont mis en
{\'e}vidence le r{\^o}le important jou{\'e} par les ``familles de
caract{\`e}res'' des groupes de Weyl concern{\'e}s. Cependant, on
s'est r{\'e}cemment rendu compte qu'il serait de grand
int{\'e}r{\^e}t de g{\'e}n{\'e}raliser la notion des familles de
caract{\`e}res aux groupes de r{\'e}flexions complexes ou, plus
pr{\'e}cis{\'e}ment, {\`a} divers types d'alg{\`e}bres de Hecke
associ{\'e}es aux groupes de r{\'e}flexions complexes.

D'une part, les groupes de r{\'e}flexions complexes et certaines
d{\'e}formations de leurs alg{\`e}bres g{\'e}n{\'e}riques (les
alg{\`e}bres cyclotomiques) interviennent naturellement pour
classifier les ``s{\'e}ries de Harish-Chandra cyclotomiques'' des
caract{\`e}res des groupes r{\'e}ductifs finis, g{\'e}n{\'e}ralisant
ainsi le r{\^o}le jou{\'e} par le groupe de Weyl et son alg{\`e}bre
de Hecke traditionnelle dans la description de la s{\'e}rie
principale. Puisque les familles de caract{\`e}res du groupe de Weyl
jouent un r{\^o}le essentiel dans la d{\'e}finition des familles de
caract{\`e}res unipotents du groupe r{\'e}ductif fini correspondant
(cf. \cite{Lu1}), on peut esp{\'e}rer que plus g{\'e}n{\'e}ralement
les familles de caract{\`e}res des alg{\`e}bres cyclotomiques jouent
un r{\^o}le-clef dans l'organisation des familles de caract{\`e}res
unipotents.

D'autre part, pour certains groupes de r{\'e}flexions complexes (et
non de Coxeter) $W$, on a des donn{\'e}es qui semblent indiquer que,
derri{\`e}re le groupe $W$ se cache un objet myst{\'e}rieux - le
\emph{Spets} (cf. \cite{BMM2}, \cite{Ma3}) - qui pourrait jouer le
r{\^o}le de ``la s{\'e}rie des groupes r{\'e}ductifs finis de groupe
de Weyl $W$''. Dans certains cas, il est possible de d{\'e}finir les
caract{\`e}res unipotents du Spets, qui sont contr{\^o}l{\'e}s par
l'alg{\`e}bre de Hecke ``spetsiale'' de $W$, une
g{\'e}n{\'e}ralisation de l'alg{\`e}bre de Hecke classique des
groupes de Weyl.

L'obstacle principal pour cette g{\'e}n{\'e}ralisation est le manque
de bases de Kazhdan-Lusztig pour les groupes de r{\'e}flexions
complexes (non de Coxeter). Cependant, des r{\'e}sultats plus
r{\'e}cents de Gyoja $\cite{Gy}$ et de Rouquier $\cite{Rou}$ ont
rendu possible la d{\'e}finition d'un substitut pour les familles
des caract{\`e}res, qui peut {\^e}tre appliqu{\'e} {\`a} tous les
groupes de r{\'e}flexions complexes. Gyoja a d{\'e}montr{\'e} (cas
par cas) que la partition en ``$p$-blocs''de l'alg{\`e}bre de
Iwahori-Hecke d'un groupe de Weyl $W$ coincide avec la partition en
familles, quand $p$ est l'unique mauvais nombre premier pour $W$.
Plus tard, Rouquier a prouv{\'e} que les familles des caract{\`e}res
d'un groupe de Weyl $W$ sont exactement les blocs de caract{\`e}res
de l'alg{\`e}bre de Iwahori-Hecke de $W$ sur un anneau de
coefficients convenable. Cette d{\'e}finition se g{\'e}n{\'e}ralise
sans probl{\`e}me {\`a} toutes les alg{\`e}bres cyclotomiques de Hecke
des groupes de r{\'e}flexions complexes. Expliquons comment.

Soit $\mu_\infty$ le groupe des racines de l'unit{\'e} de
$\mathbb{C}$ et $K$ un corps de nombres contenu dans
$\mathbb{Q}(\mu_\infty)$. On note $\mu(K)$ le groupe des racines de
l'unit{\'e} de $K$ et pour tout $d>1$, on pose
$\zeta_d:=\textrm{exp}(2\pi i/d)$. Soit $V$ un $K$-espace vectoriel
de dimension finie. Soit $W$ un sous-groupe fini de $\textrm{GL}(V)$
engendr{\'e} par des (pseudo-)r{\'e}flexions et agissant
irr{\'e}ductiblement sur $V$ et $B$ le groupe de tresses associ{\'e}
{\`a} $W$. On note $\mathcal{A}$ l'ensemble des hyperplans de
r{\'e}flexion de $W$ et $V^{\textrm{reg}}:= \mathbb{C} \otimes
V-\bigcup_{H \in \mathcal{A}}\mathbb{C} \otimes H$. Pour $x_0 \in
V^{\textrm{reg}}$, on d{\'e}finit $B:=\Pi_1(V^{\textrm{reg}}/W,x_0)$

Pour toute orbite $\mathcal{C}$ de $W$ sur $\mathcal{A}$, on note
$e_{\mathcal{C}}$ l'ordre commun des sous-groupes $W_H$, o{\`u} $H
\in \mathcal{C}$ et $W_H$ est le sous-groupe form{\'e} de 1 et
toutes les r{\'e}flexions qui fixent $H$.

On choisit un ensemble d'ind{\'e}termin{\'e}es
$\textbf{u}=(u_{\mathcal{C},j})_{(\mathcal{C} \in
\mathcal{A}/W)(0\leq j \leq e_{\mathcal{C}}-1)}$ et on d{\'e}finit
l'\emph{alg{\`e}bre de Hecke g{\'e}n{\'e}rique} $\mathcal{H}$ de $W$
comme le quotient de l'alg{\`e}bre du groupe
$\mathbb{Z}[\textbf{u},\textbf{u}^{-1}]B$ par l'id{\'e}al
engendr{\'e} par les {\'e}l{\'e}ments de la forme
$$(\textbf{s}-u_{\mathcal{C},0})(\textbf{s}-u_{\mathcal{C},1}) \ldots
(\textbf{s}-u_{\mathcal{C},e_{\mathcal{C}}-1}),$$ o{\`u}
$\mathcal{C}$ parcourt l'ensemble $\mathcal{A}/W$ et $\textbf{s}$
l'ensemble de g{\'e}n{\'e}rateurs de monodromie autour des images
dans $V^{\textrm{reg}}/W$ des {\'e}l{\'e}ments de l'orbite
d'hyperplans $\mathcal{C}$.

Si on suppose que $\mathcal{H}$ est un
$\mathbb{Z}[\textbf{u},\textbf{u}^{-1}]$-module libre de rang
$|W|$ et qu'elle est munie d'une forme sym{\'e}trisante $t$  qui
satisfait les conditions $\ref{ypo}$,  alors on a le r{\'e}sultat
suivant d{\^u} {\`a} Malle (\cite{Ma3}, 5.2) : si
$\textbf{v}=(v_{\mathcal{C},j})_{(\mathcal{C} \in
\mathcal{A}/W)(0\leq j \leq e_{\mathcal{C}}-1)}$ un ensemble de
$\sum_{\mathcal{C} \in \mathcal{A}/W}e_{\mathcal{C}}$
ind{\'e}termin{\'e}es tel que, pour tout $\mathcal{C},j$, on a
$v_{\mathcal{C},j}^{|\mu(K)|}:=\zeta_{e_\mathcal{C}}^{-j}u_{\mathcal{C},j}$,
alors l'alg{\`e}bre $K(\textbf{v})\mathcal{H}$ est semisimple
d{\'e}ploy{\'e}e. Noter bien que ces hypoth{\`e}ses sont
verifi{\'e}es pour tous les groupes de r{\'e}flexions
irr{\'e}ductibles sauf un nombre fini d'entre eux (\cite{BMM2},
remarques pr{\'e}c{\'e}dant 1.17, $\S$ 2; \cite{GIM}). Dans ce cas,
par le ``th{\'e}or{\`e}me de d{\'e}formation de Tits'', on sait que
la sp{\'e}cialisation $v_{\mathcal{C},j} \mapsto 1$ induit une
bijection $\chi \mapsto \chi_\textbf{v}$ de l'ensemble
$\mathrm{Irr}(W)$ des caract{\`e}res absolument irr{\'e}ductibles de
$W$ sur l'ensemble $\mathrm{Irr}(K(\textbf{v})\mathcal{H})$ des
caract{\`e}res absolument irr{\'e}ductibles de l'alg{\`e}bre
$K(\textbf{v})\mathcal{H}$.

Soit maintenant $y$ une ind{\'e}termin{\'e}e. La
$\mathbb{Z}_K[y,y^{-1}]$-alg{\`e}bre, not{\'e}e par
$\mathcal{H}_\phi$, obtenue comme la sp{\'e}cialisation de
$\mathcal{H}$ via le morphisme $\phi: v_{\mathcal{C},j} \mapsto
y^{n_{\mathcal{C},j}}$, o{\`u} $n_{\mathcal{C},j} \in \mathbb{Z}$
pour tous $\mathcal{C}$ et $j$, est une \emph{alg{\`e}bre de Hecke
cyclotomique}. Elle est aussi munie d'une forme sym{\'e}trisante
$t_\phi$ d{\'e}finie comme la sp{\'e}cialisation de la forme
canonique $t$. On remarque que, pour $y=1$, l'alg{\`e}bre $\mathcal{H}_\phi$ se
sp{\'e}cialise {\`a} l'alg{\`e}bre du groupe $\mathbb{Z}_K[W]$.

On appelle \emph{anneau de Rouquier} de $K$ et note par
$\mathcal{R}_K(y)$ la sous-$\mathbb{Z}_K$-alg{\`e}bre de $K(y)$
$$\mathcal{R}_K(y):=\mathbb{Z}_K[y,y^{-1},(y^n-1)^{-1}_{n\geq 1}].$$
Les \emph{blocs de Rouquier} de $\mathcal{H}_\phi$ sont les blocs de
l'alg{\`e}bre $\mathcal{R}_K(y)\mathcal{H}_\phi$. Rouquier
\cite{Rou} a montr{\'e} que si $W$ est un groupe de Weyl et
$\mathcal{H}_\phi$ est obtenue via la sp{\'e}cialisation
cyclotomique  (``spetsiale'')
$$ \phi : v_{\mathcal{C},0} \mapsto y \,\textrm{ et }\, v_{\mathcal{C},j} \mapsto 1 \textrm{ pour } j \neq 0,$$
 alors ses blocs de Rouquier coincident
avec les ``familles de caract{\`e}res'' selon Lusztig. Ainsi, les
blocs de Rouquier jouent un r{\^o}le essentiel dans le programme
``Spets'' dont l'ambition est de faire jouer {\`a} des groupes de
r{\'e}flexions complexes
le r{\^o}le de groupes de Weyl de structures encore myst{\'e}rieuses.\\

En ce qui concerne le calcul des blocs de Rouquier, le cas de la
s{\'e}rie infinie est d{\'e}j{\`a} trait{\'e}e par Brou{\'e} et Kim
dans \cite{BK} et par Kim dans \cite{Kim}. D'ailleurs, les blocs de
Rouquier de l'alg{\`e}bre de Hecke cyclotomique ``spetsiale'' de
plusieurs groupes de r{\'e}flexions complexes exceptionnels ont
{\'e}t{\'e} d{\'e}termin{\'e}s par Malle et Rouquier dans
\cite{MaRo}. G{\'e}n{\'e}ralisant les m{\'e}thodes employ{\'e}es
dans le dernier, nous avons pu calculer les blocs de Rouquier de
toutes les alg{\`e}bres de Hecke cyclotomiques de tous les groupes
de r{\'e}flexions complexes exceptionnels. De plus, nous avons
d{\'e}couvert que les blocs de Rouquier d'une alg{\`e}bre de Hecke
cyclotomique d{\'e}pendent d'une donn{\'e}e num{\'e}rique du groupe
de r{\'e}flexions complexe $W$, ses \emph{hyperplans essentiels}.\\

De fa{\c c}on plus pr{\'e}cise, les deux premiers chapitres de cette
th{\`e}se pr{\'e}sentent des r{\'e}sultats qui sont donn{\'e}s ici
pour la commodit{\'e} du lecteur. Dans le premier chapitre, qui est
consacr{\'e} {\`a} l'alg{\`e}bre commutative, nous d{\'e}montrons
des r{\'e}sultats sur la divisibilit{\'e} et
l'irr{\'e}ductibilit{\'e} qui vont {\^e}tre tr{\`e}s utiles dans le
chapitre 3. Nous introduisons aussi les notions des ``morphismes
associ{\'e}s \`{a} des mon{\^o}mes'' et des ``morphismes
adapt{\'e}s''. Le deuxi{\`e}me chapitre est l'adaptation et la
g{\'e}n{\'e}ralisation des r{\'e}sultats classiques de la
th{\'e}orie des blocs et de la th{\'e}orie des
r{\'e}pr{\'e}sentations des alg{\`e}bres sym{\'e}triques, qui
peuvent etre trouv{\'e}s dans \cite{BK} et \cite{GePf}. Par
ailleurs, nous donnons un crit{\`e}re pour qu'une
alg{\`e}bre soit semisimple d{\'e}ploy{\'e}e.\\

Dans le troisi{\`e}me chapitre, nous trouvons le coeur th{\'e}orique
de cette th{\`e}se. Son but est la d{\'e}t{\'e}rmination des blocs
de Rouquier des alg{\`e}bres de Hecke cyclotomiques des groupes de
r{\'e}flexions complexes. Nous donnons la formule explicite suivante
pour les {\'e}l{\'e}ments de Schur associ{\'e}s aux caract{\`e}res
irr{\'e}ductibles de l'alg{\`e}bre de Hecke g{\'e}n{\'e}rique d'un
groupe de r{\'e}flexions complexe $W$ : si $K$ est le corps de
d{\'e}finition de $W$ et
$\textbf{v}=(v_{\mathcal{C},j})_{(\mathcal{C} \in
\mathcal{A}/W)(0\leq j \leq e_{\mathcal{C}}-1)}$ est un ensemble
d'ind{\'e}termin{\'e}es comme ci-dessus, alors l'{\'e}l{\'e}ment de
Schur $s_\chi(\textbf{v})$ associ{\'e} au caract{\`e}re
$\chi_\textbf{v}$ de $K(\textbf{v})\mathcal{H}$ est de la forme
$$s_\chi(\textbf{v})=\xi_\chi N_\chi \prod_{i \in I_\chi} \Psi_{\chi,i}(M_{\chi,i})^{n_{\chi,i}}$$
o{\`u}
\begin{itemize}
    \item $\xi_\chi$ est un {\'e}l{\'e}ment de $\mathbb{Z}_K$,
    \item $N_\chi= \prod_{\mathcal{C},j} v_{\mathcal{C},j}^{b_{\mathcal{C},j}}$ est un mon{\^o}me dans $\mathbb{Z}_K[\textbf{v},\textbf{v}^{-1}]$
          avec $\sum_{j=0}^{e_\mathcal{C}-1}b_{\mathcal{C},j}=0$
          pour tout $\mathcal{C} \in \mathcal{A}/W$,
    \item $I_\chi$ est un ensemble d'indices,
    \item $(\Psi_{\chi,i})_{i \in I_\chi}$ est une famille de polyn{\^o}mes $K$-cyclotomiques {\`a} une
    variable
          (\ie polyn{\^o}mes minimaux sur $K$ des racines de l'unit{\'e}),
    \item $(M_{\chi,i})_{i \in I_\chi}$ est une famille de mon{\^o}mes dans $\mathbb{Z}_K[\textbf{v},\textbf{v}^{-1}]$
          et si $M_{\chi,i} = \prod_{\mathcal{C},j} v_{\mathcal{C},j}^{a_{\mathcal{C},j}}$,
          alors $\textrm{pgcd}(a_{\mathcal{C},j})=1$
          et $\sum_{j=0}^{e_\mathcal{C}-1}a_{\mathcal{C},j}=0$
          pour tout $\mathcal{C} \in \mathcal{A}/W$,
    \item ($n_{\chi,i})_{i \in I_\chi}$ est une famille d'entiers positifs.
    \end{itemize}
Cette factorisation est unique dans $K[\textbf{v},\textbf{v}^{-1}]$
et les mon{\^o}mes $(M_{\chi,i})_{i \in I_\chi}$ sont uniques {\`a}
inversion pr{\`e}s. Si $\mathfrak{p}$ est un id{\'e}al premier de
$\mathbb{Z}_K$ et $\Psi(M_{\chi, i})$ est un facteur de
$s_\chi(\textbf{v})$ tel que $\Psi_{\chi, i}(1) \in \mathfrak{p}$,
alors le mon{\^o}me $M_{\chi,i}$ s'appellera
$\mathfrak{p}$\emph{-essentiel} pour $\chi$. Nous montrons que plus
on sp{\'e}cialise notre alg{\`e}bre via des morphismes associ{\'e}s
aux mon{\^o}mes $\mathfrak{p}$-essentiels, plus la taille de ses
$\mathfrak{p}$-blocs s'agrandit.

Soit maintenant $M:=\prod_{\mathcal{C},j}
v_{\mathcal{C},j}^{a_{\mathcal{C},j}}$ un mon{\^o}me
$\mathfrak{p}$-essentiel. L'hyperplan d{\'e}fini dans
$\mathbb{C}^{\sum_{\mathcal{C} \in \mathcal{A}/W}e_{\mathcal{C}}}$
par la relation $\sum_{\mathcal{C},j}
a_{\mathcal{C},j}t_{\mathcal{C},j} = 0,$ o{\`u}
$(t_{\mathcal{C},j})_{\mathcal{C},j}$ est un ensemble de
$\sum_{\mathcal{C} \in \mathcal{A}/W}e_{\mathcal{C}}$
ind{\'e}termin{\'e}es, s'appelle \emph{hyperplan
$\mathfrak{p}$-essentiel} pour $W$. Si $\phi: v_{\mathcal{C},j}
\mapsto y^{n_{\mathcal{C},j}}$ est une sp{\'e}cialisation
cyclotomique, alors les blocs de Rouquier de $\mathcal{H}_\phi$
d{\'e}pendent des hypeplans $\mathfrak{p}$-essentiels auxquels les
$n_{\mathcal{C},j}$ appartiennent (o{\`u} $\mathfrak{p}$ parcourt
l'ensemble des id{\'e}aux premiers de $\mathbb{Z}_K$). Donc les
blocs de Rouquier d'une alg{\`e}bre de Hecke cyclotomique
d{\'e}pendent d'une
donn{\'e}e num{\'e}rique du groupe $W$.\\

Le quatri{\`e}me chapitre est la partie calculatoire de cette
th{\`e}se. Nous pr{\'e}sentons l'algorithme et les r{\'e}sultats de
la d{\'e}termination des blocs de Rouquier de toutes les
alg{\`e}bres de Hecke cyclotomiques de tous les groupes de
r{\'e}flexions complexes exceptionnels.

\tableofcontents
\newpage
 ~
\newpage
\chapter*{Introduction}
\addcontentsline{toc}{chapter}{Introduction}

The work of G.Lusztig on the irreducible characters of reductive
groups over finite fields has displayed the important role of the
``character families'' of the Weyl groups concerned. However, only
recently was it realized that it would be of great interest to
generalize the notion of character families to the complex
reflection groups, or more precisely to some types of Hecke algebras
associated with complex reflection groups.

On one hand, the complex reflection groups and their associated
``cyclotomic'' Hecke algebras appear naturally in the classification
of the ``cyclotomic Harish-Chandra series'' of the characters of the
finite reductive groups, generalizing the role of the Weyl group and
its traditional Hecke algebra in the principal series. Since the
character families of the Weyl group play an essential role in the
definition of the families of unipotent characters of the
corresponding finite reductive group (\cite{Lu1}), we can hope that
the character families of the cyclotomic Hecke algebras play a key
role in the organization of families of unipotent characters more
generally.

On the other hand, for some complex reflection groups (non-Coxeter)
$W$, some data have been gathered which seem to indicate that behind
the group $W$, there exists another mysterious object - the
\emph{Spets} (see \cite{BMM2},\cite{Ma3}) - that could play the role
of the ``series of finite reductive groups of Weyl group $W$''. In
some cases, one can define the unipotent characters of the Spets,
which are controlled by the ``spetsial'' Hecke algebra of $W$, a
generalization of the classical Hecke algebra of the Weyl groups.

The main obstacle for this generalization is the lack of
Kazhdan-Lusztig bases for the non-Coxeter complex reflection groups.
However, more recent results of Gyoja $\cite{Gy}$ and Rouquier
$\cite{Rou}$ have made possible the definition of a substitute for
families of characters which can be applied to all complex
reflection groups. Gyoja has shown (case by case) that the partition
into ``$p$-blocks'' of the Iwahori-Hecke algebra of a Weyl group $W$
coincides with the partition into families, when $p$ is the unique
bad prime number for $W$. Later, Rouquier showed that the families
of characters of a Weyl group $W$ are exactly the blocks of
characters of the Iwahori-Hecke algebra of $W$ over a suitable
coefficient ring. This definition generalizes without problem to all
the cyclotomic Hecke algebras of complex reflection groups. Let us
explain how.

Let $\mu_\infty$ be the group of all the roots of unity in
$\mathbb{C}$ and $K$ a number field contained in
$\mathbb{Q}(\mu_\infty)$. We denote by $\mu(K)$ the group of all the
roots of unity of $K$ and for all $d>1$, we put
$\zeta_d:=\textrm{exp}(2\pi i/d)$. Let $V$ be a finite dimensional
$K$-vector space. Let $W$ be a finite subgroup of GL$(V)$ generated
by (pseudo-)reflections and acting irreducibly on $V$. Denote by
$\mathcal{A}$ the set of its reflecting hyperplanes. We set
$V^{\textrm{reg}}:= \mathbb{C} \otimes V-\bigcup_{H \in
\mathcal{A}}\mathbb{C} \otimes H$. For $x_0 \in V^{\textrm{reg}}$,
we define $B:=\Pi_1(V^{\textrm{reg}}/W,x_0)$ the braid group
associated with $W$.

For every orbit $\mathcal{C}$ of $W$ on $\mathcal{A}$, we set
$e_{\mathcal{C}}$ the common order of the subgroups $W_H$, where $H$
is any element of $\mathcal{C}$ and $W_H$ the subgroup formed by 1
and all the reflections fixing the hyperplane $H$.

We choose a set of indeterminates
$\textbf{u}=(u_{\mathcal{C},j})_{(\mathcal{C} \in
\mathcal{A}/W)(0\leq j \leq e_{\mathcal{C}}-1)}$ and we denote by
$\mathbb{Z}[\textbf{u},\textbf{u}^{-1}]$ the Laurent polynomial ring
in all the indeterminates $\textbf{u}$. We define the \emph{generic
Hecke algebra} $\mathcal{H}$ of $W$ to be the quotient of the group
algebra $\mathbb{Z}[\textbf{u},\textbf{u}^{-1}]B$ by the ideal
generated by the elements of the form
$$(\textbf{s}-u_{\mathcal{C},0})(\textbf{s}-u_{\mathcal{C},1}) \ldots (\textbf{s}-u_{\mathcal{C},e_{\mathcal{C}}-1}),$$
where $\mathcal{C}$ runs over the set $\mathcal{A}/W$ and
$\textbf{s}$ over the set of monodromy generators around the images
in $V^{\textrm{reg}}/W$ of the elements of the hyperplane orbit
$\mathcal{C}$.

If we assume that $\mathcal{H}$ is a free
$\mathbb{Z}[\textbf{u},\textbf{u}^{-1}]$-module of rank $|W|$ and it
has a symmetrizing form $t$ which satisfies assumptions $\ref{ypo}$,
then we have the following result by Malle (\cite{Ma4}, 5.2): If
$\textbf{v}=(v_{\mathcal{C},j})_{(\mathcal{C} \in
\mathcal{A}/W)(0\leq j \leq e_{\mathcal{C}}-1)}$ is a set of
$\sum_{\mathcal{C} \in \mathcal{A}/W}e_{\mathcal{C}}$ indeterminates
such that, for every $\mathcal{C},j$, we have
$v_{\mathcal{C},j}^{|\mu(K)|}=\zeta_{e_\mathcal{C}}^{-j}u_{\mathcal{C},j}$,
then the $K(\textbf{v})$-algebra $K(\textbf{v})\mathcal{H}$ is split
semisimple. Note that these assumptions have been verified for all
but a finite number of irreducible complex reflection groups
(\cite{BMM2}, remarks before 1.17, $\S$ 2; \cite{GIM}). In this
case, by ``Tits' deformation theorem'', we know that the
specialization $v_{\mathcal{C},j}\mapsto 1$ induces a bijection
$\chi \mapsto \chi_{\textbf{v}}$ from the set $\mathrm{Irr}(W)$ of
absolutely irreducible characters of $W$ to the set
$\mathrm{Irr}(K(\textbf{v})\mathcal{H})$ of absolutely irreducible
characters of $K(\textbf{v})\mathcal{H}$.

Now let $y$ be a indeterminate. The $\mathbb{Z}_K[y,y^{-1}]$-algebra
obtained as a specialization of $\mathcal{H}$ via the morphism $\phi:
v_{\mathcal{C},j} \mapsto y^{n_{\mathcal{C},j}}$, where
$n_{\mathcal{C},j} \in \mathbb{Z}$ for all $\mathcal{C}$ and $j$, is
a \emph{cyclotomic Hecke algebra} and it is denoted by
$\mathcal{H}_\phi$. It also has a symmetrizing form $t_\phi$ defined
as the specialization of the canonical form $t$. We notice that,
for $y=1$, the algebra $\mathcal{H}_\phi$ specializes to the group algebra
$\mathbb{Z}_K[W]$.

We call \emph{Rouquier ring} of $K$ and denote by $\mathcal{R}_K(y)$
the $\mathbb{Z}_K$-subalgebra of $K(y)$
$$\mathcal{R}_K(y):=\mathbb{Z}_K[y,y^{-1},(y^n-1)^{-1}_{n\geq 1}].$$
The \emph{Rouquier blocks} of $\mathcal{H}_\phi$ are the blocks of
the algebra $\mathcal{R}_K(y)\mathcal{H}_\phi$. It has been shown by
Rouquier \cite{Rou}, that if $W$ is a Weyl group and
$\mathcal{H}_\phi$ is obtained via the cyclotomic specialization
(``spetsial'')
$$ \phi : v_{\mathcal{C},0} \mapsto y \,\textrm{ and }\, v_{\mathcal{C},j} \mapsto 1 \textrm{ for } j \neq 0,$$
then its Rouquier blocks coincide with the ``families of
characters'' defined by Lusztig. Thus, the Rouquier blocks play an
essential role in the program ``Spets'' whose ambition is to give to
complex reflection groups the role of Weyl groups of
as yet mysterious structures.\\

As far as the calculation of the Rouquier blocks is concerned, the
case of the infinite series has already been treated by Brou{\'e}
and Kim in \cite{BK} and Kim in \cite{Kim}. Moreover, the Rouquier
blocks of the ``spetsial'' cyclotomic Hecke algebra of many
exceptional complex reflection groups have been determined by Malle
and Rouquier in \cite{MaRo}. Generalizing the methods used in the
latter, we have been able to calculate the Rouquier blocks of all
cyclotomic Hecke algebras of all exceptional complex reflection
groups. Moreover, we discovered that the Rouquier blocks of a
cyclotomic Hecke algebra depend on a numerical datum of the complex
reflection group $W$, its
\emph{essential hyperplanes}.\\

Let us get into more details: The first two chapters of this thesis
present results which are given here for the convenience of the
reader. In the first chapter, which is dedicated on commutative
algebra, we prove some results about divisibility and irreducibility
which are going to be very useful in chapter 3. We also introduce
the notions of ``morphisms associated with monomials'' and ``adapted
morphisms''. The second chapter is the adaptation and generalization
of classic results of block theory and representation theory of
symmetric algebras, which can be found in \cite{BK} and \cite{GePf}.
We also give a criterion for an algebra to be split
semisimple.\\

In the third chapter, we find the theoretical core of this thesis.
Its aim is the determination of the Rouquier blocks of the
cyclotomic Hecke algebras of complex reflection groups. We give the
following explicit formula for the Schur elements associated with
the irreducible characters of the generic Hecke algebra of a complex
reflection group $W$: If $K$ is the field of definition of $W$ and
$\textbf{v}=(v_{\mathcal{C},j})_{(\mathcal{C} \in
\mathcal{A}/W)(0\leq j \leq e_{\mathcal{C}}-1)}$ is a set of
indeterminates like above, then the Schur element
$s_\chi(\textbf{v})$ associated with the character $\chi_\textbf{v}$
of $K(\textbf{v})\mathcal{H}$ is of the form
$$s_\chi(\textbf{v})=\xi_\chi N_\chi \prod_{i \in I_\chi} \Psi_{\chi,i}(M_{\chi,i})^{n_{\chi,i}}$$
where
\begin{itemize}
    \item $\xi_\chi$ is an element of $\mathbb{Z}_K$,
    \item $N_\chi= \prod_{\mathcal{C},j} v_{\mathcal{C},j}^{b_{\mathcal{C},j}}$ is a monomial in $\mathbb{Z}_K[\textbf{v},\textbf{v}^{-1}]$
          with $\sum_{j=0}^{e_\mathcal{C}-1}b_{\mathcal{C},j}=0$
          for all $\mathcal{C} \in \mathcal{A}/W$,
    \item $I_\chi$ is an index set,
    \item $(\Psi_{\chi,i})_{i \in I_\chi}$ is a family of $K$-cyclotomic polynomials in one variable
           (\ie minimal polynomials of the roots of unity over $K$),
    \item $(M_{\chi,i})_{i \in I_\chi}$ is a family of monomials in $\mathbb{Z}_K[\textbf{v},\textbf{v}^{-1}]$
          and if $M_{\chi,i} = \prod_{\mathcal{C},j} v_{\mathcal{C},j}^{a_{\mathcal{C},j}}$,
          then $\textrm{gcd}(a_{\mathcal{C},j})=1$
          and $\sum_{j=0}^{e_\mathcal{C}-1}a_{\mathcal{C},j}=0$
          for all $\mathcal{C} \in \mathcal{A}/W$,
    \item ($n_{\chi,i})_{i \in I_\chi}$ is a family of positive integers.
\end{itemize}
This factorization is unique in $K[\textbf{v},\textbf{v}^{-1}]$ and
the monomials $(M_{\chi,i})_{i \in I_\chi}$ are unique up to
inversion. If $\mathfrak{p}$ is a prime ideal of $\mathbb{Z}_K$ and
$\Psi(M_{\chi,i})$ is a factor of $s_\chi(\textbf{v})$ such that
$\Psi_{\chi,i}(1) \in \mathfrak{p}$, then the monomial $M_{\chi,i}$
will be called $\mathfrak{p}$\emph{-essential} for $\chi$. We show
that the more we specialize our algebra via morphisms associated
with $\mathfrak{p}$-essential monomials, the more the size of its
$\mathfrak{p}$-blocks becomes larger.

Now let $M:=\prod_{\mathcal{C},j}
v_{\mathcal{C},j}^{a_{\mathcal{C},j}}$ be a $\mathfrak{p}$-essential
monomial. The hyperplane defined in $\mathbb{C}^{\sum_{\mathcal{C}
\in \mathcal{A}/W}e_{\mathcal{C}}}$ by the relation
$\sum_{\mathcal{C},j} a_{\mathcal{C},j}t_{\mathcal{C},j} = 0,$ where
$(t_{\mathcal{C},j})_{\mathcal{C},j}$ is a set of $\sum_{\mathcal{C}
\in \mathcal{A}/W}e_{\mathcal{C}}$ indeterminates, is called
\emph{$\mathfrak{p}$-essential hyperplane} for $W$. If $\phi:
v_{\mathcal{C},j} \mapsto y^{n_{\mathcal{C},j}}$ is a cyclotomic
specialization, then the Rouquier blocks of $\mathcal{H}_\phi$
depend on which $\mathfrak{p}$-essential hyperplanes the
$n_{\mathcal{C},j}$ belong (where $\mathfrak{p}$ runs over the prime
ideals of $\mathbb{Z}_K$). Hence the Rouquier blocks of a cyclotomic
Hecke algebra depend on a numerical datum of the group $W$.\\

The fourth chapter is the calculation part of this thesis. We
present the algorithm and the results of the determination of the
Rouquier blocks of all cyclotomic Hecke algebras of all exceptional
complex reflection groups.

\chapter{On Commutative Algebra}

Throughout this chapter, all rings are assumed to be commutative
with 1. Moreover, if $R$ is a ring and $x_0,x_1,\ldots,x_m$ is a set
of indeterminates on $R$, then we denote by $R[x_0^{\pm 1},x_1^{\pm
1},\ldots,x_m^{\pm 1}]$ the Laurent polynomial ring on $m+1$
indeterminates, \ie the ring
$R[x_0,{x_0}^{-1},x_1,{x_1}^{-1},\ldots,x_m,{x_m}^{-1}]$.

\section{Localizations}

\begin{definition}\label{multiplicatively closed set}
Let $R$ be a commutative ring with $1$. We say that a subset $S$ of
$R$ is a multiplicatively closed set if $0 \notin S$, $1 \in S$ and
every finite product of elements of $S$ belongs to $S$.
\end{definition}

In the set $R \times S$, we introduce an equivalence relation such
that $(r,s)$ is equivalent to $(r',s')$ if and only if there exists
$t \in S$ such that $t(s'r-sr')=0$. We denote the equivalence class
of $(r,s)$ by $r/s$. The set of equivalence classes becomes a ring
under the operations such that the sum and the product of $r/s$ and
$r'/s'$ are given by $(s'r+sr')/ss'$ and $rr'/ss'$ respectively. We
denote this ring by $S^{-1}R$ and we call it the \emph{localization}
of $R$ at $S$. If $S$ contains no zero divisors of $R$, then any
element $r$ of $R$ can be identified with the element $r/1$ of $S^{-1}R$ and we
can regard the latter as an $R$-algebra.\\
\\
\begin{remarks}\
\emph{\begin{itemize}
  \item If $S$ is the set of all non-zero divisors of $R$, then
  $S^{-1}R$ is called \emph{the total quotient ring} of $R$. If, moreover,
  $R$ is an integral domain, the total quotient ring of $R$ is \emph{the
  field of fractions} of $R$.
  \item If $R$ is Notherian, then $S^{-1}R$ is Noetherian.
  \item If $\mathfrak{p}$ is a prime ideal of $R$, then the set
  $S:=R-\mathfrak{p}$ is a multiplicatively closed subset of $R$. Then the
  ring $S^{-1}R$ is simply denoted by $R_\mathfrak{p}$.
\end{itemize}}
\end{remarks}

The proofs for the following well known results concerning
localizations can be found in \cite{Bou2}.

\begin{proposition}\label{Bourbaki1}
Let $A$ and $B$ be two rings with multiplicative sets $S$ and $T$
respectively and $f$ an homomorphism from $A$ to $B$ such that
$f(S)$ is contained in $T$. There exists a unique homomorphism $f'$
from $S^{-1}A$ to $T^{-1}B$ such that $f'(a/1)=f(a)/1$ for every $a
\in A$. Let us suppose now that $T$ is contained in the
multiplicatively closed set of $B$ generated by $f(S)$. If $f$ is
surjective (resp. injective), then $f'$ is also surjective (resp.
injective).
\end{proposition}

\begin{corollary}\label{inclusion in localizations}
Let $A$ and $B$ be two rings with multiplicative sets $S$ and $T$
respectively such that $A \subseteq B$ and $S \subseteq T$. Then
$S^{-1}A \subseteq T^{-1}B$.
\end{corollary}

\begin{proposition}\label{different multiplicative set}
Let $A$ be a ring and $S,T$ two multiplicative sets of $A$ such that
$S \subseteq T$. We have $S^{-1}A = T^{-1}A$ if and only if every
prime ideal of $R$ that meets $T$ meets $S$.
\end{proposition}

The following proposition and its corollary give us information
about the ideals of the localization of a ring $R$ at a
multiplicatively closed subset $S$ of $R$.

\begin{proposition}\label{prime ideals of a localization}
Let $R$ be a ring and let $S$ be a multiplicatively closed subset of
$R$. Then
\begin{enumerate}
  \item Every ideal $\mathfrak{b}'$ of $S^{-1}R$ is of the form
  $S^{-1}\mathfrak{b}$ for some ideal $\mathfrak{b}$ of $R$.
  \item Let $\mathfrak{b}$ be an ideal of $R$ and let $f$ be the
  canonical surjection $R \twoheadrightarrow R/\mathfrak{b}$. Then $f(S)$
  is a multiplicatively closed subset of $R/\mathfrak{b}$ and the
  homomorphism from $S^{-1}R$ to $(f(S))^{-1}(R/\mathfrak{b})$
  canonically associated with $f$ is surjective with kernel
  $\mathfrak{b}'=S^{-1}\mathfrak{b}$. By passing to quotients, an
  isomorphism between $(S^{-1}R)/\mathfrak{b}'$ and
  $(f(S))^{-1}(R/\mathfrak{b})$ is defined.
  \item The application $\mathfrak{b}' \mapsto \mathfrak{b}$,
  restricted to the set of maximal (resp. prime) ideals of
  $S^{-1}R$, is an isomorphism (for the relation of inclusion)
  between this set and the set of maximal (resp. prime) ideals of
  $R$ that do not meet $S$.
  \item If $\mathfrak{q}'$ is a prime ideal of $S^{-1}R$ and
  $\mathfrak{q}$ is the prime ideal of $R$ such that $\mathfrak{q}'=S^{-1}\mathfrak{q}$
  (we have $\mathfrak{q} \cap S = \emptyset$), then there exists an
  isomorphism from $R_\mathfrak{q}$ to $(S^{-1}R)_{\mathfrak{q}'}$
  which sends $r/s$ to $(r/1)/(s/1)$ for $r \in R$,
  $s \in R-\mathfrak{q}$.
  \end{enumerate}
\end{proposition}

\begin{corollary}\label{an}
Let $R$ be a ring, $\mathfrak{p}$ a prime ideal of $R$ and
$S:=R-\mathfrak{p}$. For every ideal $\mathfrak{b}$ of $R$ which
does not meet $S$, let $\mathfrak{b}':= \mathfrak{b}R_\mathfrak{p}$.
Assume that $\mathfrak{b}' \neq R_\mathfrak{p}$. Then
\begin{enumerate}
  \item Let $f$ be the canonical surjection $R\twoheadrightarrow
  R/\mathfrak{b}$. The ring homomorphism from $R_\mathfrak{p}$ to
  $(R/\mathfrak{b})_{\mathfrak{p}/\mathfrak{b}}$ canonically
  associated with $f$ is surjective and its kernel is
  $\mathfrak{b}'$. Thus it defines, by passing to
  quotients, a canonical isomorphism between
  $R_\mathfrak{p}/\mathfrak{b}'$ and
  $(R/\mathfrak{b})_{\mathfrak{p}/\mathfrak{b}}$.
  \item The application $\mathfrak{b}' \mapsto \mathfrak{b}$,
  restricted to the set of prime ideals of
  $R_\mathfrak{p}$, is an isomorphism (for the relation of inclusion)
  between this set and the set of prime ideals of
  $R$ contained in $\mathfrak{p}$ (thus do not meet $S$). Therefore,
  $\mathfrak{p}R_\mathfrak{p}$ is the only maximal ideal of $R_\mathfrak{p}$.
  \item If now $\mathfrak{b}'$ is
  a prime ideal of $R_\mathfrak{p}$, then there exists an isomorphism
  from $R_\mathfrak{b}$ to $(R_\mathfrak{p})_{\mathfrak{b}'}$ which
  sends $r/s$ to $(r/1)/(s/1)$ for $r \in R$, $s \in
  R-\mathfrak{b}$.
\end{enumerate}
\end{corollary}

The notion of localization can also be extended to the modules over
the ring $R$.

\begin{definition}\label{localization of a module}
Let $R$ be a ring and $S$ a multiplicatively closed set of $R$. If
$M$ is an $R$-module, then we call localization of $M$ at $S$ and
denote by $S^{-1}M$ the $S^{-1}R$-module $M \otimes_R S^{-1}R$.
\end{definition}

\section{Integrally closed rings}

\begin{thedef}\label{integral element}
Let $R$ be a ring, $A$ an $R$-algebra and $a$ an element of $A$. The
following properties are equivalent:
\begin{description}
  \item[  (i)] The element $a$ is a root of a monic polynomial with coefficients in
  $R$.
  \item[ (ii)] The subalgebra $R[a]$ of $A$ is an $R$-module of
  finite type.
  \item[(iii)] There exists a faithful $R[a]$-module which is an
  $R$-module of finite type.
\end{description}
If $a \in A$ verifies the conditions above, we say that it is
integral over $R$.
\end{thedef}

\begin{definition}\label{integral closure}
Let $R$ be a ring and $A$ an $R$-algebra. The set of all elements of
$A$ that are integral over $R$ is an $R$-subalgebra of $A$
containing $R$; it is called the integral closure of $R$ in $A$. We
say that $R$ is integrally closed in $A$, if $R$ is an integral
domain and if it coincides with its integral closure in $A$. If now
$R$ is an integral domain and $F$ is its field of fractions, then
the integral closure of $R$ in $F$ is named simply the integral
closure of $R$, and if $R$ is integrally closed in $F$, then $R$ is
said to be integrally closed.
\end{definition}

The following proposition (\cite{Bou5}, \S 1, Prop.13) implies that
transfer theorem holds for integrally closed rings (corollary
$\ref{integrally closed polynomial ring}$).

\begin{proposition}\label{integral closure of a polynomial ring}
  If $R$ is an integral domain, let us denote by $\bar{R}$ the
  integral closure of $R$. Let $x_0, \ldots,x_m$
  be a set of indeterminates over $R$. Then the integral closure of
  $R[x_0,\ldots,x_m]$ is
  $\bar{R}[x_0,\ldots,x_m]$.
\end{proposition}

\begin{corollary}\label{integrally closed polynomial ring}
Let $R$ be an integral domain. Then $R[x_0,\ldots,x_m]$ is
integrally closed if and only if $R$ is integrally closed.
\end{corollary}

\begin{corollary} If $K$ is a field, then every polynomial ring
  $K[x_0,\ldots,x_m]$ is integrally closed.
\end{corollary}

The next proposition (\cite{Bou5}, \S 1, Prop.16) along with its
corollaries treats the
integral closures of localizations of rings.

\begin{proposition}\label{integral closure of a localization}
Let $R$ be a ring, $A$ an $R$-algebra, $\bar{R}$ the integral
closure of $R$ in $A$ and $S$ a multiplicatively closed subset of
$R$ which contains no zero divisors. Then the integral closure of
$S^{-1}R$ in $S^{-1}A$ is $S^{-1}\bar{R}$.
\end{proposition}
\begin{apod}{Let $b/s$ be an element of $S^{-1}\bar{R}$ ($s \in S, b \in \bar{R}$). Since the diagram
$$\begin{array}{ccc}
    R & \hookrightarrow & S^{-1}R \\
    \downarrow &    &  \downarrow \\
    A & \hookrightarrow & S^{-1}A
  \end{array}
$$
commutes, the element $b/1$ is integral over $S^{-1}R$. Since $1/s
\in S^{-1}R$, the element $b/s=(b/1)(1/s)$ is integral over
$S^{-1}R$.

On the other hand, let $a/t$ $(a \in A,t\in S)$ be an element of
$S^{-1}A$ integral over $S^{-1}R$. Then $a/1=(t/1)(a/t)$ is integral
over $S^{-1}R$. This means that there exist $r_i \in R$ $(1 \leq i
\leq n)$ and $s \in S$ such that
$$(a/1)^n+(r_1/s)(a/1)^{n-1}+\ldots+(r_n/s)=0.$$
The above relation can also be written as
$$(sa^n+r_1a^{n-1}+\ldots+r_n)/s=0 $$
and since $S$ contains no zero divisors of $R$, we obtain that
$$sa^n+r_1a^{n-1}+\ldots+r_n=0.$$
Multiplying the above relation with $s^{n-1}$,we deduce that
$$(sa)^n+r_1(sa)^{n-1}+\ldots+s^{n-1}r_n=0.$$
Thus, by definition, we have $sa \in \bar{R}$. Therefore, $a/1 \in
S^{-1}\bar{R}$ and hence, $a/t \in S^{-1}\bar{R}$.}
\end{apod}

\begin{corollary}Let $R$ be an integral domain, $\bar{R}$ the integral
closure of $R$ and $S$ a multiplicatively closed subset of $R$. Then
the integral closure of $S^{-1}R$ is $S^{-1}\bar{R}$.
\end{corollary}

\begin{corollary}\label{integrally closed localization}
If $R$ is an integrally closed domain and $S$ is a multiplicatively
closed subset of $R$, then $S^{-1}R$ is also integrally closed.
\end{corollary}
\subsection*{Lifting prime ideals}

\begin{definition}\label{lying over}
Let $R, R'$ be two rings and let $h:R \rightarrow R'$ be a ring
homomorphism. We say that a prime ideal $\mathfrak{a}'$ of $R'$ lies
over a prime ideal $\mathfrak{a}$ of $R$, if
$\mathfrak{a}=h^{-1}(\mathfrak{a}')$.
\end{definition}

The next result is \cite{Bou5}, \S 2, Proposition 2.

\begin{proposition}\label{primes lying over}
Let $h:R \rightarrow R'$ be a ring homomorphism such that $R'$ is
integral over $R$. Let $\mathfrak{p}$ be a prime ideal of $R$,
$S:=R-\mathfrak{p}$ and $(\mathfrak{p}_i')_{i \in I}$ the family of
all the prime ideals of $R'$ lying over $\mathfrak{p}$. If
$S'=\bigcap_{i \in I}(R'-\mathfrak{p}_i')$, then
$S^{-1}R'=S'^{-1}R'$.
\end{proposition}
\begin{apod}{By definition, we have $h(S) \subseteq S'$ and since
$h(S)^{-1}R' \simeq S^{-1}R'$, it is enough to show that if a prime
ideal $\mathfrak{q}'$ of $R'$ doesn't meet $h(S)$, then it doesn't
meet $S'$ either (see proposition $\ref{different multiplicative
set}$). Let us suppose that $\mathfrak{q}' \cap h(S) = \emptyset$
and let $\mathfrak{q}:=h^{-1}(\mathfrak{q}')$. Then we have
$\mathfrak{q} \cap S = \emptyset$, which means that $\mathfrak{q}
\subseteq \mathfrak{p}$. Since $\mathfrak{q}'$ is lying over
$\mathfrak{q}$ by definition, there exists an index $i \in I$ such
that $\mathfrak{q}' \subseteq \mathfrak{p}_i'$. Therefore,
$\mathfrak{q}' \cap S' = \emptyset$.}
\end{apod}

The following corollary deals with a case we will encounter in a
following chapter, where there exists a unique prime ideal lying
over the prime ideal $\mathfrak{p}$ of $R$. In combination with
proposition $\ref{integral closure of a localization}$, proposition
$\ref{primes lying over}$ implies that

\begin{corollary}\label{one prime lying over}
Let $R$ be an integral domain, $A$ an $R$-algebra, $\bar{R}$ the
integral closure of $R$ in $A$. Let $\mathfrak{p}$ be a prime ideal
of $R$ and $S:=R-\mathfrak{p}$. If there exists a unique prime ideal
$\bar{\mathfrak{p}}$ of $\bar{R}$ lying over $\mathfrak{p}$, then
the integral closure of $R_\mathfrak{p}$ in $S^{-1}A$ is
$\bar{R}_{\bar{\mathfrak{p}}}$.
\end{corollary}
\subsection*{Valuations}

\begin{definition}\label{valuation}
Let $R$ be a ring and $\Gamma$ a totally ordered abelian group. We
call valuation of $R$ with values in $\Gamma$ every application $v:R
\rightarrow \Gamma \cup \{\infty\}$ which satisfies the following
properties:
\begin{description}
\item[(V1)] $v(xy)=v(x)+v(y)$ for $x \in R, y \in R$.
\item[(V2)] $v(x+y) \geq \mathrm{inf}(v(x),v(y))$ for $x \in R, y \in R$.
\item[(V3)] $v(1)=0$ and $v(0)=\infty$.
\end{description}
\end{definition}
In particular, if $v(x) \neq v(y)$, property (V2) gives
$v(x+y)=\textrm{inf}(v(x),v(y))$ for $x \in R, y \in R$. Moreover,
from property (V1), we have that if $z \in R$ with $z^n=1$ for some
integer $n \geq 1$, then $nv(z)=v(z^n)=v(1)=0$ and thus $v(z)=0$.
Consequently, $v(-x)=v(-1)+v(x)=v(x)$ for all $x \in R$.\\

Now let $F$ be a field and let $v:F \rightarrow \Gamma$ be a
valuation of $F$. The set $A$ of $a \in F$ such that $v(a) \geq 0$
is a local subring of $F$. Its maximal ideal $\mathfrak{m}(A)$ is
the set of $a \in A$ such that $v(a) > 0$. For all $a \in F-A$, we
have $a^{-1} \in \mathfrak{m}(A)$. The ring $A$ is called
\emph{the ring of the valuation} $v$ on $F$.\\

We will now introduce the notion of a valuation ring. For more
information about valuation rings and their properties, see
\cite{Bou6}. Some of them will also be discussed in Chapter 2,
Section 2.4.

\begin{definition}\label{valuation ring}
Let $R$ be an integral domain contained in a field $F$. Then $R$ is
a valuation ring if for all non-zero element $x \in F$, we have $x
\in R$ or $x^{-1} \in R$. Consequently, $F$ is the field of
fractions of $R$.
\end{definition}

If $R$ is a valuation ring, then it has the following properties:
\begin{itemize}
  \item It is an integrally closed local ring.
  \item The set of the principal ideals of $R$ is totally ordered by the relation of inclusion.
  \item The set of the ideals of $R$ is totally ordered by the relation of inclusion.
\end{itemize}

Let $R$ be a valuation ring and $F$ its field of fractions. Let us
denote by $R^\times$ the set of units of $R$. Then the set
$\Gamma_R:=F^\times/R^\times$ is an abelian group, totally ordered
by the relation of inclusion of the corresponding principal ideals.
If we denote by $v_R$ the canonical homomorphism of $F^\times$ onto
$\Gamma_R$ and set $v_R(0)=\infty$, then $v_R$ is a valuation of $F$
whose ring is $R$.

The following proposition gives a characterization of integrally
closed rings in terms of valuation rings (\cite{Bou6}, \S 1, Thm.
3).

\begin{proposition}\label{intersection of valuation rings}
Let $R$ be a subring of a field $F$. The integral closure $\bar{R}$
of $R$ in $F$ is the intersection of all valuation rings in $F$
which contain $R$. Consequently, an integral domain $R$ is
integrally closed if and only if it is the intersection of a family
of valuation rings contained in its field of fractions.
\end{proposition}

This characterization helped us to prove the following result about
integrally closed rings.

\begin{proposition}\label{askhsh 12}
Let $R$ be an integrally closed ring and $f(x)=\sum_ia_ix^i$,
$g(x)=\sum_jb_jx^j$ be two polynomials in $R[x]$. If there exists an
element $c \in R$ such that all the coefficients of $f(x)g(x)$
belong to $cR$, then all the products $a_ib_j$ belong to $cR$.
\end{proposition}
\begin{apod}{Due to the proposition $\ref{intersection of valuation rings}$,
it is enough to prove
the result in the case where $R$ is a valuation ring.

From now on, let $R$ be a valuation ring. Let $v$ be a valuation of
the field of fractions of $R$ such that the ring of valuation of $v$
is $R$. Let $\kappa:=\textrm{inf}_i(v(a_i))$ and
$\lambda:=\textrm{inf}_j(v(b_j))$. Then
$\kappa+\lambda=\textrm{inf}_{i,j}(v(a_ib_j))$. We will show that
$\kappa + \lambda \geq v(c)$ and thus $c$ divides all the products
$a_ib_j$. Argue by contradiction and assume that $\kappa + \lambda <
v(c)$. Let $a_{i_1},a_{i_2},\ldots,a_{i_r}$ with
$i_1<i_2<\ldots<i_r$ be all the elements among the $a_i$ with
valuation equal to $\kappa$. Respectively, let
$b_{j_1},b_{j_2},\ldots,b_{j_s}$ with $j_1<j_2<\ldots<j_s$ be all
the elements among the $b_j$ with valuation equal to $\lambda$. We
have that $i_1+j_1<i_m+j_n$, $\forall (m,n) \neq (1,1)$. Therefore,
the coefficient $c_{i_1+j_1}$ of $x^{i_1+j_1}$ in $f(x)g(x)$ is of
the form $(a_{i_1}b_{j_1}+\sum(\textrm{terms with valuation} >
\kappa+\lambda))$ and since $v(a_{i_1}b_{j_1}) \neq
v(\sum(\textrm{terms with valuation} > \kappa+\lambda))$, we obtain
that $$v(c_{i_1+j_1})=\textrm{inf}(v(a_{i_1}b_{j_1}),
v(\sum(\textrm{terms with valuation} >
\kappa+\lambda)))=\kappa+\lambda.$$ However, since all the
coefficients of $f(x)g(x)$ are divisible by $c$, we have that
$v(c_{i_1+j_1}) \geq v(c) > \kappa + \lambda$, which is a
contradiction.}
\end{apod}

The propositions $\ref{quotient}$
and $\ref{my second lemma one variable}$ derive from the one above. We will
make use of the results in corollaries $\ref{porisma porismatos}$
and $\ref{my second lemma}$ in Chapter 3.

\begin{proposition}\label{quotient}
Let $R$ be an integrally closed domain and let $F$ be its field of
fractions. Let $\mathfrak{p}$ be a prime ideal of $R$. Then
$$(R[x])_{\mathfrak{p}R[x]} \cap F[x] = R_\mathfrak{p}[x].$$
\end{proposition}
\begin{apod}{The inclusion $R_\mathfrak{p}[x] \subseteq (R[x])_{\mathfrak{p}R[x]} \cap
F[x]$ is obvious. Now, let $f(x)$ be an element of $F[x]$. Then
$f(x)$ can be written in the form $r(x)/\xi$, where $r(x) \in R[x]$
and $\xi \in R$. If, moreover, $f(x)$ belongs to
$(R[x])_{\mathfrak{p}R[x]}$, then there exist $s(x), t(x) \in R[x]$
with $t(x) \notin \mathfrak{p}R[x]$ such that $f(x) = s(x)/t(x)$.
Thus we have $$f(x)=\frac{r(x)}{\xi}= \frac{s(x)}{t(x)}.$$ All the
coefficients of the product $r(x)t(x)$ belong to $\xi R$. Due to
proposition \ref{askhsh 12}, if $r(x)=\sum_ia_ix^i$ and
$t(x)=\sum_jb_jx^j$, then all the products $a_ib_j$ belong to $\xi
R$. Since $t(x)\notin \mathfrak{p}R[x]$, there exists $j_0$ such
that $b_{j_0} \notin \mathfrak{p}$ and $a_ib_{j_0} \in \xi R,
\forall i$. Consequently, $b_{j_0}f(x)=b_{j_0}(r(x)/\xi) \in R[x]$
and hence all the coefficients of $f(x)$ belong to
$R_\mathfrak{p}$.}
\end{apod}

\begin{corollary}\label{porisma porismatos}
Let $R$ be an integrally closed domain and let $F$ be its field of
fractions. Let $\mathfrak{p}$ be a prime ideal of $R$. Then
\begin{enumerate}
  \item $(R[x,x^{-1}])_{\mathfrak{p}R[x,x^{-1}]} \cap F[x,x^{-1}] =
R_\mathfrak{p}[x,x^{-1}].$
  \item $(R[x_0,\ldots,x_m])_{\mathfrak{p}R[x_0,\ldots,x_m]}
\cap F[x_0,\ldots,x_m] = R_\mathfrak{p}[x_0,\ldots,x_m].$
  \item $(R[x_0^{\pm1},\ldots,x_m^{\pm1}])_{\mathfrak{p}R[x_0^{\pm1},\ldots,x_m^{\pm1}]}
\cap F[x_0^{\pm1},\ldots,x_m^{\pm1}] =
R_\mathfrak{p}[x_0^{\pm1},\ldots,x_m^{\pm1}]$.
\end{enumerate}
\end{corollary}

\begin{proposition}\label{my second lemma one variable}
Let $R$ be an integrally closed domain and let $F$ be its field of
fractions. Let $r(x)$ and $s(x)$ be two elements of $R[x]$ such that
$s(x)$ divides $r(x)$ in $F[x]$. If one of the coefficients of
$s(x)$ is a unit in $R$, then $s(x)$ divides $r(x)$ in $R[x]$.
\end{proposition}
\begin{apod}{Since $s(x)$ divides $r(x)$ in $F[x]$, there exists an
element of the form $t(x)/\xi$ with $t(x) \in R[x]$ and $\xi \in R$
such that
$$r(x) = \frac{s(x)t(x)}{\xi}.$$
All the coefficients of the product $s(x)t(x)$ belong to $\xi R$.
Due to proposition \ref{askhsh 12}, if $s(x)=\sum_ia_ix^i$ and
$t(x)=\sum_jb_jx^j$, then all the products $a_ib_j$ belong to $\xi
R$. By assumption, there exists $i_0$ such that $a_{i_0}$ is a unit
in $R$ and $a_{i_0}b_j \in \xi R, \forall j$. Consequently, $b_j \in
\xi R, \forall j$ and thus $t(x)/\xi \in R[x]$.}
\end{apod}

\begin{corollary}\label{my second lemma}
Let $R$ be an integrally closed domain and let $F$ be its field of
fractions. Let $r, s$ be two elements of
$R[x_0^{\pm1},\ldots,x_m^{\pm1}]$ such that $s$ divides $r$ in
$F[x_0^{\pm1},\ldots,x_m^{\pm1}]$. If one of the coefficients of $s$
is a unit in $R$, then $s$ divides $r$ in
$R[x_0^{\pm1},\ldots,x_m^{\pm1}]$.
\end{corollary}
\subsection*{Discrete valuation rings and Krull rings}

\begin{definition}\label{discrete valuation}
Let $F$ be a field, $\Gamma$ a totally ordered abelian group and $v$
a valuation of $F$ with values in $\Gamma$. We say that the
valuation $v$ is discrete, if $\Gamma$ is isomorphic to
$\mathbb{Z}$.
\end{definition}

\begin{thedef}\label{dvr}
An integral domain $R$ is a discrete valuation ring, if it satisfies
one of the following equivalent conditions:
\begin{description}
  \item[  (i)] $\textrm{ }R$ is the ring of a discrete valuation.
  \item[ (ii)] $R$ is a local Dedekind ring.
  \item[(iii)] $R$ is a local principal ideal domain.
  \item[(iv)] $R$ is a Noetherian valuation ring.
\end{description}
\end{thedef}

By proposition $\ref{intersection of valuation rings}$, integrally
closed rings are intersections of valuation rings. Krull rings are
essentially intersections of discrete valuation rings.

\begin{definition}\label{Krull ring}
An integral domain $R$ is a Krull ring, if there exists a family of
valuations $(v_i)_{i \in I}$ of the field of fractions $F$ of $R$
with the following properties:
\begin{description}
  \item[(K1)] The valuations $(v_i)_{i \in I}$ are discrete.
  \item[(K2)] The intersection of the rings of $(v_i)_{i \in I}$ is
  $R$.
  \item[(K3)] For all $x\in F^\times$,there exists a finite number of
  $i \in I$ such that $v_i(x) \neq 0$.
\end{description}
\end{definition}

The proofs of the following results and more information about Krull
rings can be found in \cite{Bou7}, \S 1.

\begin{theorem}\label{Krull-dvr}
Let $R$ be an integral domain and let $\mathrm{Spec}_1(R)$ be the
set of its prime ideals of height $1$. Then $R$ is a Krull ring if and
only if the following properties are satisfied:
\begin{enumerate}
  \item For all $\mathfrak{p} \in \mathrm{Spec}_1(R)$, $R_\mathfrak{p}$ is a
  discrete valuation ring.
  \item $R$ is the intersection of $R_\mathfrak{p}$ for all
  $\mathfrak{p} \in \mathrm{Spec}_1(R)$.
  \item For all $r \neq 0$ in $R$, there exists a finite
  number of ideals $\mathfrak{p} \in \mathrm{Spec}_1(R)$ such that $r \in
  \mathfrak{p}$.
\end{enumerate}
\end{theorem}

Transfer theorem holds also for Krull rings.

\begin{proposition}\label{prime ideals of height 1}
Let $R$ be a Krull ring, $F$ the field of fractions of $R$ and $x$
an indeterminate. Then $R[x]$ is also a Krull ring. Moreover, its
prime ideals of height $1$ are:
\begin{itemize}
  \item the prime ideals of the form $\mathfrak{p}R[x]$, where
   $\mathfrak{p}$ is a prime ideal of height $1$ of $R$,
  \item the prime ideals of the form $\mathfrak{m} \cap R[x]$, where
   $\mathfrak{m}$ is a prime ideal of $F[x]$.
\end{itemize}
\end{proposition}

The following proposition provides us with a simple characterization
of Krull rings, when they are Noetherian.

\begin{proposition}\label{case of Krull}
Let $R$ be a Noetherian ring. Then $R$ is a Krull ring if and only
if it is integrally closed.
\end{proposition}

\begin{px} \small{\emph{Let $K$ be a finite field extension of $\mathbb{Q}$ and
$\mathbb{Z}_K$ the integral closure of $\mathbb{Z}$ in $K$. The ring
$\mathbb{Z}_K$ is a Dedekind ring and thus Noetherian and integrally
closed. Let $x_0,x_1,\ldots,x_m$ be indeterminates. Then the ring
$\mathbb{Z}_K[x_0^{\pm},x_1^{\pm},\ldots,x_m^{\pm}]$ is also
Noetherian and integrally closed and thus a Krull ring.}}
\end{px}

\section{Completions}

For all the following results concerning completions, the reader can
refer to \cite{Na}, Chapter II.

Let $I$ be an ideal of a commutative ring $R$ and let $M$ be an
$R$-module. We introduce a topology on $M$ such that the open sets
of $M$ are unions of an arbitrary number of sets of the form $m+I^nM
(m\in M)$. This topology is called the $I$\emph{-adic topology} of
$M$.

\begin{theorem}\label{16.5}
If $M$ is a Noetherian $R$-module, then for any submodule $N$ of
$M$, the $I$-adic topology of $N$ coincides with the topology of $N$
as a subspace of $M$ with the $I$-adic topology.
\end{theorem}

From now on, we will concentrate on semi-local rings and in
particular, on Noetherian semi-local rings.

\begin{definition}\label{semilocal ring}
A ring $R$ is called semi-local, if it has only a finite number of
maximal ideals. The Jacobson radical $\mathfrak{m}$ of $R$ is the
intersection of the maximal ideals of $R$.
\end{definition}

\begin{theorem}\label{17.6}
Assume that $R$ is a Noetherian semi-local ring with Jacobson
radical $\mathfrak{m}$ and let $\hat{R}$ be the completion of $R$
with respect to the $\mathfrak{m}$-adic topology. Then $\hat{R}$ is
also a semi-local ring and we have $R \subseteq \hat{R}$.
\end{theorem}

\begin{theorem}\label{17.8}
Assume that $R$ is a Noetherian semi-local ring with Jacobson
radical $\mathfrak{m}$ and that $M$ is a finitely generated
$R$-module. Let $\hat{R}$ be the completion of $R$ with respect to
the $\mathfrak{m}$-adic topology. Endow $M$ with the
$\mathfrak{m}$-adic topology. Then $M \otimes_R \hat{R}$ is the
completion of $M$ with respect to that topology.
\end{theorem}

\begin{corollary}\label{17.9}
Let $\mathfrak{a}$ be an ideal of a Noetherian semi-local ring $R$
with Jacobson radical $\mathfrak{m}$. Let $\hat{R}$ be the
completion of $R$ with respect to the $\mathfrak{m}$-adic topology.
Then the completion of $\mathfrak{a}$ is $\mathfrak{a}\hat{R}$ and
$\mathfrak{a}\hat{R}$ is isomorphic to $\mathfrak{a} \otimes_R
\hat{R}$. Furthermore, $\mathfrak{a}\hat{R} \cap R = \mathfrak{a}$
and $\hat{R}/\mathfrak{a}\hat{R}$ is the completion of
$R/\mathfrak{a}$ with respect to the $\mathfrak{m}$-adic topology.
\end{corollary}

\begin{theorem}\label{18.4}
Assume that $R$ is a Noetherian semi-local ring with Jacobson
radical $\mathfrak{m}$ and let $\hat{R}$ be the completion of $R$
with respect to the $\mathfrak{m}$-adic topology. Then
\begin{enumerate}
  \item The total quotient ring $F$ of $R$ (localization of $R$ at
  the set of non-zero divisors) is naturally a subring of the total
  quotient ring $\hat{F}$ of $\hat{R}$.
  \item For any ideal $\mathfrak{a}$ of $R$, $\mathfrak{a}\hat{R} \cap
  F =\mathfrak{a}$.
\end{enumerate}
In particular, $\hat{R} \cap F =R$.
\end{theorem}

\section{Morphisms associated with monomials}

We have the following elementary algebra result

\begin{thedef}\label{similitude}
Let $R$ be an integral domain and $M$ a free $R$-module of basis
$(e_i)_{0 \leq i \leq m}$.  Let $x=r_0e_0+r_1e_1+\ldots+r_me_m$ be a
non-zero element of $M$. We set $M^*:=\mathrm{Hom}_R(M,R)$ and
$M^*(x):= \{\varphi(x) \,|\, \varphi \in M^*\}$. Then the following
assertions are equivalent:
\begin{description}
  \item [  (i)] $\textrm{ }M^*(x)=R$.
  \item [ (ii)] $\sum_{i=0}^m Rr_i =R$.
  \item [(iii)] There exists $\varphi \in M^*$ such that $\varphi(x)=1$.
  \item [(iv)] There exists an $R$-submodule $N$ of $M$ such that $M=Rx\oplus
  N$.
\end{description}
If $x$ satisfies the conditions above, we say that it is a primitive
element of $M$.
\end{thedef}
\begin{apod}{
\begin{description}
  \item[$(i) \Leftrightarrow (ii)$]
  Let $(e_i^*)_{0 \leq i \leq m}$ be the basis of $M^*$ dual to
  $(e_i)_{0 \leq i \leq  m}$.
  Then $M^*(x)$ is generated by $(e_i^*(x))_{0 \leq i \leq m}$
  and $e_i^*(x)=r_i$.
  \item[$(ii) \Rightarrow (iii)$]
  There exist $u_0,u_1,\ldots,u_m \in R$ such that
  $\sum_{i=0}^mu_ir_i=1$. If $\varphi:=\sum_{i=0}^mu_ie_i^*$, then
  $\varphi(x)=1$.
  \item[$(iii) \Rightarrow (iv)$]
  Let $N:=\mathrm{Ker}\varphi$. For all $y \in M$, we have
  $$y = \varphi(y)x + (y -\varphi(y)x).$$
  Since $\varphi(x)=1$, we have $y \in Rx + N$. Obviously, $Rx
  \cap N =\{0\}.$
  \item[$(iv) \Rightarrow (i)$]
  Since $M$ is free and $R$ is an integral domain, $x$ is
  torsion-free (otherwise there exists $r \in R, r\neq 0$ such that
  $\sum_{i=0}^m (rr_i)e_i=0$). Therefore, the map
  $$R \rightarrow M, r \mapsto rx$$
  is an isomorphism of $R$-modules. Its inverse is a linear form on
  $Rx$ which sends $x$ to $1$. Composing it with the map
  $$M \twoheadrightarrow M/N \rightarrow R,$$
  we obtain a linear form $\varphi:M\rightarrow R$ such that
  $\varphi(x)=1$. We have $1 \in M^*(x)$ and thus $M^*(x)=R$.}
\end{description}
\end{apod}

We will apply the above result to the $\mathbb{Z}$-module
$\mathbb{Z}^{m+1}$. Let us consider $(e_i)_{0 \leq i \leq m}$ the
standard basis of $\mathbb{Z}^{m+1}$ and let
$a:=a_0e_0+a_1e_1+\ldots+a_me_m$ be an element of $\mathbb{Z}^{m+1}$
such that $\mathrm{gcd}(a_i)=1$. Then, by Bezout's theorem, there
exist $u_0,u_1,\ldots,u_m \in \mathbb{Z}$ such that
$\sum_{i=0}^mu_ia_i=1$ and hence $\sum_{i=0}^m
\mathbb{Z}a_i=\mathbb{Z}$. By theorem $\ref{similitude}$, there
exists a $\mathbb{Z}$-submodule $N_a$ of $\mathbb{Z}^{m+1}$ such
that $\mathbb{Z}^{m+1}=\mathbb{Z}a\oplus N_a$. In particular,
$N_a=\mathrm{Ker}\varphi_a$, where
$\varphi_a:=\sum_{i=0}^mu_ie_i^*$. We will denote by
$p_a:\mathbb{Z}^{m+1} \twoheadrightarrow N_a$ the projection of
$\mathbb{Z}^{m+1}$ onto $N_a$ such that
$\mathrm{Ker}p_a=\mathbb{Z}a$. We have a $\mathbb{Z}$-module
isomorphism
$i_a:N_a \tilde{\rightarrow} \mathbb{Z}^m$. Then
$f_a:=i_a \circ p_a$ is a surjective $\mathbb{Z}$-module
morphism $\mathbb{Z}^{m+1} \twoheadrightarrow \mathbb{Z}^m$
with $\mathrm{Ker}f_a=\mathbb{Z}a$.

Now let $R$ be an integral domain and let $x_0,x_1,\ldots,x_m$ be
$m+1$ indeterminates over $R$. Let $G$ be the abelian group
generated by all the monomials in $R[x_0^{\pm 1},x_1^{\pm
1},\ldots,x_m^{\pm 1}]$ with group operation the multiplication.
Then $G$ is isomorphic to the additive group $\mathbb{Z}^{m+1}$ by
the isomorphism defined as follows
$$\begin{array}{rccc}
   \theta_G: &G & \tilde{\rightarrow} & \mathbb{Z}^{m+1} \\
    &\prod_{i=0}^m x_i^{l_i} & \mapsto & (l_0,l_1,\ldots,l_m).
  \end{array}$$

\begin{lemma}\label{work on Z}
We have $R[x_0^{\pm 1},x_1^{\pm 1},\ldots,x_m^{\pm 1}]=R[G] \simeq
R[\mathbb{Z}^{m+1}]$.
\end{lemma}

Respectively, if $y_1,\ldots,y_m$ are $m$ indeterminates over $R$
and $H$ is the group generated by all the monomials in $R[y_1^{\pm
1},\ldots,y_m^{\pm1}]$, then $H \simeq \mathbb{Z}^m$ and $R[y_1^{\pm
1},\ldots,y_m^{\pm1}]=R[H] \simeq R[\mathbb{Z}^m].$

The morphism $F_a:={\theta_H}^{-1} \circ f_a \circ \theta_G: G
\twoheadrightarrow H$ induces an $R$-algebra morphism
$$\begin{array}{cccc}
    \varphi_a: & R[G] & \rightarrow & R[H] \\
        & \sum_{g \in G}a_gg & \mapsto     & \sum_{g \in G}a_gF_a(g)
  \end{array}$$
Since $F_a$ is surjective, the morphism $\varphi_a$ is also surjective.
Moreover, $\mathrm{Ker}\varphi_a$ is generated (as an $R$-module) by
the set
$$< g-1 \,|\, g \in G \textrm{ such that } \theta_G(g) \in \mathbb{Z}a >.$$

How do we translate this in the polynomial language?\\

Let $A:=R[x_0^{\pm 1},x_1^{\pm 1},\ldots,x_m^{\pm 1}]$ and
$B:=R[y_1^{\pm 1},y_2^{\pm 1},\ldots,y_m^{\pm 1}]$. The map
$\varphi_a$ is a surjective $R$-algebra morphism from $A$ to $B$
with $\mathrm{Ker}\varphi_a=(\prod_{i=0}^mx_i^{a_i}-1)A$ such that
for every monomial $N$ in $A$, $\varphi_a(N)$ is a monomial in $B$.

\begin{definition}\label{associated morphism}
Let $M:=\prod_{i=0}^mx_i^{a_i}$ be a monomial in $A$ with
$\mathrm{gcd}(a_i)=1$. An $R$-algebra morphism $\varphi_M:A
\rightarrow B$ defined as above will be called associated with the
monomial $M$.
\end{definition}

\begin{px}\label{x^5y^-3z^-2}
\small{\emph{Let $A:=R[X^{\pm 1},Y^{\pm 1},Z^{\pm 1}]$ and
$M:=X^5Y^{-3}Z^{-2}$. We have $x:=(5,-3,-2)$ and$$ (-1)\cdot
5+(-2)\cdot(-3)+0 \cdot(-2)=1.$$ Hence, with the notations of the
proof of theorem $\ref{similitude}$, we have
$$u_0=-1,u_1=-2,u_3=0.$$
The map $\varphi: \mathbb{Z}^3 \rightarrow \mathbb{Z}$ defined as
$$\varphi:=-e_0^*-2e_1^*.$$
has $\mathrm{Ker}\varphi=\{(x_0,x_1,x_2) \in \mathbb{Z}^3 \,|\, x_0
= -2x_1\} = \{(-2a,a,b) \,|\, a,b \in \mathbb{Z}\}=:N$.
\\
\\
By theorem  $\ref{similitude}$, we have $\mathbb{Z}^3=\mathbb{Z}x
\oplus N$
and the projection $p:\mathbb{Z}^3\twoheadrightarrow N$ is the map
$$(y_0,y_1,y_2)=:y \mapsto
y-\varphi(y)x=(6y_0+10y_1,-3y_0-5y_1,-2y_0-4y_1+y_2).$$ The
$\mathbb{Z}$-module $N$ is obviously isomorphic to $\mathbb{Z}^2$
via
$$i:(-2a,a,b) \mapsto (a,b).$$
Composing the two previous maps, we obtain a well defined surjection
$$\begin{array}{ccc}
    \mathbb{Z}^3 & \twoheadrightarrow & \mathbb{Z}^2 \\
    (y_0,y_1,y_2) & \mapsto & (-3y_0-5y_1,-2y_0-4y_1+y_2).
  \end{array}$$
The above surjection induces (in the way described before) an
$R$-algebra epimorphism
$$\begin{array}{cccc}
    \varphi_M: & R[X^{\pm 1},Y^{\pm 1},Z^{\pm 1}] & \twoheadrightarrow & R[X^{\pm 1},Y^{\pm 1}] \\
     & X & \mapsto & X^{-3}Y^{-2} \\
     & Y & \mapsto & X^{-5}Y^{-4} \\
     & Z & \mapsto & Y
  \end{array}$$
By straightforward calculations, we can verify that
$\mathrm{Ker}\varphi_M=(M-1)A$.}}
\end{px}

\begin{lemma}\label{primeness of q}
Let $M:=\prod_{i=0}^mx_i^{a_i}$ be a monomial in $A$ such that
$\mathrm{gcd}(a_i)=1$. Then
\begin{enumerate}
  \item The ideal $(M-1)A$ is a prime ideal of $A$.
  \item If $\mathfrak{p}$ is a prime ideal of $R$, then the ideal
$\mathfrak{q}_M:=\mathfrak{p}A+(M-1)A$ is also prime in $A$.
\end{enumerate}
\end{lemma}
\begin{apod}{
\begin{enumerate}
\item Let $\varphi_M:A \rightarrow B$ be a morphism associated with
$M$. Then $\varphi_M$ is surjective and
$\mathrm{Ker}\varphi_M=(M-1)A$. Since $B$ is an integral domain, the
ideal generated by $(M-1)$ is prime in $A$.
\item Set $R':=R / \mathfrak{p}$,
$A':=R'[x_0^{\pm 1},x_1^{\pm 1},\ldots,x_m^{\pm 1}]$, $B':=R'[
y_1^{\pm 1},\ldots,y_m^{\pm 1}]$ . Then $R'$ is an integral domain
and
$$A/\mathfrak{q}_M \simeq A'/(M-1)A' \simeq B'.$$
Since $B'$ is an integral domain, the ideal $\mathfrak{q}_M$ is
prime in $A$.}
\end{enumerate}
\end{apod}

The following assertions are now straightforward. Nevertheless, they
are stated for further reference.

\begin{proposition}\label{properties of phi}
Let $M:=\prod_{i=0}^mx_i^{a_i}$ be a monomial in $A$ with
$\mathrm{gcd}(a_i)=1$ and let $\varphi_M:A \rightarrow B$ be a
morphism associated with $M$. Let $\mathfrak{p}$ be a prime ideal of
$R$ and set $\mathfrak{q}_M:=\mathfrak{p}A+(M-1)A$. Then the
morphism $\varphi_M$ has the following properties:
\begin{enumerate}
  \item  If $f \in A$, then $\varphi_M(f) \in \mathfrak{p}B$ if and only if $f \in \mathfrak{q}_M$.
  Corollary $\ref{an}$\emph{(1)} implies that $$A_{\mathfrak{q}_M}/(M-1) A_{\mathfrak{q}_M} \simeq
   B_{\mathfrak{p}B}.$$
  \item  If $N$ is a monomial in $A$, then $\varphi_M(N)=1$ if and only if
        there exists $k \in \mathbb{Z}$ such that $N=M^k$.
\end{enumerate}
\end{proposition}

\begin{corollary}\label{q intersection zk}
$\mathfrak{q}_M \cap R = \mathfrak{p}$.
\end{corollary}
\begin{apod}{Obviously $\mathfrak{p} \subseteq \mathfrak{q}_M \cap R$.
 Let $\xi \in R$ such that $\xi \in \mathfrak{q}_M$. If
$\varphi_M$ is a morphism associated with $M$, then, by proposition
\ref{properties of phi}, $\varphi_M(\xi) \in \mathfrak{p}B$. But
$\varphi_M(\xi)=\xi$ and $\mathfrak{p}B \cap R = \mathfrak{p}$. Thus
$\xi \in \mathfrak{p}$.}
\end{apod}\\
\begin{remark}
\emph{If $m=0$ and we set $x:=x_0$, then $A:=R[x,x^{-1}]$ and
$B:=R$. The only monomials that we can associate a morphism $A
\rightarrow B$ with are $x$ and $x^{-1}$. This morphism is unique
and given by $x \mapsto 1$.}
\end{remark}\
\\

The following lemma, whose proof is straightforward when arguing by
contradiction, will be used in the proofs of propositions
$\ref{surjective is associated}$ and $\ref{surjective is adapted}$.

\begin{lemma}\label{groupsGH}
Let $G$, $H$ be two groups and $p:G \rightarrow H$ a group
homomorphism. If $R$ is an integral domain, let us denote by
$p_R:R[G] \rightarrow R[H]$ the $R$-algebra morphism induced by $p$.
If $p_R$ is surjective, then $p$ is also surjective.
\end{lemma}

\begin{proposition}\label{surjective is associated}
Let $\varphi:A\rightarrow B$ be a surjective $R$-algebra morphism
such that for every monomial $M$ in $A$, $\varphi(M)$ is a monomial
in $B$. Then $\varphi$ is associated with a monomial in $A$.
\end{proposition}
\begin{apod}{Due to the isomorphism of lemma $\ref{work on Z}$,
$\varphi$ can be considered as a surjective $R$-algebra morphism
$$\varphi: R[\mathbb{Z}^{m+1}] \rightarrow R[\mathbb{Z}^{m}].$$ The
property of $\varphi$ about the monomials implies that the above
morphism is induced by a $\mathbb{Z}$-module morphism $f:
\mathbb{Z}^{m+1} \rightarrow \mathbb{Z}^m$, which is also surjective
by lemma $\ref{groupsGH}$. Since $\mathbb{Z}^m$ is a free
$\mathbb{Z}$-module the following exact sequence sequence splits
$$0\rightarrow\mathrm{Ker}f\rightarrow\mathbb{Z}^{m+1}
\rightarrow \mathbb{Z}^m\rightarrow0$$ and we obtain that
$\mathbb{Z}^{m+1}\simeq \mathrm{Ker}f\,\oplus\,\mathbb{Z}^m$.
Therefore, $\mathrm{Ker}f$ is a $\mathbb{Z}$-module of rank $1$ and
there exists $a:=(a_0,a_1,\ldots,a_m)\in\mathbb{Z}^{m+1}$ such that
$\mathrm{Ker}f=\mathbb{Z}a$. By theorem $\ref{similitude}$, $a$ is a
primitive element of $\mathbb{Z}^{m+1}$ and we must have
$\sum_{i=0}^m \mathbb{Z}a_i=\mathbb{Z}$, whence
$\mathrm{gcd}(a_i)=1$. By definition, the morphism $\varphi$ is
associated with the monomial $\prod_{i=0}^m x_i^{a_i}$.}
\end{apod}

Now let $r \in \{1,\ldots,m+1\}$ and $\mathcal{R}:=\left\{
  \begin{array}{ll}
    R[y_r^{\pm1},\ldots,y_m^{\pm 1}], & \textrm{for } 1\leq r \leq m; \\
    R, & \textrm{for } r=m+1,
  \end{array}
\right.$ where $y_r,\ldots,y_m$ are $m-r+1$ indeterminates over $R$.

\begin{definition}\label{adapted morphism}
An $R$-algebra morphism $\varphi:A \rightarrow \mathcal{R}$ is
called adapted, if $\varphi={\varphi_r} \circ {\varphi_{r-1}} \circ
\ldots \circ {\varphi_1}$, where $\varphi_i$ is a morphism
associated with a monomial for all $i=1,\ldots,r$. The family
$\mathcal{F}:=\{\varphi_r,\varphi_{r-1},\ldots,\varphi_1\}$ is
called an adapted family for $\varphi$ whose initial morphism is
$\varphi_1$.
\end{definition}

Let us introduce the following notation: If
$M:=\prod_{i=0}^mx_i^{c_i}$ is a monomial such that
$\textrm{gcd}(c_i)=d \in \mathbb{Z}$, then
$M^\circ:=\prod_{i=0}^mx_i^{c_i/d}$.

\begin{proposition}\label{change initial}
Let $\varphi:A \rightarrow \mathcal{R}$ be an adapted morphism and
$M$ a monomial in $A$ such that $\varphi(M)=1$. Then there exists an
adapted family for $\varphi$ whose initial morphism is associated
with $M^\circ$.
\end{proposition}
\begin{apod}{Let $M:=\prod_{i=0}^mx_i^{c_i}$ be a monomial in $A$
such that $\varphi(M)=1$. Note that $\varphi(M)=1$ if and only if
$\varphi(M^\circ)=1.$ Therefore, we can assume that
$\textrm{gcd}(c_i)=1$. We will prove the desired result by induction
on $r$.
\begin{itemize}
  \item For $r=1$, due to property $\ref{properties of phi}$(2), $\varphi$ must be a morphism associated with $M$.
  \item For $r=2$, set $B:=R[z_1^{\pm1},\ldots,z_m^{\pm 1}]$.
  Let $\varphi:=\varphi_b \circ \varphi_a$,
  where \begin{itemize}
  \item $\varphi_a: A \rightarrow B$
        is a morphism associated with a monomial
        $\prod_{i=0}^mx_i^{a_i}$ in $A$ such that
        $\textrm{gcd}(a_i)=1$.
  \item $\varphi_b: B \rightarrow \mathcal{R}$
        is a morphism associated with a monomial
        $\prod_{j=1}^mz_j^{b_j}$ in $B$ such that
        $\textrm{gcd}(b_j)=1$.
\end{itemize}
By theorem $\ref{similitude}$, the element $a:=(a_0,a_1,\ldots,a_m)$
is a primitive element of $\mathbb{Z}^{m+1}$ and the element
$b:=(b_1,\ldots,b_m)$ is a primitive element of $\mathbb{Z}^{m}$.
Therefore, there exist a $\mathbb{Z}$-submodule $N_a$ of
$\mathbb{Z}^{m+1}$ and a $\mathbb{Z}$-submodule $N_b$ of
$\mathbb{Z}^{m}$ such that $\mathbb{Z}^{m+1}=\mathbb{Z}a \oplus N_a$
and $\mathbb{Z}^{m}=\mathbb{Z}b \oplus N_b$. We will denote by
$p_a:\mathbb{Z}^{m+1}\twoheadrightarrow N_a$ the projection of
$\mathbb{Z}^{m+1}$ onto $N_a$ and by $p_b:\mathbb{Z}^{m}\twoheadrightarrow
N_b$ the projection of $\mathbb{Z}^{m}$ onto $N_b$.
We have isomorphisms $i_a: N_a\tilde{\rightarrow}\mathbb{Z}^{m}$ and
$i_b: N_b \tilde{\rightarrow}\mathbb{Z}^{m-1}$.

By definition of the associated morphism, $\varphi_a$ is induced by
the morphism
$f_a:=i_a \circ p_a:\mathbb{Z}^{m+1}\twoheadrightarrow \mathbb{Z}^{m}$ and
$\varphi_b$ by
 $f_b:=i_b \circ p_b:\mathbb{Z}^{m} \twoheadrightarrow
 \mathbb{Z}^{m-1}.$
Set $f:=f_b \circ f_a$. Then $\varphi$ is the $R$-algebra
morphism induced by $f$.

The morphism $f$ is surjective. Since $\mathbb{Z}^{m-1}$ is
a free $\mathbb{Z}$-module, the following exact sequence sequence
splits
$$0\rightarrow\mathrm{Ker}f\rightarrow\mathbb{Z}^{m+1}
\rightarrow \mathbb{Z}^{m-1}\rightarrow0$$ and we obtain that
$\mathbb{Z}^{m+1}\simeq \mathrm{Ker}f\oplus\mathbb{Z}^{m-1}$.

Let $\tilde{b}:=i_a^{-1}(b)$. Then $\mathrm{Ker}f= \mathbb{Z}a
\oplus \mathbb{Z}\tilde{b}$. By assumption, we have that
$c:=(c_0,c_1,\ldots,c_m) \in \mathrm{Ker}f$. Therefore, there exist
unique $\lambda_1, \lambda_2 \in \mathbb{Z}$ such that $c=\lambda_1
a+\lambda_2 \tilde{b}$. Since $\mathrm{gcd}(c_i)=1$, we must also
have $\mathrm{gcd}(\lambda_1,\lambda_2)=1$. Hence $\sum_{i=1}^2
\mathbb{Z}\lambda_i=\mathbb{Z}$. By applying theorem
$\ref{similitude}$ to the $\mathbb{Z}$-module $\mathrm{Ker} f$, we
obtain that $c$ is a primitive element of $\mathrm{Ker}f$.
Consequently, $f=f\,' \circ f_c$, where $f_c$ is a surjective
$\mathbb{Z}$-module morphism $\mathbb{Z}^{m+1} \twoheadrightarrow
\mathbb{Z}^{m}$ with $\mathrm{Ker}f_c= \mathbb{Z}c$ and
 $f\,'$ is a surjective $\mathbb{Z}$-module
morphism $\mathbb{Z}^{m}
\twoheadrightarrow \mathbb{Z}^{m-1}$. As far as the induced
$R$-algebra morphisms are concerned, we obtain that
$\varphi=\varphi\,' \circ \varphi_c$, where $\varphi_c$ is a morphism
associated with the monomial $M=\prod_{i=0}^mx_i^{c_i}$ and
$\varphi\,'$ is also a morphism associated with a monomial (by
proposition $\ref{surjective is associated}$). Thus the assertion is proven.
  \item for $r>2$, let us suppose that the proposition is true for $1,2,\ldots,r-1$.
  If $\varphi={\varphi_r} \circ {\varphi_{r-1}} \circ \ldots \circ {\varphi_1}$, the induction hypothesis
  implies that there exist morphisms associated with monomials $\varphi_r', \varphi_{r-1}',\ldots,\varphi_2'$
  such that
  \begin{enumerate}
  \item $\varphi={\varphi_r'}\circ \ldots\circ {\varphi_2'}\circ{\varphi_1}$.
  \item $\varphi_2'$ is associated with the monomial $(\varphi_1(M))^\circ$.
  \end{enumerate}
  We have that $\varphi_2'(\varphi_1(M))=1$. Once more, by induction hypothesis we obtain that
  there exist morphisms associated with monomials $\varphi_2'',\varphi_1''$
  such that
  \begin{enumerate}
    \item ${\varphi_2'} \circ{\varphi_1}={\varphi_2''}\circ{\varphi_1''}$.
  \item $\varphi_1''$ is associated with $M^\circ$.
  \end{enumerate}
  Thus we have
  $$\varphi={\varphi_r'} \circ \ldots \circ {\varphi_3'} \circ {\varphi_2'}\circ{\varphi_1}=
  {\varphi_r'} \circ \ldots \circ{\varphi_3'} \circ{\varphi_2''} \circ{\varphi_1''}$$
  and $\varphi_1''$ is associated with $M^\circ$.}
\end{itemize}
\end{apod}

\begin{proposition}\label{surjective is adapted}
Let $\varphi:A\rightarrow \mathcal{R}$ be a surjective $R$-algebra
morphism such that for every monomial $M$ in $A$, $\varphi(M)$ is a
monomial in $\mathcal{R}$. Then $\varphi$ is an adapted morphism.
\end{proposition}
\begin{apod}{We will work again by induction on $r$.
For $r=1$, the above result is proposition $\ref{surjective is
  associated}$. For $r>1$, let us suppose that the result
is true for $1,\ldots,r-1$. Due to the isomorphism of lemma
$\ref{work on Z}$, $\varphi$ can be considered as a surjective
$R$-algebra morphism
$$\varphi: R[\mathbb{Z}^{m+1}] \rightarrow R[\mathbb{Z}^{m+1-r}].$$ The
property of $\varphi$ about the monomials implies that the above
morphism is induced by a $\mathbb{Z}$-module morphism $f:
\mathbb{Z}^{m+1} \rightarrow \mathbb{Z}^{m+1-r}$, which is also
surjective by lemma $\ref{groupsGH}$. Since $\mathbb{Z}^{m+1-r}$ is
a free $\mathbb{Z}$-module the following exact sequence sequence
splits
$$0\rightarrow\mathrm{Ker}f\rightarrow\mathbb{Z}^{m+1}
\rightarrow \mathbb{Z}^{m+1-r}\rightarrow0$$ and we obtain that
$\mathbb{Z}^{m+1}\simeq \mathrm{Ker}f\oplus\mathbb{Z}^{m+1-r}$.
Therefore, $\mathrm{Ker}f$ is a $\mathbb{Z}$-module of rank $r$, \ie
$\mathrm{Ker}f \simeq \mathbb{Z}^r$. We choose a primitive element
$a$ of $\mathrm{Ker}f$. Then there exists a $\mathbb{Z}$-submodule
$N_a$ of $\mathrm{Ker}f$ such that $\mathrm{Ker}f=\mathbb{Z}a \oplus
N_a$. Since $\mathrm{Ker}f$ is a direct summand of
$\mathbb{Z}^{m+1}$, $a$ is also a primitive element of
$\mathbb{Z}^{m+1}$ and we have
$$\mathbb{Z}^{m+1} \simeq \mathbb{Z}a \oplus N_a \oplus \mathbb{Z}^{m+1-r}.$$
Thus, by theorem $\ref{similitude}$, if $(a_0,a_1,\ldots,a_m)$ are
the coefficients of $a$ with respect to the standard basis of
$\mathbb{Z}^{m+1}$, then $\mathrm{gcd}(a_i)=1$.

Let us denote by $p_a$ the projection $\mathbb{Z}^{m+1}
\twoheadrightarrow N_a \oplus \mathbb{Z}^{m+1-r}$
 and by $p\,'$ the
projection $N_a \oplus \mathbb{Z}^{m+1-r} \twoheadrightarrow
\mathbb{Z}^{m+1-r}$. Then $f=p\,' \circ p_a$. We have a
$\mathbb{Z}$-module isomorphism
$$i: N_a \oplus \mathbb{Z}^{m+1-r} \tilde{\rightarrow}  \mathbb{Z}^{m}.$$
Set $f_a:=i \circ p_a$ and $f\,':=p\,' \circ i^{-1}$. Thus $f=f\,'
\circ f_a$. If $\varphi_a$ is the $R$-algebra morphism induced by
$f_a$, then, by definition, $\varphi_a$ is a morphism associated
with the monomial $\prod_{i=0}^m x_i^{a_i}$. The $R$-algebra
morphism $\varphi\,'$ induced by $f\,'$ is a surjective morphism
with the same property about monomials as $\varphi$ (it sends every
monomial to a monomial). By induction hypothesis, $\varphi\,'$ is an
adapted morphism. We have $\varphi=\varphi\,' \circ \varphi_a$ and
so $\varphi$ is also an adapted morphism.}
\end{apod}

\section{Irreducibility}

Let $k$ be a field and $y$ an indeterminate over $k$. We can use the
following theorem in order to determine the irreducibility of a
polynomial of the form $y^n-a$ in $k[y]$ (cf.\cite{La}, Chapter 6,
Thm. 9.1).

\begin{theorem}\label{Lang}
Let $k$ be a field, $a \in k-\{0\}$ and $n \in \mathbb{Z}$ with $n
\geq 2$. The polynomial $y^n-a$ is irreducible in $k[y]$, if for
every prime $p$ dividing $n$, we have $a \notin k^p$ and if $4$
divides $n$, we have $a \notin -4k^4$.
\end{theorem}

Let $x_0,x_1,\ldots,x_m$ be a set of $m+1$ indeterminates. We will
apply theorem $\ref{Lang}$ to the field $k(x_1,\ldots,x_m)$.

\begin{lemma}\label{irreducibility}
Let $k$ be a field. The polynomial $x_0^{a_0}-\rho\prod_{i=1}^m
x_i^{a_i}$ with $\rho \in k-\{0\}$, $a_i \in \mathbb{Z}$,
\emph{gcd}$(a_i)=1$ and $a_0 > 0$ is irreducible in
$k[x_1^{\pm1},\ldots,x_m^{\pm1}][x_0]$.
\end{lemma}
\begin{apod}{If $a_0=1$, the polynomial is of degree 1 and thus irreducible in
$k(x_1,\ldots,x_m)[x_0]$. If $a_0 \geq 2$, let us suppose that
$$\rho\prod_{i=1}^m x_i^{a_i}=(\frac{f(x_1,\ldots,x_m)}{g(x_1,\ldots,x_m)})^p$$
with $ f(x_1,\ldots,x_m),g(x_1,\ldots,x_m) \in k[x_1,\ldots,x_m]$
prime to each other,\\
$g(x_1,\ldots,x_m) \neq 0$ and $p|a_0$, $p$ prime. This relation can
be written as
$$g(x_1,\ldots,x_m)^p \rho \prod_{\{i|a_i\geq 0\}}x_i^{a_i}=f(x_1,\ldots,x_m)^p \prod_{\{i|a_i< 0\}}x_i^{-a_i}.$$
We have that
$$\textrm{gcd}(f(x_1,\ldots,x_m)^p,g(x_1,\ldots,x_m)^p)=1.$$
Since $k[x_1,\ldots,x_m]$ is a unique factorization domain and $x_i$
are primes in $k[x_1,\ldots,x_m]$, we also have that
$$\textrm{gcd}(\prod_{\{i|a_i\geq 0\}}x_i^{a_i},\prod_{\{i|a_i<
0\}}x_i^{-a_i})=1.$$ As a consequence,
$$f(x_1,\ldots,x_m)^p = \lambda\rho\prod_{\{i|a_i\geq 0\}}x_i^{a_i} \textrm{ and }
g(x_1,\ldots,x_m)^p = \lambda \prod_{\{i|a_i< 0\}}x_i^{-a_i}.$$ for
some $\lambda \in k-\{0\}$. Suppose that
$(\lambda\rho)^{1/p},\lambda^{1/p} \in k$. Once more, the fact that
$k[x_1,\ldots,x_m]$ is a unique factorization domain and $x_i$ are
primes in $k[x_1,\ldots,x_m]$ implies that
$$f(x_1,\ldots,x_m) = (\lambda\rho)^{1/p}\prod_{\{i|a_i\geq 0\}}x_i^{b_i} \textrm{ and }
g(x_1,\ldots,x_m) = \lambda^{1/p}\prod_{\{i|a_i< 0\}}x_i^{-b_i},$$
with $b_i \in \mathbb{Z}$ and $b_ip=a_i,\forall i=1,\ldots,m$. Since
$p|a_0$, this contradicts the fact that $\textrm{gcd}(a_i)=1$. In
the same way, we can show that if $4|a_0$, then $\rho\prod_{i=1}^m
x_i^{a_i} \notin -4k(x_1,\ldots, x_m)^4$. Thus, by theorem
$\ref{Lang}$, $x_0^{a_0}-\rho\prod_{i=1}^m x_i^{a_i}$ is irreducible
in $k(x_1,\ldots,x_m)[x_0]$.

Thanks to the following lemma, we obtain that
$x_0^{a_0}-\rho\prod_{i=1}^m x_i^{a_i}$ is irreducible in
$k[x_1^{\pm1},\ldots,x_m^{\pm1}][x_0]$.

\begin{lemma}\label{irreducibility in field of fractions}
Let $R$ be an integral domain with field of fractions $F$ and $f(x)$
a polynomial in $R[x]$. If $f(x)$ is irreducible in $F[x]$ and at
least one of its coefficients is a unit in $R$, then $f(x)$ is
irreducible in $R[x]$.
\end{lemma}
\begin{apod}{If $f(x)=g(x)h(x)$ for two polynomials
$g(x),h(x) \in R[x]$, then $g(x) \in R$ or $h(x) \in R$. Let us
suppose that $g(x) \in R$. Since one of the coefficients of $f(x)$
is a unit in $R$, $g(x)$ must also be a unit in $R$. Thus, $f(x)$ is
irreducible in $R[x]$.}
\end{apod}

}
\end{apod}

Lemma $\ref{irreducibility}$ implies the following proposition,
which in turn is going to be used in the proof of proposition
$\ref{second irreducible}$.

\begin{proposition}\label{M-rho}
Let $M:=\prod_{i=0}^mx_i^{a_i}$ be a monomial in
$k[x_0^{\pm1},x_1^{\pm1},\ldots,x_m^{\pm1}]$ such that
$\mathrm{gcd}(a_i)=1$ and let $\rho \in k-\{0\}$. Then $M-\rho$ is
an irreducible element of
$k[x_0^{\pm1},x_1^{\pm1},\ldots,x_m^{\pm1}]$.
\end{proposition}
\begin{apod}{Since
$\textrm{gcd}(a_i)=1$, we can suppose that $a_0 \neq 0$. Then it is
enough to show that $M-\rho$ is irreducible in the polynomial ring
$k[x_1^{\pm1},\ldots,x_m^{\pm1}][x_0]$. If $a_0
> 0$, then $M-\rho$ is irreducible in
$k[x_1^{\pm1},\ldots,x_m^{\pm1}][x_0]$ by lemma
$\ref{irreducibility}$. If now $a_0<0$, then lemma
$\ref{irreducibility}$ implies that $M^{-1}-\rho^{-1}$ is
irreducible in $k[x_1^{\pm1},\ldots,x_m^{\pm1}][x_0]$ and hence
$M-\rho$ is also irreducible in
$k[x_1^{\pm1},\ldots,x_m^{\pm1}][x_0]$.}
\end{apod}

\begin{proposition}\label{second irreducible}
Let $M:=\prod_{i=0}^mx_i^{a_i}$ be a monomial in
$k[x_0^{\pm1},x_1^{\pm1},\ldots,x_m^{\pm1}]$ such that
$\mathrm{gcd}(a_i)=1$. If $f(x)$ is an irreducible element of $k[x]$
such that $f(x) \neq x$, then $f(M)$ is irreducible in
$k[x_0^{\pm1},x_1^{\pm1},\ldots,x_m^{\pm1}]$.
\end{proposition}
\begin{apod}{Suppose that $f(M)=g \cdot h$
with $g,h \in k[x_0^{\pm1},x_1^{\pm1},\ldots,x_m^{\pm1}]$. Let
$\rho_1,\ldots,\rho_n$ be the roots of $f(x)$ in a splitting field
$k'$. Then
$$f(x)=a (x-\rho_1)\ldots (x-\rho_n)$$
for some $a \in k-\{0\}$, hence
$$f(M)=a (M-\rho_1)\ldots (M-\rho_n).$$
By proposition $\ref{M-rho}$, $M-\rho_j$ is irreducible in
$k'[x_0^{\pm1},x_1^{\pm1},\ldots,x_m^{\pm1}]$ for all $j \in
\{1,\ldots,n\}$. Since $k'[x_0^{\pm1},x_1^{\pm1},\ldots,x_m^{\pm1}]$
is a unique factorization domain, we must have
$$g=r \prod_{i=0}^m x_i^{b_i} (M-\rho_{j_1})\ldots(M-\rho_{j_s})$$
for some $r \in k$, $b_i \in \mathbb{Z}$ and $j_1,\ldots,j_s \in
\{1,\ldots,n\}$. Thus there exists $g'(x) \in k[x]$ such that
$$g=(\prod_{i=0}^m x_i^{b_i})g'(M).$$
Respectively, there exists $h'(x) \in k[x]$ such
that$$h=(\prod_{i=0}^m x_i^{-b_i}) h'(M).$$ Thus, we obtain that
$$f(M)=g'(M)h'(M). \,\,\,\,\, (\dag)$$
Since $\mathrm{gcd}(a_i)=1$, there exist integers $(u_i)_{0 \leq i
\leq m}$ such that $\sum_{i=0}^m u_ia_i=1$. Let us now consider the
$k$-algebra specialization
$$\begin{array}{cccc}
  \varphi: & k[x_0^{\pm1},x_1^{\pm1},\ldots,x_m^{\pm1}] & \rightarrow & k[x] \\
           & x_i                                        & \mapsto     & x^{u_i}.
\end{array}$$
Then $\varphi(M)=\varphi(\prod_{i=0}^m x_i^{a_i})=x^{\sum_{i=0}^m
u_ia_i}=x$. If we apply $\varphi$ to the relation $(\dag)$, we
obtain that
$$f(x)=g'(x)h'(x).$$
Since $f(x)$ is irreducible in $k[x]$, we must have that either
$g'(x) \in k$ or $h'(x) \in k$. Respectively, we deduce that either
$g$ or $h$ is a unit in
$k[x_0^{\pm1},x_1^{\pm1},\ldots,x_m^{\pm1}]$.}
\end{apod}

\chapter{On Blocks}

All the results presented in the first two sections of this chapter
have been taken from the first part of \cite{BK}.

\section{Generalities}

Let $\mathcal{O}$ be a commutative ring with a unit element and $A$
be an $\mathcal{O}$-algebra. We denote by $ZA$ the center of $A$.

An idempotent in $A$ is an element $e$ such that $e^2=e$. We say
that $e$ is a central idempotent, if it is an idempotent in $ZA$.
Two idempotents $e_1,e_2$ are orthogonal, if $e_1e_2=e_2e_1=0$.
Finally, an idempotent $e$ is primitive, if $e \neq 0$ and $e$ can
not be expressed as the sum of two non-zero orthogonal idempotents.

\begin{definition}\label{blocks}
The block-idempotents of $A$ are the central primitive idempotents
of $A$.
\end{definition}

Let $e$ be a block-idempotent of $A$. The two sided ideal $Ae$
inherits a structure of algebra, where the composition laws are
those of $A$ and the unit element is $e$. The application
$$\begin{array}{cccc}
  \pi_e: &A &\rightarrow &Ae \\
         &h &\mapsto &he
\end{array}$$
is an epimorphism of algebras. The algebra $Ae$ is called a block of
$A$. From now on, abusing the language, we will also call blocks the
block-idempotents of $A$.

\begin{lemma}\label{orthogonality of blocks}
The blocks of $A$ are mutually orthogonal.
\end{lemma}
\begin{apod}{Let $e$ be a block and $f$ a central idempotent of $A$ with $f \neq e$.
Then $ef$ and $e-ef$ are also central idempotents. We have
$e=ef+(e-ef)$ and due to the primitivity of $e$, we deduce that
either $ef=0$ or $e=ef$. If $f$ is a block too, then either $ef=0$
or $f=ef=e$. Therefore, $f$ is orthogonal to $e$.}
\end{apod}

The above lemma gives rise to the following proposition.

\begin{proposition}\label{1 sum of blocks}
Suppose that the unit element $1$ of $A$ can be expressed as a sum
of blocks: $1=\sum_{e \in E}e$. Then
\begin{enumerate}
  \item The set $E$ is the set of all the blocks of $A$.
  \item The family of morphisms $(\pi_e)_{e \in E}$ defines an
  isomorphism of algebras
  $$A \tilde{\rightarrow} \prod_{e \in E} Ae.$$
\end{enumerate}
\end{proposition}
\begin{apod}{If $f$ is a block, then $f=\sum_{e \in E}ef$. Due to
lemma $\ref{orthogonality of blocks}$, there exists $e \in E$ such
that $f=e$.}
\end{apod}

In the above context (1 is a sum of blocks), let us denote by
$\mathrm{Bl}(A)$ the set of all the blocks of $A$. Proposition
$\ref{1 sum of blocks}$ implies that the category
$_A\mathrm{\textbf{mod}}$ of $A$-modules is a direct sum of the
categories associated with the blocks:
$$_A\mathrm{\textbf{mod}} \tilde{\rightarrow} \bigoplus_{e \in \mathrm{Bl}(A)} \mathrm{}_{Ae}\mathrm{\textbf{mod}}.$$
In particular, every representation of the $\mathcal{O}$-algebra
$Ae$ defines (by composition with $\pi_e$) a representation of $A$
and we say, abusing the language, that it ``belongs to the block
$e$''.

Every indecomposable representation of $A$ belongs to one and only
one block. Thus the following partitions are defined:
$$\mathrm{Ind}(A)=\bigsqcup_{e \in \mathrm{Bl}(A)} \mathrm{Ind}(A,e) \,\,\textrm{ and
}\,\, \mathrm{Irr}(A)=\bigsqcup_{e \in \mathrm{Bl}(A)}
\mathrm{Irr}(A,e),$$ where $\mathrm{Ind}(A)$ (resp.
$\mathrm{Irr}(A)$) denotes the set of indecomposable (resp.
irreducible) representations of $A$ and $\mathrm{Ind}(A,e)$ (resp.
$\mathrm{Irr}(A,e)$) denotes the set of the elements of
$\mathrm{Ind}(A)$ (resp. $\mathrm{Irr}(A)$) which belong to $e$.\\

We will consider two situations where 1 is a sum of blocks.\\

\textbf{First case:} Suppose that 1 is a sum of orthogonal primitive
idempotents, \ie $1=\sum_{i \in P}i$, where
\begin{itemize}
  \item every $i \in P$ is a primitive idempotent,
  \item if $i,j \in P$, $i \neq j$, then $ij=ji=0$.
\end{itemize}
Let us consider the equivalence relation $\mathcal{B}$ defined on
$P$ as the symmetric and transitive closure of the relation ``$iAj
\neq \{0\}$''. Thus $(i \mathcal{B} j)$ if and only if there exist
$i_0,i_1,\ldots,i_n \in P$ with $i_0=i$ and $i_n=j$ such that for
all $k \in \{1,\ldots,n\}$, $i_{k-1}Ai_k \neq \{0\}$ or $i_kAi_{k-1}
\neq \{0\}$. To every equivalence class $B$ of $P$ with respect to
$\mathcal{B}$, we associate the idempotent $e_B:=\sum_{i \in B}i$.

\begin{proposition}\label{iAj}
The map $B \mapsto e_B$ is a bijection between the set of
equivalence classes of $\mathcal{B}$ and the set of blocks of $A$.
In particular, we have that $1=\sum_{B \in P/\mathcal{B}}e_B$ and
$1$ is sum of the blocks of $A$.
\end{proposition}
\begin{apod}{It is clear that $1=\sum_{B \in P/\mathcal{B}}e_B$. Let
$a \in A$ and let $B,B'$ be two equivalence classes of $\mathcal{B}$
with $B \neq B'$. Then, by definition of the relation $\mathcal{B}$,
$e_B a e_{B'}=0$. Since $1=\sum_{B \in P/\mathcal{B}}e_B$, we have
that $e_Ba=e_Bae_B=ae_B$. Thus $e_B \in ZA$ for all $B \in
P/\mathcal{B}$.

It remains to show that for all $B \in P/\mathcal{B}$, the central
idempotent $e_B$ is primitive. Suppose that $e_B=e+f$, where $e$ and
$f$ are two orthogonal primitive idempotents in $ZA$. Then we have a
partition $B=B_e \sqcup B_f$, where $B_e:=\{i \in B \,|\, ie=i\}$
and $B_f:=\{j \in B \,|\, jf=j\}$. For all $i \in B_e$ and $j \in
B_f$, we have $iAj=ieAfj=iAefj=\{0\}$ and so no element of $B_e$ can
be $\mathcal{B}$-equivalent to an element of $B_f$. Therefore, we
must have either $B_e=\emptyset$ or $B_f=\emptyset$, which implies
that either $e=0$ or $f=0$.}
\end{apod}

\textbf{Second case:} Suppose that $ZA$ is a subalgebra of a
commutative algebra $C$ where 1 is a sum of blocks. For example, if
$A$ is of finite type over $\mathcal{O}$, where $\mathcal{O}$ is an
integral domain with field of fractions $F$, we can choose $C$ to be
the center of the algebra $FA:=F \otimes_\mathcal{O} A$.

We set $1=\sum_{e \in E}e$, where $E$ is the set of blocks of $C$.
For all $S \subseteq E$, set $e_S:=\sum_{e \in S}e$. A subset $S$ of
$E$ is ``on $ZA$'' if $e_S \in ZA$. If $S$ and $T$ are on $ZA$, then
$S \cap T$ is on $ZA$.

\begin{proposition}\label{partition}
Let us denote by $\mathcal{P}_E(ZA)$ the set of non-empty subsets
$B$ of $E$ which are on $ZA$ and are minimal for these two
properties. Then the map $\mathcal{P}_E(ZA) \rightarrow A, B \mapsto
e_B$ induces a bijection between $\mathcal{P}_E(ZA)$ and the set of
blocks of $A$. We have $1=\sum_{B \in \mathcal{P}_E(ZA)}e_B$.
\end{proposition}
\begin{apod}{Since every idempotent in $C$ is of the form $e_S$ for
some $S \subseteq E$, it is clear that $e_B$ is a central primitive
idempotent of $A$, for all $B \in \mathcal{P}_E(ZA)$. It remains to
show that
\begin{enumerate}
  \item If $B$,$B'$ are two distinct elements of
  $\mathcal{P}_E(ZA)$, then $B \cap B' = \emptyset$.
  \item $\mathcal{P}_E(ZA)$ is a partition of $E$.
\end{enumerate}
These two properties, stated in terms of idempotents, mean:
\begin{enumerate}
  \item If $B$,$B'$ are two distinct elements of
  $\mathcal{P}_E(ZA)$, then $e_B$ and $e_{B'}$ are orthogonal.
  \item $1=\sum_{B \in \mathcal{P}_E(ZA)}e_B$.
\end{enumerate}
Let us prove them:
\begin{enumerate}
  \item We have $e_Be_{B'}=e_{B \cap B'}$ and so $B \cap B' = \emptyset$, because $B$ and $B'$ are
  minimal.
  \item Set $F:=\bigcup_{B \in \mathcal{P}_E(ZA)}B$. Then $e_F=\sum_{B \in \mathcal{P}_E(ZA)}e_B \in
  ZA$. Then $1-e_F=e_{E-F} \in ZA$, which means that $E-F$ is on
  $ZA$ . If $E-F \neq \emptyset$, then $E-F$ contains an element
  of $\mathcal{P}_E(ZA)$ in contradiction to the definition of $F$.
  Thus $F=E$ and $\mathcal{P}_E(ZA)$ is a partition of $E$.}
\end{enumerate}
\end{apod}

Let us assume that
\begin{itemize}
  \item $\mathcal{O}$ is a commutative integral domain with field of fractions
  $F$,
  \item $K$ is a field extension of $F$,
  \item $A$ is an $\mathcal{O}$-algebra, free and finitely generated as an
  $\mathcal{O}$-module.
\end{itemize}

Suppose that the $K$-algebra $KA:=K \otimes_\mathcal{O}A$ is
semisimple. Then $KA$ is isomorphic, by assumption, to a direct
product of simple algebras:
$$KA \simeq \prod_{\chi \in \mathrm{Irr}(KA)} M_\chi,$$
where $\mathrm{Irr}(KA)$ denotes the set of irreducible characters
of $KA$ and $M_\chi$ is a simple $K$-algebra.

For all $\chi \in \mathrm{Irr}(KA)$, we denote by $\pi_\chi:KA
\twoheadrightarrow M_\chi$ the projection onto the $\chi$-factor and
by $e_\chi$ the element of $KA$ such that $$\pi_{\chi'}(e_\chi)=
  \left\{
  \begin{array}{ll}
    1_{M_\chi}, & \hbox{if $\chi=\chi'$,} \\
    0, & \hbox{if $\chi \neq \chi'$.}
  \end{array}
\right.$$

The following theorem results directly from propositions $\ref{1 sum
of blocks}$ and $\ref{partition}$.

\begin{theorem}\label{minimality of blocks}\
\begin{enumerate}
  \item We have $1=\sum_{\chi \in \mathrm{Irr}(KA)}e_\chi$
   and the set $\{e_\chi\}_{\chi \in \mathrm{Irr}(KA)}$ is the set of all the blocks of the algebra $KA$.
  \item There exists a unique partition $\mathrm{Bl}(A)$ of
  $\mathrm{Irr}(KA)$ such that
  \begin{description}
    \item[(a)] For all $B \in \mathrm{Bl}(A)$, the idempotent
    $e_B:=\sum_{\chi \in B}e_\chi$ is a block of $A$.
    \item[(b)] We have $1=\sum_{B \in \mathrm{Bl}(A)}e_B$ and for
    every central idempotent $e$ of $A$, there exists a subset
    $\mathrm{Bl}(A,e)$ of $\mathrm{Bl}(A)$ such that
    $$e=\sum_{B \in \mathrm{Bl}(A,e)}e_B.$$
    In particular the set $\{e_B\}_{B \in \mathrm{Bl}(A)}$ is the set of all the blocks of $A$.
  \end{description}
\end{enumerate}
\end{theorem}

\begin{remarks}\
\emph{\begin{itemize}
  \item If $\chi \in B$ for some $B \in \mathrm{Bl}(A)$, we say that
  ``$\chi$ belongs to the block $e_B$''.
  \item For all $B \in \mathrm{Bl}(A)$, we have
  $$KAe_B \simeq \prod_{\chi \in B}M_\chi.$$
\end{itemize}}
\end{remarks}

From now on, we make the following assumptions

\begin{ypoth}\label{properties of the ring}\
\begin{description}
  \item[(int)] The ring $\mathcal{O}$ is a Noetherian and integrally
  closed domain with field of fractions $F$ and $A$ is an
  $\mathcal{O}$-algebra which is free and finitely generated as an
  $\mathcal{O}$-module.
  \item[(spl)] The field $K$ is a finite Galois extension of $F$ and
  the algebra $KA$ is split (i.e., for every simple $KA$-module $V$, $\mathrm{End}_{KA}(V) \simeq K$) semisimple.
\end{description}
\end{ypoth}

We denote by $\mathcal{O}_K$ the integral closure of $\mathcal{O}$
in $K$.

\subsection*{Blocks and integral closure}

The Galois group $\mathrm{Gal}(K/F)$ acts on $KA=K
\otimes_{\mathcal{O}} A$ (viewed as an $F$-algebra) as follows: if
$\sigma \in \mathrm{Gal}(K/F)$ and $\lambda \otimes a \in KA$, then
$\sigma(\lambda \otimes a):=\sigma(\lambda) \otimes a$.

If $V$ is a $K$-vector space and $\sigma \in \mathrm{Gal}(K/F)$, we
denote by $^\sigma V$ the $K$-vector space defined on the additive
group $V$ with multiplication $\lambda.v:=\sigma^{-1}(\lambda)v$ for
all $\lambda \in K$ and $v \in V$. If $\rho:KA \rightarrow
\mathrm{End}_K(V)$ is a representation of the $K$-algebra $KA$, then
its composition with the action of $\sigma^{-1}$ is also a
representation $^\sigma \rho: KA \rightarrow \mathrm{End}_K(^\sigma
V)$:
$$\diagram KA \rto^{\sigma^{-1}} &KA \rto^{\rho}
&\mathrm{End}_K(V). \enddiagram$$

We denote by $^\sigma \chi$ the character of $^\sigma \rho$ and we
define the action of $\mathrm{Gal}(K/F)$ on $\mathrm{Irr}(KA)$ as
follows: if $\sigma \in \mathrm{Gal}(K/F)$ and $\chi \in
\mathrm{Irr}(KA)$, then $$\sigma(\chi):={}^\sigma\!\chi = \sigma
\circ \chi \circ \sigma^{-1}.$$ This operation induces an action of
$\mathrm{Gal}(K/F)$ on the set of blocks of $KA$:
$$\sigma(e_\chi)=e_{^\sigma \chi} \textrm{ for all } \sigma \in
\mathrm{Gal}(K/F), \chi \in \mathrm{Irr}(KA).$$

Hence, the group $\mathrm{Gal}(K/F)$ acts on the set of idempotents
of $Z\mathcal{O}_KA$ and thus on the set of blocks of
$\mathcal{O}_KA$. Since $F \cap \mathcal{O}_K = \mathcal{O}$, the
idempotents of $ZA$ are the idempotents of $Z\mathcal{O}_KA$ which
are fixed by the action of $\mathrm{Gal}(K/F)$. As a consequence,
the primitive idempotents of $ZA$ are sums of the elements of the
orbits of $\mathrm{Gal}(K/F)$ on the set of primitive idempotents of
$Z\mathcal{O}_KA$. Thus, the blocks of $A$ are in bijection with the
orbits of $\mathrm{Gal}(K/F)$ on the set of blocks of
$\mathcal{O}_KA$. The following proposition is just a reformulation
of this result.

\begin{proposition}\label{Galois action on integral closure}\
\begin{enumerate}
  \item Let $B$ be a block of $A$ and $B'$ a block of
  $\mathcal{O}_KA$ contained in $B$. If $\mathrm{Gal}(K/F)_{B'}$
  denotes the stabilizer of $B'$ in $\mathrm{Gal}(K/F)$, then
  $$B=\bigcup_{\sigma \in
  \mathrm{Gal}(K/F)/\mathrm{Gal}(K/F)_{B'}}\sigma(B')
  \,\,\textrm{ i.e., }\,\,
  e_B=\sum_{\sigma \in
  \mathrm{Gal}(K/F)/\mathrm{Gal}(K/F)_{B'}}\sigma(e_{B'}).$$
  \item Two characters $\chi,\psi \in \mathrm{Irr}(KA)$ are in
  the same block of $A$ if and only if there exists $\sigma \in \mathrm{Gal}(K/F)$
  such that $\sigma(\chi)$ and $\psi$ belong to the same block of
  $\mathcal{O}_KA$.
\end{enumerate}
\end{proposition}\
\begin{remark}\emph{ For all $\chi \in B'$, we have
$\mathrm{Gal}(K/F)_\chi \subseteq \mathrm{Gal}(K/F)_{B'}.$}
\end{remark}\
\\

The assertion (2) of the proposition above allows us to transfer the
problem of the classification of the blocks of $A$ to that of the
classification of the blocks of $\mathcal{O}_KA$.

\subsection*{Blocks and prime ideals}

We denote by $\mathrm{Spec}_1(\mathcal{O})$ the set of prime ideals
of height 1 of $\mathcal{O}$. Since $\mathcal{O}$ is Noetherian and
integrally closed, it is a Krull ring and by theorem
$\ref{Krull-dvr}$, we have
$$\mathcal{O}=\bigcap_{\mathfrak{p} \in \mathrm{Spec}_1(\mathcal{O})}
\mathcal{O}_\mathfrak{p},$$ where $\mathcal{O}_\mathfrak{p}:=\{x \in
F\,|\,(\exists a \in \mathcal{O}-\mathfrak{p})(ax \in
\mathcal{O})\}$ is the localization of $\mathcal{O}$ at
$\mathfrak{p}$. More generally, if we denote by
$\mathrm{Spec}(\mathcal{O})$ the set of prime ideals of
$\mathcal{O}$, then
$$\mathcal{O}=\bigcap_{\mathfrak{p} \in \mathrm{Spec}(\mathcal{O})}
\mathcal{O}_\mathfrak{p}.$$

Let $\mathfrak{p}$ be a prime ideal of $\mathcal{O}$ and
$\mathcal{O}_\mathfrak{p}A:=\mathcal{O}_\mathfrak{p}
\otimes_{\mathcal{O}}A$. The blocks of $\mathcal{O}_\mathfrak{p}A$
are the ``$\mathfrak{p}$-blocks of $A$''. If $\chi,\psi \in
\mathrm{Irr}(KA)$ belong to the same block of
$\mathcal{O}_\mathfrak{p}A$, we write $\chi \sim_\mathfrak{p} \psi$.

\begin{proposition}\label{p-blocks}
Two characters $\chi,\psi \in \mathrm{Irr}(KA)$ belong to the same
block of $A$ if and only if there exist a finite sequence
$\chi_0,\chi_1,\ldots,\chi_n \in \mathrm{Irr}(KA)$ and a finite
sequence $\mathfrak{p}_1,\ldots,\mathfrak{p}_n \in
\mathrm{Spec}(\mathcal{O})$ such that
\begin{itemize}
  \item $\chi_0=\chi$ and $\chi_n=\psi$,
  \item for all $j$ $(1\leq j \leq n)$, $\chi_{j-1}
  \sim_{\mathfrak{p}_j} \chi_j$.
\end{itemize}
\end{proposition}
\begin{apod}{Let us denote by $\sim$ the equivalence relation on
$\mathrm{Irr}(KA)$ defined as the closure of the relation ``there
exists $\mathfrak{p} \in \mathrm{Spec}(\mathcal{O})$ such that $\chi
\sim_\mathfrak{p} \psi$''. Thus, we have to show that $\chi \sim
\psi$ if and only if $\chi$ and $\psi$ belong to the same block of
$A$.

We will first show that the equivalence relation $\sim$ is finer
than the relation ``being in the same block of A''. Let $B$ be a
block of $A$. Then $B$ is a subset of $\mathrm{Irr}(KA)$ such that
$\sum_{\chi \in B}e_\chi \in A$. Since
$\mathcal{O}=\bigcap_{\mathfrak{p} \in \mathrm{Spec}(\mathcal{O})}
\mathcal{O}_\mathfrak{p}$, we have that $\sum_{\chi \in B}e_\chi \in
\mathcal{O}_\mathfrak{p}A$ for all $\mathfrak{p} \in
\mathrm{Spec}(\mathcal{O})$. Therefore, by theorem $\ref{minimality
of blocks}$, $C$ is a union of blocks of $\mathcal{O}_\mathfrak{p}A$
for all $\mathfrak{p} \in \mathrm{Spec}(\mathcal{O})$ and, hence, a
union of equivalence classes of $\sim$.

Now we will show that the relation ``being in the same block of A''
if finer than the relation $\sim$. Let $C$ be an equivalence class
of $\sim$. Then $\sum_{\chi \in C}e_\chi \in
\mathcal{O}_\mathfrak{p}A$ for all $\mathfrak{p} \in
\mathrm{Spec}(\mathcal{O})$. Thus $\sum_{\chi \in C}e_\chi \in
\bigcap_{\mathfrak{p} \in \mathrm{Spec}(\mathcal{O})}
\mathcal{O}_\mathfrak{p}A=A$ and $C$ is a union of blocks of $A$.}
\end{apod}

\subsection*{Blocks and central morphisms}

Since $KA$ is a split semisimple $K$-algebra, we have that
$$KA \simeq \prod_{\chi \in \mathrm{Irr}(KA)} M_\chi,$$
where $M_\chi$ is a matrix algebra isomorphic to
$\mathrm{Mat}_{\chi(1)}(K)$.

Recall that $A$ is of finite type and thus integral over
$\mathcal{O}$ (\cite{Bou5}, \S1, Def.2). The map $\pi_\chi: KA
\twoheadrightarrow M_\chi$, restricted to $ZKA$, defines a map
$\omega_\chi:ZKA \twoheadrightarrow K$ (by Schur's lemma), which in
turn, restricted to $ZA$, defines the morphism
$$\omega_\chi: ZA \rightarrow \mathcal{O}_K,$$
where $\mathcal{O}_K$ denotes the integral closure of $\mathcal{O}$
in $K$.

In the case where $\mathcal{O}$ is a discrete valuation ring, we
have the following result which is proven later in this chapter,
proposition $\ref{blocks and central characters}$. For a different
approach to its proof, see \cite{BK}, Prop.1.18.

\begin{proposition}\label{omega_chi}
Suppose that $\mathcal{O}$ is a discrete valuation ring with unique
maximal ideal $\mathfrak{p}$ and $K=F$. Two characters $\chi,\chi'
\in \mathrm{Irr}(KA)$ belong to the same block of $A$ if and only if
$$\omega_\chi(a) \equiv \omega_{\chi'}(a)\,\, \mathrm{mod} \, \mathfrak{p}
\textrm{ for all } a \in ZA.$$
\end{proposition}

\section{Symmetric algebras}

Let $\mathcal{O}$ be a ring and let $A$ be an $\mathcal{O}$-algebra.
Suppose again that the assumptions $\ref{properties of the ring}$
are satisfied.

\begin{definition}\label{trace function}
A trace function on $A$ is an $\mathcal{O}$-linear map $t:A
\rightarrow \mathcal{O}$ such that $t(ab)=t(ba)$ for all $a,b \in
A$.
\end{definition}

\begin{definition}\label{symmetric algebra}
We say that a trace function $t:A \rightarrow \mathcal{O}$ is a
symmetrizing form on $A$ or that $A$ is a symmetric algebra if the
morphism
$$\hat{t}:A \rightarrow \mathrm{Hom}_\mathcal{O}(A,\mathcal{O}),\,\,
  a \mapsto (x \mapsto \hat{t}(a)(x):=t(ax))$$
is an isomorphism of $A$-modules-$A$.
\end{definition}

\begin{px}\label{symmetrizing form of the group algebra}
\small{\emph{In the case where $\mathcal{O}=\mathbb{Z}$ and
$A=\mathbb{Z}[G]$
 ($G$ a finite group), we can define the following symmetrizing form
 (``canonical'')
 on $A$
$$t:\mathbb{Z}[G] \rightarrow \mathbb{Z}, \,\,\, \sum_{g \in G}a_g g \mapsto a_1,$$
where $a_g \in \mathbb{Z}$ for all $g \in G$.}}
\end{px}

Since $A$ is a free $\mathcal{O}$-module of finite rank, we have the
following isomorphism
$$\begin{array}{ccc}
  \mathrm{Hom}_\mathcal{O}(A,\mathcal{O}) \otimes_\mathcal{O} A &\tilde{\rightarrow} &\mathrm{Hom}_\mathcal{O}(A,A)\\
  \varphi \otimes a & \mapsto & (x \mapsto \varphi(x)a).
\end{array}$$
Composing it with the isomorphism
$$\begin{array}{ccc}
  A \otimes_\mathcal{O} A &\tilde{\rightarrow} &\mathrm{Hom}_\mathcal{O}(A,\mathcal{O})\otimes_\mathcal{O} A\\
  a \otimes b & \mapsto & \hat{t}(a) \otimes b,
\end{array}$$
we obtain an isomorphism
$$ A \otimes_\mathcal{O} A \tilde{\rightarrow} \mathrm{Hom}_\mathcal{O}(A,A).$$

\begin{definition}\label{casimir}
We denote by $C_A$ and we call Casimir of $(A,t)$ the inverse image
of $\mathrm{Id}_A$ by the above isomorphism.
\end{definition}

\begin{px}\label{casimir of the group algebra}
\small{\emph{In the case where $\mathcal{O}=\mathbb{Z}$,
$A=\mathbb{Z}[G]$ ($G$ a finite group) and $t$ is the canonical
symmetrizing form, we have $C_{\mathbb{Z}[G]}=\sum_{g \in
G}g^{-1}\otimes g$.}}
\end{px}

More generally, if $(e_i)_{i \in I}$ is a basis of $A$ over
$\mathcal{O}$ and $(e_i')_{i \in I}$ is the dual basis with respect
to $t$ (\ie $t(e_ie_j')=\delta_{ij}$), then
$$C_A = \sum_{i \in I} e_i' \otimes e_i.$$

In this case, let us denote by $c_A$ the image of $C_A$ by the
multiplication $A \otimes A \rightarrow A$, \ie $c_A= \sum_{i \in I}
e_i'e_i$. It is easy to check (see also \cite{BMM2}, 7.9) the
following properties of the Casimir element:

\begin{lemma}\label{properties of the casimir}
For all $a \in A$, we have
\begin{enumerate}
  \item $\sum_i ae_i' \otimes e_i = \sum_i e_i \otimes e_i'a$.
  \item $aC_A=C_Aa$. Consequently, $c_A \in ZA$.
  \item $a=\sum_i t(ae_i')e_i=\sum_i t(ae_i)e_i'=
  \sum_i t(e_i')e_ia=\sum_i t(e_i)e_i'a.$
\end{enumerate}
\end{lemma}

If $\tau:A \rightarrow \mathcal{O}$ is a linear form, we denote by
$\tau^\vee$ its inverse image by the isomorphism $\hat{t}$, \ie
$\tau^\vee$ is the element of $A$ such that
$$t(\tau^\vee a)=\tau(a) \textrm{ for all } a \in A.$$

The element $\tau^\vee$ has the following properties:

\begin{lemma}\label{tau^vee}\
\begin{enumerate}
 \item $\tau$ is a trace function if and only if $\tau^\vee \in ZA$.
 \item We have $\tau^\vee=\sum_i \tau(e_i')e_i=\sum_i \tau(e_i)e_i'$ and
 more generally, for all $a \in A$, we have
 $\tau^\vee a=\sum_i \tau(e_i'a)e_i=\sum_i \tau(e_ia)e_i'$.
\end{enumerate}
\end{lemma}
\begin{apod}{\begin{enumerate}
  \item   Recall that $t$ is a trace function. Let $a \in A$. For all $x \in A$, we have
  $$\hat{t}(\tau^\vee a)(x) = t(\tau^\vee a x)=\tau(ax)$$
  and
  $$\hat{t}(a \tau^\vee )(x) = t(a \tau^\vee x)=t(\tau^\vee
  xa)=\tau(xa)$$
  If $\tau$ is a trace function, then $\tau(ax)=\tau(xa)$ and hence, $\hat{t}(\tau^\vee
  a)= \hat{t}(a \tau^\vee)$. Since $\hat{t}$ is an isomorphism, we
  obtain that $ \tau^\vee a = a \tau^\vee$ and thus $\tau^\vee \in ZA$.

  Now if $\tau^\vee \in ZA$ and $a,b \in A$, then $$\tau(ab)=t(\tau^\vee ab)
  =t(b \tau^\vee a) =t(ba \tau^\vee )=t(\tau^\vee ba)=\tau(ba).$$
  \item It derives from property 3 of lemma $\ref{properties of the
  casimir}$ and the definition of $\tau^\vee$.}
\end{enumerate}
\end{apod}

Let $\chi_\mathrm{reg}$ be the character of the regular
representation of $A$, \ie the linear form on $A$ defined as
$$\chi_\mathrm{reg}(a):=\mathrm{tr}_{A/\mathcal{O}}(\lambda_a),$$
where $\lambda_a:A\rightarrow A, x \mapsto ax$, is the endomorphism
of left multiplication by $a$.

\begin{proposition}\label{regular^vee}
We have $\chi_\mathrm{reg}^\vee = c_A.$
\end{proposition}
\begin{apod}{Let $a \in A$. The inverse image of $\lambda_a$ by the
isomorphism \\
$A \otimes_\mathcal{O} A \tilde{\rightarrow}
\mathrm{Hom}_\mathcal{O}(A,A)$ is $aC_A$ (by definition of the
Casimir). Hence,
$$\lambda_a = (x \mapsto \sum_i \hat{t}(e_i'a)(x)e_i) = (x \mapsto \sum_i
t(e_i'ax)e_i)$$ and
$$\mathrm{tr}_{A/\mathcal{O}}(\lambda_a)=\sum_i t(e_i'ae_i)=
t(a\sum_i e_i'e_i)=t(ac_A)=t(c_Aa).$$ Therefore, for all $a \in A$,
we have $\chi_\mathrm{reg}(a)=t(c_Aa)$, \ie $\chi_\mathrm{reg}^\vee
= c_A.$}
\end{apod}

If $A$ is a symmetric algebra with a symmetrizing form $t$, we
obtain a symmetrizing form $t^K$ on $KA$ by extension of scalars.
Every irreducible character $\chi \in \mathrm{Irr}(KA)$ is a central
function on $KA$ and thus we can define $\chi^\vee \in KA$.

\begin{definition}\label{Schur element}
For all $\chi \in \mathrm{Irr}(KA)$, we call Schur element of $\chi$
with respect to $t$ and denote by $s_\chi$ the element of $K$
defined by $$s_\chi:=\omega_\chi(\chi^\vee).$$
\end{definition}

\begin{proposition}\label{Schur element belongs to the integral closure}
For all $\chi \in \mathrm{Irr}(KA)$, $s_\chi \in \mathcal{O}_K$.
\end{proposition}

The proof of the above result will be given in proposition
$\ref{integrality of the Schur elements}$.

\begin{px}\label{Schur elements of the group algebra}
\small{\emph{Let $\mathcal{O}:=\mathbb{Z}$, $A:=\mathbb{Z}[G]$
 ($G$ a finite group) and $t$ the canonical symmetrizing form. If $K$ is an algebraically closed field of
 characteristic 0, then $KA$ is a split semisimple algebra and
 $s_\chi=|G|/\chi(1)$ for all $\chi \in \mathrm{Irr}(KA)$. Because
 of the integrality of the Schur elements, we must have
 $|G|/\chi(1) \in \mathbb{Z}=\mathbb{Z}_K \cap \mathbb{Q}$ for all $\chi \in
 \mathrm{Irr}(KA)$. Thus, we have shown that $\chi(1)$ divides $|G|$.}}
\end{px}

The following properties of the Schur elements can be derived easily
from the above (see also
\cite{Bro},\cite{Ge},\cite{GePf},\cite{GeRo},\cite{BMM2})

\begin{proposition}\label{schur elements and idempotents}\
\begin{enumerate}
  \item We have
  $$t=\sum_{\chi \in \mathrm{Irr}(KA)}\frac{1}{s_\chi}\chi.$$
  \item For all $\chi \in \mathrm{Irr}(KA)$, the central primitive
  idempotent associated with $\chi$ is
  $$e_\chi=\frac{1}{s_\chi}\chi^\vee=\frac{1}{s_\chi}\sum_{i \in I}
  \chi(e_i')e_i.$$
  \item For all $\chi \in \mathrm{Irr}(KA)$, we have
  $$s_\chi \chi(1)=\sum_{i \in I}\chi(e_i')\chi(e_i) \,\textrm{ and  }\,
  s_\chi \chi(1)^2= \chi(\sum_{i \in I}e_i'e_i)=\chi(\chi_{\mathrm{reg}}^\vee).$$
\end{enumerate}
\end{proposition}

\begin{corollary}\label{what we are searching}
The blocks of $A$ are the non-empty subsets $B$ of
$\mathrm{Irr}(KA)$ minimal for the property
$$\sum_{\chi \in B}\frac{1}{s_\chi}\chi(a) \in \mathcal{O} \textrm{ for all } a \in A.$$
\end{corollary}

\section{Twisted symmetric algebras of finite groups}

This part is an adaptation of the section ``Symmetric algebras of
finite groups'' of \cite{BK} to a more general case.\\

Let $A$ be an $\mathcal{O}$-algebra such that the assumptions
$\ref{properties of the ring}$ are satisfied with a symmetrizing
form $t$. Let $\bar{A}$ be a subalgebra of $A$ free and of finite
rank as $\mathcal{O}$-module.

We denote by $\bar{A}^\bot$ the orthogonal of $\bar{A}$ with respect
to $t$, \ie the sub-$\bar{A}$-module-$\bar{A}$ of $A$ defined as
$$\bar{A}^\bot:=\{a \in A \,|\, (\forall \bar{a} \in \bar{A})(t(a\bar{a})=0)\}.$$

\begin{proposition}\label{when symmetric subalgebra}\
\begin{enumerate}
  \item The restriction of $t$ to $\bar{A}$ is a symmetrizing form for $\bar{A}$
  if and only if $\bar{A} \oplus \bar{A}^\bot=A$.
   In this case the projection of $A$ onto $\bar{A}$ parallel to $\bar{A}^\bot$ is
  the map
  $$\mathrm{Br}_{\bar{A}}^A :A \rightarrow \bar{A}
  \textrm{ such that } t(\mathrm{Br}_{\bar{A}}^A(a)\bar{a})=t(a\bar{a})
  \textrm{ for all } a \in A \textrm{ and } \bar{a} \in \bar{A}.$$
  \item If the restriction of $t$ to $\bar{A}$ is a symmetrizing form
  for $\bar{A}$, then $\bar{A}^\bot$ is the sub-$\bar{A}$-module-$\bar{A}$
  of $A$ defined by the following two properties:
  \begin{description}
    \item[(a)] $A=\bar{A} \oplus \bar{A}^\bot$,
    \item[(b)] $\bar{A}^\bot \subseteq \mathrm{Ker}t$.
  \end{description}
\end{enumerate}
\end{proposition}
\begin{apod}{
\begin{enumerate}
  \item Let us denote by $\bar{t}$ the restriction of $t$ to
  $\bar{A}$. Suppose that $\bar{t}$ is a symmetrizing form for
  $\bar{A}$. Let $a \in A$. Then $\hat{t}(a):=(x \mapsto t(ax)) \in
  \mathrm{Hom}_\mathcal{O}(A,\mathcal{O})$. The restriction of
  $\hat{t}(a)$ to $\bar{A}$ belongs to
  $\mathrm{Hom}_\mathcal{O}(\bar{A},\mathcal{O})$ and therefore,
  there exists $\bar{a} \in \bar{A}$ such that $\bar{t}(\bar{a}\bar{x})=t(\bar{a}\bar{x})=t(a\bar{x})$
  for all $\bar{x} \in \bar{A}$. Thus $a-\bar{a} \in \bar{A}^\bot$
  and since $a=\bar{a}+(a-\bar{a})$, we obtain that $A=\bar{A}+
  \bar{A}^\bot$. If $\bar{a} \in \bar{A} \cap
  \bar{A}^\bot$, then we have $\hat{\bar{t}}(\bar{a})=0 \in
  \mathrm{Hom}_\mathcal{O}(\bar{A},\mathcal{O})$.
  Since $\hat{\bar{t}}$ is an isomorphism, we deduce that
  $\bar{a}=0$. Therefore, $A=\bar{A} \oplus \bar{A}^\bot$ and the
  definition of $\mathrm{Br}_{\bar{A}}^A$ is immediate.

  Now suppose that $A=\bar{A} \oplus \bar{A}^\bot$. We will show that the map
  $$\begin{array}{cccl}
      \hat{\bar{t}}: & \bar{A} & \rightarrow & \mathrm{Hom}_\mathcal{O}(\bar{A},\mathcal{O}) \\
                    & \bar{A} & \mapsto     & (\bar{x} \mapsto \bar{t}(\bar{a}\bar{x})=t(\bar{a}\bar{x}))
    \end{array}$$
  is an isomorphism of $\bar{A}$-modules-$\bar{A}$. The map
  $\hat{\bar{t}}$ is obviously injective, because
  $\hat{\bar{t}}(\bar{a})=0$ implies that $\bar{a} \in \bar{A} \cap \bar{A}^\bot$
  and thus $\bar{a}=0$. Now let $\bar{f}$ be an element of
  $\mathrm{Hom}_\mathcal{O}(\bar{A},\mathcal{O})$. The map $\bar{f}$ can be
  extended to a map $f \in \mathrm{Hom}_\mathcal{O}(A,\mathcal{O})$ such that
  $f(a)=\bar{f}(\mathrm{Br}_{\bar{A}}^A(a))$ for all $a \in A$, where
  $\mathrm{Br}_{\bar{A}}^A$ denotes the projection of $A$ onto $\bar{A}$ parallel to
  $\bar{A}^\bot$. Since $t$ is a symmetrizing form for $A$, there
  exists $a \in A$ such that $\hat{t}(a)=f$, \ie $t(ax)=f(x)$ for
  all $x \in A$. Consequently, if $\bar{x} \in \bar{A}$, we have
  $$t(\mathrm{Br}_{\bar{A}}^A(a)\bar{x})=t(a\bar{x})=f(\bar{x})=\bar{f}(\bar{x})$$
  and thus $\hat{\bar{t}}(\mathrm{Br}_{\bar{A}}^A(a))=f$. Hence,
  $\hat{\bar{t}}$ is surjective.
  \item Let $B$ be a sub-$\bar{A}$-module-$\bar{A}$
  of $A$ such that $A=\bar{A} \oplus B$ and $B \subseteq \mathrm{Ker}t$.
  Let $b \in B$. For all $\bar{a} \in
  \bar{A}$, we have $b\bar{a} \in B \subseteq
  \mathrm{Ker}t$ and thus $t(b\bar{a})=0$. Hence $B \subseteq \bar{A}^\bot$.
  Since the restriction of $t$ to $\bar{A}$ is a symmetrizing form
  for $\bar{A}$, we also have $A=\bar{A} \oplus \bar{A}^\bot$. Now
  let $a \in \bar{A}^\bot$. Since $A=\bar{A} \oplus B$, there exist $\bar{a} \in \bar{A}$
  and $b \in B$ such that $a=\bar{a}+b$. Since $b \in \bar{A}^\bot$,
  we must have $a=b \in B$ and therefore, $B=\bar{A}^\bot$.}
\end{enumerate}
\end{apod}

\begin{px}\label{subalgebra of group algebra}
\small{\emph{In the case where $\mathcal{O}=\mathbb{Z}$ and
$A=\mathbb{Z}[G]$
 (G a finite group), let $\bar{A}:=\mathbb{Z}[\bar{G}]$
 be the algebra of a subgroup $\bar{G}$ of $G$.
 Then the morphism $\mathrm{Br}_{\bar{A}}^A$ is the projection given by
 $$\left\{
     \begin{array}{ll}
       g \mapsto g, & \hbox{if $g \in \bar{G}$;} \\
       g \mapsto 0, & \hbox{if $g \notin \bar{G}$.}
     \end{array}
   \right.$$}}
\end{px}\

\begin{definition}\label{symmetric subalgebra}
Let $A$ be a symmetric $\mathcal{O}$-algebra with symmetrizing form
$t$. Let $\bar{A}$ be a subalgebra of $A$. We say that $\bar{A}$ is
a symmetric subalgebra of $A$, if it satisfies the following two
conditions:
\begin{enumerate}
  \item $\bar{A}$ is free (of finite rank) as an $\mathcal{O}$-module and the
  restriction $\mathrm{Res}_{\bar{A}}^A(t)$ of the form $t$ to $\bar{A}$ is a symmetrizing form
  on $\bar{A}$,
  \item $A$ is free (of finite rank) as an $\bar{A}$-module for the action
  of left multiplication by the elements of $\bar{A}$.
\end{enumerate}
\end{definition}

From now on, let us suppose that $\bar{A}$ is a symmetric subalgebra
of $A$ and set $\bar{t}:=\mathrm{Res}_{\bar{A}}^A(t)$. We denote by
$$\mathrm{Ind}_{\bar{A}}^A: _{\bar{A}}\mathrm{\textbf{mod}} \rightarrow _A\mathrm{\textbf{mod}}
\,\textrm{ and }\, \mathrm{Res}_{\bar{A}}^A: _A\mathrm{\textbf{mod}}
\rightarrow _{\bar{A}}\mathrm{\textbf{mod}} $$ the functors defined
as usual by
$$\mathrm{Ind}_{\bar{A}}^A:=A \otimes_{\bar{A}}- \textrm{ where $A$ is viewed as an $A$-module-$\bar{A}$}$$
and
$$\mathrm{Res}_{\bar{A}}^A:=A \otimes_A - \textrm{ where $A$ is viewed as an $\bar{A}$-module-$A$.}$$
Since $A$ is free as $\bar{A}$-module and as module-$\bar{A}$, the
functors $\mathrm{Res}_{\bar{A}}^A$ and $\mathrm{Ind}_{\bar{A}}^A$
are adjoint from both sides.

Moreover, let $K$ be a finite Galois extension of the field of
fractions of $\mathcal{O}$ such that the algebras $KA$ and
$K\bar{A}$ are both split semisimple.

We denote by $\langle-,-\rangle_{KA}$ the scalar product on the
$K$-vector space of trace functions for which the family
$(\chi)_{\chi \in \mathrm{Irr}(KA)}$ is orthonormal and
$\langle-,-\rangle_{K\bar{A}}$ the scalar product on the $K$-vector
space of trace functions for which the family
$(\bar{\chi})_{\bar{\chi} \in \mathrm{Irr}(K\bar{A})}$ is
orthonormal.

Since the functors $\mathrm{Res}_{\bar{A}}^A$ and
$\mathrm{Ind}_{\bar{A}}^A$ are adjoint from both sides, we obtain
the \emph{Frobenius reciprocity} formula:
$$\langle \chi, \mathrm{Ind}_{K\bar{A}}^{KA}(\bar{\chi}) \rangle_{KA} =
\langle \mathrm{Res}_{K\bar{A}}^{KA}(\chi), \bar{\chi}
\rangle_{K\bar{A}} .$$

For every element $\chi \in \mathrm{Irr}(KA)$, let
$$\mathrm{Res}_{K\bar{A}}^{KA}(\chi)=\sum_{\bar{\chi} \in
\mathrm{Irr}(K\bar{A})}m_{\chi,\bar{\chi}}\bar{\chi} \textrm{ (where
} m_{\chi,\bar{\chi}} \in \mathbb{N}).$$

Frobenius reciprocity implies that, for all $\bar{\chi} \in
\mathrm{Irr}(K\bar{A})$,

$$\mathrm{Ind}_{K\bar{A}}^{KA}(\bar{\chi})=\sum_{\chi \in
\mathrm{Irr}(KA)}m_{\chi,\bar{\chi}}\chi .$$

The following property is immediate.

\begin{lemma}\label{mxx'}
For $\chi \in \mathrm{Irr}(KA)$ and $\bar{\chi} \in
\mathrm{Irr}(K\bar{A})$, let $e(\chi)$ and $\bar{e}(\bar{\chi})$ be
respectively the block-idempotents of $KA$ and $K\bar{A}$ associated
with $\chi$ and $\bar{\chi}$. The following conditions are
equivalent:
\begin{description}
  \item[ (i)] $m_{\chi,\bar{\chi}} \neq 0$,
  \item[(ii)] $e(\chi)\bar{e}(\bar{\chi}) \neq 0$.
\end{description}
\end{lemma}

For all $\bar{\chi} \in \mathrm{Irr}(K\bar{A})$, we set
$$\mathrm{Irr}(KA,\bar{\chi}):=\{\chi \in \mathrm{Irr}(KA) \,|\, m_{\chi,\bar{\chi}} \neq 0\},$$

and for all $\chi \in \mathrm{Irr}(KA)$,
$$\mathrm{Irr}(K\bar{A},\chi):=\{\bar{\chi} \in \mathrm{Irr}(K\bar{A}) \,|\, m_{\chi,\bar{\chi}} \neq 0\}.$$

We denote respectively by $s_\chi$ and $s_{\bar{\chi}}$ the Schur
elements of $\chi$ and $\bar{\chi}$ ( with respect to the
symmetrizing forms $t$ for $A$ and $\bar{t}$ for $\bar{A}$).

\begin{lemma}\label{induction and Schur elements}
For all $\bar{\chi} \in \mathrm{Irr}(K\bar{A})$ we have
$$\frac{1}{s_{\bar{\chi}}}=\sum_{\chi \in
\mathrm{Irr}(KA,\bar{\chi})}\frac{m_{\chi,\bar{\chi}}}{s_\chi}  .$$
\end{lemma}
\begin{apod}
{It derives from the relations
$$t=\sum_{\chi \in \mathrm{Irr}(KA)}\frac{1}{s_\chi}\chi,\,\,
\bar{t}=\sum_{\bar{\chi} \in
\mathrm{Irr}(K\bar{A})}\frac{1}{s_{\bar{\chi}}}\bar{\chi},\,\,
\bar{t}= \mathrm{Res}_{\bar{A}}^{A}(t).$$}
\end{apod}

In the next chapters, we will work on the Hecke algebras of complex
reflection groups, which, under certain assumptions, are symmetric.
Sometimes the Hecke algebra of a group $W$ appears as a symmetric
subalgebra of the Hecke algebra of another group $W'$, which
contains $W$. Since we will be mostly interested in the
determination of the blocks of these algebras, it would be helpful,
if we could obtain the blocks of the former from the blocks of the
latter. This is possible with the use of a generalization of some
classical results, known as ``Clifford theory'' (see, for example,
\cite{Da}), to the twisted symmetric algebras of finite groups and
more precisely of finite cyclic groups. For the application of these
results to the Hecke algebras, the reader may refer to the Appendix
of this thesis.

\begin{definition}\label{symmetric algebra of a finite group}
We say that a symmetric $\mathcal{O}$-algebra $(A,t)$ is the twisted
symmetric algebra of a finite group $G$ over the subalgebra
$\bar{A}$, if the following conditions are satisfied:
\begin{itemize}
  \item $\bar{A}$ is a symmetric subalgebra of $A$,
  \item There exists a family $\{A_g \,|\, g \in G\}$ of
  $\mathcal{O}$-submodules of $A$ such that
  \begin{description}
    \item[(a)] $A= \bigoplus_{g \in G}A_g$,
    \item[(b)] $A_gA_h=A_{gh}$ for all $g,h \in G$,
    \item[(c)] $A_1=\bar{A}$,
    \item[(d)] $t(A_g)=0$ for all $g \in G, g \neq  1$,
    \item[(e)] $A_g \cap A^\times \neq \emptyset$ for all $g \in G$ (where
    $A^\times$ is the set of units of $A$).
    \end{description}
\end{itemize}
\end{definition}
If that is the case, then proposition $\ref{when symmetric
subalgebra}$ implies that
$$\bigoplus_{g \in G-\{1\}} A_g =
\bar{A}^\bot.$$

\begin{lemma}\label{factor group}
Let $a_g \in A_g$ such that $a_g$ is a unit in $A$. Then
$$A_g=a_g\bar{A}=\bar{A}a_g.$$
\end{lemma}
\begin{apod}{Since $a_g \in A_g$, property (b)
implies that $a_g^{-1} \in A_{g^{-1}}$. If $a \in A_g$, then $
a_g^{-1}a \in A_1=\bar{A}$. We have $a = a_g{a_g}^{-1}a \in
a_g\bar{A}$ and thus $A_g \subseteq a_g\bar{A}$. Property (b)
implies the inverse inclusion. In the same way, we show that
$A_g=\bar{A}a_g$.}
\end{apod}

From now on, let $(A,t)$ be the twisted symmetric algebra of a
finite group $G$ over the subalgebra $\bar{A}$. Due to property (e)
and the lemma above, for all $g \in G$, there exist $a_g \in A_g
\cap A^\times$ such that $A_g=a_g\bar{A}=\bar{A}a_g$.

\begin{proposition}\label{dual basis of a symmetric algebra of a finite group}
Let $(\bar{e}_i)_{i \in I}$ be a basis of $\bar{A}$ over
$\mathcal{O}$ and $(\bar{e}'_i)_{i \in I}$ its dual with respect to
the symmetrizing form $\bar{t}$.  We fix a system of representatives
$\mathrm{Rep}(A/\bar{A}):=\{a_g \,|\, g \in G\}.$ Then the families
$$ (\bar{e}_ia_g)_{i \in I, a_g \in \mathrm{Rep}(A/\bar{A})}
 \textrm{ and } (a_g^{-1}\bar{e}'_i)_{i \in I,a_g \in \mathrm{Rep}(A/\bar{A})}$$
are two $\mathcal{O}$-bases of $A$ dual to each other.
\end{proposition}

\subsection*{Action of $G$ on $Z\bar{A}$}

\begin{lemma}\label{definition of ga}
Let $\bar{a} \in Z\bar{A}$ and $g \in G$. There exists a unique
element $g(\bar{a})$ of $\bar{A}$ satisfying
$$g(\bar{a})\mathfrak{g}=\mathfrak{g}\bar{a} \textrm{ for all } \mathfrak{g} \in
A_g.\,\,\,\,\,\,\, (\dag)$$ If $a_g \in A^\times$ such that
$A_g=a_g\bar{A}$, then $$g(\bar{a})=a_g\bar{a}{a_g}^{-1}.$$
\end{lemma}
\begin{apod}{Set $g(\bar{a}):=a_g\bar{a}{a_g}^{-1}.$
 Then, for all $\mathfrak{g} \in A_g$, we have
${a_g}^{-1}\mathfrak{g} \in
  \bar{A}$ and
  $g(\bar{a})\mathfrak{g}=a_g\bar{a}{a_g}^{-1}\mathfrak{g}=a_g{a_g}^{-1}\mathfrak{g}\bar{a}=\mathfrak{g}\bar{a}$.
Now, let $y$ be another element of $A$ such that $y
\mathfrak{g}=\mathfrak{g}\bar{a} \textrm{ for all } \mathfrak{g} \in
A_g.$ Then $ya_g=a_g\bar{a}$ and hence $y=g(\bar{a})$. Therefore,
$g(\bar{a})$ is the unique element of $A$ which satisfies $(\dag)$.}
\end{apod}

\begin{proposition}\label{action of G on ZA}
The map $\bar{a} \mapsto g(\bar{a})$ defines an action of $G$ as
ring automorphism of $Z\bar{A}$.
\end{proposition}
\begin{apod}{Let $\bar{a} \in Z\bar{A}$, $g \in G$ and $a_g \in A^\times$ such that
$A_g=a_g\bar{A}$. We will show that $g(\bar{a}) \in Z\bar{A}$. If
$\bar{x} \in \bar{A}$, then $\bar{x}a_g \in A_g$ and we have
$$\bar{x}g(\bar{a})a_g=\bar{x}a_g\bar{a}=g(\bar{a})\bar{x}a_g.$$
Multiplying both sides by ${a_g}^{-1}$, we obtain that
$$\bar{x}g(\bar{a})=g(\bar{a})\bar{x}$$
and hence, $g(\bar{a}) \in Z\bar{A}$.

Since the identity $1$ of $A$ lies in $A_1$, we have
$1_G(\bar{a})=\bar{a}$. If $g_1,g_2 \in G$, then equation ($\dag$)
gives
$$g_1(g_2(\bar{a}))a_{g_1}a_{g_2}=a_{g_1}g_2(\bar{a})a_{g_2}=
a_{g_1}a_{g_2}\bar{a}.$$ Due to property (b) of the definition
$\ref{symmetric algebra of a finite group}$, the product
$a_{g_1}a_{g_2}$ generates the submodule $A_{g_1g_2}$. Therefore,
$g_1(g_2(\bar{a}))u=u\bar{a}$ for all $u \in A_{g_1g_2}$. By lemma
$\ref{definition of ga}$, we obtain that
$(g_1g_2)(\bar{a})=g_1(g_2(\bar{a})).$

Finally, let us fix $g \in G$. By definition, the map $\bar{a}
\mapsto g(\bar{a})$ is an additive automorphism of $Z\bar{A}$. If
$\bar{a}_1,\bar{a}_2 \in Z\bar{A}$, then
$$\mathfrak{g}\bar{a}_1\bar{a}_2=g(\bar{a}_2)\mathfrak{g}\bar{a}_1=g(\bar{a}_1)g(\bar{a}_2)\mathfrak{g}
\textrm{ for all } \mathfrak{g} \in A_g.$$ By lemma $\ref{definition
of ga}$, we obtain that
$g(\bar{a}_1\bar{a}_2)=g(\bar{a}_1)g(\bar{a}_2).$}
\end{apod}

Now let $\bar{b}$ be a block(-idempotent) of $\bar{A}$. If $g \in
G$, then $g(\bar{b})$ is also a block of $\bar{A}$. So we must have
either $g(\bar{b})=\bar{b}$ or $g(\bar{b})$ orthogonal to $\bar{b}$.
Set
$$\mathrm{Tr}(G,\bar{b}):=\sum_{g \in G/G_{\bar{b}}}g(\bar{b}),$$
where $G_{\bar{b}}:=\{g \in G \,|\, g(\bar{b})=\bar{b}\}$. It is
clear that
\begin{itemize}
  \item $\bar{b}$ is a central idempotent of $\bigoplus_{g \in G_{\bar{b}}}A_g=:A_{G_{\bar{b}}}$,
  \item $\mathrm{Tr}(G,\bar{b})$ is a central idempotent of $A$.
\end{itemize}
From now on, let $b:=\mathrm{Tr}(G,\bar{b})$ and $\mathfrak{g} \in
A_g$. We have
$$\bar{b}\mathfrak{g}\bar{b}= \left\{
                           \begin{array}{ll}
                             \mathfrak{g}\bar{b}=\bar{b}\mathfrak{g}, & \hbox{if $g \in G_{\bar{b}}$;} \\
                             0, & \hbox{if $g \notin G_{\bar{b}}$,}
                           \end{array}
                         \right.$$
$$\bar{b}\mathfrak{g}b=\bar{b}\mathfrak{g} \,\,\,\,\,\textrm{ and }\,\,\,\,\, b\mathfrak{g}\bar{b}=\mathfrak{g}\bar{b}.$$

\begin{proposition}\label{morita}
The applications
$$\left\{
    \begin{array}{ll}
      bA\bar{b} \otimes_{A_{G_{\bar{b}}}\bar{b}} \bar{b}Ab \rightarrow Ab \\
      ba_g\bar{a}\bar{b} \otimes \bar{b}\bar{a}'a_{g'}b \mapsto a_g\bar{a}\bar{b}\bar{a}'a_{g'}\\
      ab \mapsto \sum_{g \in G/G_{\bar{b}}} aa_g\bar{b} \otimes \bar{b}{a_g}^{-1},
    \end{array}
  \right.
 $$
and
$$\left\{
    \begin{array}{ll}
      \bar{b}Ab \otimes_{Ab} bA\bar{b} \rightarrow A_{G_{\bar{b}}}\bar{b}\\
      \bar{b}\bar{a}a_gb \otimes ba_{g'}\bar{a}'\bar{b} \mapsto \left\{
                                                            \begin{array}{ll}
                                                              \bar{b}\bar{a}a_ga_{g'}\bar{a}'\bar{b}, & \hbox{if $gg' \in G_{\bar{b}}$;} \\
                                                              0, & \hbox{if not.}
                                                            \end{array}
                                                          \right. \\
      \bar{a}a_g\bar{b} \mapsto \bar{a}a_g\bar{b} \otimes \bar{b} \textrm{ } (\textrm{ where } g \in G_{\bar{b}}),
    \end{array}
  \right.$$ define isomorphisms inverse to each other
$$ bA\bar{b} \otimes_{A_{G_{\bar{b}}}\bar{b}} \bar{b}Ab\, \tilde{\leftrightarrow}\, Ab  \,\,\,\textrm{ and
}\,\,\, \bar{b}Ab \otimes_{Ab} bA\bar{b} \,\tilde{\leftrightarrow}\,
A_{G_{\bar{b}}}\bar{b}.$$ Therefore, $bA\bar{b}$ and $\bar{b}Ab$ are
Morita equivalent. In particular, the functors
$$\mathrm{Ind}_{\bar{A}}^A=(bA\bar{b} \otimes_{A_{G_{\bar{b}}}\bar{b}}-) \,\,\,\textrm{ and
}\,\,\, \bar{b} \cdot \mathrm{Res}_{\bar{A}}^A= (\bar{b}Ab
\otimes_{Ab} -)$$ define category equivalences inverse to each other
between $_{A_{G_{\bar{b}}}\bar{b}}\emph{\textbf{mod}}$ and\\
$_{Ab}\emph{\textbf{mod}}$.
\end{proposition}

\subsection*{Multiplication of an $A$-module by an $\mathcal{O}G$-module}

Let $X$ be an $A$-module and $\rho:A \rightarrow
\mathrm{End}_\mathcal{O}(X)$ be the structural morphism. We define
an additive functor
$$ X \cdot - : _{\mathcal{O}G}\mathrm{\textbf{mod}} \rightarrow _A
\mathrm{\textbf{mod}}, Y \mapsto X \cdot Y $$ as follows:

If $Y$ is an $\mathcal{O}G$-module and $\sigma: \mathcal{O}G
\rightarrow \mathrm{End}_\mathcal{O}(Y)$ is the structural morphism,
we denote by $X \cdot Y$ the $\mathcal{O}$-module $X
\otimes_\mathcal{O} Y$. The action of $A$ on the latter is given by
the morphism
$$\rho\cdot\sigma:A\rightarrow\mathrm{End}_\mathcal{O}(X
\otimes Y), \bar{a}a_g \mapsto \rho(\bar{a}a_g) \otimes \sigma(g).$$

\begin{proposition}\label{multiplication of A and G modules}
Let $X$ be an $\bar{A}$-module. The application $$A
\otimes_{\bar{A}}X \rightarrow X \cdot \mathcal{O}G$$ defined by
$$a_g\otimes_{\bar{A}}x \mapsto \rho(a_g)(x)\otimes_\mathcal{O}g
\textrm{ (for all $x \in X$ and $g \in G$)}$$ is an isomorphism of
$A$-modules
$$\mathrm{Ind}_{\bar{A}}^A(X) \tilde{\rightarrow} X \cdot \mathcal{O}G.$$
\end{proposition}

\subsection*{Induction and restriction of $KA$-modules and
$K\bar{A}$-modules}

Let $X$ be a $KA$-module of character $\chi$ and $Y$ a $KG$-module
of character $\xi$. We denote by $\chi \cdot \xi$ the character of
the $KA$-module $X \cdot Y$. From now on, all group algebras over
$K$ will be considered split semisimple.

\begin{proposition}\label{1.39}
Let $\chi$ be an irreducible character of $KA$ which restricts to an
irreducible character $\bar{\chi}$ of $K\bar{A}$. Then
\begin{enumerate}
  \item The characters $(\chi \cdot \xi)_{\xi \in \mathrm{Irr}(KG)}$
  are distinct irreducible characters of $KA$.
  \item We have $$\mathrm{Ind}_{K\bar{A}}^{KA}(\bar{\chi})=
  \sum_{\xi \in \mathrm{Irr}(KG)} \xi(1)\chi \cdot \xi.$$
\end{enumerate}
\end{proposition}
\begin{apod}{The second relation results from proposition $\ref{multiplication of A and G
modules}$. Frobenius reciprocity now gives
$$\begin{array}{ccc}
\langle \mathrm{Ind}_{K\bar{A}}^{KA}(\bar{\chi}),
\mathrm{Ind}_{K\bar{A}}^{KA}(\bar{\chi}) \rangle_{KA} & = & \langle
\mathrm{Res}_{K\bar{A}}^{KA}(\sum_{\xi \in \mathrm{Irr}(KG)}
\xi(1)\chi \cdot \xi), \bar{\chi} \rangle_{K\bar{A}} \\
     & = & \langle
\sum_{\xi \in \mathrm{Irr}(KG)} \xi(1)^2 \bar{\chi},\bar{\chi}
\rangle_{K\bar{A}} \\
     & = & \sum_{\xi \in \mathrm{Irr}(KG)} \xi(1)^2 = |G|,
  \end{array}$$
hence from the relation in part 2 we obtain
$$|G|=\sum_{\xi,\xi' \in \mathrm{Irr}(KG)} \xi(1)\xi'(1)\langle \chi
\cdot \xi,\chi \cdot \xi'\rangle_{KA}.$$ Since $|G|=\sum_{\xi \in
\mathrm{Irr}(KG)} \xi(1)^2$, we must have $\langle \chi \cdot
\xi,\chi \cdot \xi'\rangle_{KA}=\delta_{\xi,\xi'}$ and the proof is
complete.}
\end{apod}

For all $\bar{\chi} \in \mathrm{Irr}(K\bar{A})$, we denote by
$\bar{e}(\bar{\chi})$ the block of $K\bar{A}$ associated with
$\bar{\chi}$. We have seen that if $g \in G$, then
$g(\bar{e}(\bar{\chi}))$ is also a block of $K\bar{A}$. Since
$K\bar{A}$ is split semisimple, it must be associated with an
irreducible character $g(\bar{\chi})$ of $K\bar{A}$. Thus, we can
define an action of $G$ on $\mathrm{Irr}(K\bar{A})$ such that for
all $g \in G$, $\bar{e}(g(\bar{\chi}))=g(\bar{e}(\bar{\chi}))$. We
denote by $G_{\bar{\chi}}$ the stabilizer of $\bar{\chi}$ in $G$.

\begin{proposition}\label{1.40}
Let $\bar{\chi} \in \mathrm{Irr}(K\bar{A})$ and suppose that there
exists $\tilde{\chi} \in  \mathrm{Irr}(KA_{G_{\bar{\chi}}})$ such
that
$\mathrm{Res}_{K\bar{A}}^{KA_{G_{\bar{\chi}}}}(\tilde{\chi})=\bar{\chi}$.
We set
$$\chi:=\mathrm{Ind}_{KA_{G_{\bar{\chi}}}}^{KA}(\tilde{\chi}) \,\textrm{
and }\,
\chi_\xi:=\mathrm{Ind}_{KA_{G_{\bar{\chi}}}}^{KA}(\tilde{\chi} \cdot
\xi) \textrm{ for all } \xi \in \mathrm{Irr}(KG_{\bar{\chi}}).$$
Then
\begin{enumerate}
  \item The characters $(\chi_\xi)_{\xi \in \mathrm{Irr}(KG_{\bar{\chi}})}$
  are distinct irreducible characters of $KA$.
  \item We have $$\mathrm{Ind}_{K\bar{A}}^{KA}(\bar{\chi})=
  \sum_{\xi \in \mathrm{Irr}(KG_{\bar{\chi}})} \xi(1)\chi_\xi.$$
  In particular,
  $$m_{\chi_\xi,\bar{\chi}}=\xi(1) \,\textrm{ and }\,
  \chi_\xi(1)=|G:G_{\bar{\chi}}|\bar{\chi}(1)\xi(1).$$
  \item For all $\xi \in \mathrm{Irr}(KG_{\bar{\chi}})$, we have
  $$s_{\chi_\xi}=\frac{|G_{\bar{\chi}}|}{\xi(1)}s_{\bar{\chi}}.$$
\end{enumerate}
\end{proposition}
\begin{apod}{
\begin{enumerate}
  \item By proposition $\ref{1.39}$, we obtain that the characters
   $(\tilde{\chi} \cdot \xi)_{\xi \in \mathrm{Irr}(KG_{\bar{\chi}})}$
   are distinct irreducible characters of
   $\mathrm{Irr}(KA_{G_{\bar{\chi}}})$. Now let $\bar{e}(\bar{\chi})$ be the
   block of $K\bar{A}$ associated with the irreducible character
   $\bar{\chi}$. We have seen that $\bar{e}(\bar{\chi})$ is a central
   idempotent of $KA_{G_{\bar{\chi}}}$. Proposition $\ref{morita}$
   implies that the functor
   $\mathrm{Ind}_{K\bar{A}}^{KA}$ defines a Morita equivalence
   between the category
   $_{KA_{G_{\bar{\chi}}}\bar{e}(\bar{\chi})}\mathrm{\textbf{mod}}$
   and its image. Therefore, the characters $(\mathrm{Ind}_{KA_{G_{\bar{\chi}}}}^{KA}(\tilde{\chi}
   \cdot \xi))_{\xi \in \mathrm{Irr}(KG_{\bar{\chi}})}$ are distinct
   irreducible characters of $KA$.
  \item By proposition $\ref{1.39}$, we obtain that
   $$\mathrm{Ind}_{K\bar{A}}^{KA_{G_{\bar{\chi}}}}(\bar{\chi})=
   \sum_{\xi \in \mathrm{Irr}(KG_{\bar{\chi}})} \xi(1)\tilde{\chi} \cdot \xi.$$
   Applying $\mathrm{Ind}_{KA_{G_{\bar{\chi}}}}^{KA}$ to both sides gives us the
   required relation. Obviously, $m_{\chi_\xi,\bar{\chi}}=\xi(1)$.

   Now let us calculate the value of $\chi_\xi(\bar{a})$ for any
   $\bar{a} \in \bar{A}$.
   Let $Y$ be an irreducible $KA_{G_{\bar{\chi}}}$-module of character $\psi$.Then
   $\mathrm{Ind}_{KA_{G_{\bar{\chi}}}}^{KA}(Y)=
   KA \otimes_{KA_{G_{\bar{\chi}}}} Y$ has character
   $\mathrm{Ind}_{KA_{G_{\bar{\chi}}}}^{KA}(\psi)$.
   We have $KA = \bigoplus_{g \in G/G_{\bar{\chi}}}a_g K\bar{A}$.
   Let $\bar{a} \in \bar{A}$. Then
   $$\begin{array}{ccl}
       \bar{a}\mathrm{Ind}_{KA_{G_{\bar{\chi}}}}^{KA}(Y) & = & \bigoplus_{g \in G/G_{\bar{\chi}}}\bar{a}a_g K\bar{A} \otimes_{KA_{G_{\bar{\chi}}}} Y \\
        & = & \bigoplus_{g \in G/G_{\bar{\chi}}}a_g({a_g}^{-1}\bar{a}a_g) K\bar{A} \otimes_{KA_{G_{\bar{\chi}}}} Y\\
        & = & \bigoplus_{g \in G/G_{\bar{\chi}}}a_g K\bar{A} \otimes_{KA_{G_{\bar{\chi}}}} ({a_g}^{-1}\bar{a}a_g)Y.
     \end{array}$$
   Thus,
   $\mathrm{Ind}_{KA_{G_{\bar{\chi}}}}^{KA}(\psi)(\bar{a})=
   \sum_{g \in G/G_{\bar{\chi}}}\psi({a_g}^{-1}\bar{a}a_g)$
   and
   $$\chi_\xi(\bar{a})=\sum_{g \in G/G_{\bar{\chi}}}(\tilde{\chi} \cdot
   \xi)({a_g}^{-1}\bar{a}a_g)= \sum_{g \in
   G/G_{\bar{\chi}}}\bar{\chi}({a_g}^{-1}\bar{a}a_g)\xi(1).$$
   Therefore,
   $$\chi_\xi(1)=\sum_{g \in G/G_{\bar{\chi}}}\bar{\chi}(1)\xi(1)=
   |G:G_{\bar{\chi}}|\bar{\chi}(1)\xi(1).$$
  \item Let $(\bar{e}_i)_{i \in I}$ be a basis of $\bar{A}$ as
   $\mathcal{O}$-module and let $(\bar{e}_i')_{i \in I}$ be its
   dual with respect to the symmetrizing form $\bar{t}$.
   Proposition $\ref{schur elements and idempotents}$(3), in combination with proposition
   $\ref{dual basis of a symmetric algebra of a finite group}$, gives
   $$s_{\chi_\xi}\chi_\xi(1)^2=\chi_\xi(\sum_{i \in I,g \in G}\bar{e}_i'a_g{a_g}^{-1}\bar{e}_i)=
     \chi_\xi(|G|\sum_{i \in I}\bar{e}_i'\bar{e}_i).$$
   However, $\sum_{i \in I,g \in G}\bar{e}_i'a_g{a_g}^{-1}\bar{e}_i$
   belongs to to center of $A$ (by lemma $\ref{properties of the casimir}$) and thus, for all $h \in G$,
   $${a_h}^{-1}(\sum_{i \in I}\bar{e}_i'a_g{a_g}^{-1}\bar{e}_i)a_h=
     \sum_{i \in I}\bar{e}_i'a_g{a_g}^{-1}\bar{e}_i=
     |G|\sum_{i \in I}\bar{e}_i'\bar{e}_i.$$
   Since $\sum_{i \in I}\bar{e}_i'\bar{e}_i \in \bar{A}$, by
   part 2,
   $$\begin{array}{ccl}
       \chi_\xi(|G|\sum_{i \in I}\bar{e}_i'\bar{e}_i) & = & \sum_{h \in G/G_{\bar{\chi}}}\bar{\chi}({a_h}^{-1}(|G|\sum_{i \in
   I}\bar{e}_i'\bar{e}_i)a_h)\xi(1) \\
    & & \\
        & = & \sum_{h \in G/G_{\bar{\chi}}}\bar{\chi}(|G|\sum_{i \in
   I}\bar{e}_i'\bar{e}_i)\xi(1) \\
   & & \\
        & = & |G:G_{\bar{\chi}}| |G| \xi(1) \bar{\chi}(\sum_{i \in
   I}\bar{e}_i'\bar{e}_i) \\
   & & \\
        & = & |G:G_{\bar{\chi}}|^2 |G_{\bar{\chi}}| \xi(1) s_{\bar{\chi}} \bar{\chi}(1)^2.
     \end{array}$$
   So we have
   $$s_{\chi_\xi}\chi_\xi(1)^2=|G:G_{\bar{\chi}}|^2 |G_{\bar{\chi}}| \xi(1) \bar{\chi}(1)^2 s_{\bar{\chi}}.$$
   Replacing $\chi_\xi(1)=|G:G_{\bar{\chi}}|\bar{\chi}(1)\xi(1)$
   gives
   $$s_{\chi_\xi} \xi(1) = |G_{\bar{\chi}}| s_{\bar{\chi}} .$$
\end{enumerate}}
\end{apod}

Now let $\bar{\Omega}$ be the orbit of the character $\bar{\chi} \in
\mathrm{Irr}(K\bar{A})$ under the action of $G$. We have
$|\bar{\Omega}|=|G|/|G_{\bar{\chi}}|$. Define
$$\bar{e}(\bar{\Omega})=\sum_{g \in
G/G_{\bar{\chi}}}\bar{e}(g(\bar{\chi}))= \sum_{g \in
G/G_{\bar{\chi}}}g(\bar{e}(\bar{\chi})).$$ 

If $\bar{\chi} \in \bar{\Omega}$, the set
$\mathrm{Irr}(KA,\bar{\chi})$ depends only on $\bar{\Omega}$ and we
set $\mathrm{Irr}(KA,\bar{\Omega}):=\mathrm{Irr}(KA,\bar{\chi})$.
The idempotent $\bar{e}(\bar{\Omega})$ belongs to the algebra
$(ZK\bar{A})^G$ of the elements in the center of $K\bar{A}$ fixed by
$G$ and thus to the center of $KA$ (since it commutes with all
elements of $\bar{A}$ and all $a_g$, $g \in G$). Therefore, it must
be a sum of blocks of $KA$, \ie
$$\bar{e}(\bar{\Omega}) = \sum_{\chi \in
\mathrm{Irr}(KA,\bar{\Omega})}e(\chi).$$

Let $X$ be an irreducible $KA$-module of character $\chi$ and
$\bar{X}$ an irreducible $K\bar{A}$-submodule of
$\mathrm{Res}_{K\bar{A}}^{KA}(X)$ of character $\bar{\chi}$. For $g
\in G$, the $K\bar{A}$-submodule $g(\bar{X})$ of
$\mathrm{Res}_{K\bar{A}}^{KA}(X)$ has character $g(\bar{\chi})$.
Then $\sum_{g \in G}g(\bar{X})$ is a $KA$-submodule of $X$. We
deduce that
$$\mathrm{Res}_{K\bar{A}}^{KA}(X)=( \bigoplus_{g \in
G/G_{\bar{\chi}}}g(\bar{X}))^{m_{\chi,\chi'}},$$ \ie
$$\mathrm{Res}_{K\bar{A}}^{KA}(\chi)= m_{\chi,\chi'}\sum_{g \in
G/G_{\bar{\chi}}}g(\bar{\chi}).$$ In particular, we see that
$\mathrm{Irr}(K\bar{A},\chi)$ is an orbit of $G$ on
$\mathrm{Irr}(K\bar{A})$. Notice that
$\chi(1)=m_{\chi,\chi'}|\bar{\Omega}|\bar{\chi}(1)$.\\
\\
\emph{Case where $G$ is cyclic}\\
\\
Let $G$ be a cyclic group of order $d$ and let $g$ be a generator of
$G$ (we can choose
$\mathrm{Rep}(A/\bar{A})=\{1,a_g,{a_g}^2,\ldots,{a_g}^{d-1}\}$). We
will show that the assumptions of proposition $\ref{1.40}$ are
satisfied for all irreducible characters of $K\bar{A}$.

Let $\bar{X}$ be an irreducible $K\bar{A}$-module and let
$\bar{\rho}:K\bar{A}\rightarrow \mathrm{End}_K(\bar{X})$ be the
structural morphism. Since the representation of $\bar{X}$ is
invariant by the action of $G_{\bar{\chi}}$, there exists an
automorphism $\alpha$ of the $K$-vector space $\bar{X}$ such that
$$\alpha \bar{\rho}(\bar{a}) \alpha^{-1}= g(\bar{\rho})(\bar{a}),$$
for all $g \in G_{\bar{\chi}}$.

The subgroup $G_{\bar{\chi}}$ is also cyclic. Let
$d(\bar{\chi}):=|G_{\bar{\chi}}|$. Then
$$\bar{\rho}(\bar{a})=\alpha^{d(\bar{\chi})} \bar{\rho}(\bar{a}) \alpha^{-d(\bar{\chi})}.$$
Since $\bar{X}$ is irreducible and $K\bar{A}$ is split semisimple,
$\alpha^{d(\bar{\chi})}$ must be a scalar. Instead of enlarging the
field $K$, we can assume that $K$ contains a $d(\bar{\chi})$-th root
of that scalar. By dividing $\alpha$ by that root, we reduce to the
case where $\alpha^{d(\bar{\chi})}=1$.

This allows us to extend the structural morphism
$\bar{\rho}:K\bar{A}\rightarrow \mathrm{End}_K(\bar{X})$ to a
morphism
$$\tilde{\rho}:KA_{G_{\bar{\chi}}}\rightarrow
\mathrm{End}_K(\bar{X})$$ such that
$$\tilde{\rho}(\bar{a}a_h^j):=\bar{\rho}(\bar{a})\alpha^j
 \textrm{ for } 0\leq j<d(\bar{\chi}),$$
 where $h:=g^{d/d(\bar{\chi})}$ generates $G_{\bar{\chi}}$.
The morphism $\tilde{\rho}$ defines a $KA_{G_{\bar{\chi}}}$-module
$\tilde{X}$ of character $\tilde{\chi}$.

Since the group $G$ is abelian, the set $\mathrm{Irr}(KG)$ forms a
group, which we denote by $G^\vee$. The application $\psi \mapsto
\psi \cdot \xi$, where $\psi \in \mathrm{Irr}(KA)$ and $\xi \in
G^\vee$, defines an action of $G^\vee$ on $\mathrm{Irr}(KA)$.

Let $\Omega$ be the orbit of $\tilde{\chi}$ under the action of
$(G_{\bar{\chi}})^\vee$. By proposition $\ref{1.39}$, we obtain that
$\Omega$ is a regular orbit (\ie $|\Omega|=|G_{\bar{\chi}}|$) and
that $\Omega=\mathrm{Irr}(KA_{G_{\bar{\chi}}},\bar{\chi})$.

Like in proposition $\ref{1.40}$, we introduce the notations
$$\chi:=\mathrm{Ind}_{KA_{G_{\bar{\chi}}}}^{KA}(\tilde{\chi})\textrm{ and }
\chi_\xi:=\mathrm{Ind}_{KA_{G_{\bar{\chi}}}}^{KA}(\tilde{\chi} \cdot
\xi) \textrm{ for all } \xi \in (G_{\bar{\chi}})^\vee.$$

Then
$$\mathrm{Irr}(KA,\bar{\chi})=\{\chi_\xi|\xi\in(G_{\bar{\chi}})^\vee\}
\textrm{ and } m_{\chi_\xi,\bar{\chi}}=\xi(1)=1 \textrm{ for all }
\xi \in (G_{\bar{\chi}})^\vee.$$

Recall that $d(\bar{\chi})=|G_{\bar{\chi}}|$. There exists a
surjective morphism $G \twoheadrightarrow G_{\bar{\chi}}$ defined by
$g \mapsto g^{d/d(\bar{\chi})}$, which induces an inclusion
$(G_{\bar{\chi}})^\vee \hookrightarrow G^\vee$. If $\xi \in
(G_{\bar{\chi}})^\vee$, we denote (abusing notation) by $\xi$ its
image in $G^\vee$ by the above injection. It is easy to check that
$\chi_\xi=\chi \cdot \xi$.

Hence we have proved the following result

\begin{proposition}\label{1.42}
If the group $G$ is cyclic, there exists a bijection
$$\begin{array}{ccc}
    \mathrm{Irr}(K\bar{A})/G & \tilde{\leftrightarrow} & \mathrm{Irr}(KA)/G^\vee \\
    \bar{\Omega} & \leftrightarrow & \Omega
  \end{array}$$
such that
$$\bar{e}(\bar{\Omega})=e(\Omega),\, |\bar{\Omega}||\Omega|=|G| \textrm{ and }
\left\{
  \begin{array}{ll}
    \forall \chi \in \Omega, & \mathrm{Res}_{K\bar{A}}^{KA}(\chi)=\sum_{\bar{\chi} \in \bar{\Omega}}\bar{\chi}\\
    \forall \bar{\chi} \in \bar{\Omega}, &\mathrm{Ind}_{K\bar{A}}^{KA}(\bar{\chi})=\sum_{\chi \in \Omega}\chi
  \end{array}
\right.
$$
Moreover, for all $\chi \in \Omega$ and $\bar{\chi} \in
\bar{\Omega}$, we have
$$s_\chi = |\Omega| s_{\bar{\chi}}.$$
\end{proposition}

\subsection*{Blocks of $A$ and blocks of $\bar{A}$}

Let us denote by $\mathrm{Bl}(A)$ the set of blocks of $A$ and by
$\mathrm{Bl}(\bar{A})$ the set of blocks of $\bar{A}$. For $\bar{b}
\in \mathrm{Bl}(\bar{A})$, we have set
$$\mathrm{Tr}(G,\bar{b}):=\sum_{g \in G/G_{\bar{b}}}g(\bar{b}).$$

The algebra $(Z\bar{A})^G$ is contained in both $Z\bar{A}$ and $ZA$
and the set of its blocks is
$$\mathrm{Bl}((Z\bar{A})^G)=\{\mathrm{Tr}(G,\bar{b}) \,|\, \bar{b} \in
\mathrm{Bl}(\bar{A})/G\}.$$ Moreover, $\mathrm{Tr}(G,\bar{b})$ is a
sum of blocks of $A$ and we define the subset
$\mathrm{Bl}(A,\bar{b})$ of $\mathrm{Bl}(A)$ as follows
$$\mathrm{Tr}(G,\bar{b}):=\sum_{b \in \mathrm{Bl}(A,\bar{b})}b.$$

\begin{lemma}\label{1.43}
Let $\bar{b}$ be a block of $\bar{A}$ and
$\bar{B}:=\mathrm{Irr}(K\bar{A}\bar{b})$. Then
\begin{enumerate}
  \item For all $\bar{\chi} \in \bar{B}$, we have $G_{\bar{\chi}}
  \subseteq G_{\bar{b}}$.
  \item We have
  $$\mathrm{Tr}(G,\bar{b})=\sum_{\bar{\chi} \in \bar{B}/G}
  \mathrm{Tr}(G,\bar{e}(\bar{\chi}))=
  \sum_{\{\bar{\Omega}|\bar{\Omega} \cap \bar{B} \neq \emptyset \}}\bar{e}(\bar{\Omega}).$$
\end{enumerate}
\end{lemma}
\begin{apod}{
\begin{enumerate}
  \item If $g \notin G_{\bar{b}}$, then $\bar{b}$ and $g(\bar{b})$
  are orthogonal.
  \item Note that $\bar{b}=\sum_{\bar{\chi} \in \bar{B}}\bar{e}(\bar{\chi})=
  \sum_{\bar{\chi} \in
  \bar{B}/G_{\bar{b}}}\mathrm{Tr}(G_{\bar{b}},\bar{e}(\bar{\chi}))$.
  Thus
  $$\mathrm{Tr}(G,\bar{b})=\sum_{\bar{\chi} \in \bar{B}/G}
  \mathrm{Tr}(G,\bar{e}(\bar{\chi}))=
  \sum_{\{\bar{\Omega}|\bar{\Omega} \cap \bar{B} \neq \emptyset \}}\bar{e}(\bar{\Omega}),$$
  by the definition of $\bar{e}(\bar{\Omega})$.}
\end{enumerate}
\end{apod}

Now let $G^\vee:=\mathrm{Hom}(G,K^\times)$. We suppose that $K=F$.
The multiplication of the characters of $KA$ by the characters of
$KG$ defines an action of the group $G^\vee$ on $\mathrm{Irr}(KA)$.
This action is induced by the operation of $G^\vee$ on the algebra
$A$, which is defined in the following way:
$$ \xi \cdot (\bar{a}a_g) := \xi(g)\bar{a}a_g \,\,\textrm{ for all }
\xi \in G^\vee, \bar{a} \in \bar{A}, g \in G.$$ In particular,
$G^\vee$ acts on the set of blocks of $A$. Let $b$ be a block of
$A$. Denote by $\xi \cdot b$ the product of $\xi$ and $b$ and by
$(G^\vee)_b$ the stabilizer of $b$ in $G^\vee$. We set
$$\mathrm{Tr}(G^\vee,b):=\sum_{\xi\in G^\vee/(G^\vee)_b}\xi\cdot b.$$

The set of blocks of the algebra $(ZA)^{G^\vee}$ is given by
$$\mathrm{Bl}((ZA)^{G^\vee})=\{\mathrm{Tr}(G^\vee,b) \,|\, b \in
\mathrm{Bl}(A)/G^\vee\}.$$

The following lemma is the analogue of lemma $\ref{1.43}$

\begin{lemma}\label{1.44}
Let $b$ be a block of $A$ and $B:=\mathrm{Irr}(KAb)$. Then
\begin{enumerate}
  \item For all $\chi \in B$, we have $(G^\vee)_\chi \subseteq (G^\vee)_b$.
  \item We have
  $$\mathrm{Tr}(G^\vee,b)=\sum_{\chi\in B/G^\vee} \mathrm{Tr}(G^\vee,e(\chi))=
  \sum_{\{\Omega|\Omega\cap B \neq \emptyset \}}e(\Omega).$$
\end{enumerate}
\end{lemma}
\emph{Case where $G$ is cyclic}\\
\\
For every orbit $\mathcal{Y}$ of $G^\vee$ on $\mathrm{Bl}(A)$, we
denote by $b(\mathcal{Y})$ the block of $(ZA)^{G^\vee}$ defined as
$$b(\mathcal{Y}):=\sum_{b \in \mathcal{Y}}b.$$
For every orbit $\bar{\mathcal{Y}}$ of $G$ on
$\mathrm{Bl}(\bar{A})$, we denote by $\bar{b}(\bar{\mathcal{Y}})$
the block of $(Z\bar{A})^G$ defined as
$$\bar{b}(\bar{\mathcal{Y}}):=\sum_{\bar{b} \in \bar{\mathcal{Y}}}\bar{b}.$$

The following proposition results from proposition $\ref{1.42}$ and
lemmas $\ref{1.43}$ and $\ref{1.44}$.

\begin{proposition}\label{1.45}
If the group $G$ is cyclic, there exists a bijection
$$\begin{array}{ccc}
    \mathrm{Bl}(\bar{A})/G & \tilde{\leftrightarrow} & \mathrm{Bl}(A)/G^\vee \\
    \bar{\mathcal{Y}} & \leftrightarrow & \mathcal{Y}
  \end{array}$$
such that
$$\bar{b}(\bar{\mathcal{Y}})=b(\mathcal{Y}),$$
i.e.,
$$\mathrm{Tr}(G,\bar{b})=\mathrm{Tr}(G^\vee,b) \textrm{ for all }
\bar{b} \in \bar{\mathcal{Y}} \textrm{ and } b \in \mathcal{Y}.$$ In
particular, the algebras $(Z\bar{A})^G$ and $(ZA)^{G^\vee}$ have the
same blocks.
\end{proposition}

\begin{corollary}\label{clifford}
If the blocks of $A$ are stable by the action of $G^\vee$, then the
blocks of $A$ coincide with the blocks of $(Z\bar{A})^G$.
\end{corollary}

\section{Representation theory of symmetric algebras}

For the last part of Chapter 2, except for the subsection ``A
variation for Tits' deformation theorem'', the author follows
\cite{GePf}, Chapter 7.

\subsection*{Grothendieck groups}

Let $\mathcal{O}$ be an integral domain and $K$ a field containing
$\mathcal{O}$. Let $A$ be an $\mathcal{O}$-algebra free and finitely
generated as $\mathcal{O}$-module.

Let $R_0(KA)$ be the Grothendieck group of finite-dimensional
$KA$-modules. Thus, $R_0(KA)$ is generated by expressions $[V]$, one
for each $KA$-module $V$ (up to isomorphism), with relations
$[V]=[V']+[V'']$ for each exact sequence $0 \rightarrow V'
\rightarrow V \rightarrow V'' \rightarrow 0$ of $KA$-modules. Two
$KA$-modules $V,V'$ give rise to the same element in $R_0(KA)$, if
$V$ and $V'$ have the same composition factors, counting
multiplicities. It follows that $R_0(KA)$ is free abelian with basis
given by the isomorphism classes of simple modules. Finally, let
$R_0^+(KA)$ be the subset of $R_0(KA)$ consisting of elements $[V]$,
where $V$ is a finite-dimensional $KA$-module.

\begin{definition}\label{pk}
Let $x$ be an indeterminate over $K$ and $\mathrm{Maps}(A,K[x])$ the
$K$-algebra of maps from $A$ to $K[x]$ (with pointwise
multiplication of maps as algebra multiplication). If $V$ is a
$KA$-module, let $\rho_V:KA \rightarrow \mathrm{End}_K(V)$ denote
its structural morphism. We define the map
$$\begin{array}{cccl}
    \mathfrak{p}_K: & R_0^+(KA) & \rightarrow & \emph{Maps}(A,K[x]) \\
                    & [V] & \mapsto & (a \mapsto \textrm{\emph{characteristic
                    polynomial of} } \rho_V(a)).
  \end{array}$$
Considering $\emph{Maps}(A,K[x])$ as a semigroup with respect to
multiplication, the map $\mathfrak{p}_K$ is a semigroup
homomorphism.
\end{definition}

Let $\mathrm{Irr}(KA)$ be the set of all characters $\chi_V$, where
$V$ is a simple $KA$-module.

\begin{lemma}\label{Brauer-Nesbitt}\emph{(Brauer-Nesbitt)} Assume that
$\mathrm{Irr}(KA)$ is a linearly independent set of
$\mathrm{Hom}_K(KA,K)$. Then the map $\mathfrak{p}_K$ is injective.
\end{lemma}
\begin{apod}{Let $V,V'$ be two $KA$-modules such that
$\mathfrak{p}_K([V])=\mathfrak{p}_K([V'])$. Since $[V]$, $[V']$ only
depend on the composition factors of $V,V'$, we may assume that
$V,V'$ are semisimple modules. Let
$$V=\bigoplus_{i=1}^n a_iV_i \textrm{ and } V'=\bigoplus_{i=1}^n b_iV_i,$$
where the $V_i$ are pairwise non-isomorphic simple $KA$-modules and
$a_i,b_i \geq 0$ for all $i$. We have to show that $a_i=b_i$ for all
$i$.

If, for some $i$, we have both $a_i>0$ and $b_i>0$, then we can
write $V=V_i \oplus \tilde{V}$ and $V'=V_i \oplus \tilde{V}'$. Since
$\mathfrak{p}_K$ is a semigroup homomorphism, we obtain
$$\mathfrak{p}_K([V_i]) \cdot \mathfrak{p}_K([\tilde{V}]) = \mathfrak{p}_K([V])
= \mathfrak{p}_K([V'])= \mathfrak{p}_K([V_i]) \cdot
\mathfrak{p}_K([\tilde{V}']),$$ and, dividing by
$\mathfrak{p}_K([V_i])$, we conclude that
$\mathfrak{p}_K([\tilde{V}])=\mathfrak{p}_K([\tilde{V}'])$. Thus, we
can suppose that, for all $i$, we have $a_i=0$ or $b_i=0$. Taking
characters yields that
$$\chi_V = \sum_i a_i \chi_{V_i} \textrm{ and } \chi_{V'} = \sum_i b_i \chi_{V_i}.$$
For each $a \in A$, we have that the character values $\chi_V(a)$
and $\chi_{V'}(a)$ appear as coefficients in the polynomials
$\mathfrak{p}_K([V])(a)$ and $\mathfrak{p}_K([V'])(a)$ respectively.
Hence, we have that $\sum_i (a_i-b_i)\chi_{V_i}=0 $. By assumption,
the characters $\chi_{V_i}$ are linearly independent. So we must
have $(a_i-b_i) 1_K$ for all $i$. Since for all $i$, $a_i=0$ or
$b_i=0$, this means that $a_1 1_K=0$ and $b_i 1_K=0$ for all $i$. If
the field $K$ has characteristic 0, we conclude that $a_i=b_i=0$ for
all $i$ and we are done. If $K$ has characteristic $p>0$, we
conclude that $p$ divides all $a_i$ and all $b_i$ and so
$\frac{1}{p}[V]$ and $\frac{1}{p}[V']$ exist in $R_0^+(KA)$.
Consequently, we also have
$\mathfrak{p}_K(\frac{1}{p}[V])=\mathfrak{p}_K(\frac{1}{p}[V'])$.
Repeating the above argument for $\frac{1}{p}[V]$ and
$\frac{1}{p}[V']$ yields that the multiplicity of $V_i$ in each of
these modules is still divisible by $p$. If we repeat this again and
again, we deduce that $a_i$ and $b_i$ should be divisible by
arbitrary powers of $p$. This forces $a_i=b_i=0$ for all $i$, as
desired.}
\end{apod}\\
\begin{remark}\emph{The assumption of the Brauer-Nesbitt lemma is satisfied when
(but not only when):
\begin{itemize}
  \item $KA$ is split.
  \item $K$ is a perfect field.
\end{itemize}}
\end{remark}

The following lemma implies the compatibility of the map
$\mathfrak{p}_K$ with the field extensions of $K$ (\cite{GePf},
Lemma 7.3.4).

\begin{lemma}\label{compatibility with field extensions}
Let $K \subseteq K'$ be a field extension. Then there is a canonical
map $d_K^{K'}:R_0(KA) \rightarrow R_0(K'A)$ given by $[V] \mapsto [V
\otimes_K K']$. Furthermore, we have a commutative diagram
$$\diagram R_0^+(KA) \dto^{d_K^{K'}} \rto^{\mathfrak{p}_K} &\emph{Maps}(A,K[x]) \dto^{\tau_K^{K'}} \\
           R_0^+(K'A) \rto^{\mathfrak{p}_{K'}} &\emph{Maps}(A,K'[x])
\enddiagram$$
where $\tau_K^{K'}$ is the canonical embedding. If, moreover, $KA$
is split, then $d_K^{K'}$ is an isomorphism which preserves
isomorphism classes of simple modules.
\end{lemma}
\subsection*{Integrality}

We have seen in Chapter 1 that a subring $\mathcal{R} \subseteq K$
is a valuation ring if, for each non-zero element $x \in K$, we have
$x \in \mathcal{R}$ or $x^{-1} \in \mathcal{R}$. Consequently, $K$
is the field of fractions of $\mathcal{R}$.

Such a valuation ring is a local ring whose maximal ideal we will
denote by $J(\mathcal{R})$. Valuation rings have interesting
properties, some of which are:
\begin{description}
  \item[(V1)] If $I$ is a prime ideal of $\mathcal{O}$, then there
  exists a valuation ring $\mathcal{R} \subseteq K$ such that
  $\mathcal{O} \subseteq \mathcal{R}$ and $J(\mathcal{R}) \cap
  \mathcal{O}=I$.
  \item[(V2)] Every finitely generated torsion-free module over a
  valuation ring in $K$ is free.
  \item[(V3)] The intersection of all valuation rings $\mathcal{R}
  \subseteq K$ with $\mathcal{O} \subseteq \mathcal{R}$ is the
  integral closure of $\mathcal{O}$ in $K$; each valuation ring
  itself is integrally closed in $K$ (Proposition $\ref{intersection of valuation rings}$).
\end{description}

\begin{lemma}\label{realizing modules over O}
Let $V$ be a $KA$-module. Choosing a $K$-basis of $V$, we obtain a
corresponding matrix representation $\rho:KA \rightarrow
\emph{M}_n(K)$, where $n=\emph{dim}_K(V)$. If $\mathcal{R} \subseteq
K$ is a valuation ring with $\mathcal{O} \subseteq \mathcal{R}$,
then a basis of $V$ can be chosen so that $\rho(a) \in
\emph{M}_n(\mathcal{R})$ for all $a \in A$. In that case, we say
that $V$ is realized over $\mathcal{R}$.
\end{lemma}
\begin{apod}{Let $(v_1,\ldots,v_n)$ be a $K$-basis of $V$ and $\mathcal{B}$
an $\mathcal{O}$-basis for $A$. Let $\tilde{V}$ be the
$\mathcal{O}$-submodule of $V$ spanned by the finite set $\{v_ib \,
|\, 1 \leq i \leq n, b \in \mathcal{B}\}$. Then $\tilde{V}$ is
invariant under the action of $\mathcal{R}A$ and hence a finitely
generated $\mathcal{R}A$-module. Since it is contained in a
$K$-vector space, it is also torsion-free. So (V2) implies that
$\tilde{V}$ is an $\mathcal{R}A$-lattice (a finitely generated
$\mathcal{R}A$-module which is free as $\mathcal{R}$-module) such
that $\tilde{V} \otimes_\mathcal{R} K \simeq V$. Thus any
$\mathcal{R}$-basis of $\tilde{V}$ is also a $K$-basis of $V$ with
the required property.}
\end{apod}\\
\begin{remark}
\emph{Note that the above argument only requires that $\mathcal{R}$
is a subring of $K$ such that $K$ is the field of fractions of
$\mathcal{R}$ and $\mathcal{R}$ satisfies (V2). These conditions
also hold, for example, when $\mathcal{R}$ is a principal ideal
domain with $K$ as field of fractions.}
\end{remark}\\

The following two important results derive from the above lemma.

\begin{proposition}\label{integrality of pk}
Let $V$ be a $KA$-module and $\mathcal{O}_K$ be the integral closure
of $\mathcal{O}$ in $K$. Then we have $\mathfrak{p}_K([V])(a) \in
\mathcal{O}_K[x]$ for all $a \in A$. Thus the map $\mathfrak{p}_K$
of definition $\ref{pk}$ is in fact a map $ R_0^+(KA)  \rightarrow
\emph{Maps}(A,\mathcal{O}_K[x])$.
\end{proposition}
\begin{apod}{Fix $a \in A$. Let $\mathcal{R} \subseteq K$ be a
valuation ring with $\mathcal{O} \subseteq \mathcal{R}$. By lemma
$\ref{realizing modules over O}$, there exists a basis of $V$ such
that the action of $a$ on $V$ with respect to that basis is given by
a matrix with coefficients in $\mathcal{R}$. Therefore, we have that
$\mathfrak{p}_K([V])(a) \in \mathcal{R}[x]$. Since this holds for
all valuation rings $\mathcal{R}$ in $K$ containing $\mathcal{O}$,
property (V3) implies that $\mathfrak{p}_K([V])(a) \in
\mathcal{O}_K[x]$.}
\end{apod}

Note that, in particular, proposition $\ref{integrality of pk}$
implies that $\chi_V(a)\in \mathcal{O}_K$ for all $a \in A$, where
$\chi_V$ is the character of the representation $\rho_V$.

\begin{proposition}\label{integrality of the Schur elements}
\emph{(Integrality of the Schur elements)} Assume that we have a
symmetrizing form $t$ on $A$. Let $V$ be a split simple $KA$-module
(i.e., $\mathrm{End}_{KA}(V) \simeq K$) and let $s_V$ be its Schur
element with respect to the induced form $t^K$ on $KA$. Then $s_V
\in \mathcal{O}_K$.
\end{proposition}
\begin{apod}{Let $\mathcal{R} \subseteq K$ be a
valuation ring with $\mathcal{O} \subseteq \mathcal{R}$. By lemma
$\ref{realizing modules over O}$, we can assume that $V$ affords a
representation $\rho:KA \rightarrow \mathrm{M}_n(K)$ such that
$\rho(a) \in \mathrm{M}_n(\mathcal{R})$ for all $a \in A$. Let
$\mathcal{B}$ be an $\mathcal{O}$-basis of $A$ and let
$\mathcal{B}'$ be its dual with respect to $t$. Then $s_V=\sum_{b
\in \mathcal{B}}\rho(b)_{ij}\rho(b')_{ji}$ for all $1 \leq i,j \leq
n$ (\cite{GePf}, Cor. 7.2.2). All terms in the sum lie in
$\mathcal{R}$ and so $s_V \in \mathcal{R}$. Since this holds for all
valuation rings $\mathcal{R}$ in $K$ containing $\mathcal{O}$,
property (V3) implies that $s_V \in \mathcal{O}_K$.}
\end{apod}
\subsection*{The decomposition map}

Now, we moreover assume that the ring $\mathcal{O}$ is integrally
closed in $K$. Throughout we will fix a ring homomorphism $\theta:
\mathcal{O} \rightarrow L$ into a field $L$ such that $L$ is the
field of fractions of $\theta(\mathcal{O})$. We call such a ring
homomorphism a \emph{specialization} of $\mathcal{O}$.

Let $\mathcal{R} \subseteq K$ be a valuation ring with $\mathcal{O}
\subseteq \mathcal{R}$ and $J(\mathcal{R}) \cap \mathcal{O} =
\mathrm{Ker}\theta$ (note that $\mathrm{Ker}\theta$ is a prime
ideal, since $\theta(\mathcal{O})$ is contained in a field). Let $k$
be the residue field of $\mathcal{R}$. Then the restriction of the
canonical map $\pi:\mathcal{R} \rightarrow k$ to $\mathcal{O}$ has
kernel $J(\mathcal{R}) \cap \mathcal{O} = \mathrm{Ker}\theta$. Since
$L$ is the field of fractions of $\theta(\mathcal{O})$, we may
regard $L$ as a subfield of $k$. Thus, we have a commutative diagram
$$\diagram\mathcal{O}\dto^{\theta}&\subseteq&
          \mathcal{R}\dto^{\pi}&\subseteq &K\\
           L & \subseteq & k & & \enddiagram $$

From now on, we make the following assumption:

\begin{ypothesh}\label{split assumption}
(a) $LA \textrm{ is split}$ \,\,\,\,or\,\,\,\,  (b) $L=k \textrm{
and } k \textrm{ is perfect}$.
\end{ypothesh}

The map $\theta:\mathcal{O} \rightarrow L$ induces a map $A
\rightarrow LA, a \mapsto a \otimes 1$. One consequence of the
assumption $\ref{split assumption}$ is that, due to lemma
$\ref{compatibility with field extensions}$, the map $d_L^k:R_0(LA)
\rightarrow R_0(kA)$ is an isomorphism which preserves isomorphism
classes of simple modules. Thus we can identify $R_0(LA)$ and
$R_0(kA)$. Moreover, the Brauer-Nesbitt lemma holds for $LA$, \ie
the map $\mathfrak{p}_L:R_0^+(LA) \rightarrow \mathrm{Maps}(A,L[x])$
is injective.

Let $V$ be a $KA$-module and $\mathcal{R} \subseteq K$ be a
valuation ring with $\mathcal{O} \subseteq \mathcal{R}$. By lemma
$\ref{realizing modules over O}$, there exists a $K$-basis of $V$
such that the corresponding matrix representation $\rho:KA
\rightarrow \mathrm{M}_n(K)$ ($n=\mathrm{dim}_K(V)$) has the
property that $\rho(a) \in \mathrm{M}_n(\mathcal{R})$. Then that
basis generates an $\mathcal{R}A$-lattice $\tilde{V}$ such that
$\tilde{V} \otimes_\mathcal{R} K =V$. The $k$-vector space
$\tilde{V} \otimes_\mathcal{R} k$ is a $kA$-module via
$(v\otimes1)(a\otimes1) = va\otimes1 (v \in \tilde{V},a \in A)$,
which we call the \emph{modular reduction} of $\tilde{V}$.

The matrix representation $\rho^k:kA \rightarrow \mathrm{M}_n(k)$
afforded by $k\tilde{V}$ is given by
$$\rho^k(a\otimes1)=(\pi(a_{ij})) \textrm{ where } a \in A \textrm{ and } \rho(a)=(a_{ij}).$$
To simplify notation, we shall write
$$K\tilde{V}:=\tilde{V} \otimes_\mathcal{R} K \textrm{ and }
k\tilde{V}:=\tilde{V} \otimes_\mathcal{R} k  .$$ Note that if
$\tilde{V}'$ is another $\mathcal{R}A$-lattice such that $\tilde{V}'
\otimes_\mathcal{R} K \simeq V$, then $\tilde{V}$ and $\tilde{V}'$
need not be isomorphic. The same hold for the $kA$-modules
$\tilde{V} \otimes_\mathcal{R} k $ and
$\tilde{V}'\otimes_\mathcal{R} k$.

\begin{thedef}\label{existence of decomposition maps}
Let $\theta:\mathcal{O} \rightarrow L$ be a ring homomorphism into a
field $L$ such that $L$ is the field of fractions of
$\theta(\mathcal{O})$ and $\mathcal{O}$ is integrally closed in $K$.
Assume that we have chosen a valuation ring $\mathcal{R}$ with
$\mathcal{O} \subseteq \mathcal{R} \subseteq K$ and $J(\mathcal{R})
\cap \mathcal{O} = \mathrm{Ker}\theta$ and that the assumption
$\ref{split assumption}$ is satisfied. Then
\begin{description}
  \item[(a) ] The modular reduction induces an additive map $d_\theta:R_0^+(KA) \rightarrow R_0^+(LA)$
  such that $d_\theta([K\tilde{V}])=[k\tilde{V}]$, where $\tilde{V}$
  is an $\mathcal{R}A$-lattice and $[k\tilde{V}]$ is regarded as an element of
  $R_0^+(LA)$ via the identification of $R_0(kA)$ and $R_0(LA)$.
  \item[(b) ] By Proposition $\ref{integrality of pk}$, the image of $\mathfrak{p}_K$ is
  contained in $\emph{Maps}(A,\mathcal{O}[x])$ and we have the following commutative
  diagram
  $$\diagram
  R_0^+(KA) \dto^{d_\theta}  \rto^{\mathfrak{p}_K} &\emph{Maps}(A,\mathcal{O}[x]) \dto^{\tau_\theta}\\
  R_0^+(LA) \rto^{\mathfrak{p}_L} & \mathrm{Maps}(A,L[x]) \enddiagram $$
  where $\tau_\theta:\emph{Maps}(A,\mathcal{O}[x]) \rightarrow \emph{Maps}(A,L[x])$
  is the map induced by $\theta$.
  \item[(c) ] The map $d_\theta$ is uniquely determined by the commutativity of
  the above diagram. In particular, the map $d_\theta$ depends only on
  $\theta$ and not on the choice of $\mathcal{R}$.
\end{description}
The map $d_\theta$ will be called the decomposition map associated
with the specialization $\theta:\mathcal{O} \rightarrow L$. The
matrix of that map with respect to the bases of $R_0(KA)$ and
$R_0(LA)$ consisting of the classes of the simple modules is called
the decomposition matrix with respect to $\theta$.
\end{thedef}
\begin{apod}{Let $\tilde{V}$ be an $\mathcal{R}A$-lattice and $a \in A$.
Let $(m_{ij}) \in \mathrm{M}_n(\mathcal{R})$
 be the matrix describing the action of $a$
on $\tilde{V}$ with respect to a chosen $\mathcal{R}$-basis of
$\tilde{V}$. Due to the properties of modular reduction, the action
of $a \otimes 1 \in kA$ on $k\tilde{V}$ is given by the matrix
$(\pi(m_{ij}))$. Then, by definition,
$\mathfrak{p}_L([k\tilde{V}])(a)$ is the characteristic polynomial
of $(\pi(m_{ij}))$. On the other hand, applying $\theta$ (which is
the restriction of $\pi$ to $\mathcal{O}$) to the coefficients of
the characteristic polynomial of $(m_{ij})$ returns $(\tau_\theta
\circ \mathfrak{p}_K)([K\tilde{V}])(a)$. Since the two actions
just described commute, the two polynomials obtained are equal. Thus
the following relation is established:
$$\mathfrak{p}_L([k\tilde{V}])=\tau_\theta \circ \mathfrak{p}_K([K\tilde{V}]) \textrm{ for all }
\mathcal{R}A\textrm{-lattices } \tilde{V} \textrm{    } (\dag) $$

Now let us prove (a). We have to show that the map $d_\theta$ is
well defined \ie if $\tilde{V},\tilde{V}'$ are two
$\mathcal{R}A$-lattices such that $K\tilde{V}$ and $K\tilde{V}'$
have the same composition factors (counting multiplicities), then
the classes of $k\tilde{V}$ and $k\tilde{V}'$ in $R_0(LA)$ are the
same. Moreover, for all $a \in A$, the endomorphisms
$\rho_{K\tilde{V}}(a)$ and $\rho_{K\tilde{V}'}(a)$ are conjugate. So
the equality $(\dag)$ implies that
$$\mathfrak{p}_L([k\tilde{V}])(a)=\mathfrak{p}_L([k\tilde{V}'])(a) \textrm{ for all } a\in A.$$
We have already remarked that, since the assumption $\ref{split
assumption}$ is satisfied, the Brauer-Nesbitt lemma holds for $LA$.
So we conclude $[k\tilde{V}]=[k\tilde{V}']$, as desired.

Having established the existence of $d_\theta$, we have
$[k\tilde{V}]=d_\theta([K\tilde{V}])$ for any $\mathcal{R}A$-lattice
$\tilde{V}$. Hence $(\dag)$ yields the commutativity of the diagram
in (b).

Finally, by the Brauer-Nesbitt lemma, the map $\mathfrak{p}_L$ is
injective. Hence there exists at most one map which makes the
diagram in (b) commutative. This proves (c).}
\end{apod}\
\\
\begin{remark} \emph{ Note that if $\mathcal{O}$ is a discrete valuation ring and $L$ its
residue field, we do not need the assumption $\ref{split
assumption}$ in order to define a decomposition map from $R_0^+(KA)$
to $R_0^+(LA)$ associated with the canonical map $\theta:\mathcal{O}
\rightarrow L$. For a given $KA$-module $V$, there exists an
$A$-lattice $\tilde{V}$ such that $V=\tilde{V} \otimes_\mathcal{O}
K$. The map $d_\theta:R_0^+(KA) \rightarrow
R_0^+(LA),\,[K\tilde{V}]\mapsto[\tilde{V}/L\tilde{V}]$ is well and
uniquely defined.
For the details of this construction, see \cite{CuRe}, \S16C}.
\end{remark}

Recall from proposition $\ref{integrality of pk}$ that if $V$ is a
$KA$-module, then its character $\chi_V$ restricts to a trace
function $\dot{\chi}_V: A \rightarrow \mathcal{O}$. Now, any linear
map $\lambda:A \rightarrow \mathcal{O}$ induces an $L$-linear map
$$\lambda^L:LA \rightarrow L, a\otimes1\mapsto \theta(\lambda(a))
(a\in A).$$ It is clear that if $\lambda$ is a trace function, so is
$\lambda^L$. Applying this to $\dot{\chi}_V$ shows that
$\dot{\chi}_V^L$ is a trace function on $LA$. Since character values
occur as coefficients in characteristic polynomials, theorem
$\ref{existence of decomposition maps}$ implies that
$\dot{\chi}_V^L$ is the character of $d_\theta([V])$. Moreover, for
any simple $KA$-module $V$, we have
$$\dot{\chi}_V^L = \sum_{V'} d_{VV'} \chi_{V'},$$
where the sum is over all simple $LA$-modules $V'$ (up to
isomorphism) and $D=(d_{VV'})$ is the decomposition matrix
associated with $\theta$.

The following result gives a criterion for $d_\theta$ to be trivial.
For its proof, the reader may refer, for example, to \cite{GePf},
Thm. 7.4.6.

\begin{theorem}\label{Tits}\emph{(Tits' deformation theorem)}
Assume that $KA$ and $LA$ are split. If $LA$ is semisimple, then
$KA$ is also semisimple and the decomposition map $d_\theta$ is an
isomorphism which preserves isomorphism classes of simple modules.
In particular, the map $\emph{Irr}(KA) \rightarrow \emph{Irr}(LA),
\chi \mapsto \dot{\chi}^L$ is a bijection.
\end{theorem}

Finally, if $A$ is symmetric, we can check whether the assumption of
Tits' deformation theorem is satisfied, using the following theorem
(cf. \cite{GePf}, Thm. 7.4.7).

\begin{theorem}\label{semisimplicity}\emph{(Semisimplicity
criterion)} Assume that $KA$ and $LA$ are split and that $A$ is
symmetric with symmetrizing form $t$. For any simple $KA$-module
$V$, let $s_V \in \mathcal{O}$ be the Schur element with respect to
the induced symmetrizing form $t^K$ on $KA$. Then $LA$ is semisimple
if and only if $\theta(s_V) \neq 0$ for all $V$.
\end{theorem}

\begin{corollary}\label{injective preserves splitness}
Let $K$ be the field of fractions of $\mathcal{O}$.
Assume that $KA$ is split semisimple and that $A$ is
symmetric with symmetrizing form $t$. If the map $\theta$ is
injective, then $LA$ is split semisimple.
\end{corollary}

\subsection*{A variation for Tits' deformation
theorem}

Let us suppose that $\mathcal{O}$ is a Krull ring and
$\theta:\mathcal{O} \rightarrow L$ is a ring homomorphism like
above. We will give a criterion for $LA$ to be split semisimple.

\begin{theorem}\label{Lehrer}
Let $K$ be the field of fractions of $\mathcal{O}$. Assume that $KA$ is split semisimple and that $A$ is symmetric with
symmetrizing form $t$. For any simple $KA$-module $V$, let $s_V \in
\mathcal{O}$ be the Schur element with respect to the induced
symmetrizing form $t^K$ on $KA$. If $\mathrm{Ker}\theta$ is a prime
ideal of $\mathcal{O}$ of height $1$, then $LA$ is split semisimple if
and only if $\theta(s_V) \neq 0$ for all $V$.
\end{theorem}
\begin{apod}{If $LA$ is split semisimple, then theorem
$\ref{semisimplicity}$ implies that $\theta(s_V) \neq 0$ for all
$V$. Now let us denote by $\mathrm{Irr}(KA)$ the set of irreducible
characters of $KA$. If $\chi$ is the character afforded by a simple
$KA$-module $V_\chi$, then $s_\chi:=s_{V_\chi}$. We set
$\mathfrak{q}:=\mathrm{Ker}\theta$ and suppose that $s_\chi \notin
\mathfrak{q}$ for all $\chi \in \mathrm{Irr}(KA)$. Since $KA$ is
split semisimple, it is isomorphic to a product of matrix algebras
over $K$:
$$KA \simeq \prod_{\chi \in \mathrm{Irr}(KA)}\mathrm{End}_{K}(V_\chi)$$
Let us denote by $\pi_\chi:KA \twoheadrightarrow
\mathrm{End}_{K}(V_\chi)$ the projection onto the $\chi$-factor,
such that $\pi:=\prod_{\chi \in \mathrm{Irr}(KA)}\pi_\chi$ is the
above isomorphism. Then $\chi=\mathrm{tr}_{V_\chi} \circ \pi_\chi$,
where $\mathrm{tr}_{V_\chi}$ denotes the standard trace on
$\mathrm{End}_{K}(V_\chi)$.

Let $\mathcal{B},\mathcal{B}'$ be two dual bases of $A$ with respect
to the symmetrizing form $t$. By lemma $\ref{tau^vee}$, for all $a
\in KA$ and $\chi \in \mathrm{Irr}(KA)$, we have
$$\chi^\vee a= \sum_{b \in \mathcal{B}} \chi(b'a)b.$$
Applying $\pi$ to both sides yields
$$\pi(\chi^\vee) \pi(a)= \sum_{b \in \mathcal{B}} \chi(b'a)\pi(b).$$
By definition of the Schur element,
$\pi(\chi^\vee)=\pi_\chi(\chi^\vee)=\omega_\chi(\chi^\vee)=s_\chi$.
Thus, if $\alpha \in \mathrm{End}_{K}(V_\chi)$, then
$$\pi^{-1}(\alpha)=\frac{1}{s_\chi}\sum_{b \in \mathcal{B}}\mathrm{tr}_{V_\chi}(\pi_\chi(b')\alpha)b. \,\,\,\,\,\,\,\,\,\,(\dag)$$

Since $\mathcal{O}$ is a Krull ring and $\mathfrak{q}$ is a prime
ideal of height $1$ of $\mathcal{O}$, the ring
$\mathcal{O}_{\mathfrak{q}}$ is, by theorem $\ref{Krull-dvr}$, a
discrete valuation ring. Thanks to lemma $\ref{realizing modules
over O}$, there exists an $\mathcal{O}_{\mathfrak{q}}A$-lattice
$\tilde{V}_\chi$ such that $K \otimes_{\mathcal{O}_{\mathfrak{q}}}
\tilde{V}_\chi \simeq V_\chi$.

Moreover, $1/s_\chi \in \mathcal{O}_\mathfrak{q}$ for all ${\chi \in
\mathrm{Irr}(KA)}$. Due to the relation $(\dag)$, the map $\pi$
induces an isomorphism
$$\mathcal{O}_{\mathfrak{q}}A \simeq
\prod_{\chi \in
\mathrm{Irr}(KA)}\mathrm{End}_{\mathcal{O}_{\mathfrak{q}}}(\tilde{V}_\chi),$$
\ie $\mathcal{O}_{\mathfrak{q}}A$ is the product of matrix algebras
over $\mathcal{O}_{\mathfrak{q}}$. Since
$\mathrm{Ker}\theta=\mathfrak{q}$, the above isomorphism remains
after applying $\theta$. Therefore, we obtain that $LA$ is a product
of matrix algebras over $L$ and thus split semisimple.}
\end{apod}

If that is the case, then the assumption of Tits' deformation
theorem is satisfied and there exists a bijection $\mathrm{Irr}(KA)
\leftrightarrow \mathrm{Irr}(LA)$.

\subsection*{Symmetric algebras over discrete valuation rings}

From now on, we assume that the following conditions are satisfied:
\begin{itemize}
  \item $\mathcal{O}$ is a discrete valuation ring in $K$ and $K$ is
  perfect; let $v:K\rightarrow \mathbb{Z} \cup \{\infty\}$ be the
  corresponding valuation.
  \item $KA$ is split semisimple.
  \item $\theta:\mathcal{O} \rightarrow L$ is the canonical map onto
  the residue field $L$ of $\mathcal{O}$.
  \item $A$ is a symmetric algebra with symmetrizing form $t$.
\end{itemize}

We have already seen that we have a well-defined decomposition map
$d_\theta:R_0^+(KA) \rightarrow R_0^+(LA)$. The decomposition matrix
associated with $d_\theta$ is the $|\mathrm{Irr}(KA)| \times
|\mathrm{Irr}(LA)|$ matrix $D=(d_{\chi\phi})$ with non-negative
integer entries such that
$$d_\theta([V_\chi])=\sum_{\phi \in
\mathrm{Irr}(LA)}d_{\chi\phi}[V_\phi'] \,\,\textrm{ for } \chi \in
\mathrm{Irr}(KA),$$ where $V_\chi$ is a simple $KA$-module with
character $\chi$ and $V_\phi'$ is a simple $LA$-module with
character $\phi$. We sometimes call the characters of $KA$
``ordinary'' and the characters of $LA$ ``modular''. We say that
$\phi \in \mathrm{Irr}(LA)$ is a \emph{modular constituent} of $\chi
\in \mathrm{Irr}(KA)$, if $d_{\chi\phi} \neq 0$.

The rows of $D$ describe the decomposition of $d_\theta([V_\chi])$
in the standard basis of $R_0(LA)$. An interpretation of the columns
is given by the following result (cf. \cite{GePf}, Thm. 7.5.2),
which is part of Brauer's classical theory of modular
representations.

\begin{theorem}\label{Brauer reciprocity}\emph{(Brauer reciprocity)}
For each $\phi \in \mathrm{Irr}(LA)$, there exists some primitive
idempotent $e_\phi \in A$ such that
$$[e_\phi KA]=\sum_{\chi \in \mathrm{Irr}(KA)}d_{\chi\phi}[V_\chi]
\in R_0^+(KA).$$
\end{theorem}
Let $\phi \in \mathrm{Irr}(LA)$. Consider the map $\psi(\phi):ZKA
\rightarrow K$ defined by
$$\psi(\phi):=\sum_{\chi \in \mathrm{Irr}(KA)}
\frac{d_{\chi\phi}}{s_\chi}\omega_\chi,$$ where $\omega_\chi: ZKA
\twoheadrightarrow K$ is the central morphism associated with $\chi
\in \mathrm{Irr}(KA)$, as defined at the end of section 2.1.

\begin{theorem}\label{Geck-Rouquier}
The map $\psi(\phi)$ restricts to a map $ZA \rightarrow
\mathcal{O}$. In particular,
$$\psi(\phi)(1)=\sum_{\chi \in \mathrm{Irr}(KA)}
\frac{d_{\chi\phi}}{s_\chi} \in \mathcal{O}.$$
\end{theorem}
\begin{apod}{Let us denote by $t^K$ the induced symmetrizing
form on $KA$. If $e_\phi$ is an idempotent as in theorem
$\ref{Brauer reciprocity}$, then we can define a $K$-linear map
$\lambda_\phi:ZKA \rightarrow K,\,z \mapsto t^K(ze_\phi)$. We claim
that $\lambda_\phi=\psi(\phi)$. Since $KA$ is split semisimple, the
elements $\{\chi^\vee\,|\, \chi \in \mathrm{Irr}(KA)\}$ form a basis
of $ZKA$. It is, therefore, sufficient to show that
$$\lambda_\phi(\chi^\vee)=\psi(\phi)(\chi^\vee) \textrm{ for all } \chi
\in \mathrm{Irr}(KA).$$ We have
$\psi(\phi)(\chi^\vee)=d_{\chi\phi}1_K$. Now consider the left-hand
side.
$$\begin{array}{rl}
    \lambda_\phi(\chi^\vee) & =t^K(\chi^\vee e_\phi)=\chi(e_\phi) =\mathrm{dim}_K(V_\chi e_\phi)1_K \\
     & \\
     & =\mathrm{dim}_K(\mathrm{Hom}_K(e_\phi
KA,V_\chi)1_K=d_{\chi\phi}1_K.
  \end{array}$$
Hence the above claim is established.

Finally, it remains to observe that since $e_\phi \in A$, the
function $\lambda_\phi$ takes values in $\mathcal{O}$ on all
elements of $A$.}
\end{apod}

Finally, we will treat the block distribution of characters. For
this purpose, we introduce the following notions.

\begin{definition}\label{brauer graph}\
\begin{enumerate}
  \item The Brauer graph associated with $A$ has vertices labeled by
  the irreducible characters of $KA$ and an edge joining
  $\chi,\chi' \in \emph{Irr}(KA)$ if $\chi \neq \chi'$ and there
  exists some $\phi \in \emph{Irr}(LA)$ such that $d_{\chi\phi} \neq 0 \neq
  d_{\chi'\phi}$. A connected component of a Brauer graph is called
  a block.
  \item Let $\chi \in \mathrm{Irr}(KA)$. Recall that $0 \neq s_\chi
  \in \mathcal{O}$. Let $\delta_\chi:=v(s_\chi)$, where $v$
  is the given valuation. Then $\delta_\chi$ is called the defect of
  $\chi$ and we have $\delta_\chi \geq 0$ for all $\chi \in \mathrm{Irr}(KA)$. If $B$ is a block, then
  $\delta_B:=\mathrm{max}\{\delta_\chi \,|\, \chi \in B\}$ is called the
  defect of $B$.
\end{enumerate}
\end{definition}

By \cite{Fe}, 17.9, each block $B$ of $A$ corresponds to a central
primitive idempotent (block-idempotent, by definition
$\ref{blocks}$) $e_B$ of $A$. If $\chi \in B$ and $e_\chi$ is its
corresponding central primitive idempotent in $KA$, then $e_Be_\chi
\neq 0$.

Every $\chi \in \mathrm{Irr}(KA)$ determines a central morphism
$\omega_\chi:ZKA \rightarrow K$. Since $\mathcal{O}$ is integrally
closed, we have $\omega_\chi(z) \in \mathcal{O}$ for all $z \in ZA$.
We have the following standard results relating blocks with central
morphisms.

\begin{proposition}\label{blocks and central characters}
Let $\chi,\chi' \in \emph{Irr}(KA)$. Then $\chi$ and $\chi'$ belong
to the same block of $A$ if and only if
$$\theta(\omega_\chi(z))=\theta(\omega_{\chi'}(z)) \textrm{ for all
} z \in ZA.$$ i.e.,
$$\omega_\chi(z) \equiv \omega_{\chi'}(z) \,\,\emph{mod}\,J(\mathcal{O})  \textrm{ for all
} z \in ZA.$$
\end{proposition}
\begin{apod}{First assume that $\chi,\chi'$ belong to the same
block of $A$, \ie they belong to a connected component of the Brauer
graph. It is sufficient to consider the case where $\chi,\chi'$ are
directly linked on the Brauer graph, \ie there exists some $\phi \in
\mathrm{Irr}(LA)$ such that $d_{\chi\phi} \neq 0 \neq
d_{\chi'\phi}$. Let $\tilde{V}_\chi$ be an $A$-lattice such that
$K\tilde{V}_\chi$ affords $\chi$. Let $z \in ZA$. Then $z \otimes 1$
acts by the scalar $\theta(\omega_\chi(z))$ on every modular
constituent of $k\tilde{V}_\chi$. Similarly, $z \otimes 1$ acts by
the scalar $\theta(\omega_{\chi'}(z))$ on every modular constituent
of $k\tilde{V}_{\chi'}$, where $\tilde{V}_{\chi'}$ is an $A$-lattice
such that $K\tilde{V}_{\chi'}$ affords $\chi'$. Since, by
assumption, $K\tilde{V}_\chi$ and $K\tilde{V}_{\chi'}$ have a
modular constituent in common, we have
$\theta(\omega_\chi(z))=\theta(\omega_{\chi'}(z))$, as desired.

Now assume that $\chi$ belongs to the block $B$ and $\chi'$ to the
block $B'$, with $B \neq B'$. Let $e_B$, $e_{B'}$ be the
corresponding central primitive idempotents. Then
$\omega_\chi(e_B)=1$ and $\omega_{\chi'}(e_B)=0$. Consequently,
$\theta(\omega_\chi(e_B)) \neq \theta(\omega_{\chi'}(e_B))$.}
\end{apod}

\begin{theorem}\label{blocks of defect 0}\emph{(Blocks of defect 0)}
Let $\chi \in \emph{Irr}(KA)$ with $\theta(s_\chi) \neq 0$. Then
$\chi$ is an isolated vertex in the Brauer graph and the
corresponding decomposition matrix is just \emph{(1)}.
\end{theorem}
\begin{apod}{Let $t^K$ be the induced symmetrizing
form on $KA$ and $\hat{t}^K$ the isomorphism from $KA$ to
$\mathrm{Hom}_K(KA,K)$ induced by $t^K$. The irreducible character
$\chi \in \mathrm{Irr}(KA)$ is a trace function on $KA$ and thus we
can define $\chi^\vee:=(\hat{t}^K)^{-1}(\chi) \in KA$. Since $\chi$
restricts to a trace function $A \rightarrow \mathcal{O}$, we have
in fact $\chi^\vee \in ZA$. By definition, we have that
$\omega_\chi(\chi^\vee)=s_\chi$ and $\omega_{\chi'}(\chi^\vee)=0$
for any $\chi' \in \mathrm{Irr}(KA), \chi' \neq \chi$. Now assume
that there exists some character $\chi'$ which is linked to $\chi$
in the Brauer graph. Proposition $\ref{blocks and central
characters}$ implies that $0 \neq
\theta(s_\chi)=\theta(\omega_\chi(\chi^\vee))=\theta(\omega_{\chi'}(\chi^\vee))=0$,
which is absurd.

It remains to show that $d_\theta([V_\chi])$ is the class of a
simple module in $R_0^+(LA)$.  By lemma $\ref{realizing modules over
O}$, there exists a basis of $V_\chi$ and a corresponding
representation $\rho:KA \rightarrow \mathrm{M}_n(K)$ afforded by
$V_\chi$ such that $\rho(a) \in \mathrm{M}_n(\mathcal{O})$ for all
$a \in A$. Let $\mathcal{B}$ be an $\mathcal{O}$-basis of $A$ and
let $\mathcal{B}'$ be its dual with respect to $t$. We have seen in
proposition $\ref{integrality of the Schur elements}$ that
$s_\chi=\sum_{b \in \mathcal{B}}\rho(b)_{ij}\rho(b')_{ji}$ for all
$1 \leq i,j \leq n$. All terms in this relation lie in
$\mathcal{O}$. So we can apply the map $\theta$ and obtain a similar
relation for $\theta(s_\chi)$ with respect to the module
$L\tilde{V}_\chi$, where $\tilde{V}_\chi \subseteq V_\chi$ is the
$A$-lattice spanned by the above basis of $V_\chi$. Since
$\theta(s_\chi) \neq 0$, the module $L\tilde{V}_\chi$ is simple
(\cite{GePf}, Lemma 7.2.3).}
\end{apod}

The following result is an immediate consequence of theorems
$\ref{Geck-Rouquier}$ and $\ref{blocks of defect 0}$.

\begin{proposition}\label{Malle-Rouquier}
Let $\chi \in \mathrm{Irr}(KA)$. Then $\chi$ is a block by itself if
and only if $\theta(s_\chi) \neq 0.$
\end{proposition}

\chapter{On Hecke algebras}

\section{Complex reflection groups and associated braid groups}

Let $\mu_\infty$ be the group of all the roots of unity in
$\mathbb{C}$ and $K$ a number field contained in
$\mathbb{Q}(\mu_\infty)$. We denote by $\mu(K)$ the group of all the
roots of unity of $K$. For every integer $d>1$, we set
$\zeta_d:=\mathrm{exp}(2\pi i/d)$ and denote by $\mu_d$ the group of
all the $d$-th roots of unity. Let $V$ be a $K$-vector space of
finite dimension $r$.

\begin{definition}\label{pseudo-reflection}
A pseudo-reflection of $\mathrm{GL}(V)$ is a non-trivial element $s$
of $\mathrm{GL}(V)$ which acts trivially on a hyperplane, called the
reflecting hyperplane of $s$.

If $W$ is a finite subgroup of $\mathrm{GL}(V)$ generated by
pseudo-reflections, then $(V,W)$ is called a $K$-reflection group of
rank $r$.
\end{definition}

We have the following classification of complex reflection groups,
also known as ``Shephard-Todd classification''. For more details
about the classification, one may refer to \cite{ShTo}.

\begin{theorem} Let $(V,W)$ be an irreducible complex
reflection group (i.e., $W$ acts irreducibly on $V$). Then one of
the following assertions is true:
\begin{itemize}
  \item There exist non-zero integers $d,e,r$ such
  that $(V,W) \simeq G(de,e,r)$, where $G(de,e,r)$ is the group of all monomial
  $r \times r$ matrices with entries in $\mu_{de}$ and product of all non-zero
  entries in $\mu_d$.
  \item $(V,W)$ is isomorphic to one of the 34 exceptional groups
  $G_n$ $(n=4,\ldots,37)$.
\end{itemize}
\end{theorem}

The following theorem has been proved (using a case by case
analysis) by Benard \cite{Ben} and Bessis \cite{Bes1} and
generalizes a well known result for Weyl groups.

\begin{theorem}\label{field of definition}
Let $(V,W)$ be a reflection group. Let $K$ be the field generated by
the traces on $V$ of all the elements of $W$. Then all irreducible
$KW$-representations are absolutely irreducible i.e., $K$ is a
splitting field for $W$. The field $K$ is called the field of
definition of the reflection group $W$.
\end{theorem}

\begin{itemize}
  \item If $K \subseteq \mathbb{R}$, then $W$ is a (finite) Coxeter
  group.
  \item If $K=\mathbb{Q}$, then $W$ is a Weyl group.
\end{itemize}

For the following definitions and results about braid groups we
follow \cite{BMR}.\\

Let $X$ be a topological space. Given a point $x_0 \in X$, we denote
by $\Pi_1(X,x_0)$ the fundamental group with base point $x_0$.

Let $V$ be a $K$-vector space as before. Let $W$ be a finite
subgroup of $\mathrm{GL}(V)$ generated by pseudo-reflections and
acting irreducibly on $V$. We denote by $\mathcal{A}$ the set of its
reflecting hyperplanes. We define the \emph{regular variety}
$V^{\textrm{reg}}:= \mathbb{C} \otimes V-\bigcup_{H \in
\mathcal{A}}\mathbb{C} \otimes H$. For $x_0 \in V^{\textrm{reg}}$,
we define $P:=\Pi_1(V^{\textrm{reg}},x_0)$ the \emph{pure braid
group} (at $x_0$) associated with $W$. If now $p:V^{\textrm{reg}}
\rightarrow V^{\textrm{reg}}/W$ denotes the canonical surjection, we
define $B:=\Pi_1(V^{\textrm{reg}}/W,p(x_0))$ the \emph{braid group}
(at $x_0$) associated with $W$.

The projection $p$ induces a surjective map $B\twoheadrightarrow W,
\sigma \mapsto \bar{\sigma}$ as follows: Let
$\tilde{\sigma}:[0,1]\rightarrow V^{\textrm{reg}}$ be a path in
$V^{\textrm{reg}}$ such that $\tilde{\sigma}(0)=x_0$, which lifts
$\sigma$. Then $\bar{\sigma}$ is defined by the equality
$\bar{\sigma}(x_0)=\tilde{\sigma}(1)$. Note that the map $\sigma
\mapsto \bar{\sigma}$ is an anti-morphism.

Denoting by $W^\mathrm{op}$ the group opposite to $W$, we have the
following short exact sequence
$$1\rightarrow P\rightarrow B\rightarrow
W^\mathrm{op}\rightarrow1,$$ where the map $B\rightarrow
W^\mathrm{op}$ is defined by $\sigma \mapsto \bar{\sigma}.$
\\

Now, for every hyperplane $H \in \mathcal{A}$, we set $e_H$ the
order of the group $W_H$, where $W_H$ is the subgroup of $W$ formed
by 1 and all the reflections fixing the hyperplane $H$. The group
$W_H$ is cyclic: if $s_H$ denotes an element of $W_H$ with
determinant $\zeta_H:=\zeta_{e_H}$, then $W_H=<s_H>$ and $s_H$ is
called a \emph{distinguished reflection} in $W$.

Let $L_H:=\mathrm{Im}(s-\mathrm{id}_V)$. Then, for all $x \in V$, we
have $x=\mathrm{pr}_H(x)+\mathrm{pr}_{L_H}(x)$ with
$\mathrm{pr}_H(x) \in H$ and $\mathrm{pr}_{L_H}(x) \in L_H$. Thus,
$s_H(x)=\mathrm{pr}_H(x)+\zeta_H\mathrm{pr}_{L_H}(x)$.

If $t \in \mathbb{R}$, we set $\zeta_H^t:=\mathrm{exp}(2\pi it/e_H)$
and we denote by $s_H^t$ the element of $\mathrm{GL}(V)$ (a
pseudo-reflection if $t\neq 0$) defined by
$$s_H^t(x):=\mathrm{pr}_H(x)+\zeta_H^t\mathrm{pr}_{L_H}(x).$$

For $x \in V$, we denote by $\sigma_{H,x}$ the path in $V$ from $x$
to $s_H(x)$ defined by
$$\sigma_{H,x}:[0,1] \rightarrow V,\,\, t \mapsto s_H^t(x).$$

Let $\gamma$ be a path in $V^\mathrm{reg}$ with initial point $x_0$
and terminal point $x_H$. Then $\gamma^{-1}$ is the path in
$V^\mathrm{reg}$ with initial point $x_H$ and terminal point $x_0$
such that
$$\gamma^{-1}(t)=\gamma(1-t) \textrm{ for all } t \in [0,1].$$
Thus, we can define the path $s_H(\gamma^{-1}):t \mapsto
s_H(\gamma^{-1}(t))$, which goes from $s_H(x_H)$ to $s_H(x_0)$ and
lies also in $V^\mathrm{reg}$, since for all $x \in V^\mathrm{reg}$,
$s_H(x) \in V^\mathrm{reg}$ (If $s_H(x) \notin V^\mathrm{reg}$, then
$s_H(x)$ must belong to a hyperplane $H'$. If $s_{H'}$ is a
distinguished pseudo-reflection with reflecting hyperplane $H'$,
then $s_{H'}(s_H(x))=s_H(x)$ and ${s_H}^{-1}(s_{H'}(s_H(x)))=x$.
However, ${s_H}^{-1}s_{H'}s_H$ is a reflection and $x$ belongs to
its reflecting hyperplane, $s_H^{-1}(H')$. This contradicts the fact
that $x$ belongs to $V^\mathrm{reg}$.). Now we define a path from
$x_0$ to $s_H(x_0)$ as follows:
$$\sigma_{H,\gamma}:=s_H(\gamma^{-1}(t)) \cdot \sigma_{H,x_H} \cdot
\gamma$$

If $x_H$ is chosen ``close to $H$ and far from the other reflecting
hyperplanes'', the path $\sigma_{H,\gamma}$ lies in $V^\mathrm{reg}$
and its homotopy class doesn't depend on the choice of $x_H$. The
element it induces in the braid group $B$, $\textbf{s}_{H,\gamma}$,
is a distinguished braid reflection around the image of $H$ in
$V^\mathrm{reg}/W$.

\begin{proposition}\label{braid reflections}\
\begin{enumerate}
  \item The braid group $B$ is generated by the distinguished braid
   reflections around the images of the hyperplanes $H \in
   \mathcal{A}$ in $V^\mathrm{reg}/W$.
  \item The image of $\emph{\textbf{s}}_{H,\gamma}$ in $W$ is $s_H$.
  \item Whenever $\gamma'$ is a path in $V^\mathrm{reg}$ from $x_0$
   to $x_H$, if $\tau$ denotes the loop in $V^\mathrm{reg}$ defined
   by $\tau:=\gamma'^{-1}\gamma$, then
   $$\sigma_{H,\gamma'}=s_H(\tau) \cdot \sigma_{H,\gamma} \cdot
   \tau^{-1}$$
   and in particular $\emph{\textbf{s}}_{H,\gamma}$ and
  $\emph{\textbf{s}}_{H,\gamma}$ are conjugate in $P$.
  \item The path $\prod_{j=e_H-1}^{j=0}\sigma_{H,s_H^j(\gamma)}$,
  a loop in $V^\mathrm{reg}$, induces the element
  $\emph{\textbf{s}}_{H,\gamma}^{e_H}$ in the braid group $B$ and belongs to the pure braid
  group $P$. It is a distinguished braid reflection around $H$ in
  $P$.
\end{enumerate}
\end{proposition}

\begin{definition}\label{s-distinguished braid reflection}
Let $s$ be a distinguished pseudo-reflection in $W$ with reflecting
hyperplane $H$. An $s$-distinguished braid reflection or monodromy
generator is a distinguished braid reflection $\emph{\textbf{s}}$
around the image of $H$ in $V^\mathrm{reg}/W$ such that
$\bar{\emph{\textbf{s}}}=s$.
\end{definition}

\begin{definition}\label{pi}
Let $x_0 \in V^\mathrm{reg}$ as before. We denote by \textbf{$\pi$}
the element of $P$ defined by the loop  $t \mapsto
x_0\mathrm{exp}(2\pi it).$
\end{definition}

\begin{lemma}\label{pi in ZP}
We have \textbf{$\pi$}$\in ZP$.
\end{lemma}

\begin{thedef}\label{length function}
Given $\mathcal{C} \in \mathcal{A}/W$, there is a unique length
function $l_\mathcal{C}:B \rightarrow \mathbb{Z}$ defined as
follows: if $b=\emph{\textbf{s}}_1^{n_1} \cdot
\emph{\textbf{s}}_2^{n_2} \cdot \cdot \cdot
\emph{\textbf{s}}_m^{n_m}$ where (for all $j$) $n_j \in \mathbb{Z}$
and $\emph{\textbf{s}}_j$ is a distinguished braid reflection around
an element of $\mathcal{C}_j$, then
$$l_\mathcal{C}(b)=\sum_{\{j\,|\,\mathcal{C}_j=\mathcal{C}\}}n_j.$$
Thus, the length function $l:B \rightarrow \mathbb{Z}$ is defined,
for all $b \in B$, as
$$l(b)=\sum_{\mathcal{C} \in \mathcal{A}/W}l_\mathcal{C}(b).$$
\end{thedef}

We say that $B$ has an \emph{Artin-like} presentation (\cite{Op2},
5.2), if it has a presentation of the form
$$<\textbf{s} \in \textbf{S} \,|\, \{\textbf{v}_i=\textbf{w}_i\}_{i \in
I}>,$$ where $\textbf{S}$ is a finite set of distinguished braid
reflections and $I$ is a finite set of relations which are
multi-homogeneous, \ie such that for all $i$, $\textbf{v}_i$ and
$\textbf{w}_i$ are positive words in elements of $\textbf{S}$ (and
hence, for each $\mathcal{C} \in \mathcal{A}/W$, we have
$l_\mathcal{C}(\textbf{v}_i)=l_\mathcal{C}(\textbf{w}_i)$).

The following result by Bessis (\cite{Bes3}, Thm.0.1) shows that any
braid group has an Artin-like presentation.

\begin{theorem}\label{bessis}
Let $W$ be a complex reflection group with associated braid group
$B$. Then there exists a subset
$\emph{\textbf{S}}=\{\emph{\textbf{s}}_1,\ldots,\emph{\textbf{s}}_n\}$
of $B$ such that
\begin{enumerate}
  \item The elements $\emph{\textbf{s}}_1,\ldots,\emph{\textbf{s}}_n$ are
  distinguished braid reflection and therefore, their images $s_1,\ldots,s_n$ in
  $W$ are distinguished reflections.
  \item The set $\emph{\textbf{S}}$ generates $B$ and therefore,
  $S:=\{s_1,\ldots,s_n\}$ generates $W$.
  \item There exists a set $\mathcal{R}$ of relations of the form
  $\emph{\textbf{w}}_1=\emph{\textbf{w}}_2$, where $\emph{\textbf{w}}_1$ and
  $\emph{\textbf{w}}_2$ are positive words of equal length in the elements
  of $\emph{\textbf{S}}$, such that $<\emph{\textbf{S} }\,|\, \mathcal{R}>$ is a
  presentation of $B$.
  \item Viewing now $\mathcal{R}$ as a set of relations in $S$, the
  group $W$ is presented by
  $$<S \,|\, \mathcal{R}; (\forall s \in S)(s^{e_s}=1)>$$
  where $e_s$ denotes the order of $s$ in $W$.
\end{enumerate}
\end{theorem}

\section{Generic Hecke algebras}

Let $K,V,W,\mathcal{A},P,B$ be defined as in the previous section.
For every orbit $\mathcal{C}$ of $W$ on $\mathcal{A}$, we set
$e_{\mathcal{C}}$ the common order of the subgroups $W_H$, where $H$
is any element of $\mathcal{C}$ and $W_H$ the subgroup formed by 1
and all the reflections fixing the hyperplane $H$.

We choose a set of indeterminates
$\textbf{u}=(u_{\mathcal{C},j})_{(\mathcal{C} \in
\mathcal{A}/W)(0\leq j \leq e_{\mathcal{C}}-1)}$ and we denote by
$\mathbb{Z}[\textbf{u},\textbf{u}^{-1}]$ the Laurent polynomial ring
in all the indeterminates $\textbf{u}$. We define the \emph{generic
Hecke algebra} $\mathcal{H}$ of $W$ to be the quotient of the group
algebra $\mathbb{Z}[\textbf{u},\textbf{u}^{-1}]B$ by the ideal
generated by the elements of the form
$$(\textbf{s}-u_{\mathcal{C},0})(\textbf{s}-u_{\mathcal{C},1}) \ldots (\textbf{s}-u_{\mathcal{C},e_{\mathcal{C}}-1}),$$
where $\mathcal{C}$ runs over the set $\mathcal{A}/W$ and
$\textbf{s}$ runs over the set of monodromy generators around the
images in $V^{\textrm{reg}}/W$ of the elements of the hyperplane
orbit $\mathcal{C}$.

\begin{px}
\small{\emph{Let $W:=G_4=<s,t \,|\, sts=tst, s^3=t^3=1>$. Then $s$
and $t$ are conjugate in $W$ and their reflecting hyperplanes belong
to the same orbit in $\mathcal{A}/W$. The generic Hecke algebra of
$W$ can be presented as follows
$$\begin{array}{rll}
   \mathcal{H}=<S,T \,\,|&STS=TST, &(S-u_0)(S-u_1)(S-u_2)=0, \\
                                &         &(T-u_0)(T-u_1)(T-u_2)=0>.
  \end{array}$$}}
\end{px}

We make some assumptions for the algebra $\mathcal{H}$. Note that
they have been verified for all but a finite number of irreducible
complex reflection groups (\cite{BMM2}, remarks before 1.17, $\S$ 2;
\cite{GIM}).

\begin{ypoth}\label{ypo}
The algebra $\mathcal{H}$ is a free
$\mathbb{Z}[\textbf{\emph{u}},\textbf{\emph{u}}^{-1}]$-module of
rank $|W|$. Moreover, there exists a linear form
$t:\mathcal{H}\rightarrow
\mathbb{Z}[\textbf{\emph{u}},\textbf{\emph{u}}^{-1}]$ with the
following properties:
\begin{enumerate}
    \item $t$ is a symmetrizing form for $\mathcal{H}$, i.e.,
     $t(hh')=t(h'h)$ for all $h,h' \in \mathcal{H}$ and the map
     $$\begin{array}{cccc}
     \hat{t}: & \mathcal{H} & \rightarrow & \textrm{\emph{Hom}}(\mathcal{H},\mathbb{Z}[\textbf{\emph{u}},\textbf{\emph{u}}^{-1}]) \\
              & h & \mapsto & (h' \mapsto t(hh'))
     \end{array}$$
     is an isomorphism.
    \item Via the specialization $u_{\mathcal{C},j} \mapsto
     \zeta_{e_\mathcal{C}}^j$, the form $t$ becomes the canonical
     symmetrizing form on the group algebra $\mathbb{Z}W$.
    \item If we denote by $\alpha \mapsto \alpha^*$ the automorphism of
     $\mathbb{Z}[\emph{\textbf{u}},\emph{\textbf{u}}^{-1}]$ consisting of the
     simultaneous inversion of the indeterminates, then for all $b \in B$, we
     have
          $$t(b^{-1})^*=\frac{t(b\pi)}{t(\pi)},$$
     where $\pi$ is the central element of $P$ defined in
     $\ref{pi}$.
\end{enumerate}
\end{ypoth}

We know that the form $t$ is unique (\cite{BMM2}, 2.1). From now on,
let us suppose that the assumptions $\ref{ypo}$ are satisfied. Then
we have the following result by G.Malle (\cite{Ma4}, 5.2).

\begin{theorem}\label{Semisimplicity Malle}
Let $\textbf{\emph{v}}=(v_{\mathcal{C},j})_{(\mathcal{C} \in
\mathcal{A}/W)(0\leq j \leq e_{\mathcal{C}}-1)}$ be a set of
$\sum_{\mathcal{C} \in \mathcal{A}/W}e_{\mathcal{C}}$ indeterminates
such that, for every $\mathcal{C},j$, we have
$v_{\mathcal{C},j}^{|\mu(K)|}=\zeta_{e_\mathcal{C}}^{-j}u_{\mathcal{C},j}$.
Then the $K(\textbf{\emph{v}})$-algebra
$K(\textbf{\emph{v}})\mathcal{H}$ is split semisimple.
\end{theorem}

By ``Tits' deformation theorem'' (theorem $\ref{Tits}$), it follows
that the specialization $v_{\mathcal{C},j}\mapsto 1$ induces a
bijection $\chi \mapsto \chi_{\textbf{v}}$ from the set
$\mathrm{Irr}(W)$ of absolutely irreducible characters of $W$ to the
set $\mathrm{Irr}(K(\textbf{v})\mathcal{H})$ of absolutely
irreducible characters of $K(\textbf{v})\mathcal{H}$, such that the
following diagram is commutative $$\begin{array}{rccc}
  \chi_\textbf{v} : & \mathcal{H} & \rightarrow & \mathbb{Z}_K[\textbf{v},\textbf{v}^{-1}] \\
  & \downarrow &  & \downarrow \\
  \chi: & \mathbb{Z}_KW &\rightarrow &\mathbb{Z}_K .
\end{array}$$

Since the assumptions $\ref{ypo}$ are satisfied and the algebra
$K(\textbf{v})\mathcal{H}$ is split semisimple, we can define the
Schur element $s_{\chi_{\textbf{v}}}$ for every irreducible
character $\chi_{\textbf{v}}$ of $K(\textbf{v})\mathcal{H}$ with
respect to the symmetrizing form $t$. The bijection $\chi \mapsto
\chi_{\textbf{v}}$ from $\mathrm{Irr}(W)$ to
$\mathrm{Irr}(K(\textbf{v})\mathcal{H})$ implies that the
specialization $v_{\mathcal{C},j}\mapsto 1$ sends
$s_{\chi_{\textbf{v}}}$ to $|W|/\chi(1)$ (which is the Schur element
of $\chi$ in the group algebra with respect to the canonical
symmetrizing form). The following result is simply the application
of proposition $\ref{schur elements and idempotents}$ to this case.

\begin{proposition}\label{recall form of blocks}\
\begin{enumerate}
  \item We have
  $$t=\sum_{\chi_{\textbf{\emph{v}}} \in \mathrm{Irr}(K(\textbf{\emph{v}})\mathcal{H})}
  \frac{1}{s_{\chi_{\textbf{\emph{v}}}}}\chi_{\textbf{\emph{v}}}.$$
  \item For all $\chi_{\textbf{\emph{v}}} \in \mathrm{Irr}(K(\textbf{\emph{v}})\mathcal{H})$, the
  block-idempotent of $K(\textbf{\emph{v}})\mathcal{H}$ associated with $\chi_{\textbf{\emph{v}}}$ is
  $e_{\chi_{\textbf{\emph{v}}}}=\chi_{\textbf{\emph{v}}}^\vee/s_{\chi_{\textbf{\emph{v}}}}.$
\end{enumerate}
\end{proposition}
\begin{remark}
\emph{ The bijection $\mathrm{Irr}(W) \leftrightarrow
\mathrm{Irr}(K(\textbf{v})\mathcal{H})$, $\chi \mapsto
\chi_\textbf{v}$ allows us to write $\mathrm{Irr}(W)$ instead of
$\mathrm{Irr}(K(\textbf{v})\mathcal{H})$ and $\chi$ instead of
$\chi_\textbf{v}$ in all the relations above.}
\end{remark}\\

Our first result concerns the form of the Schur elements associated
with the irreducible characters of $K(\textbf{v})\mathcal{H}$
(always assuming that the asumptions $\ref{ypo}$ are satisfied). We
will see later that this result plays a crucial role in the
determination of the blocks of Hecke algebras.

\begin{theorem}\label{Schur element generic}
The Schur element $s_\chi(\textbf{\emph{v}})$ associated with the
character $\chi_{\textbf{\emph{v}}}$ of
$K(\textbf{\emph{v}})\mathcal{H}$ is an element of
$\mathbb{Z}_K[\textbf{\emph{v}},\textbf{\emph{v}}^{-1}]$ of the form
$$s_\chi({\textbf{\emph{v}}})=\xi_\chi N_\chi \prod_{i \in I_\chi} \Psi_{\chi,i}(M_{\chi,i})^{n_{\chi,i}}$$
where
\begin{itemize}
    \item $\xi_\chi$ is an element of $\mathbb{Z}_K$,
    \item $N_\chi= \prod_{\mathcal{C},j} v_{\mathcal{C},j}^{b_{\mathcal{C},j}}$ is a monomial in $\mathbb{Z}_K[\textbf{\emph{v}},\textbf{\emph{v}}^{-1}]$
          with $\sum_{j=0}^{e_\mathcal{C}-1}b_{\mathcal{C},j}=0$
          for all $\mathcal{C} \in \mathcal{A}/W$,
    \item $I_\chi$ is an index set,
    \item $(\Psi_{\chi,i})_{i \in I_\chi}$ is a family of $K$-cyclotomic polynomials in one variable
           (i.e., minimal polynomials of the roots of unity over $K$),
    \item $(M_{\chi,i})_{i \in I_\chi}$ is a family of monomials in $\mathbb{Z}_K[\textbf{\emph{v}},\textbf{\emph{v}}^{-1}]$
          and if $M_{\chi,i} = \prod_{\mathcal{C},j} v_{\mathcal{C},j}^{a_{\mathcal{C},j}}$,
          then $\textrm{\emph{gcd}}(a_{\mathcal{C},j})=1$
          and $\sum_{j=0}^{e_\mathcal{C}-1}a_{\mathcal{C},j}=0$
          for all $\mathcal{C} \in \mathcal{A}/W$,
    \item ($n_{\chi,i})_{i \in I_\chi}$ is a family of positive integers.
\end{itemize}
\end{theorem}
\begin{apod}{By proposition $\ref{Schur element belongs to the integral closure}$,
we have that $s_\chi(\textbf{v}) \in
\mathbb{Z}_K[\textbf{v},\textbf{v}^{-1}]$. The rest is a case by
case analysis. For $W$ an irreducible complex reflection group, we
will denote by $\mathcal{H}(W)$ its generic Hecke algebra defined
over the splitting field of theorem $\ref{Semisimplicity Malle}$.

Let us first consider the group $G(d,1,r)$ for $d \geq 1, r>2$. By
\cite{Mat}, Cor. 6.5, the Schur elements of $\mathcal{H}(G(d,1,r))$
are of the desired form. In \cite{Kim} it was shown that
$\mathcal{H}(G(de,1,r))$ for a specific choice of parameters becomes
the twisted symmetric algebra of the cyclic group of order $e$ over
a symmetric subalgebra which is isomorphic to
$\mathcal{H}(G(de,e,r))$ (it had been already shown in \cite{BK} for
$d=1$). Thus, by proposition $\ref{1.42}$, the Schur elements of
$\mathcal{H}(G(de,1,r))$ are multiples by some integer of the Schur
elements of $\mathcal{H}(G(de,e,r))$. Therefore, the assertion is
established for all the groups of the infinite series $G(de,e,r)$
with $r>2$.

For the groups $G(de,e,2)$, $G_7$, $G_{11}$ and $G_{19}$, the
generic Schur elements are determined by Malle in \cite{Ma2}. In the
same article, we can find the specializations of parameters which
permit us to calculate, using again proposition $\ref{1.42}$, the
Schur elements of
\begin{itemize}
\item $\mathcal{H}(G_4)$, $\mathcal{H}(G_5)$, $\mathcal{H}(G_6)$ from
$\mathcal{H}(G_7)$.
\item $\mathcal{H}(G_8)$, $\mathcal{H}(G_9)$, $\mathcal{H}(G_{10})$,
$\mathcal{H}(G_{12})$, $\mathcal{H}(G_{13})$, $\mathcal{H}(G_{14})$,
$\mathcal{H}(G_{15})$ from $\mathcal{H}(G_{11})$.
\item $\mathcal{H}(G_{16})$, $\mathcal{H}(G_{17})$, $\mathcal{H}(G_{18})$,
$\mathcal{H}(G_{20})$, $\mathcal{H}(G_{21})$, $\mathcal{H}(G_{22})$
from $\mathcal{H}(G_{19})$.
\end{itemize}
For more details, the reader may refer to the Appendix, where the
above specializations are given explicitly.

The generic Schur elements for the remaining non-Coxeter exceptional
complex reflection groups, \ie the groups $G_{24}$, $G_{25}$,
$G_{26}$, $G_{27}$, $G_{29}$, $G_{31}$, $G_{32}$, $G_{33}$,
$G_{34}$, have been also calculated by Malle in \cite{Ma5}.

As far as the exceptional real reflection groups are concerned, \ie
the groups
$G_{23}=H_3,G_{28}=F_4,G_{30}=H_4,G_{35}=E_6,G_{36}=E_7,G_{37}=E_8$,
the Schur elements have been
calculated
\begin{itemize}
  \item for $E_6$ and $E_7$ by Surowski (\cite{Sur78}),
  \item for $E_8$ by Benson (\cite{Ben79}),
  \item for $F_4$ and $H_3$ by Lusztig (\cite{Lu79b} and
  \cite{Lu82} respectively),
  \item for $H_4$ by Alvis and Lusztig (\cite{AlLu82}).
\end{itemize}

To obtain the desired formula from the data given in the
above articles, we used the GAP Package CHEVIE (where some mistakes in these
articles have been corrected).}
\end{apod}

In the Appendix of this thesis, we give the factorization of the
generic Schur elements of the groups $G_7$, $G_{11}$, $G_{19}$,
$G_{26}$, $G_{28}$ and $G_{32}$, so that the reader may verify the
above result. The Schur elements for $G_{25}$ are also obtained as
specializations of the Schur elements of $G_{26}$. The groups
$G_{23}$, $G_{24}$, $G_{27}$, $G_{29}$, $G_{30}$ $G_{31}$, $G_{33}$,
$G_{34}$, $G_{35}$, $G_{36}$ and $G_{37}$ are all generated by
reflections of order 2 whose reflecting hyperplanes belong to one
single orbit. Therefore, the splitting field of their generic Hecke
algebra is of the form $K(v_0,v_1)$, where $K$ is the field of
definition of the group. In these cases, the generic Hecke algebra
is essentially one-parametered and it is easy to check that the
irreducible factors of the generic Schur elements over
$K[v_0^\pm,v_1^\pm]$ are $K$-cyclotomic polynomials taking values on
$v:=v_0v_1^{-1}$.\\
\\
\begin{remark}\emph{ It is a consequence of \cite{Raph}, Thm. 3.4, that the
irreducible factors of the generic Schur elements over
$\mathbb{C}[\textbf{v},\textbf{v}^{-1}]$ are divisors of Laurent
polynomials of the form $M(\textbf{v})^n-1$, where
\begin{itemize}
\item $M(\textbf{v})$ is a monomial in $\mathbb{C}[\textbf{v},\textbf{v}^{-1}]$,
\item $n$ is a positive integer.
\end{itemize}}
\end{remark}\

Thanks to proposition $\ref{second irreducible}$, the factorization
of the proposition $\ref{Schur element generic}$ is unique in
$K[\textbf{v},\textbf{v}^{-1}]$. However, this does not mean that
the monomials $M_{\chi,i}$ appearing in it are unique. Let
$$\Psi_{\chi,i}(M_{\chi,i})= u \Psi_{\chi,j}(M_{\chi,j}),$$
where $i, j \in I$, $\Psi_{\chi,i}, \Psi_{\chi,j}$ are two
$K$-cyclotomic polynomials, $M_{\chi,i},M_{\chi,j}$ are two
monomials in $K[\textbf{v},\textbf{v}^{-1}]$ with the properties
described in $\ref{Schur element generic}$ and $u$ is a unit element
in $K[\textbf{v},\textbf{v}^{-1}]$. Let $\varphi_i$ be a morphism
associated with the monomial $M_{\chi,i}$ (see definition
$\ref{associated morphism}$) from
$\mathbb{Z}_K[\textbf{v},\textbf{v}^{-1}]$ to a Laurent polynomial
ring with one indeterminate less. If we apply $\varphi_i$ to the
above equality, we obtain
$$\Psi_{\chi,i}(1)=\varphi_i(u)\varphi_i(\Psi_{\chi,j}(M_{\chi,j})).$$
Since $\Psi_{\chi,i}(1) \in \mathbb{Z}_K$ and $\varphi_i$ sends
$M_{\chi,j}$ to a monomial, we deduce that
$\varphi_i(M_{\chi,j})=1$. By proposition $\ref{properties of
phi}$(2) and the fact that $M_{\chi,j}$ satisfies the conditions
described in proposition $\ref{Schur element generic}$, we must have
$$M_{\chi,j} = M_{\chi,i}^{\pm 1}.$$
Hence, the monomials $M_{\chi,i}$ appearing in the factorization of
the generic Schur element are unique up to inversion.

If $M_{\chi,j} = M_{\chi,i}$, then $\Psi_{\chi,j} = \Psi_{\chi,i}$
and $u=1$. If $M_{\chi,j} = M_{\chi,i}^{- 1}$, then $\Psi_{\chi,j}$
is conjugate to $\Psi_{\chi,i}$ and $u$ is of the form $\zeta
M_{\chi,i}^{\mathrm{deg}(\Psi_{\chi,i})}$, where
$\zeta$ is a root of unity. Therefore, the coefficient $\xi_\chi$ is
unique up to a root of unity.\\
\\
\begin{remark}\emph{ The first cyclotomic polynomial never appears in the
factorization of a Schur element $s_\chi(\textbf{v})$. Otherwise the
specialization $v_{\mathcal{C},j} \mapsto 1$ would send
$s_\chi(\textbf{v})$ to 0 and not to $|W|/\chi(1)$ as it should.}
\end{remark}
\\

Now let $\mathfrak{p}$ be a prime ideal of $\mathbb{Z}_K$ and let us
denote by $A$ the ring $\mathbb{Z}_K[\textbf{v},\textbf{v}^{-1}]$.
If $\Psi(M_{\chi,i})$ is a factor of $s_\chi(\textbf{v})$ and
$\Psi_{\chi,i}(1) \in \mathfrak{p}$, then the monomial $M_{\chi,i}$
is called $\mathfrak{p}$\emph{-essential} for $\chi$ in $A$. Due to
proposition $\ref{properties of phi}$(1), we have
$$\Psi_{\chi,i}(1) \in \mathfrak{p} \Leftrightarrow
  \Psi_{\chi,i}(M_{\chi,i}) \in \mathfrak{q}_i,\,\,\,\,\,(\dag)$$
where $\mathfrak{q}_i:=(M_{\chi,i}-1)A +\mathfrak{p}A.$ Recall that
$\mathfrak{q}_i$ is a prime ideal of $A$ (lemma $\ref{primeness of
q}$). Due to the primeness of $\mathfrak{q}_i$, the following
proposition is an immediate consequence of $(\dag)$.

\begin{proposition}\label{p-essential}
Let $M:=\prod_{\mathcal{C},j}v_{\mathcal{C},j}^{a_{\mathcal{C},j}}$
be a monomial in $A$ with $\mathrm{gcd}(a_{\mathcal{C},j})=1$ and
$\mathfrak{q}_M:=(M-1)A +\mathfrak{p}A$. Then  $M$ is
$\mathfrak{p}$-essential for $\chi$ in $A$ if and only if
$s_\chi/\xi_\chi \in \mathfrak{q}_M$, where $\xi_\chi$ denotes the
coefficient of $s_\chi$ in $\ref{Schur element generic}$.
\end{proposition}

From now on, for reasons of convenience and only for this section,
we will substitute the set of indeterminates
$\textbf{v}=(v_{\mathcal{C},j})_{(\mathcal{C} \in
\mathcal{A}/W)(0\leq j \leq e_{\mathcal{C}}-1)}$ with the set
$\{x_0,x_1,\ldots,x_m\}$, where $m:=(\sum_{\mathcal{C} \in
\mathcal{A}/W}e_{\mathcal{C}})-1$. Hence, the algebra $\mathcal{H}$
will be considered over the ring $A=\mathbb{Z}_K[x_0^{\pm
1},x_1^{\pm 1},\ldots,x_m^{\pm 1}]$.\\

Let $M:=\prod_{i=0}^m x_i^{a_i}$ be a monomial in $A$, such that
$a_i \in \mathbb{Z}$ and gcd$(a_i)=1$. Let $B:=\mathbb{Z}_K[y_1^{\pm
1},\ldots,y_m^{\pm 1}]$ and consider $\varphi_M: A \rightarrow B$ a
$\mathbb{Z}_K$-algebra morphism associated with $M$. Let us denote
by $\mathcal{H}_{\varphi_M}$ the algebra obtained by $\mathcal{H}$
via the specialization $\varphi_M$.

\begin{proposition}\label{associated morphism preserves splitness}
The algebra $K(y_1,\ldots,y_m)\mathcal{H}_{\varphi_M}$ is split
semisimple.
\end{proposition}
\begin{apod}{By theorem $\ref{Semisimplicity Malle}$, the algebra
$K(x_0,x_1,\ldots,x_m)\mathcal{H}$ is split semisimple. The ring $A$
is a Krull ring and $\mathrm{Ker}\varphi_M=(M-1)A$ is a prime ideal
of height 1 of $A$. Due to the form of the generic Schur elements,
given in theorem $\ref{Schur element generic}$, the fact that the
morphism $\varphi_M$ sends every monomial in $A$ to a monomial in
$B$ implies that $\varphi_M(s_\chi) \neq 0 $ for all $\chi \in
\mathrm{Irr}(K(x_0,x_1,\ldots,x_m)\mathcal{H})$ (we have already
explained why the first cyclotomic polynomial never appears in the
factorization of the generic Schur elements). Thus we can apply
theorem $\ref{Lehrer}$ and obtain that the algebra
$K(y_1,\ldots,y_m)\mathcal{H}_{\varphi_M}$ is split semisimple.}
\end{apod}

By ``Tits'deformation theorem'', the specialization $y_i \mapsto 1$
induces a bijection from the set
$\mathrm{Irr}(K(y_1,\ldots,y_m)\mathcal{H}_{\varphi_M})$ of
absolutely irreducible characters of
$K(y_1,\ldots,y_m)\mathcal{H}_{\varphi_M}$ to the set
$\mathrm{Irr}(W)$. The Schur elements of the former are the
specializations of the Schur elements of
$K(x_0,x_1,\ldots,x_m)\mathcal{H}$ via $\varphi_M$ and thus of the
form described in $\ref{Schur element generic}$.

From now on, whenever we refer to irreducible characters, we mean
irreducible characters of the group $W$. Due to the existing
bijections
$$\mathrm{Irr}(K(x_0,x_1,\ldots,x_m)\mathcal{H}) \leftrightarrow
  \mathrm{Irr}(K(y_1,\ldots,y_m)\mathcal{H}_{\varphi_M}) \leftrightarrow
  \mathrm{Irr}(W),$$
it makes sense to compare the blocks of $\mathcal{H}$ and
$\mathcal{H}_{\varphi_M}$ (in terms of partitions of $
\mathrm{Irr}(W)$) over suitable rings.\\

Let $\mathfrak{p}$ be a prime ideal of $\mathbb{Z}_K$ and
$\mathfrak{q}_M:=(M-1)A + \mathfrak{p}A$.

\begin{theorem}\label{Aq B}
The  blocks of $B_{\mathfrak{p}B}\mathcal{H}_{\varphi_M}$ coincide
with the blocks of $A_{\mathfrak{q}_M}\mathcal{H}$.
\end{theorem}
\begin{apod}{Let us denote by $\mathfrak{n}_M$ the kernel of $\varphi_M$, \ie
$\mathfrak{n}_M:=(M-1)A$. By proposition $\ref{properties of
phi}$(1), we have that $A_{\mathfrak{q}_M}/\mathfrak{n}_M
A_{\mathfrak{q}_M} \simeq B_{\mathfrak{p}B}$. Therefore, it is
enough to show that the canonical surjection
$A_{\mathfrak{q}_M}\mathcal{H} \twoheadrightarrow
(A_{\mathfrak{q}_M}/\mathfrak{n}_M A_{\mathfrak{q}_M})\mathcal{H}$
induces a block bijection between these two algebras.

From now on, the symbol $ $ $\widehat{}$ $ $ will stand for
$\mathfrak{q}_M$-adic completion. It is immediate that, via the
canonical surjection, a block of $A_{\mathfrak{q}_M}\mathcal{H}$ is
a sum of blocks of $(A_{\mathfrak{q}_M}/\mathfrak{n}_M
A_{\mathfrak{q}_M})\mathcal{H}$. Now let $\bar{e}$ be a block of
$(A_{\mathfrak{q}_M}/\mathfrak{n}_M A_{\mathfrak{q}_M})\mathcal{H}$.
By theorem $\ref{17.6}$, all Noetherian local rings are contained in
their completions. Thus $\bar{e}$ can be written as a sum of blocks
of $\widehat{(A_{\mathfrak{q}_M}/\mathfrak{n}_M
A_{\mathfrak{q}_M})}\mathcal{H}$. Due to corollary $\ref{17.9}$, the
last algebra is isomorphic to
$(\hat{A}_{\mathfrak{q}_M}/\mathfrak{n}_M\hat{A}_{\mathfrak{q}_M})\mathcal{H}$,
which is in turn isomorphic to the quotient algebra
$\hat{A}_{\mathfrak{q}_M}\mathcal{H}/\mathfrak{n}_M\hat{A}_{\mathfrak{q}_M}\mathcal{H}$.
By the theorems of lifting idempotents (see \cite{The}, Thm.3.2) and
the following lemma, $\bar{e}$ is lifted to a sum of central
primitive idempotents in $\hat{A}_{\mathfrak{q}_M}\mathcal{H}$.
However, by the fact that $K(x_0,x_1,\ldots,x_m)\mathcal{H}$ is
split semisimple, we have that the blocks of
$\hat{A}_{\mathfrak{q}_M}\mathcal{H}$ belong to
$K(x_0,x_1,\ldots,x_m)\mathcal{H}$. But $K(x_0,x_1,\ldots,x_m) \cap
\hat{A}_{\mathfrak{q}_M}=A_{\mathfrak{q}_M}$ (theorem $\ref{18.4}$)
and $A_{\mathfrak{q}_M}\mathcal{H} \cap
Z(\hat{A}_{\mathfrak{q}_M}\mathcal{H}) \subseteq
Z(A_{\mathfrak{q}_M}\mathcal{H})$. Therefore, $\bar{e}$ is lifted to
a sum of blocks in $A_{\mathfrak{q}_M}\mathcal{H}$ and this provides
the block bijection.
\end{apod}

\begin{lemma}\label{central lifts to central}
Let $\mathcal{O}$ be a Noetherian ring and $\mathfrak{q}$ a prime
ideal of $\mathcal{O}$. Let $H$ be an
$\mathcal{O}_\mathfrak{q}$-algebra free and of finite rank as a
$\mathcal{O}_\mathfrak{q}$-module. Let $\mathfrak{p}$ be a prime
ideal of $\mathcal{O}$ such that $\mathfrak{p} \subseteq
\mathfrak{q}$ and $e$ an idempotent of $H$ whose image $\bar{e}$ in
$(\mathcal{O}_\mathfrak{q}/\mathfrak{p}_\mathfrak{q})H:=\mathcal{O}_\mathfrak{q}/\mathfrak{p}\mathcal{O}_\mathfrak{q}
\otimes_{\mathcal{O}_\mathfrak{q}}H$ is central. Then $e$ is
central.
\end{lemma}
\begin{apod}{We set $P:=\mathcal{O}_\mathfrak{q}/\mathfrak{p}_\mathfrak{q}$.
Since $\bar{e}$ is central, we have
$$\bar{e}PH(1-\bar{e})=(1-\bar{e})PH\bar{e}=\{0\},$$ \ie
$$eH(1-e) \subseteq \mathfrak{p}\mathcal{O}_\mathfrak{q}H \,\,\textrm{
and }\,\, (1-e)He \subseteq \mathfrak{p}\mathcal{O}_\mathfrak{q}H.$$
Since $e$ and $(1-e)$ are idempotents, we get
$$eH(1-e) \subseteq \mathfrak{p}\mathcal{O}_\mathfrak{q}eH(1-e) \,\,\textrm{ and
}\,\, (1-e)He \subseteq
\mathfrak{p}\mathcal{O}_\mathfrak{q}(1-e)He.$$ However,
$\mathfrak{p}\mathcal{O}_\mathfrak{q} \subseteq
\mathfrak{q}\mathcal{O}_\mathfrak{q}$, the latter being the maximal
ideal of $\mathcal{O}_\mathfrak{q}$. By Nakayama's lemma,
$$eH(1-e)=(1-e)He=\{0\}.$$
Thus, from
$$H=eHe \oplus eH(1-e) \oplus (1-e)He \oplus (1-e)H(1-e)$$
we deduce that
$$H=eHe \oplus (1-e)H(1-e)$$
and consequently, $e$ is central.}
\end{apod}

}
\\
\begin{remark}\emph{ Since $\mathfrak{q}_M=\mathfrak{q}_{M^{-1}}$,
proposition $\ref{Aq B}$ implies that the blocks of
$B_{\mathfrak{p}B}\mathcal{H}_{\varphi_M}$ coincide with the blocks
of $B_{\mathfrak{p}B}\mathcal{H}_{\varphi_{M^{-1}}}$.}
\end{remark}

\begin{proposition}\label{simple inclusion}
If two irreducible characters $\chi$ and $\psi$ are in the same
block of $A_{\mathfrak{p}A}\mathcal{H}$, then they are in the same
block of $A_{\mathfrak{q}_M}\mathcal{H}$.
\end{proposition}
\begin{apod}{Let $C$ be a block of $A_{\mathfrak{q}_M}\mathcal{H}$. Then $\sum_{\chi \in C} e_\chi
\in A_{\mathfrak{q}_M}\mathcal{H} \subset
A_{\mathfrak{p}A}\mathcal{H}$. Thus $C$ is a union of blocks of
$A_{\mathfrak{p}A}\mathcal{H}$.}
\end{apod}

\begin{corollary}\label{union of blocks}
If two irreducible characters $\chi$ and $\psi$ are in the same
block of $A_{\mathfrak{p}A}\mathcal{H}$, then they are in the same
block of $B_{\mathfrak{p}B}\mathcal{H}_{\varphi_M}$.
\end{corollary}

The corollary above implies that the size of $\mathfrak{p}$-blocks
grows larger as the number of indeterminates becomes smaller.
However, we will now see that the size of blocks remains the same,
if our specialization is not associated with a
$\mathfrak{p}$-essential monomial.

\begin{proposition}\label{not essential for block}
Let $C$ be a block of $A_{\mathfrak{p}A}\mathcal{H}$. If $M$ is not
a $\mathfrak{p}$-essential monomial for any $\chi \in C$, then $C$
is a block of $A_{\mathfrak{q}_M}\mathcal{H}$ (and thus of
$B_{\mathfrak{p}B}\mathcal{H}_{\varphi_M}$).
\end{proposition}
\begin{apod}{ Using the notations of proposition $\ref{Schur element generic}$,
we have that, for all $\chi \in C$, $s_\chi/\xi_\chi \notin
\mathfrak{q}_M$.  Since $C$ is a block of
$A_{\mathfrak{p}A}\mathcal{H}$, we have  $$\sum_{\chi \in C} e_\chi
=\sum_{\chi \in C} \frac{\chi^\vee}{s_\chi} \in
A_{\mathfrak{p}A}\mathcal{H}.$$ If $\mathcal{B},\mathcal{B}'$ are
two bases of $\mathcal{H}$ dual to each other, then
$\chi^\vee=\sum_{b \in \mathcal{B}}\chi(b)b'$ and the above relation
implies that $$\sum_{\chi \in C}\frac{\chi(b)}{s_\chi} \in
A_{\mathfrak{p}A}, \forall b \in \mathcal{B}.$$ Let $f_b:=\sum_{\chi
\in C}(\chi(b)/s_\chi) \in A_{\mathfrak{p}A}$. Then $f_b$ is of the
form $r_b/(\xi s)$, where $\xi:=\prod_{\chi \in C}\xi_\chi \in
\mathbb{Z}_K$ and $s:=\prod_{\chi \in C} s_\chi/\xi_\chi \in A$.
Since $\mathfrak{q}_M$ is a prime ideal of $A$, the element $s$, by
assumption, doesn't belong to $\mathfrak{q}_M$. We also have that
$r_b/\xi \in A_{\mathfrak{p}A}$. By corollary $\ref{porisma
porismatos}$, there exists $\xi' \in \mathbb{Z}_K-\mathfrak{p}$ such
that $r_b/\xi=r_b'/\xi'$ for some $r_b' \in A$. Since
$\mathfrak{q}_M \cap \mathbb{Z}_K = \mathfrak{p}$ (corollary $\ref{q
intersection zk}$), the element $\xi'$ doesn't belong to the ideal
$\mathfrak{q}_M$ either. Therefore, $f_b=r_b'/(\xi' s) \in
A_{\mathfrak{q}_M}, \forall b \in \mathcal{B}$ and hence,
$\sum_{\chi \in C} e_\chi \in A_{\mathfrak{q}_M}\mathcal{H}$. Thus
$C$ is a union of blocks of $A_{\mathfrak{q}_M}\mathcal{H}$. Since
the blocks of $A_{\mathfrak{q}_M}\mathcal{H}$ are unions of blocks
of $A_{\mathfrak{p}A}\mathcal{H}$, we eventually obtain that $C$ is
a block of $A_{\mathfrak{q}_M}\mathcal{H}$.}
\end{apod}

\begin{corollary}\label{not essential for all}
If $M$ is not a $\mathfrak{p}$-essential monomial for any $\chi \in
\mathrm{Irr}(W)$, then the blocks of $A_{\mathfrak{q}_M}\mathcal{H}$
coincide with the blocks of $A_{\mathfrak{p}A}\mathcal{H}$.
\end{corollary}

Of course all the above results hold for $B$ in the place of $A$, if
we further specialize $B$ (and $\mathcal{H}_{\varphi_M}$) via a
morphism associated with a monomial
in $B$ (always assuming that the assumptions $\ref{ypo}$ are satisfied).\\

Now let $r \in \{1,\ldots,m+1\}$ and $R:=\left\{
  \begin{array}{ll}
    \mathbb{Z}_K[y_r^{\pm1},\ldots,y_m^{\pm 1}], & \textrm{for } 1\leq r \leq m; \\
    \mathbb{Z}_K, & \textrm{for } r=m+1,
  \end{array}
\right.$ where $y_r,\ldots,y_m$ are $m-r+1$ indeterminates over
$\mathbb{Z}_K$. We shall recall some definitions given in Chapter 1.

\begin{definition}\label{adapted}
A $\mathbb{Z}_K$-algebra morphism $\varphi:A \rightarrow R$ is
called adapted, if $\varphi={\varphi_r} \circ {\varphi_{r-1}} \circ
\ldots \circ {\varphi_1}$, where $\varphi_i$ is a morphism
associated with a monomial for all $i=1,\ldots,r$. The family
$\mathcal{F}:=\{\varphi_r,\varphi_{r-1},\ldots,\varphi_1\}$ is
called an adapted family for $\varphi$ whose initial morphism is
$\varphi_1$.
\end{definition}

Let $\varphi:A \rightarrow R$ be an adapted morphism and let $F$ be
the field of fractions of $R$. Let us denote by
$\mathcal{H}_\varphi$ the algebra obtained as the specialization of
$\mathcal{H}$ via $\varphi$. Applying Proposition $\ref{associated
morphism preserves splitness}$ $r$ times, we obtain that the algebra
$F\mathcal{H}_\varphi$ is split semisimple. Again, by
``Tits'deformation theorem'', the specialization $y_i \mapsto 1$
induces a bijection from the set
$\mathrm{Irr}(F\mathcal{H}_{\varphi})$ of absolutely irreducible
characters of $F\mathcal{H}_{\varphi}$ to the set $\mathrm{Irr}(W)$.
Therefore, whenever we refer to irreducible characters, we mean
irreducible characters of the group $W$.

We shall repeat here proposition $\ref{change initial}$, proved in
Chapter 1. Recall that if $M:=\prod_{i=0}^mx_i^{b_i}$ is a monomial
such that $\textrm{gcd}(b_i)=d \in \mathbb{Z}$, then
$M^\circ:=\prod_{i=0}^mx_i^{b_i/d}$.

\begin{proposition}\label{apo panw}
Let $\varphi:A \rightarrow R$ be an adapted morphism and $M$ a
monomial in $A$ such that $\varphi(M)=1$. Then there exists an
adapted family for $\varphi$ whose initial morphism is associated
with $M^\circ$.
\end{proposition}

\begin{proposition}\label{blocks of initial monomial}
Let $\varphi:A \rightarrow R$ be an adapted morphism and
$\mathcal{H}_\varphi$ the algebra obtained as the specialization of
$\mathcal{H}$ via $\varphi$. If $M$ is a monomial in $A$ such that
$\varphi(M)=1$ and
$\mathfrak{q}_{M^\circ}:=(M^\circ-1)A+\mathfrak{p}A$, then the
blocks of $R_{\mathfrak{p}R}\mathcal{H}_\varphi$ are unions of
blocks of $A_{\mathfrak{q}_{M^\circ}}\mathcal{H}$.
\end{proposition}
\begin{apod}{Let $M$ be a monomial
in $A$ such that $\varphi(M)=1$. Due to proposition $\ref{apo
panw}$, there exists an adapted family for $\varphi$ whose initial
morphism $\varphi_1$ is associated with $M^\circ$. Let us denote by
$B$ the image of $\varphi_1$ and by $\mathcal{H}_{\varphi_1}$ the
algebra obtained as the specialization of $\mathcal{H}$ via
$\varphi_1$. Thanks to theorem $\ref{Aq B}$, the blocks of
$B_{\mathfrak{p}B}\mathcal{H}_{\varphi_1}$ coincide with the blocks
of $A_{\mathfrak{q}_{M^\circ}}\mathcal{H}$. Now, by corollary
$\ref{union of blocks}$, if two irreducible characters belong to the
same $\mathfrak{p}$-block of a Hecke algebra, then they belong to
the same $\mathfrak{p}$-block of its specialization via a morphism
associated with a monomial. Inductively, we obtain that the blocks
of $R_{\mathfrak{p}R}\mathcal{H}_\varphi$ are unions of blocks of
$B_{\mathfrak{p}B}\mathcal{H}_{\varphi_1}$ and thus of
$A_{\mathfrak{q}_{M^\circ}}\mathcal{H}$.}
\end{apod}

We will now state and prove our main result concerning the
$\mathfrak{p}$-blocks of Hecke algebras. Let us recall the
factorization of the generic Schur element $s_\chi$ associated with
the irreducible character $\chi$ in $\ref{Schur element generic}$.
We have defined a monomial $M:=\prod_{i=0}^mx_i^{a_i}$ with
$\mathrm{gcd}(a_i)=1$ to be $\mathfrak{p}$-essential for $\chi$ if
$s_\chi$ has a factor of the form $\Psi(M)$, where $\Psi$ is a
$K$-cyclotomic polynomial such that $\Psi(1) \in \mathfrak{p}$. We
have seen (proposition $\ref{p-essential}$) that $M$ is
$\mathfrak{p}$-essential for $\chi$ if and only if $s_\chi/\xi_\chi
\in \mathfrak{q}_M$. A monomial of that form is called generally
$\mathfrak{p}$\emph{-essential for} $W$ if it is
$\mathfrak{p}$-essential for some irreducible character $\chi$ of
$W$. We can easily find all $\mathfrak{p}$-essential monomials for
$W$ by looking at the unique factorization of the generic Schur
elements in $K[x_0^{\pm},x_1^{\pm},\ldots,x_m^{\pm}]$.

Let $\varphi:A \rightarrow R$ be an adapted morphism and
$\mathcal{H}_\varphi$ the algebra obtained as the specialization of
$\mathcal{H}$ via $\varphi$. Let $M_1,\ldots,M_k$ be the
$\mathfrak{p}$-essential monomials for $W$ such that $\phi(M_j)=1$
for all $j=1,\ldots,k$. We have $M_j^\circ=M_j$ for all
$j=1,\ldots,k$. Set $\mathfrak{q}_0:=\mathfrak{p}A$,
$\mathfrak{q}_j:=\mathfrak{p}A+(M_j-1)A$ for $j=1,\ldots,k$ and
$\mathcal{Q}:=\{\mathfrak{q}_0,\mathfrak{q}_1,\ldots,\mathfrak{q}_k\}$.

Now let $\mathfrak{q} \in \mathcal{Q}$. If two irreducible
characters $\chi,\psi$ belong to the same block of
$A_\mathfrak{q}\mathcal{H}$, we write $\chi \sim_{\mathfrak{q}}
\psi$.

\begin{theorem}\label{main theorem}
Two irreducible characters $\chi,\psi \in \textrm{\emph{Irr}}(W)$
are in the same block of $R_{\mathfrak{p}R}\mathcal{H}_\varphi$ if
and only if there exist a finite sequence
$\chi_0,\chi_1,\ldots,\chi_n \in \textrm{\emph{Irr}}(W)$ and a
finite sequence $\mathfrak{q}_{j_1},\ldots,\mathfrak{q}_{j_n} \in
\mathcal{Q}$ such that
\begin{itemize}
  \item $\chi_0=\chi$ and $\chi_n=\psi$,
  \item for all $i$ $(1\leq i \leq n)$, $\chi_{i-1} \sim_{\mathfrak{q}_{j_i}}\chi_i$.
\end{itemize}
\end{theorem}
\begin{apod}{Let us denote by $\sim$ the equivalence relation on
$\textrm{Irr}(W)$ defined as the closure of the relation ``there
exists $\mathfrak{q} \in \mathcal{Q}$ such that $\chi
\sim_\mathfrak{q} \psi$''. Therefore, we have to show that $\chi$
and $\psi$ are in the same block of
$R_{\mathfrak{p}R}\mathcal{H}_\varphi$ if and only if $\chi \sim
\psi$.

If $\chi \sim \psi$, then proposition $\ref{blocks of initial
monomial}$ implies that $\chi$ and $\psi$ are in the same block of
$R_{\mathfrak{p}R}\mathcal{H}_\varphi$. Now let $C$ be an
equivalence class of $\sim$. We have that $C$ is a union of blocks
of $A_\mathfrak{q}\mathcal{H}, \forall \mathfrak{q} \in \mathcal{Q}$
and thus
$$\sum_{\theta \in C}\frac{\theta^\vee}{s_\theta} \in
A_\mathfrak{q}\mathcal{H}, \forall \mathfrak{q} \in \mathcal{Q}.$$
If $\mathcal{B},\mathcal{B}'$ are two dual bases of $\mathcal{H}$
with respect to the symmetrizing form $t$, then
$\theta^\vee=\sum_{b\in \mathcal{B}}\theta(b)b'$ and hence
$$\sum_{\theta \in C}\frac{\theta(b)}{s_\theta} \in
A_\mathfrak{q}, \forall \mathfrak{q} \in \mathcal{Q}, \forall b \in
\mathcal{B}.$$ Let us recall the form of the Schur element
$s_\theta$ in $\ref{Schur element generic}$. Set
$\xi_C:=\prod_{\theta \in C}\xi_\theta$ and  $s_C:=\prod_{\theta \in
C}(s_\theta/\xi_\theta)$. Then, for all $b \in \mathcal{B}$, there
exists an element $r_{C,b} \in A$ such that
$$\sum_{\theta \in C}\frac{\theta(b)}{s_\theta} = \frac{r_{C,b}}{\xi_C s_C}.$$

The element $s_C \in A$ is product of terms ($K$-cyclotomic
polynomials taking values on monomials) which are irreducible in
$K[x_0^{\pm 1},x_1^{\pm 1},\ldots,x_m^{\pm 1}]$, due to proposition
$\ref{second irreducible}$. We also have $s_C \notin \mathfrak{p}A$.

Fix $b \in \mathcal{B}$. The ring $K[x_0^{\pm 1},x_1^{\pm
1},\ldots,x_m^{\pm 1}]$ is a unique factorization domain and thus
the quotient $r_{C,b}/s_C$ can be written uniquely in the form
$r/\alpha s$ where
\begin{itemize}
\item $r,s \in A$,
\item $\alpha \in \mathbb{Z}_K$,
\item $s|s_C$ in $A$,
\item $\mathrm{gcd}(r,s)=1$ in $K[x_0^{\pm 1},x_1^{\pm 1},\ldots,x_m^{\pm 1}]$
\end{itemize}
 and for $\xi:= \alpha \xi_C $ we have
$$\frac{r_{C,b}}{\xi_C s_C} = \frac{r}{\xi s}\in
A_\mathfrak{q}, \forall \mathfrak{q} \in \mathcal{Q}.$$ Thus, for
all $\mathfrak{q} \in \mathcal{Q}$, there exist
$r_\mathfrak{q},s_\mathfrak{q} \in A$ with $s_\mathfrak{q} \notin
\mathfrak{q}$ such that
$$\frac{r}{\xi s}=\frac{r_\mathfrak{q}}{s_\mathfrak{q}}.$$

Since $(r,s)=1$, we obtain that $s|s_\mathfrak{q}$ in
$K[x_0^{\pm1},x_1^{\pm 1},\ldots,x_m^{\pm 1}]$. However, $s|s_C$ in
$A$ and thus $s$ is a product of $K$-cyclotomic polynomials taking
values on monomials. Consequently, at least one of the coefficients
of $s$ is a unit in $A$. Corollary $\ref{my second lemma}$ implies
that $s|s_\mathfrak{q}$ in $A$. Therefore, $s \notin \mathfrak{q}$
for all $\mathfrak{q} \in \mathcal{Q}$.

Moreover, we have that $r/\xi \in A_{\mathfrak{q}_0} =
A_{\mathfrak{p}A}$. By corollary $\ref{porisma porismatos}$, there
exist $r' \in A$ and $\xi' \in \mathbb{Z}_K-\mathfrak{p}$ such that
$r/\xi=r'/\xi'$. Then
$$\frac{r}{\xi s}=\frac{r'}{\xi' s}\in
A_\mathfrak{q}, \forall \mathfrak{q} \in \mathcal{Q}.$$

Now let us suppose that $\varphi(\xi' s)=\xi' \varphi(s)$ belongs to
$\mathfrak{p}R$. Since $\xi' \notin \mathfrak{p}$, we must have
$\varphi(s) \in \mathfrak{p}R$. However, the morphism $\varphi$
always sends monomial to monomial. Since $s \notin \mathfrak{p}A$
and $s|s_C$, $s$ must have a factor of the form $\Psi(M)$, where
\begin{itemize}
  \item $M:=\prod_{i=0}^mx_i^{a_i}$ is a monomial in $A$ such that $\mathrm{gcd}(a_i)=1$ and $\varphi(M)=1$,
  \item $\Psi$ is a $K$-cyclotomic polynomial such that $\Psi(1) \in \mathfrak{p}$.
\end{itemize}
Thus $s \in \mathfrak{q}_M:=(M-1)A+\mathfrak{p}A$. Since $\Psi(M) |
s_C$, $M$ is a $\mathfrak{p}$-essential monomial for some
irreducible character $\theta \in C$, \ie $M \in
\{M_1,\ldots,M_k\}$. This contradicts the fact that $s \notin
\mathfrak{q}$ for all $\mathfrak{q} \in \mathcal{Q}$. Therefore,
$\varphi(\xi' s) \notin \mathfrak{p}R$.

So we have
$$\frac{\varphi(r_{C,b})}{\varphi(\xi_C s_C)}=
\frac{\varphi(r')}{\varphi(\xi' s)} \in R_{\mathfrak{p}R}$$ and this
holds for all $b \in \mathcal{B}$. Consequently,
$$\sum_{\theta \in C}\frac{\varphi(\theta^\vee)}{\varphi(s_\theta)} \in
R_{\mathfrak{p}R}\mathcal{H}_\varphi.$$ Thus $C$ is a union of
blocks of $R_{\mathfrak{p}R}\mathcal{H}_\varphi$.}
\end{apod}\\
\begin{remark}
\emph{We can obtain corollary $\ref{not essential for all}$ as an
application of the above theorem for
$\mathcal{Q}=\{\mathfrak{q}_0\}$.}
\end{remark}
\\

To summarize: Theorem $\ref{main theorem}$ allows us to calculate
the blocks of $R_{\mathfrak{p}R}\mathcal{H}_\varphi$ for all adapted
morphisms $\varphi:A \rightarrow R$, if we know the blocks of
$A_{\mathfrak{p}A}\mathcal{H}$ and the blocks of
$A_{\mathfrak{q}_M}\mathcal{H}$ for all $\mathfrak{p}$-essential
monomials $M$. Thus the study of the blocks of the generic Hecke
algebra $\mathcal{H}$ in a finite number of cases suffices to
calculate the
$\mathfrak{p}$-blocks of all Hecke algebras obtained via such specializations.\\

The following result will be only used in the section about
cyclotomic Hecke algebras. However, it is also a result on generic
Hecke algebras, so it will be stated now.

Let $n$ be an integer, $n \neq 0$. Define $I^n:A \rightarrow
A':=\mathbb{Z}_K[y_0^{\pm 1},y_1^{\pm 1},\ldots,y_m^{\pm 1}]$ to be
the $\mathbb{Z}_K$-algebra morphism $x_i \mapsto y_i^n$. Obviously,
$I^n$ is injective.

\begin{lemma}\label{power}
Let $n$ be an integer, $n \neq 0$. Let $I^n:A \rightarrow A'$ be the
$\mathbb{Z}_K$-algebra morphism defined above and $\mathcal{H}'$ the
algebra obtained as the specialization of $\mathcal{H}$ via $I^n$.
Then the blocks of $A'_{\mathfrak{p}A'}\mathcal{H}'$ coincide with
the blocks of $A_{\mathfrak{p}A}\mathcal{H}$.
\end{lemma}
\begin{apod}{Since the map $I^n$ is injective, we can consider $A$ as a subring of $A'$ via the
identification $x_i = y_i^n$ for all $i=0,1,\ldots,m$. By corollary
$\ref{inclusion in localizations}$, we obtain that
$A_{\mathfrak{p}A}$ is contained in $A'_{\mathfrak{p}A'}$ and hence,
the blocks of $A_{\mathfrak{p}A}\mathcal{H}$ are unions of blocks of
$A'_{\mathfrak{p}A'}\mathcal{H}'$.

Now let $C$ be a block of $A'_{\mathfrak{p}A'}\mathcal{H}'$. Since
the field of fractions of $A$ is a splitting field for $\mathcal{H}$
(and thus for $\mathcal{H}'$), we obtain that
$$\sum_{\chi \in C} e_\chi \in (A'_{\mathfrak{p}A'} \cap
K(x_0,x_1,\ldots,x_m))\mathcal{H}'.$$ If $A'_{\mathfrak{p}A'} \cap
K(x_0,x_1,\ldots,x_m)=A_{\mathfrak{p}A}$, then $C$ is also a union
of blocks of $A_{\mathfrak{p}A}\mathcal{H}$ and we obtain the
desired result.

In order to prove that $A'_{\mathfrak{p}A'} \cap
K(x_0,x_1,\ldots,x_m)=A_{\mathfrak{p}A}$, it suffices to show that:
\begin{description}
  \item[(a)] The ring $A_{\mathfrak{p}A}$ is integrally closed.
  \item[(b)] The ring $A'_{\mathfrak{p}A'}$ is integral over
  $A_{\mathfrak{p}A}$.
\end{description}

Since the ring $A$ is integrally closed, part (a) is immediate by
corollary $\ref{integrally closed localization}$ which states that
any localization of an integrally closed ring is also integrally
closed.

For part (b), we have that $A'$ is integral over $A$, since
$y_i^n-x_i=0$, for all $i=0,1,\ldots,m$. Moreover, $A'$ is
integrally closed and thus the integral closure of $A$ in
$K(y_0,y_1,\ldots,y_m)$. The only prime ideal of $A'$ lying over
$\mathfrak{p}A$ is, obviously, $\mathfrak{p}A'$. Following corollary
$\ref{one prime lying over}$, we obtain that the integral closure of
$A_{\mathfrak{p}A}$ in $K(y_0,y_1,\ldots,y_m)$ is
$A'_{\mathfrak{p}A'}$. Thus $A'_{\mathfrak{p}A'}$ is integral over
$A_{\mathfrak{p}A}$.}
\end{apod}

We can consider $I^n$ as an endomorphism of $A$ and denote it by
$I^n_A$. If $k$ is another integer, $k \neq 0$, then ${I^k_A} \circ
I^n_A = {I^n_A}\circ I^k_A = I^{kn}_A$. If now $\varphi:A
\rightarrow R$ is an adapted morphism, we can easily check that
$\varphi \circ I^n_A = {I^n_R} \circ \varphi$. Abusing notation, we
write $\varphi \circ I^n = {I^n} \circ \varphi$.

\begin{corollary}\label{i gives same blocks}
Let $\varphi:A \rightarrow R$ be an adapted morphism and
$\mathcal{H}_\varphi$ the algebra obtained as the specialization of
$\mathcal{H}$ via $\varphi$. Let $\phi: A \rightarrow R$ be a
$\mathbb{Z}_K$-algebra morphism such that ${I^\alpha}\circ
\varphi={I^\beta}\circ \phi$ for some $\alpha,\beta \in
\mathbb{Z}-\{0\}$. If $\mathcal{H}_\phi$ is the algebra obtained as
the specialization of $\mathcal{H}$ via $\phi$, then the blocks of
$R_{\mathfrak{p}R}\mathcal{H}_\phi$ coincide with the blocks of
$R_{\mathfrak{p}R}\mathcal{H}_\varphi$ and we can use theorem
$\ref{main theorem}$ to calculate them.
\end{corollary}

\section{Cyclotomic Hecke algebras}

Let $y$ be an indeterminate. We set $x:=y^{|\mu(K)|}.$

\begin{definition}\label{specialization}
A cyclotomic specialization of $\mathcal{H}$ is a
$\mathbb{Z}_K$-algebra morphism $\phi:
\mathbb{Z}_K[\textbf{\emph{v}},\textbf{\emph{v}}^{-1}]\rightarrow
\mathbb{Z}_K[y,y^{-1}]$ with the following properties:
\begin{itemize}
  \item $\phi: v_{\mathcal{C},j} \mapsto y^{n_{\mathcal{C},j}}$ where
  $n_{\mathcal{C},j} \in \mathbb{Z}$ for all $\mathcal{C}$ and $j$.
  \item For all $\mathcal{C} \in \mathcal{A}/W$, and if $z$ is another
  indeterminate, the element of $\mathbb{Z}_K[y,y^{-1},z]$ defined by
  $$\Gamma_\mathcal{C}(y,z):=\prod_{j=0}^{e_\mathcal{C}-1}(z-\zeta_{e_\mathcal{C}}^jy^{n_{\mathcal{C},j}})$$
  is invariant by the action of $\textrm{\emph{Gal}}(K(y)/K(x))$.
\end{itemize}
\end{definition}

If $\phi$ is a cyclotomic specialization of $\mathcal{H}$,
the corresponding \emph{cyclotomic Hecke algebra} is the
$\mathbb{Z}_K[y,y^{-1}]$-algebra, denoted by $\mathcal{H}_\phi$,
which is obtained as the specialization of the
$\mathbb{Z}_K[\textbf{v},\textbf{v}^{-1}]$-algebra $\mathcal{H}$ via
the morphism $\phi$. It also has a symmetrizing form $t_\phi$
defined as the specialization of the canonical form $t$.\\
\\
\begin{remark} \emph{Sometimes we describe the morphism $\phi$ by the
formula}
$$u_{\mathcal{C},j} \mapsto \zeta_{e_\mathcal{C}}^j x^{n_{\mathcal{C},j}}.$$
\emph{If now we set $q:=\zeta x$ for some root of unity $\zeta \in
\mu(K)$, then the cyclotomic specialization $\phi$ becomes a
$\zeta$-\emph{cyclotomic specialization} and $\mathcal{H}_\phi$ can
be also considered over $\mathbb{Z}_K[q,q^{-1}]$.}
\end{remark}
\begin{px}\label{spetsial}
\small{\emph{The spetsial Hecke algebra $\mathcal{H}_q^s(W)$ is the
1-cyclotomic algebra obtained by the specialization
$$u_{\mathcal{C},0} \mapsto q,\,\, u_{\mathcal{C},j} \mapsto \zeta_{e_\mathcal{C}}^j \textrm{
for } 1 \leq j \leq e_\mathcal{C}-1, \textrm{ for all } \mathcal{C}
\in \mathcal{A}/W.$$ For example, if $W:=G_4$, then
$$\mathcal{H}_q^s(W)=<S,T \textrm{ }|\, STS=TST,
(S-q)(S^2+S+1)=(T-q)(T^2+T+1)=0>.$$}}
\end{px}

Set $A:=\mathbb{Z}_K[\textbf{v},\textbf{v}^{-1}]$ and
$\Omega:=\mathbb{Z}_K[y,y^{-1}]$. Let $\phi: A \rightarrow \Omega$
be a cyclotomic specialization such that $\phi(v_{\mathcal{C},j})=
y^{n_{\mathcal{C},j}}$. Recall that, for $\alpha \in
\mathbb{Z}-\{0\}$, we denote by $I^\alpha:\Omega \rightarrow \Omega$
the monomorphism $y \mapsto y^\alpha$.

\begin{theorem}\label{cyclotomic}
Let $\phi:A \rightarrow \Omega$ be a cyclotomic specialization like
above. Then there exist an adapted $\mathbb{Z}_K$-algebra morphism
$\varphi:A \rightarrow \Omega$ and $\alpha \in \mathbb{Z}-\{0\}$
such that
$$\phi={I^\alpha}\circ\varphi.$$
\end{theorem}
\begin{apod}{We set $d:=\mathrm{gcd}(n_{\mathcal{C},j})$ and consider the
cyclotomic specialization $\varphi: v_{\mathcal{C},j} \mapsto
y^{n_{\mathcal{C},j}/d}$. We have $\phi=I^d \circ\varphi$. Since
$\mathrm{gcd}(n_{\mathcal{C},j}/d)=1$,
there exist $a_{\mathcal{C},j} \in \mathbb{Z}$ such that
$$\sum_{\mathcal{C},j}a_{\mathcal{C},j}(n_{\mathcal{C},j}/d)=1.$$
We have
$y=\varphi(\prod_{\mathcal{C},j}v_{\mathcal{C},j}^{a_{\mathcal{C},j}})$
an hence, $\varphi$ is surjective.  Then, by proposition
$\ref{surjective is adapted}$, $\varphi$ is adapted.}
\end{apod}

Let $\varphi$ be defined as in theorem $\ref{cyclotomic}$ and
$\mathcal{H}_\varphi$ the corresponding cyclotomic Hecke algebra.
Proposition $\ref{associated morphism preserves splitness}$ implies
that the algebra $K(y)\mathcal{H}_\varphi$ is split semisimple. Due
to corollary $\ref{injective preserves splitness}$ and the theorem
above, we deduce that the algebra $K(y)\mathcal{H}_\phi$ is also
split semisimple. For $y=1$ this algebra specializes to the group
algebra $KW$ (the form $t_\phi$ becoming the canonical form on the
group algebra). Thus, by ``Tits' deformation theorem'', the
specialization $v_{\mathcal{C},j} \mapsto 1$ defines the following bijections
$$\begin{array}{ccccc}
    \textrm{Irr}(W) & \leftrightarrow & \textrm{Irr}(K(y)\mathcal{H}_\phi) & \leftrightarrow & \textrm{Irr}(K(\textbf{v})\mathcal{H}) \\
    \chi & \mapsto & \chi_{\phi} & \mapsto & \chi_{\textbf{v}}.
  \end{array}$$

The following result is an immediate consequence of the proposition
$\ref{Schur element generic}$.

\begin{proposition}\label{Schur element cyclotomic}
The Schur element $s_{\chi_\phi}(y)$ associated with the irreducible
character $\chi_\phi$ of $K(y)\mathcal{H}_\phi$ is a Laurent
polynomial in $y$ of the form
$$s_{\chi_\phi}(y)=\psi_{\chi,\phi} y^{a_{\chi,\phi}} \prod_{\Phi \in
C_K}\Phi(y)^{n_{\chi,\phi}}$$ where $\psi_{\chi,\phi} \in
\mathbb{Z}_K$, $a_{\chi,\phi} \in \mathbb{Z}$, $n_{\chi,\phi} \in
\mathbb{N}$ and $C_K$ is a set of $K$-cyclotomic polynomials.
\end{proposition}

Let $\mathfrak{p}$ be a prime ideal of $\mathbb{Z}_K$. Theorem
$\ref{cyclotomic}$ allows us to use theorem $\ref{main theorem}$ for
the calculation of the blocks of
$\Omega_{\mathfrak{p}\Omega}\mathcal{H}_\phi$, since they coincide
with the blocks of $\Omega_{\mathfrak{p}\Omega}\mathcal{H}_\varphi$
by corollary $\ref{i gives same blocks}$. Therefore, we need to know
which $\mathfrak{p}$-essential monomials are sent to 1 by $\phi$.

Let $M:=\prod_{\mathcal{C},j}v_{\mathcal{C},j}^{a_{\mathcal{C},j}}$
be a $\mathfrak{p}$-essential monomial for $W$ in $A$. Then
$$\phi(M)=1 \Leftrightarrow \sum_{\mathcal{C},j}a_{\mathcal{C},j}n_{\mathcal{C},j}=0.$$
Set $m:=\sum_{\mathcal{C}\in \mathcal{A}/W}e_\mathcal{C}$. The
hyperplane defined in $\mathbb{C}^m$ by the relation
$$\sum_{\mathcal{C},j}a_{\mathcal{C},j}t_{\mathcal{C},j}=0,$$ where
$(t_ {\mathcal{C},j})_{ \mathcal{C},j}$ is a set of $m$
indeterminates, is called \emph{$\mathfrak{p}$-essential hyperplane}
for $W$. A hyperplane in $\mathbb{C}^m$ is called \emph{essential}
for $W$, if it is $\mathfrak{p}$-essential for some prime ideal
$\mathfrak{p}$ of $\mathbb{Z}_K$.

In order to calculate the blocks of
$\Omega_{\mathfrak{p}\Omega}\mathcal{H}_\phi$, we check to which
$\mathfrak{p}$-essential hyperplanes the $n_{\mathcal{C},j}$ belong:
\begin{itemize}
  \item If the $n_{\mathcal{C},j}$ belong to no $\mathfrak{p}$-essential hyperplane, then the
blocks of $\Omega_{\mathfrak{p}\Omega}\mathcal{H}_\phi$ coincide
with the blocks of $A_{\mathfrak{p}A}\mathcal{H}$.
  \item If the
 $n_{\mathcal{C},j}$ belong to exactly one $\mathfrak{p}$-essential
 hyperplane, corresponding to the $\mathfrak{p}$-essential monomial
 $M$, then the blocks of
 $\Omega_{\mathfrak{p}\Omega}\mathcal{H}_\phi$ coincide
 with the blocks of $A_{\mathfrak{q}_M}\mathcal{H}$.
  \item If the $n_{\mathcal{C},j}$ belong to more than one $\mathfrak{p}$-essential
hyperplane, then we use theorem $\ref{main theorem}$ to calculate
the blocks of $\Omega_{\mathfrak{p}\Omega}\mathcal{H}_\phi$.
\end{itemize}

If now $n_{\mathcal{C},j}=n \in \mathbb{Z}$ for all $\mathcal{C},j$,
then $\Omega_{\mathfrak{p}\Omega}\mathcal{H}_\phi \simeq
\Omega_{\mathfrak{p}\Omega}W$ and the $n_{\mathcal{C},j}$ belong to
all $\mathfrak{p}$-essential hyperplanes. Due to theorem $\ref{main
theorem}$, we obtain the following proposition

\begin{proposition}\label{group algebra}
Let $\mathfrak{p}$ be a prime ideal of $\mathbb{Z}_K$ lying over a
prime number $p$. If two irreducible characters $\chi$ and $\psi$
are in the same block of
$\Omega_{\mathfrak{p}\Omega}\mathcal{H}_\phi$, then they are in the
same $p$-block of $W$.
\end{proposition}
\begin{apod}{The blocks of $\Omega_{\mathfrak{p}\Omega}{H}_\phi$ are unions of the blocks of
$A_{\mathfrak{q}_M}\mathcal{H}$ for all $\mathfrak{p}$-essential
monomials $M$ such that $\phi(M)=1$, whereas the $p$-blocks of $W$
are unions of the blocks of $A_{\mathfrak{q}_M}\mathcal{H}$ for all
$\mathfrak{p}$-essential monomials $M$.}
\end{apod}\
\\
\begin{remark}
\emph{It is well known that, since the ring
$\Omega_{\mathfrak{p}\Omega}$ is a discrete valuation ring (by
theorem $\ref{Krull-dvr}$), the blocks of
$\Omega_{\mathfrak{p}\Omega}W$ are the $p$-blocks of $W$ as
determined by Brauer theory. The reason is the following:}

\emph{Let $\hat{\Omega}_\mathfrak{p}$ be the $\mathfrak{p}$-adic
completion of $\Omega_\mathfrak{p}$. Then the $p$-blocks of $W$
correspond to the central primitive idempotents of
$\hat{\Omega}_\mathfrak{p}W$. Since $K(y)W$ is a split semisimple
algebra and $\Omega_\mathfrak{p}$ is a local Noetherian ring, by
theorem $\ref{18.4}$, the central primitive idempotents of
$\hat{\Omega}_\mathfrak{p}W$ belong to
$$K(y)W \cap \hat{\Omega}_\mathfrak{p}W=\Omega_\mathfrak{p}W.$$
Since we are working with group algebras, we also have that
$$Z(\hat{\Omega}_\mathfrak{p}W) \cap \Omega_\mathfrak{p}W =
Z(\Omega_\mathfrak{p}W).$$ Thus the central primitive idempotents of
$\hat{\Omega}_\mathfrak{p}W$ coincide with the central primitive
idempotents of $\Omega_\mathfrak{p}W$.}
\end{remark}
\\

However, we know from Brauer theory that if the order of the group
$W$ is prime to $p$, then every character of $W$ is a $p$-block by
itself (see, for example, \cite{Se}, 15.5, Prop.43). It is an
immediate consequence of proposition $\ref{group algebra}$ that

\begin{proposition}\label{p prime to the order of the group}
If $\mathfrak{p}$ is a prime ideal of $\mathbb{Z}_K$ lying over a
prime number $p$ which doesn't divide the order of the group $W$,
then the blocks of $\Omega_{\mathfrak{p}\Omega}\mathcal{H}_\phi$ are
singletons.
\end{proposition}

\section{Rouquier blocks of the cyclotomic Hecke algebras}

\begin{definition}\label{Rouquier ring}
We call Rouquier ring of $K$ and denote by $\mathcal{R}_K(y)$ the
$\mathbb{Z}_K$-subalgebra of $K(y)$
$$\mathcal{R}_K(y):=\mathbb{Z}_K[y,y^{-1},(y^n-1)^{-1}_{n\geq 1}]$$
\end{definition}

Let $\phi: v_{\mathcal{C},j} \mapsto y^{n_{\mathcal{C},j}}$ be a
cyclotomic specialization and $\mathcal{H}_\phi$ the corresponding
cyclotomic Hecke algebra. The \emph{Rouquier blocks} of
$\mathcal{H}_\phi$ are the blocks of the algebra
$\mathcal{R}_K(y)\mathcal{H}_\phi$.\\
\\
\begin{remark}\emph{ It has been shown by Rouquier
\cite{Rou}, that if $W$ is a Weyl group and $\mathcal{H}_\phi$ is
obtained via the ``spetsial'' cyclotomic specialization (see example
$\ref{spetsial}$), then its Rouquier blocks coincide with the
``families of characters'' defined by Lusztig. Thus, the Rouquier
blocks play an essential role in the program ``Spets'' (see
\cite{BMM2}) whose ambition is to give to complex reflection groups
the role of Weyl groups of as yet mysterious structures.}
\end{remark}

\begin{proposition}\label{Some properties of the Rouquier
ring}\emph{(Some properties of the Rouquier ring)}
\begin{enumerate}
  \item The group of units $\mathcal{R}_K(y)^\times$ of the Rouquier ring $\mathcal{R}_K(y)$
  consists of the elements of the form
  $$u y^n \prod_{\Phi \in \mathrm{Cycl}(K)} \Phi(y)^{n_\phi},$$
  where $u \in \mathbb{Z}_K^\times$, $n, n_\phi \in \mathbb{Z}$, $\mathrm{Cycl}(K)$ is the set
  of $K$-cyclotomic polynomials and
  $n_\phi=0$ for all but a finite number of $\Phi$.
  \item The prime ideals of $\mathcal{R}_K(y)$ are
  \begin{itemize}
    \item the zero ideal $\{0\}$,
    \item the ideals of the form $\mathfrak{p}\mathcal{R}_K(y)$,
    where $\mathfrak{p}$ is a prime ideal of $\mathbb{Z}_K$,
    \item the ideals of the form $P(y)\mathcal{R}_K(y)$, where
    $P(y)$ is an irreducible element of $\mathbb{Z}_K[y]$ of degree
    at least $1$, prime to $y$ and to $\Phi(y)$ for all $\Phi \in
    \mathrm{Cycl}(K)$.
  \end{itemize}
  \item The Rouquier ring $\mathcal{R}_K(y)$ is a Dedekind ring.
\end{enumerate}
\end{proposition}
\begin{apod}{
\begin{enumerate}
   \item This part is immediate from the definition of
   $K$-cyclotomic polynomials.
   \item Since $\mathcal{R}_K(y)$ is an integral domain, the zero
   ideal is prime.

   The ring $\mathbb{Z}_K$ is a Dedekind ring
   and thus a Krull ring, by proposition $\ref{case of
   Krull}$. Proposition $\ref{prime ideals of height 1}$ implies that the ring
   $\mathbb{Z}_K[y]$ is also a Krull ring whose prime ideals of height 1
    are of the form $\mathfrak{p}\mathbb{Z}_K[y]$
   ($\mathfrak{p}$ prime in $\mathbb{Z}_K$) and
   $P(y)\mathbb{Z}_K[y]$ ($P(y)$ irreducible in $\mathbb{Z}_K[y]$ of degree at least 1).
   Moreover, $\mathbb{Z}_K$ has an infinite
   number of non-zero prime ideals whose intersection is the zero ideal. Since
   all non-zero prime ideals of $\mathbb{Z}_K$ are maximal, we
   obtain that every prime ideal of $\mathbb{Z}_K$ is the intersection of
   maximal ideals. Thus $\mathbb{Z}_K$ is, by definition, a Jacobson
   ring (cf. \cite{Eis}, \S 4.5). The general form of the
   Nullstellensatz (\cite{Eis}, Thm.4.19) implies that
   for every maximal ideal $\mathfrak{m}$ of $\mathbb{Z}_K[y]$, the
   ideal $\mathfrak{m} \cap \mathbb{Z}_K$ is a maximal ideal of
   $\mathbb{Z}_K$. We deduce that the maximal ideals of
   $\mathbb{Z}_K[y]$ are of the form $\mathfrak{p}\mathbb{Z}_K[y]+P(y)\mathbb{Z}_K[y]$
   ($\mathfrak{p}$ prime in $\mathbb{Z}_K$ and $P(y)$ of degree at least 1 irreducible
   modulo $\mathfrak{p}$). Since
   $\mathbb{Z}_K[y]$ has Krull dimension 2, we have now described
   all its prime ideals.

   The ring $\mathcal{R}_K(y)$ is a localization of $\mathbb{Z}_K[y]$. Therefore,
   in order to prove that the non-zero prime ideals of
   $\mathcal{R}_K(y)$ are the ones described above, it is enough to
   show that $\mathfrak{m}\mathcal{R}_K(y)=\mathcal{R}_K(y)$ for all
   maximal ideals $\mathfrak{m}$ of $\mathbb{Z}_K[y]$. For this, it
   suffices to show that $\mathfrak{p}\mathcal{R}_K(y)$ is a maximal
   ideal of $\mathcal{R}_K(y)$ for all prime ideals $\mathfrak{p}$
   of $\mathbb{Z}_K$.

   Let $\mathfrak{p}$ be a prime ideal of $\mathbb{Z}_K$. Then
   $$\mathcal{R}_K(y)/\mathfrak{p}\mathcal{R}_K(y) \simeq \mathbb{F}_\mathfrak{p}[y,y^{-1},(y^n-1)^{-1}_{n\geq 1}],$$
   where $\mathbb{F}_\mathfrak{p}$ denotes the finite field
   $\mathbb{Z}_K/\mathfrak{p}$. Since $\mathbb{F}_\mathfrak{p}$ is finite,
   every polynomial in $\mathbb{F}_\mathfrak{p}[y]$ is a product of
   elements which divide $y$ or $y^n-1$ for some $n \in \mathbb{N}$.
   Thus every element of $\mathbb{F}_\mathfrak{p}[y]$ is invertible in
   $\mathcal{R}_K(y)/\mathfrak{p}\mathcal{R}_K(y)$. Consequently, we
   obtain that
   $$\mathcal{R}_K(y)/\mathfrak{p}\mathcal{R}_K(y) \simeq \mathbb{F}_\mathfrak{p}(y)$$
   and thus $\mathfrak{p}$ generates a maximal ideal in
   $\mathcal{R}_K(y)$.
   \item  The ring $\mathcal{R}_K(y)$ is the localization of a
   Noetherian integrally closed ring and thus Noetherian and
   integrally closed itself. Moreover, following the description of
   its prime ideals in part 2, it has Krull dimension 1.}
\end{enumerate}
\end{apod}
\begin{remark}
 \emph{If $P(y)$ is an irreducible element of $\mathbb{Z}_K[y]$ of degree
 at least 1, prime to $y$ and to $\Phi(y)$ for all $\Phi \in
 \mathrm{Cycl}(K)$, then the field $\mathcal{R}_K(y)/P(y)\mathcal{R}_K(y)$
 is isomorphic to the field of fractions of the ring
 $\mathbb{Z}_K[y]/P(y)\mathbb{Z}_K[y]$}.
\end{remark}
\\

Now let us recall the form of the Schur elements of the cyclotomic
Hecke algebra $\mathcal{H}_\phi$ given in proposition $\ref{Schur
element cyclotomic}$. If $\chi_\phi$ is an irreducible character of
$K(y)\mathcal{H}_\phi$, then its Schur element $s_{\chi_\phi}(y)$ is
of the form
$$s_{\chi_\phi}(y)=\psi_{\chi,\phi} y^{a_{\chi,\phi}} \prod_{\Phi \in
C_K}\Phi(y)^{n_{\chi,\phi}}$$ where $\psi_{\chi,\phi} \in
\mathbb{Z}_K$, $a_{\chi,\phi} \in \mathbb{Z}$, $n_{\chi,\phi} \in
\mathbb{N}$ and $C_K$ is a set of $K$-cyclotomic polynomials.

\begin{definition}\label{bad}
A prime ideal $\mathfrak{p}$ of $\mathbb{Z}_K$ lying over a prime
number $p$ is $\phi$-bad for $W$, if there exists $\chi_\phi \in
\textrm{\emph{Irr}}(K(y)\mathcal{H}_\phi)$ with $\psi_{\chi,\phi}
\in \mathfrak{p}$. If $\mathfrak{p}$ is $\phi$-bad for $W$, we say
that $p$ is a $\phi$-bad prime number for $W$.
\end{definition}
\begin{remark}\emph{ If $W$ is a Weyl group and $\phi$ is the
``spetsial'' cyclotomic specialization, then the $\phi$-bad prime
ideals are the ideals generated by the bad prime numbers (in the
``usual'' sense) for $W$ (see \cite{GeRo}, 5.2).}
\end{remark}\

Note that if $\mathfrak{p}$ is $\phi$-bad for $W$, then $p$ must
divide the order of the group (since
$s_{\chi_\phi}(1)=|W|/\chi(1)$).\\

Let us denote by $\mathcal{O}$ the Rouquier ring. By proposition
$\ref{p-blocks}$, the Rouquier blocks of $\mathcal{H}_\phi$ are
unions of the blocks of $\mathcal{O}_\mathcal{P}\mathcal{H}_\phi$
for all prime ideals $\mathcal{P}$ of $\mathcal{O}$. However, in all
of the following cases, due to the form of the Schur elements, the
blocks of $\mathcal{O}_\mathcal{P}\mathcal{H}_\phi$ are singletons
(\ie $e_{\chi_\phi}=\chi_\phi^\vee /s_{\chi_\phi} \in
\mathcal{O}_\mathcal{P}\mathcal{H}_\phi$ for all $\chi_\phi \in
\mathrm{Irr}(K(y)\mathcal{H}_\phi)$):
\begin{itemize}
  \item $\mathcal{P}$ is the zero ideal $\{0\}$.
  \item $\mathcal{P}$ is of the form $P(y)\mathcal{O}$, where
$P(y)$ is an irreducible element of $\mathbb{Z}_K[y]$ of degree at
least $1$, prime to $y$ and to $\Phi(y)$ for all $\Phi \in
\mathrm{Cycl}(K)$.
  \item $\mathcal{P}$ is of the form $\mathfrak{p}\mathcal{O}$, where
$\mathfrak{p}$ is a prime ideal of $\mathbb{Z}_K$ which is not
$\phi$-bad for $W$.
\end{itemize}
Therefore, the blocks of $\mathcal{O}\mathcal{H}_\phi$ are, simply,
unions of the blocks of
$\mathcal{O}_{\mathfrak{p}\mathcal{O}}\mathcal{H}_\phi$ for all
$\phi$-bad prime ideals $\mathfrak{p}$ of $\mathbb{Z}_K$. By
proposition $\ref{prime ideals of a localization}$(4), we obtain
that $\mathcal{O}_{\mathfrak{p}\mathcal{O}} \simeq
\Omega_{\mathfrak{p}\Omega}$, where
$\Omega:=\mathbb{Z}_K[y,y^{-1}]$. In the previous section we saw how
we can use theorem $\ref{main theorem}$ to calculate the blocks of
$\Omega_{\mathfrak{p}\Omega}\mathcal{H}_\phi$ and thus obtain
the Rouquier blocks of $\mathcal{H}_\phi$.\\

The following description of the Rouquier blocks results from
proposition $\ref{omega_chi}$ and the description of $\phi$-bad
prime ideals for $W$.

\begin{proposition}\label{Rouquier blocks and central characters}
Let $\chi,\psi \in \emph{Irr}(W)$. The characters $\chi_\phi$ and
$\psi_\phi$ are in the same Rouquier block of $\mathcal{H}_\phi$ if
and only if there exists a finite sequence
$\chi_0,\chi_1,\ldots,\chi_n \in \emph{Irr}(W)$ and a finite
sequence $\mathfrak{p}_1,\ldots,\mathfrak{p}_n$ of $\phi$-bad prime
ideals of $\mathbb{Z}_K$ such that
\begin{itemize}
  \item $(\chi_0)_\phi=\chi_\phi$ and $(\chi_n)_\phi=\psi_\phi$,
  \item for all $j$ $(1\leq j \leq n)$,\,\,
        $\omega_{(\chi_{j-1})_\phi} \equiv \omega_{(\chi_j)_\phi}
        \emph{ mod } \mathfrak{p}_j\mathcal{O}.$
\end{itemize}
\end{proposition}
Following the notations in \cite{BMM2}, 6B, for every element $P(y)
\in \mathbb{C}(y)$, we call
\begin{itemize}
  \item \emph{valuation of $P(y)$ at $y$} and denote by $\mathrm{val}_y(P)$ the order of $P(y)$
  at 0 (we have $\mathrm{val}_y(P)<0$ if 0 is a pole of $P(y)$ and $\mathrm{val}_y(P)>0$ if 0 is a zero of $P(y)$),
  \item \emph{degree of $P(y)$ at $y$} and denote by $\mathrm{deg}_y(P)$ the opposite of the
  valuation of $P(1/y)$.
\end{itemize}
Moreover, if $y^n=x$, then
$$\mathrm{val}_x(P(y)):=\frac{\mathrm{val}_y(P)}{n} \,\textrm{ and
}\, \mathrm{deg}_x(P(y)):=\frac{\mathrm{deg}_y(P)}{n}.$$ For $\chi
\in \mathrm{Irr}(W)$, we define
$$a_{\chi_\phi}:=\mathrm{val}_x(s_{\chi_\phi}(y)) \,\textrm{ and }\,
A_{\chi_\phi}:=\mathrm{deg}_x(s_{\chi_\phi}(y)).$$ The following
result is proven in \cite{BK}, Prop.2.9.

\begin{proposition}\label{aA}
Let $x:=y^{|\mu(K)|}$.
\begin{enumerate}
  \item For all $\chi \in \emph{Irr}(W)$, we have
        $$\omega_{\chi_\phi}(\pi)=t_\phi(\pi)x^{a_{\chi_\phi}+A_{\chi_\phi}},$$
        where $\pi$ is the central element of the pure braid group defined in
        $\ref{pi}$.
  \item Let $\chi,\psi \in \mathrm{Irr}(W)$. If $\chi_\phi$ and
        $\psi_\phi$ belong to the same Rouquier block, then
        $$a_{\chi_\phi}+A_{\chi_\phi}=a_{\psi_\phi}+A_{\psi_\phi}.$$
\end{enumerate}
\end{proposition}
\begin{apod}
{\begin{enumerate}
      \item If $P(y) \in \mathbb{C}[y,y^{-1}]$, we denote by
      $P(y)^*$ the polynomial whose coefficients are the complex
      conjugates of those of $P(y)$. By \cite{BMM2}, 2.8, we know
      that the Schur element $s_{\chi_\phi}(y)$ is semi-palindromic and satisfies
      $$s_{\chi_\phi}(y^{-1})^*=\frac{t_\phi(\pi)}{\omega_{\chi_\phi}(\pi)}s_{\chi_\phi}(y).$$
      We deduce (\cite{BMM2}, 6.5, 6.6) that
      $$\frac{t_\phi(\pi)}{\omega_{\chi_\phi}(\pi)}=\xi
      x^{-(a_{\chi_\phi}+A_{\chi_\phi})},$$
      for some $\xi \in \mathbb{C}$.
      For $y=x=1$, the first equation gives
      $t_\phi(\pi)=\omega_{\chi_\phi}(\pi)$ and the second
      one $\xi=1$. Thus we obtain
      $$\omega_{\chi_\phi}(\pi)=t_\phi(\pi)x^{a_{\chi_\phi}+A_{\chi_\phi}}.$$
      \item Suppose that $\chi_\phi$ and $\psi_\phi$ belong to the same Rouquier
      block. Due to proposition $\ref{Rouquier blocks and central
      characters}$, it is enough to show that if there exists a
      $\phi$-bad prime ideal $\mathfrak{p}$ of $\mathbb{Z}_K$ such that
      $\omega_{\chi_\phi} \equiv \omega_{\psi_\phi}
      \textrm{ mod } \mathfrak{p}\mathcal{O}$,
      then $a_{\chi_\phi}+A_{\chi_\phi}=a_{\psi_\phi}+A_{\psi_\phi}.$
      If $\omega_{\chi_\phi} \equiv \omega_{\psi_\phi}
      \textrm{ mod } \mathfrak{p}\mathcal{O}$, then, in
      particular, $\omega_{\chi_\phi}(\pi) \equiv
      \omega_{\psi_\phi}(\pi)
      \textrm{ mod } \mathfrak{p}\mathcal{O}$. Part 1 implies
      that
      $$t_\phi(\pi)x^{a_{\chi_\phi}+A_{\chi_\phi}} \equiv
      t_\phi(\pi)x^{a_{\psi_\phi}+A_{\psi_\phi}} \textrm{ mod } \mathfrak{p}\mathcal{O}.$$
      We know by \cite{BMM2}, 2.1 that $t_\phi(\pi)$ is of the
      form $\xi x^M$, where $\xi$ is a root of unity and $M \in
      \mathbb{Z}$. Thus $t_\phi(\pi) \notin
      \mathfrak{p}\mathcal{O}$ and the above congruence gives
      $$x^{a_{\chi_\phi}+A_{\chi_\phi}} \equiv
        x^{a_{\psi_\phi}+A_{\psi_\phi}} \textrm{ mod } \mathfrak{p}\mathcal{O},$$
      whence
      $$a_{\chi_\phi}+A_{\chi_\phi}=a_{\psi_\phi}+A_{\psi_\phi}.$$}
\end{enumerate}
\end{apod}
\begin{remark}\emph{ For all Coxeter groups, Kazhdan-Lusztig theory states that if $\chi_\phi$ and
$\psi_\phi$ belong to the same Rouquier block, then
$a_{\chi_\phi}=a_{\psi_\phi}$ and $A_{\chi_\phi}=A_{\psi_\phi}$ (cf.
\cite{Lu82}). The same assertion has been proven
\begin{itemize}
  \item for the imprimitive complex reflection groups in \cite{BK}.
  \item for the ``spetsial'' complex reflection groups in
  \cite{MaRo}.
\end{itemize}
The results of the next chapter prove that it holds for the groups
$G_{12}$, $G_{22}$ and $G_{31}$. We conjecture that it is true for
all the remaining exceptional complex reflection groups (\ie the
groups $G_5,G_7,G_9,G_{10},G_{11},G_{13},G_{15\ldots21}$).}
\end{remark}

\chapter{Rouquier blocks of the cyclotomic Hecke algebras of
the exceptional complex reflection groups}

All the notations used in this chapter have been explained in
Chapter 3.

\section{General principles}

Let $W$ be a complex reflection group such that the assumptions
$\ref{ypo}$ are satisfied. Let $\mathcal{H}$ be its generic Hecke
algebra defined over the ring
$A:=\mathbb{Z}_K[\textbf{v},\textbf{v}^{-1}]$. Let $\mathfrak{p}$ be
a prime ideal of $\mathbb{Z}_K$ lying over a prime number $p$ which
divides the order of the group $W$. We can determine the
$\mathfrak{p}$-essential hyperplanes for each character $\chi \in
\mathrm{Irr}(W)$ by looking at the factorization of its generic
Schur element in $K[\textbf{v},\textbf{v}^{-1}]$ (see theorem
$\ref{Schur element generic}$).

Let $\phi:v_{\mathcal{C},j} \mapsto y^{n_{\mathcal{C},j}}$ be a
cyclotomic specialization and let $\mathcal{H}_\phi$ be the
cyclotomic Hecke algebra obtained by $\mathcal{H}$ via $\phi$. Let
us denote by $\mathcal{O}$ the Rouquier ring. We can distinguish
three cases.
\begin{itemize}
  \item If the $n_{\mathcal{C},j}$ belong to no
  $\mathfrak{p}$-essential hyperplane, then the blocks of
  $\mathcal{O}_{\mathfrak{p}\mathcal{O}}\mathcal{H}_\phi$ coincide with the
  blocks of $A_{\mathfrak{p}A}\mathcal{H}$.
  \item If the $n_{\mathcal{C},j}$ belong to exactly one
  $\mathfrak{p}$-essential hyperplane, corresponding to the $\mathfrak{p}$-essential monomial
  $M$, then the blocks of $\mathcal{O}_{\mathfrak{p}\mathcal{O}}\mathcal{H}_\phi$ coincide with the
  blocks of $A_{\mathfrak{q}_M}\mathcal{H}$, where
  $\mathfrak{q}_M:=\mathfrak{p}A+(M-1)A$.
  \item If the $n_{\mathcal{C},j}$ belong to more than one
  $\mathfrak{p}$-essential hyperplane, then we use theorem $\ref{main
  theorem}$ in order to calculate the blocks of
  $\mathcal{O}_{\mathfrak{p}\mathcal{O}}\mathcal{H}_\phi$.
\end{itemize}

Now recall that the Rouquier blocks of $\mathcal{H}_\phi$ are unions
of the blocks of
$\mathcal{O}_{\mathfrak{p}\mathcal{O}}\mathcal{H}_\phi$ for all
$\phi$-bad prime ideals $\mathfrak{p}$ of $\mathbb{Z}_K$. We
distinguish  again three cases.
\begin{itemize}
  \item If the $n_{\mathcal{C},j}$ belong to no
  essential hyperplane, then the Rouquier blocks of
  $\mathcal{H}_\phi$ are unions of the blocks of
  $A_{\mathfrak{p}A}\mathcal{H}$ for all $\phi$-bad prime ideals
  $\mathfrak{p}$. We say that these are \emph{the Rouquier blocks
  associated with no essential hyperplane for} $W$.
  \item If the $n_{\mathcal{C},j}$ belong to exactly one
  essential hyperplane, corresponding to the essential monomial
  $M$, then the Rouquier blocks of $\mathcal{H}_\phi$
  are unions of the blocks of $A_{\mathfrak{q}_M}\mathcal{H}$, where
  $\mathfrak{q}_M:=\mathfrak{p}A+(M-1)A$, for all $\phi$-bad prime ideals
  $\mathfrak{p}$ (if $M$ is not $\mathfrak{p}$-essential,
  then, by corollary $\ref{not essential for all}$, the blocks of $A_{\mathfrak{q}_M}\mathcal{H}$
  coincide with the blocks of $A_{\mathfrak{p}A}\mathcal{H}$).
  We say that these are \emph{the Rouquier blocks
  associated with that essential hyperplane}.
  \item If the $n_{\mathcal{C},j}$ belong to more than one
  essential hyperplane, then the Rouquier blocks of
  $\mathcal{H}_\phi$ are unions of the Rouquier blocks
  associated with the essential hyperplanes to which the $n_{\mathcal{C},j}$ belong.
\end{itemize}

Therefore, if we know the blocks of $A_{\mathfrak{p}A}\mathcal{H}$
and $A_{\mathfrak{q}_M}\mathcal{H}$ for all $\mathfrak{p}$-essential
monomials $M$, for all $\mathfrak{p}$, we know the Rouquier blocks
of $\mathcal{H}_\phi$ for any cyclotomic specialization $\phi$.

In order to calculate the blocks of $A_{\mathfrak{p}A}\mathcal{H}$
(resp. of $A_{\mathfrak{q}_M}\mathcal{H}$), we find a cyclotomic
specialization $\phi:v_{\mathcal{C},j} \mapsto
y^{n_{\mathcal{C},j}}$ such that the $n_{\mathcal{C},j}$ belong to
no $\mathfrak{p}$-essential hyperplane (resp. the
$n_{\mathcal{C},j}$ belong to the $\mathfrak{p}$-essential
hyperplane corresponding to $M$ and no other) and we calculate the
blocks of $\mathcal{O}_{\mathfrak{p}\mathcal{O}}\mathcal{H}_\phi$.

The algorithm presented in the next section uses some theorems
proved in previous chapters in order to form a partition of
$\mathrm{Irr}(W)$ into sets which are unions of blocks of
$\mathcal{O}_{\mathfrak{p}\mathcal{O}}\mathcal{H}_\phi$. These
theorems are
\begin{description}
  \item[$\ref{Malle-Rouquier}$]An irreducible character $\chi$ is a block by
  itself in $\mathcal{O}_{\mathfrak{p}\mathcal{O}}\mathcal{H}_\phi$
  if and only if $\textrm{ }\textrm{ }s_{\chi_\phi} \notin \mathfrak{p}\mathbb{Z}_K[y,y^{-1}]$.
  \item[$\ref{group algebra}$ ]If $\chi,\psi$ belong to the same block of
  $\mathcal{O}_{\mathfrak{p}\mathcal{O}}\mathcal{H}_\phi$, then they are in the
  same $\textrm{ }p$-block of $W$.
  \item[$\ref{aA}$ ]If $\chi,\psi$ are in the same block of
  $\mathcal{O}_{\mathfrak{p}\mathcal{O}}\mathcal{H}_\phi$, then
  $a_{\chi_\phi}+A_{\chi_\phi}=a_{\psi_\phi}+A_{\psi_\phi}.$
  \item[$\ref{not essential for block}$]Let $C$ be a block of $A_{\mathfrak{p}A}\mathcal{H}$. If $M$ is not
  a $\mathfrak{p}$-essential monomial for any $\textrm{ }\chi \in C$ , then $C$
  is a block of $A_{\mathfrak{q}_M}\mathcal{H}$.
\end{description}
If the partition obtained is minimal, then it represents the blocks
of $\mathcal{O}_{\mathfrak{p}\mathcal{O}}\mathcal{H}_\phi$.\\

With the help of the package CHEVIE of GAP, we created a program
that follows this algorithm to obtain the Rouquier blocks of all
cyclotomic Hecke algebras of the groups $G_7$, $G_{11}$, $G_{19}$,
$G_{26}$, $G_{28}$ and $G_{32}$. We used Clifford theory (for more
details, see Appendix) in order to obtain the Rouquier blocks for
\begin{itemize}
  \item $G_4$, $G_5$, $G_6$ from $G_7$,
  \item $G_8$, $G_9$, $G_{10}$, $G_{12}$, $G_{13}$, $G_{14}$, $G_{15}$ from $G_{11}$,
  \item $G_{16}$, $G_{17}$, $G_{18}$, $G_{20}$, $G_{21}$, $G_{22}$ from $G_{19}$,
  \item $G_{25}$ from $G_{26}$.
\end{itemize}
In all of the above cases (except for the ``spetsial'' case for
$G_{32}$), we can determine that the partition into
$\mathfrak{p}$-blocks obtained is minimal. This is done either by
using again the above theorems or by applying the results of
Clifford theory.

For all remaining groups, the Rouquier blocks of the ``spetsial''
cyclotomic Hecke algebra have already been calculated in \cite{MaRo}
(along with those of $G_{32}$), where more criteria for the
partition of $\mathrm{Irr}(W)$ into $\mathfrak{p}$-blocks are given.
Since they are groups generated by reflections of order 2 whose
reflecting hyperplanes belong to one single orbit, their generic
Hecke algebras are defined over a ring of the form
$\mathbb{Z}[x_0^\pm,x_1^\pm]$ and the only essential monomial is
$x_0x_1^{-1}$. If $x_i \mapsto y^{a_i}$ is a cyclotomic
specialization and $a_0=a_1$, then the specialized algebra is the
group algebra, whose Rouquier blocks are known (they are unions of
the group's $p$-blocks for all primes $p$ dividing the order of the
group). According to the above algorithm, it is enough to study one
case where $a_0 \neq a_1$ and thus the ``spetsial'' case covers our
needs.
\section{Algorithm}

Let $\mathfrak{p}$ be a prime ideal of $\mathbb{Z}_K$ lying over a
prime number $p$ which divides the order of the group $W$. As we saw
in the previous section, we need to calculate the blocks of
$A_{\mathfrak{p}A}\mathcal{H}$ and the blocks of
$A_{\mathfrak{q}_M}\mathcal{H}$ for all $\mathfrak{p}$-essential
monomials $M$.

Together with Jean Michel, we have programmed into GAP the
factorized generic Schur elements for all exceptional complex
reflection groups, verifying thus theorem $\ref{Schur element
generic}$. These data have been stored under the name ``SchurData''
and correspond to the following presentation of the Schur element of
an irreducible character $\chi$:
$$s_\chi=\xi_\chi N_\chi \prod_{i \in I_\chi} \Psi_{\chi,i}(M_{\chi,i})^{n_{\chi,i}}\,\,\,\,\,(\dag)$$
(for the notations, the reader should refer to theorem $\ref{Schur
element generic}$). Firstly, we determine the
$\mathfrak{p}$-essential monomials/$\mathfrak{p}$-essential
hyperplanes for $W$. Given the prime ideal $\mathfrak{p}$, GAP
provides us with a way to determine whether an element of
$\mathbb{Z}_K$ belongs to $\mathfrak{p}$. In the above formula, if
$\Psi_{\chi,i}(1) \in \mathfrak{p}$, then $M_{\chi,i}$ is a
$\mathfrak{p}$-essential monomial.

If now we are interested in calculating the blocks of
$A_{\mathfrak{p}A}\mathcal{H}$, we follow the steps below:
\begin{enumerate}
  \item We select the characters $\chi \in \mathrm{Irr}(W)$ whose Schur element
  has its coefficient $\xi_\chi$ in $\mathfrak{p}$.
  The remaining ones will be blocks of $A_{\mathfrak{p}A}\mathcal{H}$ by themselves, thanks to proposition
  $\ref{Malle-Rouquier}$. Thus we form a first partition $\lambda_1$
  of $\mathrm{Irr}(W)$; one part formed by the selected characters,
  each remaining character forming a part by itself.
  \item We calculate the $p$-blocks of $W$.
  By proposition $\ref{group algebra}$, if two irreducible characters aren't
  in the same $p$-block of $W$, then they can not be in the same block of $A_{\mathfrak{p}A}\mathcal{H}$.
  We intersect the partition $\lambda_1$ with the partition obtained
  by the $p$-blocks of $W$ and we obtain a finer partition, named
  $\lambda_2$.
  \item We find a cyclotomic specialization $\phi:v_{\mathcal{C},j} \mapsto y^{n_{\mathcal{C},j}}$
  such that the $n_{\mathcal{C},j}$ belong to no
  $\mathfrak{p}$-essential hyperplane. This is done by trying and
  checking random values for the $n_{\mathcal{C},j}$.
  The blocks of $A_{\mathfrak{p}A}\mathcal{H}$ coincide with the
  blocks of
  $\mathcal{O}_{\mathfrak{p}\mathcal{O}}\mathcal{H}_\phi$.
  Following proposition $\ref{aA}$, we take the intersection of the
  partition we already have with the subsets of $\mathrm{Irr}(W)$,
  where the sum $a_{\chi_\phi}+A_{\chi_\phi}$ remains constant. This
  procedure is repeated several times, because sometimes the
  partition becomes finer after some repetitions. Finally, we obtain
  the partition $\lambda_3$, which is the finest of all.
\end{enumerate}

If we are interested in calculating the blocks of
$A_{\mathfrak{q}_M}\mathcal{H}$ for some $\mathfrak{p}$-essential
monomial $M$, the procedure is more or less the same:
\begin{enumerate}
  \item We select the characters $\chi \in \mathrm{Irr}(W)$ for which $M$ is a
  $\mathfrak{p}$-essential monomial.
  We form a first partition $\lambda_1$ of $\mathrm{Irr}(W)$; one part formed by the selected characters,
  each remaining character forming a part by itself. The idea is
  that, by proposition $\ref{not essential for block}$, if $M$ is not
  $\mathfrak{p}$-essential for any character in a block $C$ of
  $A_{\mathfrak{p}A}\mathcal{H}$, then $C$ is a block of
  $A_{\mathfrak{q}_M}\mathcal{H}$. This explains step 4.
  \item We calculate the $p$-blocks of $W$.
  By proposition $\ref{group algebra}$, if two irreducible characters aren't
  in the same $p$-block of $W$, then they can not be in the same block of $A_{\mathfrak{q}_M}\mathcal{H}$.
  We intersect the partition $\lambda_1$ with the partition obtained
  by the $p$-blocks of $W$ and we obtain a finer partition, named
  $\lambda_2$.
  \item We find a cyclotomic specialization $\phi:v_{\mathcal{C},j} \mapsto y^{n_{\mathcal{C},j}}$
  such that the $n_{\mathcal{C},j}$ belong to the
  $\mathfrak{p}$-essential hyperplane defined by $M$ and to no other, (again by trying and
  checking random values for the $n_{\mathcal{C},j}$). We repeat the
  third step as described for $A_{\mathfrak{p}A}\mathcal{H}$ to obtain partition $\lambda_3$.
  \item We take the union of $\lambda_3$ and the partition defined by the blocks
  of $A_{\mathfrak{p}A}\mathcal{H}$.
\end{enumerate}

The above algorithm is, due to step 3, heuristic. However, as we
have said in the previous section, for the cases we have used it
($G_7$, $G_{11}$, $G_{19}$, $G_{26}$, $G_{28}$, $G_{32}$), we have
been able to determine (using the criteria explained also in the
previous section) that the partition obtained at the end is minimal
and corresponds to the blocks we are looking for.\\

The Rouquier blocks associated with no essential hyperplane (resp.
with  the essential hyperplane corresponding to some essential
monomial $M$) are unions of the blocks of
$A_{\mathfrak{p}A}\mathcal{H}$ (resp. of
$A_{\mathfrak{q}_M}\mathcal{H}$) for all $\mathfrak{p}$ lying over
primes which divide the order of the group $W$. We have observed
that the above algorithm provides us with the correct Rouquier
blocks for all exceptional complex reflection groups in all cases,
except for the ``spetsial'' case of $G_{34}$.

\section{Results}

Using the algorithm of the previous section, we have been able to
calculate the Rouquier blocks associated with all essential
hyperplanes for all exceptional complex reflection groups.

We will give here the example of $G_7$ and show how we obtain the
blocks of $G_6$ from those of $G_7$. Nevertheless, let us first
explain the notations of characters used by the CHEVIE package.

Let $W$ be an exceptional irreducible complex reflection group. For
$\chi \in \mathrm{Irr}(W)$, we set $d(\chi):=\chi(1)$ and we denote
by $b(\chi)$ the valuation of the fake degree of $\chi$ (for the
definition of the fake degree see \cite{Brou}, 1.20). The
irreducible characters $\chi$ of $W$ are determined by the
corresponding pairs $(d(\chi),b(\chi))$ and we write
$\chi=\phi_{d,b}$, where $d:=d(\chi)$ and $b:=b(\chi)$. If two
irreducible characters $\chi$ and $\chi'$ have $d(\chi)=d(\chi')$
and $b(\chi)=b(\chi')$, we use primes `` $'$ '' to distinguish them
(following \cite{Ma3},\cite{MaRo}).

\begin{px}\label{example g7}
\emph{\small{The generic Hecke algebra of $G_7$ is
$$\begin{array}{rccl}
    \mathcal{H}(G_7) & = & <S,T,U \,\,| &  STU=TUS=UST  \\
     &  &  & (S-x_0)(S-x_1)=0 \\
     &  &  & (T-y_0)(T-y_1)(T-y_2)=0 \\
     &  &  & (U-z_0)(U-z_1)(U-z_2)=0>
  \end{array}$$
Let $\phi$ be a cyclotomic specialization of $\mathcal{H}(G_7)$ with
$$\phi(x_i)=\zeta_2^i x^{a_i} ,
\phi(y_j)=\zeta_3^j x^{b_j}, \phi(z_k)=\zeta_3^k x^{c_k}.$$ The
$\phi$-bad prime numbers are 2 and 3. We will now give all the
essential hyperplanes for $G_7$ and the non-trivial Rouquier blocks
associated with each.
\begin{description}\scriptsize
\item[No essential hyperplane]\hfil\break
$\{\phi_{2,9}',\phi_{2,15}\}$, $\{\phi_{2,7}',\phi_{2,13}'\}$,
$\{\phi_{2, 11}',\phi_{2,5}'\}$, $\{\phi_{2,7}'',\phi_{2,13}''\}$,
$\{\phi_{2,11}'', \phi_{2,5}''\}$, $\{\phi_{2,9}'',\phi_{2,3}'\}$,
$\{\phi_{2,11}''',\phi_{2, 5}'''\}$,
$\{\phi_{2,9}''',\phi_{2,3}''\}$, $\{\phi_{2,7}''',\phi_{2,1}\}$,
 $\{\phi_{3,6},\phi_{3,10},\phi_{3,2}\}$, $\{\phi_{3,4},\phi_{3,8},\phi_{3,
12}\}$\item[$c_1-c_2=0$]\hfil\break $\{\phi_{1,4}',\phi_{1,8}'\}$,
$\{\phi_{1,8}'',\phi_{1,12}'\}$, $\{\phi_{1, 12}'',\phi_{1,16}\}$,
$\{\phi_{1,10}',\phi_{1,14}'\}$, $\{\phi_{1,14}'', \phi_{1,18}'\}$,
$\{\phi_{1,18}'',\phi_{1,22}\}$, $\{\phi_{2,9}',\phi_{2, 15}\}$,
$\{\phi_{2,7}',\phi_{2,11}',\phi_{2,13}',\phi_{2,5}'\}$, $\{\phi_{2,
7}'',\phi_{2,13}''\}$,
$\{\phi_{2,11}'',\phi_{2,9}'',\phi_{2,5}'',\phi_{2, 3}'\}$,
$\{\phi_{2,11}''',\phi_{2,5}'''\}$, $\{\phi_{2,9}''',\phi_{2,7}''',
\phi_{2,3}'',\phi_{2,1}\}$, $\{\phi_{3,6},\phi_{3,10},\phi_{3,2}\}$,
 $\{\phi_{3,4},\phi_{3,8},\phi_{3,12}\}$\item[$c_0-c_1=0$]\hfil\break
$\{\phi_{1,0},\phi_{1,4}'\}$, $\{\phi_{1,4}'',\phi_{1,8}''\}$,
$\{\phi_{1, 8}''',\phi_{1,12}''\}$, $\{\phi_{1,6},\phi_{1,10}'\}$,
$\{\phi_{1,10}'', \phi_{1,14}''\}$,
$\{\phi_{1,14}''',\phi_{1,18}''\}$,\\ $\{\phi_{2,9}',\phi_{2,
7}',\phi_{2,15},\phi_{2,13}'\}$, $\{\phi_{2,11}',\phi_{2,5}'\}$,
$\{\phi_{2, 7}'',\phi_{2,11}'',\phi_{2,13}'',\phi_{2,5}''\}$,
$\{\phi_{2,9}'',\phi_{2, 3}'\}$,\\
$\{\phi_{2,11}''',\phi_{2,9}''',\phi_{2,5}''',\phi_{2,3}''\}$,
 $\{\phi_{2,7}''',\phi_{2,1}\}$, $\{\phi_{3,6},\phi_{3,10},\phi_{3,2}\}$,
 $\{\phi_{3,4},\phi_{3,8},\phi_{3,12}\}$\item[$c_0-c_2=0$]\hfil\break
$\{\phi_{1,0},\phi_{1,8}'\}$, $\{\phi_{1,4}'',\phi_{1,12}'\}$,
$\{\phi_{1, 8}''',\phi_{1,16}\}$, $\{\phi_{1,6},\phi_{1,14}'\}$,
$\{\phi_{1,10}'',\phi_{1, 18}'\}$,
$\{\phi_{1,14}''',\phi_{1,22}\}$,\\ $\{\phi_{2,9}',\phi_{2,11}',
\phi_{2,15},\phi_{2,5}'\}$, $\{\phi_{2,7}',\phi_{2,13}'\}$,
$\{\phi_{2,7}'', \phi_{2,9}'',\phi_{2,13}'',\phi_{2,3}'\}$,
$\{\phi_{2,11}'',\phi_{2,5}''\}$,\\
 $\{\phi_{2,11}''',\phi_{2,7}''',\phi_{2,5}''',\phi_{2,1}\}$, $\{\phi_{2,
9}''',\phi_{2,3}''\}$, $\{\phi_{3,6},\phi_{3,10},\phi_{3,2}\}$,
$\{\phi_{3,4},
\phi_{3,8},\phi_{3,12}\}$\item[$b_1-b_2=0$]\hfil\break
$\{\phi_{1,4}'',\phi_{1,8}'''\}$, $\{\phi_{1,8}'',\phi_{1,12}''\}$,
 $\{\phi_{1,12}',\phi_{1,16}\}$, $\{\phi_{1,10}'',\phi_{1,14}'''\}$,
 $\{\phi_{1,14}'',\phi_{1,18}''\}$,\\ $\{\phi_{1,18}',\phi_{1,22}\}$,
 $\{\phi_{2,9}',\phi_{2,15}\}$, $\{\phi_{2,7}',\phi_{2,13}'\}$, $\{\phi_{2,
11}',\phi_{2,5}'\}$,
$\{\phi_{2,7}'',\phi_{2,11}''',\phi_{2,13}'',\phi_{2, 5}'''\}$,\\
$\{\phi_{2,11}'',\phi_{2,9}''',\phi_{2,5}'',\phi_{2,3}''\}$,
 $\{\phi_{2,9}'',\phi_{2,7}''',\phi_{2,3}',\phi_{2,1}\}$, $\{\phi_{3,6},
\phi_{3,10},\phi_{3,2}\}$,
$\{\phi_{3,4},\phi_{3,8},\phi_{3,12}\}$\item[$b_0-b_1=0$]\hfil\break
$\{\phi_{1,0},\phi_{1,4}''\}$, $\{\phi_{1,4}',\phi_{1,8}''\}$,
$\{\phi_{1,8}', \phi_{1,12}'\}$, $\{\phi_{1,6},\phi_{1,10}''\}$,
$\{\phi_{1,10}',\phi_{1, 14}''\}$,
$\{\phi_{1,14}',\phi_{1,18}'\}$,\\
$\{\phi_{2,9}',\phi_{2,7}'', \phi_{2,15},\phi_{2,13}''\}$,
$\{\phi_{2,7}',\phi_{2,11}'',\phi_{2,13}', \phi_{2,5}''\}$,
$\{\phi_{2,11}',\phi_{2,9}'',\phi_{2,5}',\phi_{2,3}'\}$,\\
 $\{\phi_{2,11}''',\phi_{2,5}'''\}$, $\{\phi_{2,9}''',\phi_{2,3}''\}$,
 $\{\phi_{2,7}''',\phi_{2,1}\}$, $\{\phi_{3,6},\phi_{3,10},\phi_{3,2}\}$,
 $\{\phi_{3,4},\phi_{3,8},\phi_{3,12}\}$\item[$b_0-b_2=0$]\hfil\break
$\{\phi_{1,0},\phi_{1,8}'''\}$, $\{\phi_{1,4}',\phi_{1,12}''\}$,
$\{\phi_{1, 8}',\phi_{1,16}\}$, $\{\phi_{1,6},\phi_{1,14}'''\}$,
$\{\phi_{1,10}',\phi_{1, 18}''\}$, $\{\phi_{1,14}',\phi_{1,22}\}$,\\
$\{\phi_{2,9}',\phi_{2,11}''', \phi_{2,15},\phi_{2,5}'''\}$,
$\{\phi_{2,7}',\phi_{2,9}''',\phi_{2,13}', \phi_{2,3}''\}$,
$\{\phi_{2,11}',\phi_{2,7}''',\phi_{2,5}',\phi_{2,1}\}$,
 $\{\phi_{2,7}'',\phi_{2,13}''\}$,\\ $\{\phi_{2,11}'',\phi_{2,5}''\}$,
 $\{\phi_{2,9}'',\phi_{2,3}'\}$, $\{\phi_{3,6},\phi_{3,10},\phi_{3,2}\}$,
$\{\phi_{3,4},\phi_{3,8},\phi_{3,12}\}$\item[$a_0-a_1-2b_0+b_1+b_2-2c_0+c_1+c_2=0$]\hfil\break
$\{\phi_{1,6},\phi_{2,9}',\phi_{2,15},\phi_{3,4},\phi_{3,8},\phi_{3,12}\}$,
 $\{\phi_{2,7}',\phi_{2,13}'\}$, $\{\phi_{2,11}',\phi_{2,5}'\}$, $\{\phi_{2,
7}'',\phi_{2,13}''\}$, $\{\phi_{2,11}'',\phi_{2,5}''\}$,\\
$\{\phi_{2,9}'', \phi_{2,3}'\}$, $\{\phi_{2,11}''',\phi_{2,5}'''\}$,
$\{\phi_{2,9}''',\phi_{2, 3}''\}$, $\{\phi_{2,7}''',\phi_{2,1}\}$,
$\{\phi_{3,6},\phi_{3,10},\phi_{3,
2}\}$\item[$a_0-a_1-2b_0+b_1+b_2+c_0-2c_1+c_2=0$]\hfil\break
$\{\phi_{1,10}',\phi_{2,7}',\phi_{2,13}',\phi_{3,4},\phi_{3,8},\phi_{3,12}\}$,
 $\{\phi_{2,9}',\phi_{2,15}\}$, $\{\phi_{2,11}',\phi_{2,5}'\}$, $\{\phi_{2,
7}'',\phi_{2,13}''\}$, $\{\phi_{2,11}'',\phi_{2,5}''\}$,\\
$\{\phi_{2,9}'', \phi_{2,3}'\}$, $\{\phi_{2,11}''',\phi_{2,5}'''\}$,
$\{\phi_{2,9}''',\phi_{2, 3}''\}$, $\{\phi_{2,7}''',\phi_{2,1}\}$,
$\{\phi_{3,6},\phi_{3,10},\phi_{3,
2}\}$\item[$a_0-a_1-2b_0+b_1+b_2+c_0+c_1-2c_2=0$]\hfil\break
$\{\phi_{1,14}',\phi_{2,11}',\phi_{2,5}',\phi_{3,4},\phi_{3,8},\phi_{3,12}\}$,
 $\{\phi_{2,9}',\phi_{2,15}\}$, $\{\phi_{2,7}',\phi_{2,13}'\}$, $\{\phi_{2,
7}'',\phi_{2,13}''\}$, $\{\phi_{2,11}'',\phi_{2,5}''\}$,\\
$\{\phi_{2,9}'', \phi_{2,3}'\}$, $\{\phi_{2,11}''',\phi_{2,5}'''\}$,
$\{\phi_{2,9}''',\phi_{2, 3}''\}$, $\{\phi_{2,7}''',\phi_{2,1}\}$,
$\{\phi_{3,6},\phi_{3,10},\phi_{3,
2}\}$\item[$a_0-a_1-b_0-b_1+2b_2-c_0-c_1+2c_2=0$]\hfil\break
$\{\phi_{1,16},\phi_{2,7}''',\phi_{2,1},\phi_{3,6},\phi_{3,10},\phi_{3,2}\}$,
 $\{\phi_{2,9}',\phi_{2,15}\}$, $\{\phi_{2,7}',\phi_{2,13}'\}$, $\{\phi_{2,
11}',\phi_{2,5}'\}$, $\{\phi_{2,7}'',\phi_{2,13}''\}$,\\
$\{\phi_{2,11}'', \phi_{2,5}''\}$, $\{\phi_{2,9}'',\phi_{2,3}'\}$,
$\{\phi_{2,11}''',\phi_{2, 5}'''\}$,
$\{\phi_{2,9}''',\phi_{2,3}''\}$, $\{\phi_{3,4},\phi_{3,8},\phi_{3,
12}\}$\item[$a_0-a_1-b_0-b_1+2b_2-c_0+2c_1-c_2=0$]\hfil\break
$\{\phi_{1,12}'',\phi_{2,9}''',\phi_{2,3}'',\phi_{3,6},\phi_{3,10},\phi_{3,
2}\}$, $\{\phi_{2,9}',\phi_{2,15}\}$,
$\{\phi_{2,7}',\phi_{2,13}'\}$,
 $\{\phi_{2,11}',\phi_{2,5}'\}$, $\{\phi_{2,7}'',\phi_{2,13}''\}$,\\ $\{\phi_{2,
11}'',\phi_{2,5}''\}$, $\{\phi_{2,9}'',\phi_{2,3}'\}$,
$\{\phi_{2,11}''', \phi_{2,5}'''\}$, $\{\phi_{2,7}''',\phi_{2,1}\}$,
$\{\phi_{3,4},\phi_{3,8},
\phi_{3,12}\}$\item[$a_0-a_1-b_0-b_1+2b_2+2c_0-c_1-c_2=0$]\hfil\break
$\{\phi_{1,8}''',\phi_{2,11}''',\phi_{2,5}''',\phi_{3,6},\phi_{3,10},\phi_{3,
2}\}$, $\{\phi_{2,9}',\phi_{2,15}\}$,
$\{\phi_{2,7}',\phi_{2,13}'\}$,
 $\{\phi_{2,11}',\phi_{2,5}'\}$, $\{\phi_{2,7}'',\phi_{2,13}''\}$,\\ $\{\phi_{2,
11}'',\phi_{2,5}''\}$, $\{\phi_{2,9}'',\phi_{2,3}'\}$,
$\{\phi_{2,9}''', \phi_{2,3}''\}$, $\{\phi_{2,7}''',\phi_{2,1}\}$,
$\{\phi_{3,4},\phi_{3,8},
\phi_{3,12}\}$\item[$a_0-a_1-b_0+b_2-c_0+c_1=0$]\hfil\break
$\{\phi_{1,12}'',\phi_{1,6},\phi_{2,9}'',\phi_{2,3}'\}$,
$\{\phi_{2,9}', \phi_{2,15}\}$, $\{\phi_{2,7}',\phi_{2,13}'\}$,
$\{\phi_{2,11}',\phi_{2, 5}'\}$, $\{\phi_{2,7}'',\phi_{2,13}''\}$,
$\{\phi_{2,11}'',\phi_{2,5}''\}$,\\
 $\{\phi_{2,11}''',\phi_{2,5}'''\}$, $\{\phi_{2,9}''',\phi_{2,3}''\}$,
 $\{\phi_{2,7}''',\phi_{2,1}\}$, $\{\phi_{3,6},\phi_{3,10},\phi_{3,2}\}$,
$\{\phi_{3,4},\phi_{3,8},\phi_{3,12}\}$\item[$a_0-a_1-b_0+b_2-c_1+c_2=0$]\hfil\break
$\{\phi_{1,16},\phi_{1,10}',\phi_{2,7}'',\phi_{2,13}''\}$,
$\{\phi_{2,9}', \phi_{2,15}\}$, $\{\phi_{2,7}',\phi_{2,13}'\}$,
$\{\phi_{2,11}',\phi_{2, 5}'\}$, $\{\phi_{2,11}'',\phi_{2,5}''\}$,
$\{\phi_{2,9}'',\phi_{2,3}'\}$,\\
 $\{\phi_{2,11}''',\phi_{2,5}'''\}$, $\{\phi_{2,9}''',\phi_{2,3}''\}$,
 $\{\phi_{2,7}''',\phi_{2,1}\}$, $\{\phi_{3,6},\phi_{3,10},\phi_{3,2}\}$,
$\{\phi_{3,4},\phi_{3,8},\phi_{3,12}\}$\item[$a_0-a_1-b_0+b_2+c_0-c_2=0$]\hfil\break
$\{\phi_{1,8}''',\phi_{1,14}',\phi_{2,11}'',\phi_{2,5}''\}$,
$\{\phi_{2,9}', \phi_{2,15}\}$, $\{\phi_{2,7}',\phi_{2,13}'\}$,
$\{\phi_{2,11}',\phi_{2, 5}'\}$, $\{\phi_{2,7}'',\phi_{2,13}''\}$,
$\{\phi_{2,9}'',\phi_{2,3}'\}$,\\
 $\{\phi_{2,11}''',\phi_{2,5}'''\}$, $\{\phi_{2,9}''',\phi_{2,3}''\}$,
 $\{\phi_{2,7}''',\phi_{2,1}\}$, $\{\phi_{3,6},\phi_{3,10},\phi_{3,2}\}$,
$\{\phi_{3,4},\phi_{3,8},\phi_{3,12}\}$\item[$a_0-a_1-b_0+b_1-c_0+c_2=0$]\hfil\break
$\{\phi_{1,12}',\phi_{1,6},\phi_{2,9}''',\phi_{2,3}''\}$,
$\{\phi_{2,9}', \phi_{2,15}\}$, $\{\phi_{2,7}',\phi_{2,13}'\}$,
$\{\phi_{2,11}',\phi_{2, 5}'\}$, $\{\phi_{2,7}'',\phi_{2,13}''\}$,
$\{\phi_{2,11}'',\phi_{2,5}''\}$,\\
 $\{\phi_{2,9}'',\phi_{2,3}'\}$, $\{\phi_{2,11}''',\phi_{2,5}'''\}$,
 $\{\phi_{2,7}''',\phi_{2,1}\}$, $\{\phi_{3,6},\phi_{3,10},\phi_{3,2}\}$,
$\{\phi_{3,4},\phi_{3,8},\phi_{3,12}\}$\item[$a_0-a_1-b_0+b_1+c_1-c_2=0$]\hfil\break
$\{\phi_{1,8}'',\phi_{1,14}',\phi_{2,11}''',\phi_{2,5}'''\}$,
$\{\phi_{2,9}', \phi_{2,15}\}$, $\{\phi_{2,7}',\phi_{2,13}'\}$,
$\{\phi_{2,11}',\phi_{2, 5}'\}$, $\{\phi_{2,7}'',\phi_{2,13}''\}$,
$\{\phi_{2,11}'',\phi_{2,5}''\}$,\\
 $\{\phi_{2,9}'',\phi_{2,3}'\}$, $\{\phi_{2,9}''',\phi_{2,3}''\}$, $\{\phi_{2,
7}''',\phi_{2,1}\}$, $\{\phi_{3,6},\phi_{3,10},\phi_{3,2}\}$,
$\{\phi_{3,4},
\phi_{3,8},\phi_{3,12}\}$\item[$a_0-a_1-b_0+b_1+c_0-c_1=0$]\hfil\break
$\{\phi_{1,4}'',\phi_{1,10}',\phi_{2,7}''',\phi_{2,1}\}$,
$\{\phi_{2,9}', \phi_{2,15}\}$, $\{\phi_{2,7}',\phi_{2,13}'\}$,
$\{\phi_{2,11}',\phi_{2, 5}'\}$, $\{\phi_{2,7}'',\phi_{2,13}''\}$,
$\{\phi_{2,11}'',\phi_{2,5}''\}$,\\
 $\{\phi_{2,9}'',\phi_{2,3}'\}$, $\{\phi_{2,11}''',\phi_{2,5}'''\}$,
 $\{\phi_{2,9}''',\phi_{2,3}''\}$, $\{\phi_{3,6},\phi_{3,10},\phi_{3,2}\}$,
$\{\phi_{3,4},\phi_{3,8},\phi_{3,12}\}$\item[$a_0-a_1-b_0+2b_1-b_2-c_0-c_1+2c_2=0$]\hfil\break
$\{\phi_{1,12}',\phi_{2,9}'',\phi_{2,3}',\phi_{3,6},\phi_{3,10},\phi_{3,2}\}$,
 $\{\phi_{2,9}',\phi_{2,15}\}$, $\{\phi_{2,7}',\phi_{2,13}'\}$, $\{\phi_{2,
11}',\phi_{2,5}'\}$, $\{\phi_{2,7}'',\phi_{2,13}''\}$,\\
$\{\phi_{2,11}'', \phi_{2,5}''\}$,
$\{\phi_{2,11}''',\phi_{2,5}'''\}$, $\{\phi_{2,9}''',\phi_{2,
3}''\}$, $\{\phi_{2,7}''',\phi_{2,1}\}$,
$\{\phi_{3,4},\phi_{3,8},\phi_{3,
12}\}$\item[$a_0-a_1-b_0+2b_1-b_2-c_0+2c_1-c_2=0$]\hfil\break
$\{\phi_{1,8}'',\phi_{2,11}'',\phi_{2,5}'',\phi_{3,6},\phi_{3,10},\phi_{3,
2}\}$, $\{\phi_{2,9}',\phi_{2,15}\}$,
$\{\phi_{2,7}',\phi_{2,13}'\}$,
 $\{\phi_{2,11}',\phi_{2,5}'\}$, $\{\phi_{2,7}'',\phi_{2,13}''\}$,\\ $\{\phi_{2,
9}'',\phi_{2,3}'\}$, $\{\phi_{2,11}''',\phi_{2,5}'''\}$,
$\{\phi_{2,9}''', \phi_{2,3}''\}$, $\{\phi_{2,7}''',\phi_{2,1}\}$,
$\{\phi_{3,4},\phi_{3,8},
\phi_{3,12}\}$\item[$a_0-a_1-b_0+2b_1-b_2+2c_0-c_1-c_2=0$]\hfil\break
$\{\phi_{1,4}'',\phi_{2,7}'',\phi_{2,13}'',\phi_{3,6},\phi_{3,10},\phi_{3,
2}\}$, $\{\phi_{2,9}',\phi_{2,15}\}$,
$\{\phi_{2,7}',\phi_{2,13}'\}$,
 $\{\phi_{2,11}',\phi_{2,5}'\}$, $\{\phi_{2,11}'',\phi_{2,5}''\}$,\\ $\{\phi_{2,
9}'',\phi_{2,3}'\}$, $\{\phi_{2,11}''',\phi_{2,5}'''\}$,
$\{\phi_{2,9}''', \phi_{2,3}''\}$, $\{\phi_{2,7}''',\phi_{2,1}\}$,
$\{\phi_{3,4},\phi_{3,8},
\phi_{3,12}\}$\item[$a_0-a_1-b_1+b_2-c_0+c_2=0$]\hfil\break
$\{\phi_{1,16},\phi_{1,10}'',\phi_{2,7}',\phi_{2,13}'\}$,
$\{\phi_{2,9}', \phi_{2,15}\}$, $\{\phi_{2,11}',\phi_{2,5}'\}$,
$\{\phi_{2,7}'',\phi_{2, 13}''\}$, $\{\phi_{2,11}'',\phi_{2,5}''\}$,
$\{\phi_{2,9}'',\phi_{2,3}'\}$,\\
 $\{\phi_{2,11}''',\phi_{2,5}'''\}$, $\{\phi_{2,9}''',\phi_{2,3}''\}$,
 $\{\phi_{2,7}''',\phi_{2,1}\}$, $\{\phi_{3,6},\phi_{3,10},\phi_{3,2}\}$,
$\{\phi_{3,4},\phi_{3,8},\phi_{3,12}\}$\item[$a_0-a_1-b_1+b_2+c_1-c_2=0$]\hfil\break
$\{\phi_{1,12}'',\phi_{1,18}',\phi_{2,9}',\phi_{2,15}\}$,
$\{\phi_{2,7}', \phi_{2,13}'\}$, $\{\phi_{2,11}',\phi_{2,5}'\}$,
$\{\phi_{2,7}'',\phi_{2, 13}''\}$, $\{\phi_{2,11}'',\phi_{2,5}''\}$,
$\{\phi_{2,9}'',\phi_{2,3}'\}$,\\
 $\{\phi_{2,11}''',\phi_{2,5}'''\}$, $\{\phi_{2,9}''',\phi_{2,3}''\}$,
 $\{\phi_{2,7}''',\phi_{2,1}\}$, $\{\phi_{3,6},\phi_{3,10},\phi_{3,2}\}$,
$\{\phi_{3,4},\phi_{3,8},\phi_{3,12}\}$\item[$a_0-a_1-b_1+b_2+c_0-c_1=0$]\hfil\break
$\{\phi_{1,8}''',\phi_{1,14}'',\phi_{2,11}',\phi_{2,5}'\}$,
$\{\phi_{2,9}', \phi_{2,15}\}$, $\{\phi_{2,7}',\phi_{2,13}'\}$,
$\{\phi_{2,7}'',\phi_{2, 13}''\}$, $\{\phi_{2,11}'',\phi_{2,5}''\}$,
$\{\phi_{2,9}'',\phi_{2,3}'\}$,\\
 $\{\phi_{2,11}''',\phi_{2,5}'''\}$, $\{\phi_{2,9}''',\phi_{2,3}''\}$,
 $\{\phi_{2,7}''',\phi_{2,1}\}$, $\{\phi_{3,6},\phi_{3,10},\phi_{3,2}\}$,
 $\{\phi_{3,4},\phi_{3,8},\phi_{3,12}\}$\item[$a_0-a_1=0$]\hfil\break
$\{\phi_{1,0},\phi_{1,6}\}$, $\{\phi_{1,4}',\phi_{1,10}'\}$,
$\{\phi_{1,8}', \phi_{1,14}'\}$, $\{\phi_{1,4}'',\phi_{1,10}''\}$,
$\{\phi_{1,8}'',\phi_{1, 14}''\}$, $\{\phi_{1,12}',\phi_{1,18}'\}$,
$\{\phi_{1,8}''',\phi_{1,14}'''\}$,
 $\{\phi_{1,12}'',\phi_{1,18}''\}$, $\{\phi_{1,16},\phi_{1,22}\}$, $\{\phi_{2,
9}',\phi_{2,15}\}$, $\{\phi_{2,7}',\phi_{2,13}'\}$,
$\{\phi_{2,11}',\phi_{2, 5}'\}$, $\{\phi_{2,7}'',\phi_{2,13}''\}$,
$\{\phi_{2,11}'',\phi_{2,5}''\}$,
 $\{\phi_{2,9}'',\phi_{2,3}'\}$, $\{\phi_{2,11}''',\phi_{2,5}'''\}$,
 $\{\phi_{2,9}''',\phi_{2,3}''\}$, $\{\phi_{2,7}''',\phi_{2,1}\}$, $\{\phi_{3,
6},\phi_{3,4},\phi_{3,10},\phi_{3,8},\phi_{3,2},\phi_{3,12}\}$\item[$a_0-a_1+b_1-b_2-c_0+c_1=0$]\hfil\break
$\{\phi_{1,8}'',\phi_{1,14}''',\phi_{2,11}',\phi_{2,5}'\}$,
$\{\phi_{2,9}', \phi_{2,15}\}$, $\{\phi_{2,7}',\phi_{2,13}'\}$,
$\{\phi_{2,7}'',\phi_{2, 13}''\}$, $\{\phi_{2,11}'',\phi_{2,5}''\}$,
$\{\phi_{2,9}'',\phi_{2,3}'\}$,\\
 $\{\phi_{2,11}''',\phi_{2,5}'''\}$, $\{\phi_{2,9}''',\phi_{2,3}''\}$,
 $\{\phi_{2,7}''',\phi_{2,1}\}$, $\{\phi_{3,6},\phi_{3,10},\phi_{3,2}\}$,
$\{\phi_{3,4},\phi_{3,8},\phi_{3,12}\}$\item[$a_0-a_1+b_1-b_2-c_1+c_2=0$]\hfil\break
$\{\phi_{1,12}',\phi_{1,18}'',\phi_{2,9}',\phi_{2,15}\}$,
$\{\phi_{2,7}', \phi_{2,13}'\}$, $\{\phi_{2,11}',\phi_{2,5}'\}$,
$\{\phi_{2,7}'',\phi_{2, 13}''\}$, $\{\phi_{2,11}'',\phi_{2,5}''\}$,
$\{\phi_{2,9}'',\phi_{2,3}'\}$,\\
 $\{\phi_{2,11}''',\phi_{2,5}'''\}$, $\{\phi_{2,9}''',\phi_{2,3}''\}$,
 $\{\phi_{2,7}''',\phi_{2,1}\}$, $\{\phi_{3,6},\phi_{3,10},\phi_{3,2}\}$,
$\{\phi_{3,4},\phi_{3,8},\phi_{3,12}\}$\item[$a_0-a_1+b_1-b_2+c_0-c_2=0$]\hfil\break
$\{\phi_{1,4}'',\phi_{1,22},\phi_{2,7}',\phi_{2,13}'\}$,
$\{\phi_{2,9}', \phi_{2,15}\}$, $\{\phi_{2,11}',\phi_{2,5}'\}$,
$\{\phi_{2,7}'',\phi_{2, 13}''\}$, $\{\phi_{2,11}'',\phi_{2,5}''\}$,
$\{\phi_{2,9}'',\phi_{2,3}'\}$,\\
 $\{\phi_{2,11}''',\phi_{2,5}'''\}$, $\{\phi_{2,9}''',\phi_{2,3}''\}$,
 $\{\phi_{2,7}''',\phi_{2,1}\}$, $\{\phi_{3,6},\phi_{3,10},\phi_{3,2}\}$,
$\{\phi_{3,4},\phi_{3,8},\phi_{3,12}\}$\item[$a_0-a_1+b_0-2b_1+b_2-2c_0+c_1+c_2=0$]\hfil\break
$\{\phi_{1,10}'',\phi_{2,7}'',\phi_{2,13}'',\phi_{3,4},\phi_{3,8},\phi_{3,
12}\}$, $\{\phi_{2,9}',\phi_{2,15}\}$,
$\{\phi_{2,7}',\phi_{2,13}'\}$,
 $\{\phi_{2,11}',\phi_{2,5}'\}$, $\{\phi_{2,11}'',\phi_{2,5}''\}$,\\ $\{\phi_{2,
9}'',\phi_{2,3}'\}$, $\{\phi_{2,11}''',\phi_{2,5}'''\}$,
$\{\phi_{2,9}''', \phi_{2,3}''\}$, $\{\phi_{2,7}''',\phi_{2,1}\}$,
$\{\phi_{3,6},\phi_{3,10},
\phi_{3,2}\}$\item[$a_0-a_1+b_0-2b_1+b_2+c_0-2c_1+c_2=0$]\hfil\break
$\{\phi_{1,14}'',\phi_{2,11}'',\phi_{2,5}'',\phi_{3,4},\phi_{3,8},\phi_{3,
12}\}$, $\{\phi_{2,9}',\phi_{2,15}\}$,
$\{\phi_{2,7}',\phi_{2,13}'\}$,
 $\{\phi_{2,11}',\phi_{2,5}'\}$, $\{\phi_{2,7}'',\phi_{2,13}''\}$,\\ $\{\phi_{2,
9}'',\phi_{2,3}'\}$, $\{\phi_{2,11}''',\phi_{2,5}'''\}$,
$\{\phi_{2,9}''', \phi_{2,3}''\}$, $\{\phi_{2,7}''',\phi_{2,1}\}$,
$\{\phi_{3,6},\phi_{3,10},
\phi_{3,2}\}$\item[$a_0-a_1+b_0-2b_1+b_2+c_0+c_1-2c_2=0$]\hfil\break
$\{\phi_{1,18}',\phi_{2,9}'',\phi_{2,3}',\phi_{3,4},\phi_{3,8},\phi_{3,12}\}$,
 $\{\phi_{2,9}',\phi_{2,15}\}$, $\{\phi_{2,7}',\phi_{2,13}'\}$, $\{\phi_{2,
11}',\phi_{2,5}'\}$, $\{\phi_{2,7}'',\phi_{2,13}''\}$,\\
$\{\phi_{2,11}'', \phi_{2,5}''\}$,
$\{\phi_{2,11}''',\phi_{2,5}'''\}$, $\{\phi_{2,9}''',\phi_{2,
3}''\}$, $\{\phi_{2,7}''',\phi_{2,1}\}$,
$\{\phi_{3,6},\phi_{3,10},\phi_{3,
2}\}$\item[$a_0-a_1+b_0-b_1-c_0+c_1=0$]\hfil\break
$\{\phi_{1,4}',\phi_{1,10}'',\phi_{2,7}''',\phi_{2,1}\}$,
$\{\phi_{2,9}', \phi_{2,15}\}$, $\{\phi_{2,7}',\phi_{2,13}'\}$,
$\{\phi_{2,11}',\phi_{2, 5}'\}$, $\{\phi_{2,7}'',\phi_{2,13}''\}$,
$\{\phi_{2,11}'',\phi_{2,5}''\}$,\\
 $\{\phi_{2,9}'',\phi_{2,3}'\}$, $\{\phi_{2,11}''',\phi_{2,5}'''\}$,
 $\{\phi_{2,9}''',\phi_{2,3}''\}$, $\{\phi_{3,6},\phi_{3,10},\phi_{3,2}\}$,
$\{\phi_{3,4},\phi_{3,8},\phi_{3,12}\}$\item[$a_0-a_1+b_0-b_1-c_1+c_2=0$]\hfil\break
$\{\phi_{1,8}',\phi_{1,14}'',\phi_{2,11}''',\phi_{2,5}'''\}$,
$\{\phi_{2,9}', \phi_{2,15}\}$, $\{\phi_{2,7}',\phi_{2,13}'\}$,
$\{\phi_{2,11}',\phi_{2, 5}'\}$, $\{\phi_{2,7}'',\phi_{2,13}''\}$,
$\{\phi_{2,11}'',\phi_{2,5}''\}$,\\
 $\{\phi_{2,9}'',\phi_{2,3}'\}$, $\{\phi_{2,9}''',\phi_{2,3}''\}$, $\{\phi_{2,
7}''',\phi_{2,1}\}$, $\{\phi_{3,6},\phi_{3,10},\phi_{3,2}\}$,
$\{\phi_{3,4},
\phi_{3,8},\phi_{3,12}\}$\item[$a_0-a_1+b_0-b_1+c_0-c_2=0$]\hfil\break
$\{\phi_{1,0},\phi_{1,18}',\phi_{2,9}''',\phi_{2,3}''\}$,
$\{\phi_{2,9}', \phi_{2,15}\}$, $\{\phi_{2,7}',\phi_{2,13}'\}$,
$\{\phi_{2,11}',\phi_{2, 5}'\}$, $\{\phi_{2,7}'',\phi_{2,13}''\}$,
$\{\phi_{2,11}'',\phi_{2,5}''\}$,\\
 $\{\phi_{2,9}'',\phi_{2,3}'\}$, $\{\phi_{2,11}''',\phi_{2,5}'''\}$,
 $\{\phi_{2,7}''',\phi_{2,1}\}$, $\{\phi_{3,6},\phi_{3,10},\phi_{3,2}\}$,
$\{\phi_{3,4},\phi_{3,8},\phi_{3,12}\}$\item[$a_0-a_1+b_0-b_2-c_0+c_2=0$]\hfil\break
$\{\phi_{1,8}',\phi_{1,14}''',\phi_{2,11}'',\phi_{2,5}''\}$,
$\{\phi_{2,9}', \phi_{2,15}\}$, $\{\phi_{2,7}',\phi_{2,13}'\}$,
$\{\phi_{2,11}',\phi_{2, 5}'\}$, $\{\phi_{2,7}'',\phi_{2,13}''\}$,
$\{\phi_{2,9}'',\phi_{2,3}'\}$,\\
 $\{\phi_{2,11}''',\phi_{2,5}'''\}$, $\{\phi_{2,9}''',\phi_{2,3}''\}$,
 $\{\phi_{2,7}''',\phi_{2,1}\}$, $\{\phi_{3,6},\phi_{3,10},\phi_{3,2}\}$,
$\{\phi_{3,4},\phi_{3,8},\phi_{3,12}\}$\item[$a_0-a_1+b_0-b_2+c_1-c_2=0$]\hfil\break
$\{\phi_{1,4}',\phi_{1,22},\phi_{2,7}'',\phi_{2,13}''\}$,
$\{\phi_{2,9}', \phi_{2,15}\}$, $\{\phi_{2,7}',\phi_{2,13}'\}$,
$\{\phi_{2,11}',\phi_{2, 5}'\}$, $\{\phi_{2,11}'',\phi_{2,5}''\}$,
$\{\phi_{2,9}'',\phi_{2,3}'\}$,\\
 $\{\phi_{2,11}''',\phi_{2,5}'''\}$, $\{\phi_{2,9}''',\phi_{2,3}''\}$,
 $\{\phi_{2,7}''',\phi_{2,1}\}$, $\{\phi_{3,6},\phi_{3,10},\phi_{3,2}\}$,
$\{\phi_{3,4},\phi_{3,8},\phi_{3,12}\}$\item[$a_0-a_1+b_0-b_2+c_0-c_1=0$]\hfil\break
$\{\phi_{1,0},\phi_{1,18}'',\phi_{2,9}'',\phi_{2,3}'\}$,
$\{\phi_{2,9}', \phi_{2,15}\}$, $\{\phi_{2,7}',\phi_{2,13}'\}$,
$\{\phi_{2,11}',\phi_{2, 5}'\}$, $\{\phi_{2,7}'',\phi_{2,13}''\}$,
$\{\phi_{2,11}'',\phi_{2,5}''\}$,\\
 $\{\phi_{2,11}''',\phi_{2,5}'''\}$, $\{\phi_{2,9}''',\phi_{2,3}''\}$,
 $\{\phi_{2,7}''',\phi_{2,1}\}$, $\{\phi_{3,6},\phi_{3,10},\phi_{3,2}\}$,
$\{\phi_{3,4},\phi_{3,8},\phi_{3,12}\}$\item[$a_0-a_1+b_0+b_1-2b_2-2c_0+c_1+c_2=0$]\hfil\break
$\{\phi_{1,14}''',\phi_{2,11}''',\phi_{2,5}''',\phi_{3,4},\phi_{3,8},\phi_{3,
12}\}$, $\{\phi_{2,9}',\phi_{2,15}\}$,
$\{\phi_{2,7}',\phi_{2,13}'\}$,
 $\{\phi_{2,11}',\phi_{2,5}'\}$, $\{\phi_{2,7}'',\phi_{2,13}''\}$,\\ $\{\phi_{2,
11}'',\phi_{2,5}''\}$, $\{\phi_{2,9}'',\phi_{2,3}'\}$,
$\{\phi_{2,9}''', \phi_{2,3}''\}$, $\{\phi_{2,7}''',\phi_{2,1}\}$,
$\{\phi_{3,6},\phi_{3,10},
\phi_{3,2}\}$\item[$a_0-a_1+b_0+b_1-2b_2+c_0-2c_1+c_2=0$]\hfil\break
$\{\phi_{1,18}'',\phi_{2,9}''',\phi_{2,3}'',\phi_{3,4},\phi_{3,8},\phi_{3,
12}\}$, $\{\phi_{2,9}',\phi_{2,15}\}$,
$\{\phi_{2,7}',\phi_{2,13}'\}$,
 $\{\phi_{2,11}',\phi_{2,5}'\}$, $\{\phi_{2,7}'',\phi_{2,13}''\}$,\\ $\{\phi_{2,
11}'',\phi_{2,5}''\}$, $\{\phi_{2,9}'',\phi_{2,3}'\}$,
$\{\phi_{2,11}''', \phi_{2,5}'''\}$, $\{\phi_{2,7}''',\phi_{2,1}\}$,
$\{\phi_{3,6},\phi_{3,10},
\phi_{3,2}\}$\item[$a_0-a_1+b_0+b_1-2b_2+c_0+c_1-2c_2=0$]\hfil\break
$\{\phi_{1,22},\phi_{2,7}''',\phi_{2,1},\phi_{3,4},\phi_{3,8},\phi_{3,12}\}$,
 $\{\phi_{2,9}',\phi_{2,15}\}$, $\{\phi_{2,7}',\phi_{2,13}'\}$, $\{\phi_{2,
11}',\phi_{2,5}'\}$, $\{\phi_{2,7}'',\phi_{2,13}''\}$,\\
$\{\phi_{2,11}'', \phi_{2,5}''\}$, $\{\phi_{2,9}'',\phi_{2,3}'\}$,
$\{\phi_{2,11}''',\phi_{2, 5}'''\}$,
$\{\phi_{2,9}''',\phi_{2,3}''\}$, $\{\phi_{3,6},\phi_{3,10},\phi_{3,
2}\}$\item[$a_0-a_1+2b_0-b_1-b_2-c_0-c_1+2c_2=0$]\hfil\break
$\{\phi_{1,8}',\phi_{2,11}',\phi_{2,5}',\phi_{3,6},\phi_{3,10},\phi_{3,2}\}$,
 $\{\phi_{2,9}',\phi_{2,15}\}$, $\{\phi_{2,7}',\phi_{2,13}'\}$, $\{\phi_{2,
7}'',\phi_{2,13}''\}$, $\{\phi_{2,11}'',\phi_{2,5}''\}$,\\
$\{\phi_{2,9}'', \phi_{2,3}'\}$, $\{\phi_{2,11}''',\phi_{2,5}'''\}$,
$\{\phi_{2,9}''',\phi_{2, 3}''\}$, $\{\phi_{2,7}''',\phi_{2,1}\}$,
$\{\phi_{3,4},\phi_{3,8},\phi_{3,
12}\}$\item[$a_0-a_1+2b_0-b_1-b_2-c_0+2c_1-c_2=0$]\hfil\break
$\{\phi_{1,4}',\phi_{2,7}',\phi_{2,13}',\phi_{3,6},\phi_{3,10},\phi_{3,2}\}$,
 $\{\phi_{2,9}',\phi_{2,15}\}$, $\{\phi_{2,11}',\phi_{2,5}'\}$, $\{\phi_{2,
7}'',\phi_{2,13}''\}$, $\{\phi_{2,11}'',\phi_{2,5}''\}$,\\
$\{\phi_{2,9}'', \phi_{2,3}'\}$, $\{\phi_{2,11}''',\phi_{2,5}'''\}$,
$\{\phi_{2,9}''',\phi_{2, 3}''\}$, $\{\phi_{2,7}''',\phi_{2,1}\}$,
$\{\phi_{3,4},\phi_{3,8},\phi_{3,
12}\}$\item[$a_0-a_1+2b_0-b_1-b_2+2c_0-c_1-c_2=0$]\hfil\break
$\{\phi_{1,0},\phi_{2,9}',\phi_{2,15},\phi_{3,6},\phi_{3,10},\phi_{3,2}\}$,
 $\{\phi_{2,7}',\phi_{2,13}'\}$, $\{\phi_{2,11}',\phi_{2,5}'\}$, $\{\phi_{2,
7}'',\phi_{2,13}''\}$, $\{\phi_{2,11}'',\phi_{2,5}''\}$,\\
$\{\phi_{2,9}'', \phi_{2,3}'\}$, $\{\phi_{2,11}''',\phi_{2,5}'''\}$,
$\{\phi_{2,9}''',\phi_{2, 3}''\}$, $\{\phi_{2,7}''',\phi_{2,1}\}$,
$\{\phi_{3,4},\phi_{3,8},\phi_{3, 12}\}$
\end{description} Now, the
generic Hecke algebra of $G_6$ is
$$\begin{array}{rccl}
    \mathcal{H}(G_6) & = & < V,W \,\,| &  VWVWVW=WVWVWV, \\
     &  &  & (V-x_0)(V-x_1)=0 \\
     &  &  & (W-z_0)(W-z_1)(W-z_2)=0>
  \end{array}$$
As we can see in lemma $\ref{g7}$ of the Appendix,
$\mathcal{H}(G_6)$ is isomorphic to the subalgebra $\bar{A}:=<S,U>$
of the following specialization of $\mathcal{H}(G_7)$
$$\begin{array}{rccl}
    A & := & <S,T,U \,\,| &  STU=TUS=UST, T^3=1  \\
     &  &  & (S-x_0)(S-x_1)=0 \\
     &  &  & (U-z_0)(U-z_1)(U-z_2)=0>
  \end{array}$$
The algebra $A$ is the twisted symmetric algebra of the cyclic group
$C_3$ over the symmetric subalgebra $\bar{A}$ and the
block-idempotents of $A$ and $\bar{A}$ coincide for all further
specializations of the parameters. If we denote by $\phi$ the
characters of $\bar{A}$ and by $\psi$ the characters of $A$, we have
$$\begin{array}{llllll}\scriptsize
    \mathrm{Ind}_{\bar{A}}^A(\phi_{1,0}) & = & \psi_{1,0}+\psi_{1,4''}+\psi_{1,8'''}   &  \mathrm{Ind}_{\bar{A}}^A(\phi_{1,4}) & = & \psi_{1,4'}+\psi_{1,8''}+\psi_{1,12''}\\
    \mathrm{Ind}_{\bar{A}}^A(\phi_{1,8}) & = & \psi_{1,8'}+\psi_{1,12'}+\psi_{1,16}   &  \mathrm{Ind}_{\bar{A}}^A(\phi_{1,6}) & = & \psi_{1,6}+\psi_{1,10''}+\psi_{1,14'''}\\
    \mathrm{Ind}_{\bar{A}}^A(\phi_{1,10}) & = & \psi_{1,10'}+\psi_{1,14''}+\psi_{1,18''}   &  \mathrm{Ind}_{\bar{A}}^A(\phi_{1,14}) & = & \psi_{1,14'}+\psi_{1,18'}+\psi_{1,22}\\
    \mathrm{Ind}_{\bar{A}}^A(\phi_{2,5''}) & = & \psi_{2,9'}+\psi_{2,13''}+\psi_{2,5'''}  &  \mathrm{Ind}_{\bar{A}}^A(\phi_{2,3''}) & = &\psi_{2,7'}+\psi_{2,11''}+\psi_{2,3''}\\
    \mathrm{Ind}_{\bar{A}}^A(\phi_{2,3'})& = & \psi_{2,11'}+\psi_{2,7'''}+\psi_{2,3'}  & \mathrm{Ind}_{\bar{A}}^A(\phi_{2,7}) & = &\psi_{2,7''}+\psi_{2,11'''}+\psi_{2,15}\\
    \mathrm{Ind}_{\bar{A}}^A(\phi_{2,1}) & = & \psi_{2,9''}+\psi_{2,5'}+\psi_{2,1}&   \mathrm{Ind}_{\bar{A}}^A(\phi_{2,5'}) & = & \psi_{2,9'''}+\psi_{2,13'}+\psi_{2,5''}\\
    \mathrm{Ind}_{\bar{A}}^A(\phi_{3,2}) & = & \psi_{3,6}+\psi_{3,10}+\psi_{3,2} &   \mathrm{Ind}_{\bar{A}}^A(\phi_{3,4}) & = & \psi_{3,4}+\psi_{3,8}+\psi_{3,12}
  \end{array}$$
Let $\phi$ be a cyclotomic specialization of $\mathcal{H}(G_7)$ with
$$\phi(x_i)=\zeta_2^i x^{a_i} ,
  \phi(z_k)=\zeta_3^k x^{c_k}.$$ It corresponds to the cyclotomic
specialization $\phi'$ of $\mathcal{H}(G_7)$ with
$$\phi'(x_i)=\zeta_2^i x^{a_i} ,
\phi'(y_j)=\zeta_3^j , \phi(z_k)=\zeta_3^k x^{c_k}.$$
 Therefore, the
essential hyperplanes for $G_6$ are obtained from the essential
hyperplanes for $G_7$ by setting $b_0=b_1=b_2=0$ and the non-trivial
Rouquier blocks associated with each are:
\begin{description}\scriptsize
\item[No essential hyperplane]\hfil\break
$\{\phi_{2,5}'',\phi_{2,7}\}$, $\{\phi_{2,3}'',\phi_{2,5}'\}$,
$\{\phi_{2,3}', \phi_{2,1}\}$\item[$c_1-c_2=0$]\hfil\break
$\{\phi_{1,4},\phi_{1,8}\}$, $\{\phi_{1,10},\phi_{1,14}\}$,
$\{\phi_{2,5}'', \phi_{2,7}\}$,
$\{\phi_{2,3}'',\phi_{2,3}',\phi_{2,1},\phi_{2,5}'\}$\item[$c_0-c_1=0$]\hfil\break
$\{\phi_{1,0},\phi_{1,4}\}$, $\{\phi_{1,6},\phi_{1,10}\}$,
$\{\phi_{2,5}'', \phi_{2,3}'',\phi_{2,7},\phi_{2,5}'\}$,
$\{\phi_{2,3}',\phi_{2,1}\}$\item[$c_0-c_2=0$]\hfil\break
$\{\phi_{1,0},\phi_{1,8}\}$, $\{\phi_{1,6},\phi_{1,14}\}$,
$\{\phi_{2,5}'', \phi_{2,3}',\phi_{2,7},\phi_{2,1}\}$,
$\{\phi_{2,3}'',\phi_{2,5}'\}$\item[$a_0-a_1-2c_0+c_1+c_2=0$]\hfil\break
$\{\phi_{1,6},\phi_{2,5}'',\phi_{2,7},\phi_{3,4}\}$,
$\{\phi_{2,3}'',\phi_{2, 5}'\}$,
$\{\phi_{2,3}',\phi_{2,1}\}$\item[$a_0-a_1+c_0-2c_1+c_2=0$]\hfil\break
$\{\phi_{1,10},\phi_{2,3}'',\phi_{2,5}',\phi_{3,4}\}$,
$\{\phi_{2,5}'', \phi_{2,7}\}$,
$\{\phi_{2,3}',\phi_{2,1}\}$\item[$a_0-a_1+c_0+c_1-2c_2=0$]\hfil\break
$\{\phi_{1,14},\phi_{2,3}',\phi_{2,1},\phi_{3,4}\}$,
$\{\phi_{2,5}'',\phi_{2, 7}\}$,
$\{\phi_{2,3}'',\phi_{2,5}'\}$\item[$a_0-a_1-c_0-c_1+2c_2=0$]\hfil\break
$\{\phi_{1,8},\phi_{2,3}',\phi_{2,1},\phi_{3,2}\}$,
$\{\phi_{2,5}'',\phi_{2, 7}\}$,
$\{\phi_{2,3}'',\phi_{2,5}'\}$\item[$a_0-a_1-c_0+2c_1-c_2=0$]\hfil\break
$\{\phi_{1,4},\phi_{2,3}'',\phi_{2,5}',\phi_{3,2}\}$,
$\{\phi_{2,5}'',\phi_{2, 7}\}$,
$\{\phi_{2,3}',\phi_{2,1}\}$\item[$a_0-a_1+2c_0-c_1-c_2=0$]\hfil\break
$\{\phi_{1,0},\phi_{2,5}'',\phi_{2,7},\phi_{3,2}\}$,
$\{\phi_{2,3}'',\phi_{2, 5}'\}$,
$\{\phi_{2,3}',\phi_{2,1}\}$\item[$a_0-a_1-c_0+c_1=0$]\hfil\break
$\{\phi_{1,4},\phi_{1,6},\phi_{2,3}',\phi_{2,1}\}$,
$\{\phi_{2,5}'',\phi_{2, 7}\}$,
$\{\phi_{2,3}'',\phi_{2,5}'\}$\item[$a_0-a_1-c_1+c_2=0$]\hfil\break
$\{\phi_{1,8},\phi_{1,10},\phi_{2,5}'',\phi_{2,7}\}$,
$\{\phi_{2,3}'',\phi_{2, 5}'\}$,
$\{\phi_{2,3}',\phi_{2,1}\}$\item[$a_0-a_1+c_0-c_2=0$]\hfil\break
$\{\phi_{1,0},\phi_{1,14},\phi_{2,3}'',\phi_{2,5}'\}$,
$\{\phi_{2,5}'', \phi_{2,7}\}$,
$\{\phi_{2,3}',\phi_{2,1}\}$\item[$a_0-a_1-c_0+c_2=0$]\hfil\break
$\{\phi_{1,8},\phi_{1,6},\phi_{2,3}'',\phi_{2,5}'\}$,
$\{\phi_{2,5}'',\phi_{2, 7}\}$,
$\{\phi_{2,3}',\phi_{2,1}\}$\item[$a_0-a_1+c_1-c_2=0$]\hfil\break
$\{\phi_{1,4},\phi_{1,14},\phi_{2,5}'',\phi_{2,7}\}$,
$\{\phi_{2,3}'',\phi_{2, 5}'\}$,
$\{\phi_{2,3}',\phi_{2,1}\}$\item[$a_0-a_1+c_0-c_1=0$]\hfil\break
$\{\phi_{1,0},\phi_{1,10},\phi_{2,3}',\phi_{2,1}\}$,
$\{\phi_{2,5}'',\phi_{2, 7}\}$,
$\{\phi_{2,3}'',\phi_{2,5}'\}$\item[$a_0-a_1=0$]\hfil\break
$\{\phi_{1,0},\phi_{1,6}\}$, $\{\phi_{1,4},\phi_{1,10}\}$,
$\{\phi_{1,8}, \phi_{1,14}\}$, $\{\phi_{2,5}'',\phi_{2,7}\}$,
$\{\phi_{2,3}'',\phi_{2,5}'\}$,
 $\{\phi_{2,3}',\phi_{2,1}\}$,
 $\{\phi_{3,2},\phi_{3,4}\}$
\end{description}}}
\end{px}

Since it will take too many pages to describe the Rouquier blocks
associated with all essential hyperplanes of all exceptional complex
reflection groups, we have stored these data in a computer file and
created two GAP functions which display them. These functions are
called ``AllBlocks'' and ``DisplayAllBlocks'' and they can be found
on my webpage
$$\textrm{\emph{http://www.math.jussieu.fr/$\sim$chlouveraki}}$$
Let us give an example of their use for the group $G_4$.

\begin{px}`` \verb|gap>| '' \emph{\small is the GAP prompt}
\begin{verbatim}
gap> W:=ComplexReflectionGroup(4);
\end{verbatim}
\emph{\small \#creates the group $W$}
\begin{verbatim}
gap> DisplayAllBlocks(W);
No essential hyperplane
[["phi{1,0}"],["phi{1,4}"],["phi{1,8}"], ["phi{2,5}"],
["phi{2,3}"],["phi{2,1}"],["phi{3,2}"]]
c_1-c_2=0
[["phi{1,0}"],["phi{1,4}","phi{1,8}","phi{2,5}"],
["phi{2,3}","phi{2,1}"],["phi{3,2}"]]
c_0-c_1=0
[["phi{1,0}","phi{1,4}","phi{2,1}"],["phi{1,8}"],
["phi{2,5}","phi{2,3}"],["phi{3,2}"]]
c_0-c_2=0
[["phi{1,0}","phi{1,8}","phi{2,3}"],["phi{1,4}"],
["phi{2,5}","phi{2,1}"],["phi{3,2}"]]
2c_0-c_1-c_2=0
[["phi{1,0}","phi{2,5}","phi{3,2}"],
["phi{1,4}"],["phi{1,8}"], ["phi{2,3}"],["phi{2,1}"]]
c_0-2c_1+c_2=0
[["phi{1,0}"],["phi{1,4}","phi{2,3}","phi{3,2}"], ["phi{1,8}"],
["phi{2,5}"],["phi{2,1}"]]
c_0+c_1-2c_2=0
[["phi{1,0}"],["phi{1,4}"],["phi{1,8}","phi{2,1}","phi{3,2}"],
["phi{2,5}"],["phi{2,3}"]]
\end{verbatim}\emph{\small \# displays all essential hyperplanes for $W$ and the
Rouquier blocks associated with each}
\begin{verbatim}
gap> AllBlocks(W);
[rec( cond:=[ ],
block:=[[1],[2],[3],[4],[5],[6],[7]]),
rec( cond:=[0,1,-1],
block:=[[1],[2,3,4],[5,6],[7]]),
rec( cond:=[1,-1,0],
block:=[[1,2,6],[3],[4,5],[7]]),
rec( cond:=[1,0,-1],
block:=[[1,3,5],[2],[4,6],[7]]),
rec( cond:=[2,-1,-1],
block:=[[1,4,7],[2],[3],[5],[6]]),
rec( cond:=[1,-2,1],
block:=[[1],[2,5,7],[3],[4],[6]]),
rec( cond:=[1,1,-2],
block:=[[1],[2],[3,6,7],[4],[5]])]
\end{verbatim}
\emph{\small \# displays the same data in a way easy to work with:
the essential hyperplanes are represented by the vectors cond (such
that cond$*[c_0,c_1,c_2]=0$ is the corresponding essential
hyperplane
 and cond$:=[\,\,]$ means ``No essential hyperplane'') and the characters are
 defined by their indexes in the list ``CharNames($W$)'':}
\begin{verbatim}
gap> CharNames(W);
[ "phi{1,0}", "phi{1,4}", "phi{1,8}", "phi{2,5}",
  "phi{2,3}", "phi{2,1}", "phi{3,2}" ]
\end{verbatim}
\end{px}\

Let $W$ be an exceptional irreducible complex reflection group.
Since we have the Rouquier blocks associated with all essential
hyperplanes for $W$, we have created the function ``RouquierBlocks''
which calculates the Rouquier blocks of any cyclotomic Hecke algebra
associated with $W$. Given a cyclotomic specialization
$u_{\mathcal{C},j} \mapsto \zeta_{e_\mathcal{C}}^j
x^{n_{\mathcal{C},j}}$, this function checks to which essential
hyperplanes the $n_{\mathcal{C},j}$ belong and returns the blocks
obtained as unions of the Rouquier blocks associated with these
hyperplanes.

The function ``RouquierBlocks'' along with the function
``DisplayRouquier Blocks'' (the first returns the characters' index
in the list ``CharNames($W$)'', whereas the second returns their
name) can be also found on my webpage. Before we give an example of
their use, let us explain how we create a cyclotomic Hecke algebra
in GAP with the use of the package CHEVIE (we copy here the relative
part in the GAP manual, which can be found on J.Michel's webpage
\emph{http://www.math.jussieu.fr/$\sim$jmichel}):

The command ``Hecke($W$, \emph{para})'' returns the cyclotomic Hecke
algebra associated with the complex reflection group $W$. The
following forms are accepted for \emph{para}: if \emph{para} is a
single value, it is replicated to become a list of same length as
the number of generators of $W$. Otherwise, \emph{para} should be a
list of the same length as the number of generators of $W$, with
possibly unbound entries (which means it can also be a list of
lesser length). There should be at least one entry bound for each
orbit of reflections, and if several entries are bound for one
orbit, they should all be identical. Now again, an entry for a
reflection of order $e$ can be either a single value or a list of
length $e$. If it is a list, it is interpreted as the list
$[u_0,...,u_{e-1}]$ of parameters for that reflection. If it is a
single value $q$, it is interpreted as the partly specialized list
of parameters $[q,\mathrm{E}(e),...,\mathrm{E}(e)^{e-1}]$ (in GAP,
$\mathrm{E}(e)$ represents $\zeta_e$).

Let us now give an example of the definition of a cyclotomic Hecke
algebra and the use of the functions ``RouquierBlocks'' and
``DisplayRouquier Blocks'' on it.

\begin{px}
\emph{\small The generic Hecke algebra of $G_4$ is
$$\begin{array}{rccll}
    \mathcal{H}(G_4) & = & <S,T \,\,| &  STS=TST, & (S-x_0)(S-x_1)(S-x_2)=0 \\
     &  &  &  & (T-x_0)(T-x_1)(T-x_2)=0>
  \end{array}$$
If we want to calculate the Rouquier blocks of the cyclotomic Hecke
algebra
$$\begin{array}{rccll}
    \mathcal{H}_\phi & = & <S,T \,\,| &  STS=TST, & (S-1)(S-\zeta_3 x)(S-\zeta_3^2x^2)=0 \\
     &  &  &  & (T-1)(T-\zeta_3x)(T-\zeta_3^2x^2)=0>
  \end{array}$$
we use the following commands:}
\begin{verbatim}
gap> W:=ComplexReflectionGroup(4);
gap> H:=Hecke(W,[[1,E(3)*x,E(3)^2*x^2]]);
\end{verbatim}
\emph{\small \# here the single value
$[1,\mathrm{E}(3)*x,\mathrm{E}(3)^2*x^2]$ is interpreted, according
to the rules, as
$[\,[1,\mathrm{E}(3)*x,\mathrm{E}(3)^2*x^2],[1,\mathrm{E}(3)*x,\mathrm{E}(3)^2*x^2]\,]$}
\begin{verbatim}
gap> RouquierBlocks(H);
[ [ 1 ], [ 2, 5, 7 ], [ 3 ], [ 4 ], [ 6 ] ]
gap> DisplayRouquierBlocks(H);
[["phi{1,0}"],["phi{1,4}","phi{2,3}","phi{3,2}"],
["phi{1,8}"],["phi{2,5}"],[ "phi{2,1}"]]
\end{verbatim}
\end{px}

\section{Essential hyperplanes}

We have checked for all exceptional complex reflection groups that the
$\mathfrak{p}$-Rouquier blocks associated with no or some
$\mathfrak{p}$-essential hyperplane (\ie the blocks of
$A_{\mathfrak{p}A}\mathcal{H}$ or $A_{\mathfrak{q}_M}\mathcal{H}$
respectively) are fixed by the action of the Galois group
$\mathrm{Gal}(K/\mathbb{Q})$. This implies that if a
hyperplane is $\mathfrak{p}'$-essential for $W$ for some prime ideal
$\mathfrak{p}'$ lying over a prime number $p$, then it is
$\mathfrak{p}$-essential for all prime ideals $\mathfrak{p}$ lying
over $p$. Therefore, we can talk about determining the $p$-essential
hyperplanes for $W$, where $p$ is a prime number dividing the order
of the group.

\begin{px} \small{\emph{The prime numbers which divide the order of the group
$G_7$ are 2 and 3. The essential hyperplanes for $G_7$ are already
given in example $\ref{example g7}$ (note that different letters
represent different hyperplane orbits).\\
\\
The only 3-essential hyperplanes for $G_7$ are:\\
$$\begin{array}{ccc}
    c_1-c_2=0, & c_0-c_1=0, & c_0-c_2=0  \\
    b_1-b_2=0, & b_0-b_1=0, & b_0-b_2=0
  \end{array}$$
All its remaining essential hyperplanes are strictly 2-essential.\\
\\
From these, we can obtain the $p$-essential hyperplanes (where
$p=2,3$)
\begin{itemize}
  \item for $G_6$ by setting $b_0=b_1=b_2=0$,
  \item for $G_5$ by setting $a_0=a_1=0$,
  \item for $G_4$ by setting $a_0=a_1=b_0=b_1=b_2=0$.
\end{itemize}
}}
\end{px}

For the $p$-essential hyperplanes of the other groups, the reader
may refer to my webpage and use the function
``EssentialHyperplanes'' which is applied as follows
\begin{verbatim}
gap> EssentialHyperplanes(W,p);
\end{verbatim}
and returns
\begin{itemize}
  \item the essential hyperplanes for $W$, if $p=0$.
  \item the $p$-essential hyperplanes for $W$, if $p$ divides the order of
  $W$.
  \item error, if $p$ doesn't divide the order of $W$.
\end{itemize}

\begin{px}\
\begin{verbatim}
gap> W:=ComplexReflectionGroup(4);
gap> EssentialHyperplanes(W,0);
c_1-c_2=0
c_0-c_1=0
c_0-c_2=0
2c_0-c_1-c_2=0
c_0-2c_1+c_2=0
c_0+c_1-2c_2=0
gap> EssentialHyperplanes(W,2);
2c_0-c_1-c_2=0
c_0-2c_1+c_2=0
c_0+c_1-2c_2=0
c_0-c_1=0
c_1-c_2=0
c_0-c_2=0
gap> EssentialHyperplanes(W,3);
c_1-c_2=0
c_0-c_1=0
c_0-c_2=0
gap> EssentialHyperplanes(W,5);
Error, The number p should divide the order of the group
\end{verbatim}
\end{px}
\begin{remark} \emph{For the groups $G_{12}$, $G_{22}$, $G_{23}$, $G_{24}$, $G_{27}$,
$G_{29}$, $G_{30}$, $G_{31}$, $G_{33}$, $G_{34}$, $G_{35}$,
$G_{36}$, $G_{37}$ the only essential hypeplane is $a_0=a_1$, which
is $p$-essential for all the prime numbers $p$ dividing the order of
the group.}
\end{remark}

\chapter*{Appendix}
\addcontentsline{toc}{chapter}{Appendix}

Let $W$ be a complex reflection group and let us denote by
$\mathcal{H}(W)$ its generic Hecke algebra. Suppose that the
assumptions $\ref{ypo}$ are satisfied. Let $W'$ be another complex
reflection group such that, for some specialization of the
parameters, $\mathcal{H}(W)$ is the twisted symmetric algebra of a
finite cyclic group $G$ over the symmetric subalgebra
$\mathcal{H}(W')$. Then, if we know the blocks of $\mathcal{H}(W)$,
we can use propositions $\ref{1.42}$ and $\ref{1.45}$ in order to
calculate the blocks of $\mathcal{H}(W')$.

Moreover, in all the cases that will be studied below, if we denote
by $\chi'$ the (irreducible) restriction to $\mathcal{H}(W')$ of an
irreducible character $\chi \in \mathrm{Irr}(\mathcal{H}(W))$, then
the Schur elements verify
$$s_\chi = |W:W'| s_{\chi'}.$$ Therefore, if the Schur elements of
$\mathcal{H}(W)$ verify theorem $\ref{Schur element generic}$, so do
the Schur elements of $\mathcal{H}(W')$.

\subsection*{The groups $G_4$, $G_5$, $G_6$, $G_7$}

The following table gives the specializations of the parameters of
the generic Hecke algebra $\mathcal{H}(G_7)$,
$(x_0,x_1;y_0,y_1,y_2;z_0,z_1,z_2)$, which give the generic Hecke
algebras of the groups $G_4$, $G_5$ and $G_6$ (\cite{Ma2}, Table
4.6).

\begin{center}
\begin{tabular}{|c|c|c|c|c|}
  \hline
  Group & Index & S & T & U \\
  \hline
  $G_7$ & 1 & $x_0,x_1$ & $y_0,y_1,y_2$         & $z_0,z_1,z_2$\\
  $G_5$ & 2 & $1,-1$    & $y_0,y_1,y_2$         & $z_0,z_1,z_2$\\
  $G_6$ & 3 & $x_0,x_1$ & $1,\zeta_3,\zeta_3^2$ & $z_0,z_1,z_2$\\
  $G_4$ & 6 & $1,-1$    & $1,\zeta_3,\zeta_3^2$ & $z_0,z_1,z_2$\\
  \hline
\end{tabular}
\end{center}
$$\textrm{Specializations of the parameters for }\mathcal{H}(G_7)$$
\\

\begin{limma}\label{g7}\
\begin{itemize}
  \item The algebra $\mathcal{H}(G_7)$ specialized via
  $$(x_0,x_1;y_0,y_1,y_2;z_0,z_1,z_2) \mapsto
  (1,-1;y_0,y_1,y_2;z_0,z_1,z_2)$$
  is the twisted symmetric algebra of the
  cyclic group $C_2$ over the symmetric subalgebra $\mathcal{H}(G_5)$ with
  parameters $(y_0,y_1,y_2;z_0,z_1,z_2)$. The
  block-idempotents of the two algebras coincide.
  \item The algebra $\mathcal{H}(G_7)$ specialized via
  $$(x_0,x_1;y_0,y_1,y_2;z_0,z_1,z_2) \mapsto
  (x_0,x_1;1,\zeta_3,\zeta_3^2;z_0,z_1,z_2)$$
  is the twisted symmetric algebra of the
  cyclic group $C_3$ over the symmetric subalgebra $\mathcal{H}(G_6)$ with
  parameters $(x_0,x_1;z_0,z_1,z_2)$. The
  block-idempotents of the two algebras coincide.
  \item The algebra $\mathcal{H}(G_6)$ specialized via
  $$(x_0,x_1;z_0,z_1,z_2) \mapsto
  (1,-1;z_0,z_1,z_2)$$
  is the twisted symmetric algebra of the
  cyclic group $C_2$ over the symmetric subalgebra $\mathcal{H}(G_4)$ with
  parameters $(z_0,z_1,z_2)$. The
  block-idempotents of the two algebras coincide.
\end{itemize}
\end{limma}
\begin{apod}{\small We have
$$\begin{array}{rccl}
    \mathcal{H}(G_7) & = & <S,T,U\,\, | &  STU=TUS=UST  \\
     &  &  & (S-x_0)(S-x_1)=0 \\
     &  &  & (T-y_0)(T-y_1)(T-y_2)=0 \\
     &  &  & (U-z_0)(U-z_1)(U-z_2)=0>
  \end{array}$$
\begin{itemize}
\item Let
  $$\begin{array}{rccl}
     A & := & <S,T,U \,\,| &  STU=TUS=UST, S^2=1  \\
     &  &  & (T-y_0)(T-y_1)(T-y_2)=0 \\
     &  &  & (U-z_0)(U-z_1)(U-z_2)=0>
  \end{array}$$
   and $$\bar{A}:=<T,U>.$$
   Then $$A=\bar{A} \oplus S \bar{A}  \textrm{ and } \bar{A} \simeq \mathcal{H}(G_5).$$
  \item Let
   $$\begin{array}{rccl}
     A & := & <S,T,U \,\,| &  STU=TUS=UST, T^3=1  \\
     &  &  & (S-x_0)(S-x_1)=0 \\
     &  &  & (U-z_0)(U-z_1)(U-z_2)=0>
  \end{array}$$
  and $$\bar{A}:=<S,U>.$$
  Then $$A=\bigoplus_{i=0}^2 T^i \bar{A}  \textrm{ and } \bar{A} \simeq \mathcal{H}(G_6).$$
  \item Let
  $$\begin{array}{rccl}
     A & := & <S,U \,\,| &  SUSUSU=USUSUS, S^2=1  \\
     &  &  & (U-z_0)(U-z_1)(U-z_2)=0>
  \end{array}$$
   and $$\bar{A}:=<U,SUS>.$$
   Then $$A=\bar{A} \oplus S \bar{A} \textrm{ and } \bar{A} \simeq \mathcal{H}(G_4).$$}
\end{itemize}
\end{apod}

The Schur elements of all irreducible characters of
$\mathcal{H}(G_7)$ are calculated in \cite{Ma2} and they are
obtained by permutation of the parameters from the following ones
(for an explanation concerning the notations of characters, see
section 4.3): {\footnotesize
\\
\\
$s_{\phi_{1,0}}=\Phi_{1}(x_0/x_1)\cdot\allowbreak\Phi_{1}(x_0y_0^2z_0^2/x_1y_1y_2z_1z_2)
\cdot\allowbreak\Phi_{1}(y_0/y_1)\cdot\allowbreak\Phi_{1}(y_0/y_2)\cdot
\allowbreak\Phi_{1}(z_0/z_1)\cdot\allowbreak\Phi_{1}(z_0/z_2)\cdot\allowbreak
\Phi_{1}(x_0y_0z_0/x_1y_1z_1)\cdot\allowbreak\Phi_{1}(x_0y_0z_0/x_1y_1z_2)\cdot
\allowbreak\Phi_{1}(x_0y_0z_0/x_1y_2z_1)\cdot\allowbreak
\Phi_{1}(x_0y_0z_0/x_1y_2z_2)$
\\
\\
$s_{\phi_{2,9'}}=2y_2/y_0\Phi_{1}(y_0/y_1)\cdot\allowbreak\Phi_{1}(y_2/y_0)\cdot\allowbreak
\Phi_{1}(z_1/z_0)\cdot\allowbreak\Phi_{1}(z_2/z_0)\cdot\allowbreak\
\Phi_{1}(r/x_0y_0z_0)\cdot\allowbreak\Phi_{1}(r/x_0y_2z_1)\cdot\allowbreak
\Phi_{1}(r/x_0y_2z_2)\cdot\allowbreak\Phi_{1}(r/x_1y_0z_0)\cdot\allowbreak
\Phi_{1}(r/x_1y_2z_1)\cdot\allowbreak\Phi_{1}(r/x_1y_2z_2)$\\
where $r=\root 2\of{x_0x_1y_1y_2z_1z_2}$
\\
\\
$s_{\phi_{3,6}}=3\Phi_{1}(x_1/x_0)\cdot\allowbreak\Phi_{1}(x_0y_0z_0/r)\cdot
\allowbreak\Phi_{1}(x_0y_0z_1/r)\cdot\allowbreak\Phi_{1}(x_0y_0z_2/r)\cdot
\allowbreak\Phi_{1}(x_0y_1z_0/r)\cdot\allowbreak\Phi_{1}(x_0y_1z_1/r)\cdot
\allowbreak\Phi_{1}(x_0y_1z_2/r)\cdot\allowbreak\Phi_{1}(x_0y_2z_0/r)\cdot
\allowbreak\Phi_{1}(x_0y_2z_1/r)\cdot\allowbreak\Phi_{1}(x_0y_2z_2/r)$\\
where $\ r=\root 3\of{x_0^2x_1y_0y_1y_2z_0z_1z_2}$}
\\
\\
Following theorem $\ref{Semisimplicity Malle}$ and \cite{Ma4}, Table
8.1, if we set
$$\begin{array}{cl}
    X_i^{12}:=(\zeta_2)^{-i}x_i & (i=0,1), \\
    Y_j^{12}:=(\zeta_3)^{-j}y_j & (j=0,1,2), \\
    Z_k^{12}:=(\zeta_3)^{-k}z_k & (k=0,1,2),
  \end{array}$$
then $\mathbb{Q}(\zeta_{12})(X_0,X_1,Y_0,Y_1,Y_2,Z_0,Z_1,Z_2)$ is a
splitting field for $\mathcal{H}(G_7)$. Hence the factorization of
the Schur elements over that field is as described by theorem
$\ref{Schur element generic}$.

\subsection*{The groups $G_8$, $G_9$, $G_{10}$, $G_{11}$, $G_{12}$, $G_{13}$, $G_{14}$, $G_{15}$}

The following table gives the specializations of the parameters of
the generic Hecke algebra $\mathcal{H}(G_{11})$,
$(x_0,x_1;y_0,y_1,y_2;z_0,z_1,z_2,z_3)$, which give the generic
Hecke algebras of the groups $G_8,\ldots,G_{15}$ (\cite{Ma2}, Table
4.9).

\begin{center}
\begin{tabular}{|c|c|c|c|c|}
  \hline
  Group & Index & S & T & U \\
  \hline
  $G_{11}$ & 1  & $x_0,x_1$ & $y_0,y_1,y_2$         & $z_0,z_1,z_2,z_3$\\
  $G_{10}$ & 2  & $1,-1$    & $y_0,y_1,y_2$         & $z_1,z_1,z_2,z_3$\\
  $G_{15}$ & 2  & $x_0,x_1$ & $y_0,y_1,y_2$         & $\sqrt{u_0},\sqrt{u_1},-\sqrt{u_0},-\sqrt{u_1}$\\
  $G_9$    & 3  & $x_0,x_1$ & $1,\zeta_3,\zeta_3^2$ & $z_0,z_1,z_2,z_3$\\
  $G_{14}$ & 4  & $x_0,x_1$ & $y_0,y_1,y_2$         & $1,i,-1,-i$\\
  $G_8$    & 6  & $1,-1$    & $1,\zeta_3,\zeta_3^2$ & $z_0,z_1,z_2,z_3$\\
  $G_{13}$ & 6  & $x_0,x_1$ & $1,\zeta_3,\zeta_3^2$ & $\sqrt{u_0},\sqrt{u_1},-\sqrt{u_0},-\sqrt{u_1}$\\
  $G_{12}$ & 12 & $x_0,x_1$ & $1,\zeta_3,\zeta_3^2$ & $1,i,-1,-i$\\
  \hline
\end{tabular}
\end{center}

$$\textrm{Specializations of the parameters for }\mathcal{H}(G_{11})$$

\begin{limma}\label{g11}\
\begin{itemize}
 \item The algebra $\mathcal{H}(G_{11})$ specialized via
  $$(x_0,x_1;y_0,y_1,y_2;z_0,z_1,z_2,z_3) \mapsto
  (1,-1;y_0,y_1,y_2;z_0,z_1,z_2,z_3)$$
  is the twisted symmetric algebra of the
  cyclic group $C_2$ over the symmetric subalgebra $\mathcal{H}(G_{10})$ with
  parameters $(y_0,y_1,y_2;z_0,z_1,z_2,z_3)$. The
  block-idempo-tents of the two algebras coincide.
 \item The algebra $\mathcal{H}(G_{11})$ specialized via
  $$(x_0,x_1;y_0,y_1,y_2;z_0,z_1,z_2,z_3) \mapsto
  (x_0,x_1;1,\zeta_3,\zeta_3^2;z_0,z_1,z_2,z_3)$$
  is the twisted symmetric algebra of the
  cyclic group $C_3$ over the symmetric subalgebra $\mathcal{H}(G_9)$ with
  parameters $(x_0,x_1;z_0,z_1,z_2,z_3)$. The
  block-idempotents of the two algebras coincide.
 \item The algebra $\mathcal{H}(G_9)$ specialized via
  $$(x_0,x_1;z_0,z_1,z_2,z_3) \mapsto
  (1,-1;z_0,z_1,z_2,z_3)$$
  is the twisted symmetric algebra of the
  cyclic group $C_2$ over the symmetric subalgebra $\mathcal{H}(G_8)$ with
  parameters $(z_0,z_1,z_2,z_3)$. The
  block-idempotents of the two algebras coincide.
 \item The algebra $\mathcal{H}(G_{11})$ specialized via
  $$(x_0,x_1;y_0,y_1,y_2;z_0,z_1,z_2,z_3) \mapsto
  (x_0,x_1;y_0,y_1,y_2;1,i,-1,-i)$$
  is the twisted symmetric algebra of the
  cyclic group $C_4$ over the symmetric subalgebra $\mathcal{H}(G_{14})$ with
  parameters $(x_0,x_1;y_0,y_1,y_2)$. The
  block-idempotents of the two algebras coincide.
 \item The algebra $\mathcal{H}(G_{14})$ specialized via
  $$(x_0,x_1;y_0,y_1,y_2) \mapsto
  (x_0,x_1;1,\zeta_3,\zeta_3^2)$$
  is the twisted symmetric algebra of the
  cyclic group $C_3$ over the symmetric subalgebra $\mathcal{H}(G_{12})$ with
  parameters $(x_0,x_1)$. The
  block-idempotents of the two algebras coincide.
 \item The algebra $\mathcal{H}(G_{11})$ specialized via
  $$(x_0,x_1;y_0,y_1,y_2;z_0,z_1,z_2,z_3) \mapsto
  (x_0,x_1;y_0,y_1,y_2;\sqrt{u_0},\sqrt{u_1},-\sqrt{u_0},-\sqrt{u_1})$$
  is the twisted symmetric algebra of the
  cyclic group $C_2$ over the symmetric subalgebra $\mathcal{H}(G_{15})$ with
  parameters $(x_0,x_1;y_0,y_1,y_2;u_0,u_1)$. The
  block-idempotents of the two algebras coincide.
 \item The algebra $\mathcal{H}(G_{15})$ specialized via
  $$(x_0,x_1;y_0,y_1,y_2;u_0,u_1) \mapsto
  (x_0,x_1;1,\zeta_3,\zeta_3^2;u_0,u_1)$$
  is the twisted symmetric algebra of the
  cyclic group $C_3$ over the symmetric subalgebra $\mathcal{H}(G_{13})$ with
  parameters $(x_0,x_1;u_0,u_1)$. The
  block-idempotents of the two algebras coincide.
\end{itemize}
\end{limma}
\begin{apod}{\small We have $$\begin{array}{rccl}
    \mathcal{H}(G_{11}) & = & <S,T,U \,\,| &  STU=TUS=UST  \\
     &  &  & (S-x_0)(S-x_1)=0 \\
     &  &  & (T-y_0)(T-y_1)(T-y_2)=0 \\
     &  &  & (U-z_0)(U-z_1)(U-z_2)(U-z_3)=0>
  \end{array}$$
\begin{itemize}
\item Let
  $$\begin{array}{rccl}
     A & := & <S,T,U \,\,| &  STU=TUS=UST, S^2=1  \\
     &  &  & (T-y_0)(T-y_1)(T-y_2)=0 \\
     &  &  & (U-z_0)(U-z_1)(U-z_2)(U-z_3)=0>
  \end{array}$$
   and $$\bar{A}:=<T,U>.$$
   Then $$ A=\bar{A} \oplus S \bar{A} \textrm{ and } \bar{A} \simeq \mathcal{H}(G_{10}).$$
\item Let
   $$\begin{array}{rccl}
     A & := & <S,T,U \,\,| &  STU=TUS=UST, T^3=1  \\
     &  &  & (S-x_0)(S-x_1)=0 \\
     &  &  & (U-z_0)(U-z_1)(U-z_2)(U-z_3)=0>
  \end{array}$$
  and $$\bar{A}:=<S,U>.$$ Then
  $$A=\bigoplus_{i=0}^2 T^i \bar{A} \textrm{ and } \bar{A}  \simeq \mathcal{H}(G_9).$$
\item Let
  $$\begin{array}{rccl}
     A & := & <S,U \,\,| &  SUSUSU=USUSUS, S^2=1  \\
     &  &  & (U-z_0)(U-z_1)(U-z_2)(U-z_3)=0>
  \end{array}$$
   and $$\bar{A}:=<U,SUS>.$$
   Then $$A=\bar{A} \oplus S \bar{A} \textrm{ and } \bar{A}  \simeq \mathcal{H}(G_8).$$
\item Let
   $$\begin{array}{rccl}
    A& := & <S,T,U \,\,| &  STU=TUS=UST, U^4=1  \\
     &  &  & (S-x_0)(S-x_1)=0 \\
     &  &  & (T-y_0)(T-y_1)(T-y_2)=0>
  \end{array}$$
  and $$\bar{A}:=<S,T>.$$ Then
  $$A=\bigoplus_{i=0}^3 U^i \bar{A} \textrm{ and } \bar{A}\simeq \mathcal{H}(G_{14}).$$
\item Let
   $$\begin{array}{rccl}
    A& := & <S,T \,\,| &  STSTSTST=TSTSTSTS, T^3=1  \\
     &  &  & (S-x_0)(S-x_1)=0 >
  \end{array}$$
  and $$\bar{A}:=<S,TST^2,T^2ST>.$$
  Then
  $$A=\bigoplus_{i=0}^2 T^i \bar{A}\textrm{ and } \bar{A}  \simeq \mathcal{H}(G_{12}).$$
 \item Let
  $$\begin{array}{rccl}
     A & := & <S,T,U \,\,| &  STU=TUS=UST  \\
     &  &  & (S-x_0)(S-x_1)=0 \\
     &  &  & (T-y_0)(T-y_1)(T-y_2)=0 \\
     &  &  & (U^2-u_0)(U^2-u_1)=0>
  \end{array}$$
   and $$\bar{A}:=<S,T,U^2>.$$
   Then $$ A=\bar{A} \oplus U \bar{A} \textrm{ and } \bar{A} \simeq \mathcal{H}(G_{15}).$$
 \item Let
   $$\begin{array}{rccl}
     A & := & <U^2,S,T \,\,| & STU^2=U^2ST, U^2STST=TU^2STS, T^3=1  \\
     &  &  & (S-x_0)(S-x_1)=0 \\
     &  &  & (U^2-u_0)(U^2-u_1)=0>
  \end{array}$$
  and $$\bar{A}:=<U^2,S,T^2ST>.$$ Then
  $$ A=\bigoplus_{i=0}^2 T^i \bar{A}\textrm{ and } \bar{A} \simeq \mathcal{H}(G_{13}).$$}
\end{itemize}
\end{apod}

The Schur elements of all irreducible characters of
$\mathcal{H}(G_{11})$ are calculated in \cite{Ma2} and they are
obtained by permutation of the parameters from the following ones:
\\
\\
{\footnotesize
$s_{\phi_{1,0}}=\Phi_{1}(x_0/x_1)\cdot\allowbreak\Phi_{1}(y_0/y_1)\cdot\allowbreak
\Phi_{1}(y_0/y_2)\cdot\allowbreak\Phi_{1}(z_0/z_1)\cdot\allowbreak\Phi_{1}(
z_0/z_2)\cdot\allowbreak\Phi_{1}(z_0/z_3)\cdot\allowbreak\Phi_{1}(x_0y_0
z_0/x_1y_1z_1)\cdot\allowbreak\Phi_{1}(x_0y_0z_0/x_1y_1z_2)\cdot\allowbreak
\Phi_{1}(x_0y_0z_0/x_1y_1z_3)\cdot\allowbreak\Phi_{1}(x_0y_0z_0/x_1y_2z_
1)\cdot\allowbreak \Phi_{1}(x_0y_0z_0/x_1y_2z_2)\cdot\allowbreak
\\ \Phi_{1}(x_0y_0z_0/x_1y_2z_3)\cdot\allowbreak\Phi_{1}(x_0y_0^2z_0^2/x_1y_1y_2z_1z_2)
\cdot\allowbreak
\Phi_{1}(x_0y_0^2z_0^2/x_1y_1y_2z_1z_3)\cdot\allowbreak
\Phi_{1}(x_0y_0^2z_0^2/x_1y_1y_2z_2z_3)\cdot\allowbreak\Phi_{1}(x_0^2y_0^2z_0^3/x_1^2y_1y_2z_1z_2z_3)$
\\
\\
$s_{\phi_{2,1}}=-2z_1/z_0\Phi_{1}(y_0/y_2)\cdot\allowbreak\Phi_{1}(y_1/y_2)\cdot
\allowbreak\Phi_{1}(z_0/z_2)\cdot\allowbreak\Phi_{1}(z_0/z_3)\cdot\allowbreak
\Phi_{1}(z_1/z_2)\cdot\allowbreak\Phi_{1}(z_1/z_3)\cdot\allowbreak\Phi_{1}
(y_0z_0z_1/y_2z_2z_3)\cdot\allowbreak\Phi_{1}(y_1z_0z_1/y_2z_2z_3)\cdot
\allowbreak\Phi_{1}(r/x_0y_2z_2)\cdot\allowbreak\Phi_{1}(r/x_0y_2z_3)\cdot
\allowbreak\Phi_{1}(r/x_1y_2z_2)\cdot\allowbreak\\
\Phi_{1}(r/x_1y_2z_3)\cdot
\allowbreak\Phi_{1}(r/x_0y_0z_1)\cdot\allowbreak\Phi_{1}(r/x_0y_1z_1)\
\cdot\allowbreak\Phi_{1}(r/x_1y_0z_1)\cdot\allowbreak\Phi_{1}(r/x_1y_1z_1)$\\
where $r=\root 2\of{x_0x_1y_0y_1z_0z_1}$
\\
\\
$s_{\phi_{3,2}}=3\Phi_{1}(x_1/x_0)\cdot\allowbreak\Phi_{1}(z_1/z_3)\cdot\allowbreak
\Phi_{1}(z_2/z_3)\cdot\allowbreak\Phi_{1}(z_0/z_3)\cdot\allowbreak\Phi_{1}
(r/x_1y_0z_3)\cdot\allowbreak\Phi_{1}(r/x_1y_1z_3)\cdot\allowbreak\Phi_1
(r/x_1y_2z_3)\cdot\allowbreak\Phi_{1}(x_0y_0z_0/r)\cdot\allowbreak\Phi
_{1}(x_0y_0z_1/r)\cdot\allowbreak\Phi_{1}(x_0y_0z_2/r)\cdot\allowbreak\Phi_{1}
(x_0y_1z_0/r)\cdot\allowbreak
\\ \Phi_{1}(x_0y_1z_1/r)\cdot\allowbreak\
\Phi_{1}(x_0y_1z_2/r)\cdot\allowbreak\Phi_{1}(x_0y_2z_0/r)\cdot\allowbreak
\Phi_{1}(x_0y_2z_1/r)\cdot\allowbreak\Phi_{1}(x_0y_2z_2/r)$\\
where $r=\root 3\of{x_0^2x_1y_0y_1y_2z_0z_1z_2}$
\\
\\
$s_{\phi_{4,21}}=-4\Phi_{1}(y_0/y_1)\cdot\allowbreak\Phi_{1}(y_0/y_2)\cdot\allowbreak\
\Phi_{1}(r/x_0y_0z_0)\cdot\allowbreak\Phi_{1}(r/x_1y_0z_0)\cdot\allowbreak
\Phi_{1}(x_0y_0z_1/r)\cdot\allowbreak\Phi_{1}(x_0y_0z_2/r)\cdot\allowbreak
\Phi_{1}(x_0y_0z_3/r)\cdot\allowbreak\Phi_{1}(x_1y_0z_1/r)\cdot\allowbreak
\Phi_{1}(x_1y_0z_2/r)\cdot\allowbreak\Phi_{1}(x_1y_0z_3/r)\cdot\allowbreak
\\
\Phi_{1}(x_0x_1y_0y_1z_0z_1/r^2)\cdot\allowbreak\Phi_{1}(x_0x_1y_0y_1
z_0z_2/r^2)\cdot\allowbreak\Phi_{1}(x_0x_1y_0y_1z_0z_3/r^2)\cdot\allowbreak
\\ \Phi_{1}(x_0x_1y_0y_2z_0z_1/r^2)\cdot\allowbreak\Phi_{1}(x_0x_1y_0y_2z_0z_2/r^2)
\cdot\allowbreak\Phi_{1}(x_0x_1y_0y_2z_0z_3/r^2)$\\
where $r=\root 4\of{x_0^2x_1^2y_0^2y_1y_2z_0z_1z_2z_3}$}
\\
\\
Following theorem $\ref{Semisimplicity Malle}$ and \cite{Ma4}, Table
8.1, if we set
$$\begin{array}{cl}
    X_i^{24}:=(\zeta_2)^{-i}x_i & (i=0,1), \\
    Y_j^{24}:=(\zeta_3)^{-j}y_j & (j=0,1,2), \\
    Z_k^{24}:=(\zeta_4)^{-k}z_k & (k=0,1,2,3),
  \end{array}$$
then $\mathbb{Q}(\zeta_{24})(X_0,X_1,Y_0,Y_1,Y_2,Z_0,Z_1,Z_2,Z_3)$
is a splitting field for $\mathcal{H}(G_{11})$. Hence the
factorization of the Schur elements over that field is as described
by theorem $\ref{Schur element generic}$.

\subsection*{The groups $G_{16}$, $G_{17}$, $G_{18}$, $G_{19}$, $G_{20}$, $G_{21}$, $G_{22}$}

The following table gives the specializations of the parameters of
the generic Hecke algebra $\mathcal{H}(G_{19})$,
$(x_0,x_1;y_0,y_1,y_2;z_0,z_1,z_2,z_3,z_4)$, which give the generic
Hecke algebras of the groups $G_{16},\ldots,G_{22}$ (\cite{Ma2},
Table 4.12).

\begin{center}
\begin{tabular}{|c|c|c|c|c|}
  \hline
  Group & Index & S & T & U \\
  \hline
  $G_{19}$ & 1  & $x_0,x_1$ & $y_0,y_1,y_2$         & $z_0,z_1,z_2,z_3,z_4$\\
  $G_{18}$ & 2  & $1,-1$    & $y_0,y_1,y_2$         & $z_0,z_1,z_2,z_3,z_4$\\
  $G_{17}$ & 3  & $x_0,x_1$ & $1,\zeta_3,\zeta_3^2$ & $z_0,z_1,z_2,z_3,z_4$\\
  $G_{21}$ & 5  & $x_0,x_1$ & $y_0,y_1,y_2$         & $1,\zeta_5,\zeta_5^2,\zeta_5^3,\zeta_5^4$\\
  $G_{16}$ & 6  & $1,-1$    & $1,\zeta_3,\zeta_3^2$ & $z_0,z_1,z_2,z_3,z_4$\\
  $G_{20}$ & 10 & $1,-1$    & $y_0,y_1,y_2$         & $1,\zeta_5,\zeta_5^2,\zeta_5^3,\zeta_5^4$\\
  $G_{22}$ & 15 & $x_0,x_1$ & $1,\zeta_3,\zeta_3^2$ & $1,\zeta_5,\zeta_5^2,\zeta_5^3,\zeta_5^4$\\
  \hline
\end{tabular}
\end{center}

$$\textrm{Specializations of the parameters for }\mathcal{H}(G_{19})$$

\begin{limma}\label{g19}\
\begin{itemize}
 \item The algebra $\mathcal{H}(G_{19})$ specialized via
  $$(x_0,x_1;y_0,y_1,y_2;z_0,z_1,z_2,z_3,z_4) \mapsto
  (1,-1;y_0,y_1,y_2;z_0,z_1,z_2,z_3,z_4)$$
  is the twisted symmetric algebra of the
  cyclic group $C_2$ over the symmetric subalgebra $\mathcal{H}(G_{18})$ with
  parameters $(y_0,y_1,y_2;z_0,z_1,z_2,z_3,z_4)$. The
  block-idempotents of the two algebras coincide.
 \item The algebra $\mathcal{H}(G_{19})$ specialized via
  $$(x_0,x_1;y_0,y_1,y_2;z_0,z_1,z_2,z_3,z_4) \mapsto
  (x_0,x_1;1,\zeta_3,\zeta_3^2;z_0,z_1,z_2,z_3,z_4)$$
  is the twisted symmetric algebra of the
  cyclic group $C_3$ over the symmetric subalgebra $\mathcal{H}(G_{17})$ with
  parameters $(x_0,x_1;z_0,z_1,z_2,z_3,z_4)$. The
  block-idempotents of the two algebras coincide.
\item The algebra $\mathcal{H}(G_{17})$ specialized via
  $$(x_0,x_1;z_0,z_1,z_2,z_3,z_4) \mapsto
  (1,-1;z_0,z_1,z_2,z_3,z_4)$$
  is the twisted symmetric algebra of the
  cyclic group $C_2$ over the symmetric subalgebra $\mathcal{H}(G_{16})$ with
  parameters $(z_0,z_1,z_2,z_3,z_4)$. The
  block-idempotents of the two algebras coincide.
\item The algebra $\mathcal{H}(G_{19})$ specialized via
  $$(x_0,x_1;y_0,y_1,y_2;z_0,z_1,z_2,z_3,z_4) \mapsto
  (x_0,x_1;y_0,y_1,y_2;1,\zeta_5,\zeta_5^2,\zeta_5^3,\zeta_5^4)$$
  is the twisted symmetric algebra of the
  cyclic group $C_5$ over the symmetric subalgebra $\mathcal{H}(G_{21})$ with
  parameters $(x_0,x_1;y_0,y_1,y_2)$. The
  block-idempotents of the two algebras coincide.
\item The algebra $\mathcal{H}(G_{21})$ specialized via
  $$(x_0,x_1;y_0,y_1,y_2) \mapsto
  (1,-1;y_0,y_1,y_2)$$
  is the twisted symmetric algebra of the
  cyclic group $C_2$ over the symmetric subalgebra $\mathcal{H}(G_{20})$ with
  parameters $(y_0,y_1,y_2)$. The
  block-idempotents of the two algebras coincide.
\item The algebra $\mathcal{H}(G_{21})$ specialized via
  $$(x_0,x_1;y_0,y_1,y_2) \mapsto
  (x_0,x_1;1,\zeta_3,\zeta_3^2)$$
  is the twisted symmetric algebra of the
  cyclic group $C_3$ over the symmetric subalgebra $\mathcal{H}(G_{22})$ with
  parameters $(x_0,x_1)$. The
  block-idempotents of the two algebras coincide.
\end{itemize}
\end{limma}
\begin{apod}{\small We have $$\begin{array}{rccl}
    \mathcal{H}(G_{19}) & = & <S,T,U \,\,| &  STU=TUS=UST  \\
     &  &  & (S-x_0)(S-x_1)=0 \\
     &  &  & (T-y_0)(T-y_1)(T-y_2)=0 \\
     &  &  & (U-z_0)(U-z_1)(U-z_2)(U-z_3)(U-z_4)=0>
  \end{array}$$
\begin{itemize}
\item Let
  $$\begin{array}{rccl}
     A & := & <S,T,U \,\,| &  STU=TUS=UST, S^2=1  \\
     &  &  & (T-y_0)(T-y_1)(T-y_2)=0 \\
     &  &  & (U-z_0)(U-z_1)(U-z_2)(U-z_3)(U-z_4)=0>
  \end{array}$$
   and $$\bar{A}:=<T,U>.$$ Then
   $$A=\bar{A} \oplus S \bar{A} \textrm{ and } \bar{A} \simeq \mathcal{H}(G_{18}).$$
  \item Let
   $$\begin{array}{rccl}
     A & := & <S,T,U \,\,| &  STU=TUS=UST, T^3=1  \\
     &  &  & (S-x_0)(S-x_1)=0 \\
     &  &  & (U-z_0)(U-z_1)(U-z_2)(U-z_3)(U-z_4)=0>
  \end{array}$$
  and $$\bar{A}:=<S,U>.$$
  Then $$A=\bigoplus_{i=0}^2 T^i \bar{A} \textrm{ and } \bar{A} \simeq \mathcal{H}(G_{17}).$$
  \item Let
  $$\begin{array}{rccl}
     A & := & <S,U \,\,| &  SUSUSU=USUSUS, S^2=1  \\
     &  &  & (U-z_0)(U-z_1)(U-z_2)(U-z_3)(U-z_4)=0>
  \end{array}$$
   and $$\bar{A}:=<U,SUS>.$$ Then
   $$A=\bar{A} \oplus S \bar{A} \textrm{ and } \bar{A} \simeq
   \mathcal{H}(G_{16}).$$
   \item Let
   $$\begin{array}{rccl}
    A& := & <S,T,U \,\,| &  STU=TUS=UST, U^5=1  \\
     &  &  & (S-x_0)(S-x_1)=0 \\
     &  &  & (T-y_0)(T-y_1)(T-y_2)=0>
  \end{array}$$
  and $$\bar{A}:=<S,T>.$$ Then
  $$A=\bigoplus_{i=0}^4 U^i \bar{A}  \textrm{ and } \bar{A} \simeq \mathcal{H}(G_{21}).$$
\item Let
   $$\begin{array}{rccl}
    A& := & <S,T \,\,| &  STSTSTSTST=TSTSTSTSTS, S^2=1  \\
     &  &  & (T-y_0)(T-y_1)(T-y_2)=0>
  \end{array}$$
  and $$\bar{A}:=<T,STS>.$$ Then
  $$ A=\bar{A} \oplus S\bar{A} \textrm{ and } \bar{A} \simeq \mathcal{H}(G_{20}).$$
\item Let
   $$\begin{array}{rccl}
    A& := & <S,T \,\,| &  STSTSTSTST=TSTSTSTSTS, T^3=1  \\
     &  &  & (S-x_0)(S-x_1)=0 >
  \end{array}$$
  and $$\bar{A}:=<S,TST^2,T^2ST>.$$  Then
  $$A=\bigoplus_{i=0}^2 T^i \bar{A} \textrm{ and } \bar{A}\simeq
  \mathcal{H}(G_{22}).$$}
\end{itemize}
\end{apod}

The Schur elements of all irreducible characters of
$\mathcal{H}(G_{19})$ are calculated in \cite{Ma2} and they are
obtained by permutation of the parameters from the following ones:
\\
\\
{\footnotesize
$s_{\phi_{1,0}}=\Phi_{1}(x_0/x_1)\cdot\allowbreak\Phi_{1}(y_0/y_1)\cdot\allowbreak
\Phi_{1}(y_0/y_2)\cdot\allowbreak\Phi_{1}(z_0/z_1)\cdot\allowbreak\Phi_{1}(
z_0/z_2)\cdot\allowbreak\Phi_{1}(z_0/z_3)\cdot\allowbreak\Phi_{1}(z_0/z_4)
\cdot\allowbreak\Phi_{1}(x_0y_0z_0/x_1y_1z_1)\cdot\allowbreak\Phi_{1}(
x_0y_0z_0/x_1y_1z_2)\cdot\allowbreak\Phi_{1}(x_0y_0z_0/x_1y_1z_3)\cdot
\allowbreak\Phi_{1}(x_0y_0z_0/x_1y_1z_4)\cdot\allowbreak\\
\Phi_{1}(x_0y_0z_0/x_1y_2z_1)
\cdot\allowbreak\Phi_{1}(x_0y_0z_0/x_1y_2z_2)\cdot\allowbreak
\Phi_{1}(x_0y_0z_0/x_1y_2z_3)\cdot\allowbreak\Phi_{1}(x_0y_0z_0/x_1y_2z_4)\cdot
\allowbreak\\
\Phi_{1}(x_0y_0^2z_0^2/x_1y_1y_2z_1z_2)\cdot\allowbreak\Phi_1
(x_0y_0^2z_0^2/x_1y_1y_2z_1z_3)\cdot\allowbreak
\Phi_{1}(x_0y_0^2z_0^2/x_1y_1y_2z_1z_4) \cdot\allowbreak\\
\Phi_{1}(x_0y_0^2z_0^2/x_1y_1y_2z_2z_3)\cdot
\allowbreak\Phi_{1}(x_0y_0^2z_0^2/x_1y_1y_2z_2z_4)\cdot\allowbreak\Phi_{1}
(x_0y_0^2z_0^2/x_1y_1y_2z_3z_4)\cdot\allowbreak\\
\Phi_{1}(x_0^2y_0^2z_0^3/x_1^2y_1y_2z_1z_2z_3)\cdot\allowbreak\Phi_{1}(x_0^2y_0^2z_0^3/x_1^2y_1y_2z_1z_2z_4)
\cdot\allowbreak\Phi_{1}(x_0^2y_0^2z_0^3/x_1^2y_1y_2z_1z_3z_4)\cdot\allowbreak\\
\Phi_{1}(x_0^2y_0^2z_0^3/x_1^2y_1y_2z_2z_3z_4)\cdot\allowbreak
\Phi_{1}(x_0^2y_0^3z_0^4/x_1^2y_1^2y_2z_1z_2z_3z_4)\cdot\allowbreak\Phi_{1}
(x_0^2y_0^3z_0^4/x_1^2y_1y_2^2z_1z_2z_3z_4)\cdot\allowbreak\\
\Phi_{1}(x_0^3y_0^4z_0^4/x_1^3y_1^2y_2^2z_1z_2z_3z_4)$
\\
\\
$s_{\phi_{2,31'}}=-2\Phi_{1}(y_0/y_2)\cdot\allowbreak\Phi_{1}(y_1/y_2)\cdot\allowbreak
\Phi_{1}(z_0/z_2)\cdot\allowbreak\Phi_{1}(z_0/z_3)\cdot\allowbreak\Phi_{1}
(z_0/z_4)\cdot\allowbreak\Phi_{1}(z_1/z_2)\cdot\allowbreak\Phi_{1}(z_1/
z_3)\cdot\allowbreak\Phi_{1}(z_1/z_4)\cdot\allowbreak\Phi_{1}(y_0z_0z_1/
y_2z_2z_3)\cdot\allowbreak\Phi_{1}(y_0z_0z_1/y_2z_2z_4)\cdot\allowbreak
\Phi_{1}(y_0z_0z_1/y_2z_3z_4)\cdot\allowbreak\Phi_{1}(y_1z_0z_1/y_2z_2z_3)
\cdot\allowbreak\Phi_{1}(y_1z_0z_1/y_2z_2z_4)\cdot\allowbreak\Phi_{1}(y_1z_0z_1/y_2z_3z_4)
\cdot\allowbreak\Phi_{1}(y_0y_1z_0z_1^2/y_2^2z_2z_3z_4)\cdot
\allowbreak\Phi_{1}(y_0y_1z_0^2z_1/y_2^2z_2z_3z_4)\cdot\allowbreak\Phi_{1}
(r/x_0y_0z_0)\cdot\allowbreak\Phi_{1}(x_0y_0z_1/r)\cdot\allowbreak\Phi_{1}
(x_1y_0z_0/r)\cdot\allowbreak\Phi_{1}(r/x_1y_0z_1)\cdot\allowbreak\
\Phi_{1}(r/x_1y_2z_2)\cdot\allowbreak\Phi_{1}(r/x_1y_2z_3)\cdot\allowbreak
\Phi_{1}(r/x_1y_2z_4)\cdot\allowbreak\Phi_{1}(r/x_0y_2z_2)\cdot\allowbreak
\Phi_{1}(r/x_0y_2z_3)\cdot\allowbreak\Phi_{1}(r/x_0y_2z_4)\cdot\allowbreak
\Phi_{1}(rz_0z_1/x_0y_2z_2z_3z_4)\cdot\allowbreak\\ \Phi_{1}(rz_0z_1/x_1y_2z_2z_3z_4)$\\
where $r=\root 2\of{x_0x_1y_0y_1z_0z_1}$
\\
\\
$s_{\phi_{3,22'}}=3\Phi_{1}(x_1/x_0)\cdot\allowbreak\Phi_{1}(z_0/z_3)\cdot\allowbreak
\Phi_{1}(z_0/z_4)\cdot\allowbreak\Phi_{1}(z_1/z_3)\cdot\allowbreak\Phi_{1}
(z_1/z_4)\cdot\allowbreak\Phi_{1}(z_2/z_3)\cdot\allowbreak\Phi_{1}(z_2/z_4)
\cdot\allowbreak\Phi_{1}(x_0z_0z_1/x_1z_3z_4)\cdot\allowbreak\Phi_{1}
(x_0z_0z_2/x_1z_3z_4)\cdot\allowbreak\Phi_{1}(x_0z_1z_2/x_1z_3z_4)\cdot
\allowbreak\Phi_{1}(r/x_1y_0z_3)\cdot\allowbreak\Phi_{1}(r/x_1y_0z_4)\cdot
\allowbreak\Phi_{1}(r/x_1y_1z_3)\cdot\allowbreak\Phi_{1}(r/x_1y_1z_4)\cdot
\allowbreak\Phi_{1}(r/x_1y_2z_3)\cdot\allowbreak\Phi_{1}(r/x_1y_2z_4)\cdot
\allowbreak\Phi_{1}(x_0y_0z_0/r)\cdot\allowbreak\Phi_{1}(x_0y_0z_1/r)\
\cdot\allowbreak\Phi_{1}(x_0y_0z_2/r)\cdot\allowbreak\Phi_{1}(x_0y_1z_0/r)
\cdot\allowbreak\Phi_{1}(x_0y_1z_1/r)\cdot\allowbreak\Phi_{1}(x_0y_1z_2/r)
\cdot\allowbreak\Phi_{1}(x_0y_2z_0/r)\cdot\allowbreak\Phi_{1}(x_0y_2z_1/r)
\cdot\allowbreak\Phi_{1}(x_0y_2z_2/r)\cdot\allowbreak\Phi_{1}(r^2/x_0x_1y_0y_1z_3z_4)
\cdot\allowbreak\Phi_{1}(r^2/x_0x_1y_0y_2z_3z_4)\cdot\allowbreak
\Phi_{1}(r^2/x_0x_1y_1y_2z_3z_4)$\\
where $r=\root 3\of{x_0^2x_1y_0y\ _1y_2z_0z_1z_2}$
\\
\\
$s_{\phi_{4,18}}=-4\Phi_{1}(y_1/y_0)\cdot\allowbreak\Phi_{1}(y_0/y_2)\cdot
\allowbreak\Phi_{1}(z_0/z_4)\cdot\allowbreak\Phi_{1}(z_1/z_4)\cdot\allowbreak
\Phi_{1}(z_2/z_4)\cdot\allowbreak\Phi_{1}(z_3/z_4)\cdot\allowbreak\Phi_{1}
(x_0y_0z_0/r)\cdot\allowbreak\Phi_{1}(x_0y_0z_1/r)\cdot\allowbreak\Phi_{1}
(x_0y_0z_2/r)\cdot\allowbreak\Phi_{1}(x_0y_0z_3/r)\cdot\allowbreak\Phi_{1}(
x_1y_0z_0/r)\cdot\allowbreak\Phi_{1}(x_1y_0z_1/r)\cdot\allowbreak
\Phi_{1}(x_1y_0z_2/r)\cdot\allowbreak\Phi_{1}(x_1y_0z_3/r)\cdot\allowbreak
\Phi_{1}(r/x_0y_1z_4)\cdot\allowbreak\Phi_{1}(r/x_1y_1z_4)\cdot\allowbreak
\Phi_{1}(r/x_0y_2z_4)\cdot\allowbreak\Phi_{1}(r/x_1y_2z_4)\cdot\allowbreak\\
\Phi_{1}(r^2/x_0x_1y_0y_1z_0z_1)\cdot\allowbreak\Phi_{1}(r^2/x_0x_1y_0y_1z_0z_2)
\cdot\allowbreak\Phi_{1}(x_0x_1y_0y_1z_0z_3/r^2)\cdot\allowbreak\Phi_{1}
(x_0x_1y_0y_1z_1z_2/r^2)\cdot\allowbreak\Phi_{1}(r^2/x_0x_1y_0y_1z_1z_3
)\cdot\allowbreak\Phi_{1}(r^2/x_0x_1y_0y_1z_2z_3)\cdot\allowbreak
\Phi_{1}(r^2/x_0x_1y_1y_2z_0z_4)\cdot\allowbreak\Phi_{1}(r^2/x_0x_1y_1y_2z_1z_4)
\cdot\allowbreak\Phi_{1}(r^2/x_0x_1y_1y_2z_2z_4)\cdot\allowbreak\Phi_{1}
(r ^2/x_0x_1y_1y_2z_3z_4)$\\
where $r=\root 4\of{x_0^2x_1^2y_0^2y_1y_2z_0z_1z_2z_3}$
\\
\\
$s_{\phi_{5,16}}=5\Phi_{1}(x_0/x_1)\cdot\allowbreak\Phi_{1}(y_2/y_0)\cdot\allowbreak
\Phi_{1}(y_2/y_1)\cdot\allowbreak\Phi_{1}(x_0y_0z_0/r)\cdot\allowbreak
\Phi_{1}(x_0y_0z_1/r)\cdot\allowbreak\Phi_{1}(x_0y_0z_2/r)\cdot\allowbreak
\Phi_{1}(x_0y_0z_3/r)\cdot\allowbreak\Phi_{1}(x_0y_0z_4/r)\cdot\allowbreak
\ \Phi_{1}(x_0y_1z_0/r)\cdot\allowbreak\Phi_{1}(x_0y_1z_1/r)\cdot
\allowbreak\Phi_{1}(x_0y_1z_2/r)\cdot\allowbreak\Phi_{1}(x_0y_1z_3/r)\cdot
\allowbreak\Phi_{1}(x_0y_1z_4/r)\cdot\allowbreak\Phi_{1}(r/x_1y_2z_0)\cdot
\allowbreak\Phi_{1}(r/x_1y_2z_1)\cdot\allowbreak\Phi_{1}(r/x_1y_2z_2)\cdot
\allowbreak\Phi_{1}(r/x_1y_2z_3)\cdot\allowbreak\Phi_{1}(r/x_1y_2z_4)\cdot
\allowbreak\Phi_{1}(x_0x_1y_0y_1z_0z_1/r^2)\cdot\allowbreak\Phi_{1}
(x_0x_1y_0y_1z_0z_2/r^2)\cdot\allowbreak\\
\Phi_{1}(x_0x_1y_0y_1z_0z_3/r^2)\cdot
\allowbreak\Phi_{1}(x_0x_1y_0y_1z_0z_4/r^2)\cdot\allowbreak\Phi_{1}
(x_0x_1y_0y_1z_1z_2/r^2)\cdot\allowbreak\Phi_{1}(x_0x_1y_0y_1z_1z_3/r^2)
\cdot\allowbreak\Phi_{1}(x_0x_1y_0y_1z_1z_4/r^2)\cdot\allowbreak\Phi_{1}
(x_0x_1y_0y_1z_2z_3/r^2)\cdot\allowbreak\Phi_{1}(x_0x_1y_0y_1z_2z_4/r^2)\cdot\allowbreak
\Phi_{1}(x_0x_1y_0y_1z_3z_4/r^2)$\\
where $r=\root 5\of{x_0^3x_1^2y_0^2y_1^2y _2z_0z_1z_2z_3z_4}$
\\
\\
$s_{\phi_{6,15}}=-6\Phi_{1}(z_0/z_1)\cdot\allowbreak\Phi_{1}(z_0/z_2)\cdot\allowbreak
\Phi_{1}(z_0/z_3)\cdot\allowbreak\Phi_{1}(z_0/z_4)\cdot\allowbreak\Phi_{1}
(r/x_0y_0z_0)\cdot\allowbreak\Phi_{1}(r/x_0y_1z_0)\cdot\allowbreak\Phi_
{1}(r/x_0y_2z_0)\cdot\allowbreak\Phi_{1}(x_1y_0z_0/r)\cdot\allowbreak\Phi_{1}
(x_1y_1z_0/r)\cdot\allowbreak\Phi_{1}(x_1y_2z_0/r)\cdot\allowbreak
\Phi_{1}(x_0x_1y_0y_1z_0z_1/r^2)\cdot\allowbreak\\
\Phi_{1}(x_0x_1y_0y_1z_0z_2/r^2)\cdot\allowbreak\Phi_{1}(x_0x_1y_0y_1z_0z_3/r^2)\cdot\allowbreak
\Phi_{1}(x_0x_1y_0y_1z_0z_4/r^2)\cdot\allowbreak\Phi_{1}(x_0x_1y_0y_2z_0z_1/r^2)
\cdot\allowbreak\Phi_{1}(x_0x_1y_0y_2z_0z_2/r^2)\cdot\allowbreak\Phi_{1}
(x_0x_1y_0y_2z_0z_3/r^2)\cdot\allowbreak\Phi_{1}(x_0x_1y_0y_2z_0z_4/r^2)
\cdot\allowbreak\Phi_{1}(x_0x_1y_1y_2z_0z_1/r^2)\cdot\allowbreak\Phi_{1}
(x_0x_1y_1y_2z_0z_2/r^2)\cdot\allowbreak\Phi_{1}(x_0x_1y_1y_2z_0z_3/r^2)
\cdot\allowbreak\Phi_{1}(x_0x_1y_1y_2z_0z_4/r^2)\cdot\allowbreak\\
\Phi_{1} (x_0^2x_1y_0y_1y_2z_0z_1z_2/r^3)\cdot\allowbreak
\Phi_{1}(x_0^2x_1y_0y_1y_2z_0z_1
z_3/r^3)\cdot\allowbreak\Phi_{1}(x_0^2x_1y_0y_1y_2z_0z_1z_4/r^3)\cdot
\allowbreak\\
\Phi_{1}(x_0^2x_1y_0y_1y_2z_0z_2z_3/r^3)\cdot\allowbreak
\Phi_{1} (x_0^2x_1y_0y_1y_2z_0z_2z_4/r^3)\cdot\allowbreak
\Phi_{1}(x_0^2x_1y_0y_1y_2z_0z _3z_4/r^3)$\\ where
$r=\root6\of{x_0^3x_1^3y_0^2y_1^2y_2^2z_0^2z_1z_2z_3z_4}$}
\\
\\
Following theorem $\ref{Semisimplicity Malle}$ and \cite{Ma4}, Table
8.1, if we set
$$\begin{array}{cl}
    X_i^{60}:=(\zeta_2)^{-i}x_i & (i=0,1), \\
    Y_j^{60}:=(\zeta_3)^{-j}y_j & (j=0,1,2), \\
    Z_k^{60}:=(\zeta_5)^{-k}z_k & (k=0,1,2,3,4),
  \end{array}$$
then
$\mathbb{Q}(\zeta_{60})(X_0,X_1,Y_0,Y_1,Y_2,Z_0,Z_1,Z_2,Z_3,Z_4)$ is
a splitting field for $\mathcal{H}(G_{19})$. Hence the factorization
of the Schur elements over that field is as described by theorem
$\ref{Schur element generic}$.

\subsection*{The groups $G_{25}$, $G_{26}$}

The following table gives the specialization of the parameters of
the generic Hecke algebra $\mathcal{H}(G_{26})$,
$(x_0,x_1;y_0,y_1,y_2)$, which give the generic Hecke algebra of the
group $G_{25}$ (\cite{Ma5}, Theorem 6.3).

\begin{center}
\begin{tabular}{|c|c|c|c|}
  \hline
  Group & Index & S & T \\
  \hline
  $G_{26}$ & 1  & $x_0,x_1$ & $y_0,y_1,y_2$\\
  $G_{25}$ & 2  & $1,-1$    & $y_0,y_1,y_2$\\
  \hline
\end{tabular}
\end{center}

$$\textrm{Specialization of the parameters for }\mathcal{H}(G_{26})$$

\begin{limma}\label{g26}\
The algebra $\mathcal{H}(G_{26})$ specialized via
$$(x_0,x_1;y_0,y_1,y_2) \mapsto
(1,-1;y_0,y_1,y_2)$$ is the twisted symmetric algebra of the cyclic
group $C_2$ over the symmetric subalgebra $\mathcal{H}(G_{25})$ with
parameters $(y_0,y_1,y_2)$. The block-idempotents of the two
algebras coincide.
\end{limma}
\begin{apod}
{\small We have $$\begin{array}{rccl}
    \mathcal{H}(G_{26}) & = & <S,T,U \,\,| &  STST=TSTS, UTU=TUT, SU=US  \\
     &  &  & (S-x_0)(S-x_1)=0 \\
     &  &  & (T-y_0)(T-y_1)(T-y_2)=0 \\
     &  &  & (U-y_0)(U-y_1)(U-y_2)=0>
  \end{array}$$
  Let
$$\begin{array}{rccl}
     A & := & <S,T,U \,\,| & STST=TSTS, UTU=TUT, SU=US, S^2=1  \\
     &  &  & (T-y_0)(T-y_1)(T-y_2)=0 \\
     &  &  & (U-y_0)(U-y_1)(U-y_2)=0>
  \end{array}$$
   and $$\bar{A}:=<SUS,T,U>.$$ Then
   $$A=\bar{A} \oplus S \bar{A} \textrm{ and } \bar{A} \simeq \mathcal{H}(G_{25}).$$}
\end{apod}

The Schur elements of all irreducible characters of
$\mathcal{H}(G_{26})$ are calculated in \cite{Ma5} and they are
obtained by permutation of the parameters from the following ones:
\\
\\
{\footnotesize
$s_{\phi_{1,0}}=-\Phi_{1}(x_0/x_1)\cdot\allowbreak\Phi_{1}(y_0/y_1)\cdot\allowbreak
\Phi_{1}(y_0/y_2)\cdot\allowbreak\Phi_{2}(x_0y_0/x_1y_1)\cdot\allowbreak
\Phi_{2}(x_0y_0/x_1y_2)\cdot\allowbreak\Phi_{1}(x_0y_0^2/x_1y_1^2)\cdot
\allowbreak\Phi_{1}(x_0y_0^2/x_1y_2^2)\cdot\allowbreak\Phi_{2}(x_0y_0^3/x_1y_1^2y_2)
\cdot\allowbreak\Phi_{2}(x_0y_0^3/x_1y_1y_2^2)\cdot\allowbreak
\Phi_{6}(x_0y_0^2/x_1y_1y_2)\cdot\allowbreak\Phi_{2}(y_0^2/y_1y_2)\cdot
\allowbreak\Phi_{6}(y_0/y_1)\cdot\allowbreak\Phi_{6}(y_0/y_2)$
\\
\\
$s_{\phi_{2,3}}=y_1/y_0\Phi_{1}(x_0/x_1)\cdot\allowbreak\Phi_{1}(y_0/y_2)\cdot
\allowbreak\Phi_{1}(y_1/y_2)\cdot\allowbreak\Phi_{1}(x_0y_0/x_1y_2)\cdot
\allowbreak\Phi_{1}(x_0y_1/x_1y_2)\cdot\allowbreak\Phi_{2}(x_0y_0/x_1y_2)\cdot
\allowbreak\Phi_{2}(x_0y_1/x_1y_2)\cdot\allowbreak\Phi_{2}(x_0y_0/x_1y_1)
\cdot\allowbreak\Phi_{2}(x_0y_1/x_1y_0)\cdot\allowbreak\Phi_{6}
(x_0y_0y_1/x_1y_2^2)\cdot\allowbreak\Phi_{2}(y_0y_1/y_2^2)\cdot\allowbreak\Phi_{6}
(y_0/y_1)$
\\
\\
$s_{\phi_{3,6}}=-\Phi_{1}(x_0/x_1)\cdot\allowbreak\Phi_{3}(x_0/x_1)\cdot
\allowbreak\Phi_{2}(x_0y_0/x_1y_1)\cdot\allowbreak\Phi_{2}(x_0y_0/x_1y_2)\cdot
\allowbreak\Phi_{2}(x_0y_1/x_1y_0)\cdot\allowbreak\Phi_{2}(x_0y_1/x_1y_2)\cdot
\allowbreak\Phi_{2}(x_0y_2/x_1y_0)\cdot\allowbreak\Phi_{2}(x_0y_2/x_1y_1)
\cdot\allowbreak\Phi_{2}(y_0y_1/y_2^2)\cdot\allowbreak\Phi_{2}(y_0y_2/y_1^2
)\cdot\allowbreak\Phi_{2}(y_1y_2/y_0^2)$
\\
\\
$s_{\phi_{3,1}}=-\Phi_{1}(x_1/x_0)\cdot\allowbreak\Phi_{1}(y_0/y_1)\cdot
\allowbreak\Phi_{1}(y_0/y_2)\cdot\allowbreak\Phi_{2}(y_0/y_2)\cdot\allowbreak\Phi_{1}
(y_1/y_2)\cdot\allowbreak\Phi_{2}(y_0y_1/y_2^2)\cdot\allowbreak\Phi_{2}(
y_0^2/y_1y_2)\cdot\allowbreak\Phi_{6}(y_0/y_2)\cdot\allowbreak\Phi_{2}
(x_0y_0/x_1y_2)\cdot\allowbreak\Phi_{2}(x_0y_1/x_1y_0)\cdot\allowbreak
\Phi_{1}(x_0y_0^2/x_1y_1^2)\cdot\allowbreak\\
\Phi_{2}(x_0y_0^2y_1/x_1y_2^3)$
\\
\\
$s_{\phi_{6,2}}=\Phi_{1}(x_0/x_1)\cdot\allowbreak\Phi_{1}(y_1/y_0)\cdot\allowbreak
\Phi_{1}(y_0/y_2)\cdot\allowbreak\Phi_{1}(y_1/y_2)\cdot\allowbreak\Phi_{2}(
y_2/y_0)\cdot\allowbreak\Phi_{6}(y_0/y_2)\cdot\allowbreak\Phi_{2}(y_0y_2
/y_1^2)\cdot\allowbreak\Phi_{1}(x_0y_1/x_1y_2)\cdot\allowbreak\Phi_{2}
(x_1y_0/x_0y_2)\cdot\allowbreak\Phi_{2}(x_0y_1/x_1y_2)\cdot\allowbreak
\Phi_{2}(x_0y_0^3/x_1y_1^2y_2)$
\\
\\
$s_{\phi_{8,3}}=2\Phi_{1}(y_0/y_1)\cdot\allowbreak\Phi_{1}(y_0/y_2)\cdot\allowbreak
\Phi_{2}(x_1y_2/x_0y_1)\cdot\allowbreak\Phi_{2}(x_1y_1/x_0y_2)\cdot
\allowbreak\Phi_{2}(ry_0/x_1y_2^2)\cdot\allowbreak\Phi_{2}(ry_0/x_1y_1^2)\cdot
\allowbreak\Phi_{1}(ry_2/x_1y_0y_1)\cdot\allowbreak\Phi_{1}(ry_1/x_1y_0y_2)
\cdot\allowbreak\Phi_{3}(ry_0/x_1y_1y_2)\cdot\allowbreak\Phi_{3}(ry_0/x_0
y_1y_2)$\\
where $r=\root 2\of{-x_0x_1y_1y_2}$
\\
\\
$s_{\phi_{9,7}}=\Phi_{1}(\zeta_3^2)\cdot\allowbreak\Phi_{6}(y_0/y_1)\cdot\allowbreak
\ \Phi_{6}(y_2/y_0)\cdot\allowbreak\Phi_{6}(y_1/y_2)\cdot\allowbreak
\Phi_{2}(\zeta_3x_0y_1y_2/x_1y_0^2)\cdot\allowbreak
\Phi_{2}(\zeta_3x_0y_0y_2/x_1y_1^2)\cdot\allowbreak\Phi_{2}(\zeta_3x_0y_0y_1/x_1y_2^2)\cdot
\allowbreak\Phi_{1}(x_1/x_0)\cdot\allowbreak\Phi_{1}(\zeta_3x_0/x_1)$}
\\
\\
Following theorem $\ref{Semisimplicity Malle}$ and \cite{Ma4}, Table
8.2, if we set
$$\begin{array}{cl}
    X_i^6:=(\zeta_2)^{-i}x_i & (i=0,1), \\
    Y_j^6:=(\zeta_3)^{-j}y_j & (j=0,1,2),
  \end{array}$$
then $\mathbb{Q}(\zeta_3)(X_1,X_2,Y_1,Y_2,Y_2)$ is a splitting field
for $\mathcal{H}(G_{26})$. Hence the factorization of the Schur
elements over that field is as described by theorem
$\ref{Schur element generic}$.\\

\subsection*{The group $G_{28}$ (``$F_4$'')}

Let $\mathcal{H}(G_{28})$ be the generic Hecke algebra of the real
reflection group $G_{28}$ over the ring
$\mathbb{Z}[x_0^\pm,x_1^\pm,y_0^\pm,y_1^\pm]$. We have
$$\begin{array}{rccl}
     \mathcal{H}(G_{28}) & = & <S_1,S_2,T_1,T_2 \,\,| &S_1S_2S_1=S_2S_1S_2,\,\, T_1T_2T_1=T_2T_1T_2 \\
     &  &  &  S_1T_1=T_1S_1,\,\, S_1T_2=T_2S_1,\,\, S_2T_2=T_2S_2, \\
     &  &  &  S_2T_1S_2T_1=T_1S_2T_1S_2,\\
     &  &  & (S_i-x_0)(S_i-x_1)=(T_i-y_0)(T_i-y_1)=0>
  \end{array}$$
The Schur elements of all irreducible characters of
$\mathcal{H}(G_{28})$ have been calculated in \cite{Lu79b} and they
are obtained by permutation of the parameters from the following
ones:
\\
\\
{\footnotesize
$s_{\phi_{1,0}}=\Phi_{1}(y_0/y_1)\cdot\allowbreak\Phi_{6}(y_0/y_1)\cdot\allowbreak
\Phi_{1}(x_0/x_1)\cdot\allowbreak\Phi_{6}(x_0/x_1)\cdot\allowbreak\Phi_{1}(
x_0y_0^2/x_1y_1^2)\cdot\allowbreak
\Phi_{6}(x_0y_0/x_1y_1)\cdot\allowbreak
\Phi_{1}(x_0^2y_0/x_1^2y_1)\cdot\allowbreak\Phi_{4}(x_0y_0/x_1y_1)\cdot
\allowbreak\Phi_{2}(x_0y_0/x_1y_1)\cdot\allowbreak\Phi_{2}(x_0y_0/x_1y_1)$
\\
\\
$s_{\phi_{2,4''}}=-y_1/y_0\Phi_{6}(y_0/y_1)\cdot\allowbreak\Phi_{3}(x_0/x_1)\cdot
\allowbreak\Phi_{6}(x_0/x_1)\cdot\allowbreak\Phi_{1}(x_0/x_1)\cdot\allowbreak\
\Phi_{1}(x_0/x_1)\cdot\allowbreak\Phi_{1}(x_0^2y_0/x_1^2y_1)\cdot
\allowbreak\Phi_{2}(x_0y_0/x_1y_1)\cdot\allowbreak\Phi_{2}(x_0y_1/x_1y_0)\cdot
\allowbreak\Phi_{1}(x_0^2y_1/x_1^2y_0)$
\\
\\
$s_{\phi_{4,8}}=2\Phi_{6}(y_0/y_1)\cdot\allowbreak\Phi_{6}(x_1/x_0)\cdot\allowbreak
\Phi_{2}(x_0y_1/x_1y_0)\cdot\allowbreak\Phi_{2}(x_0y_1/x_1y_0)\cdot
\allowbreak\Phi_{2}(x_1y_1/x_0y_0)\cdot\allowbreak\Phi_{2}(x_0y_0/x_1y_1)$
\\
\\
$s_{\phi_{4,1}}=\Phi_{1}(y_0/y_1)\cdot\allowbreak\Phi_{6}(y_0/y_1)\cdot\allowbreak
\Phi_{1}(x_1/x_0)\cdot\allowbreak\Phi_{6}(x_0/x_1)\cdot\allowbreak\Phi_{2}(
x_0y_1/x_1y_0)\cdot\allowbreak
\Phi_{6}(x_0y_0/x_1y_1)\cdot\allowbreak
\Phi_{2}(x_0y_0/x_1y_1)\cdot\allowbreak\Phi_{2}(x_0y_0/x_1y_1)$
\\
\\
$s_{\phi_{6,6''}}=3\Phi_{1}(y_1/y_0)\cdot\allowbreak\Phi_{1}(y_1/y_0)\cdot\allowbreak
\Phi_{1}(x_1/x_0)\cdot\allowbreak\Phi_{1}(x_1/x_0)\cdot\allowbreak\Phi_{6}
(x_0y_0/x_1y_1)\cdot\allowbreak
\Phi_{2}(x_0y_1/x_1y_0)\cdot\allowbreak \Phi_{2}(x_1y_0/x_0y_1)$
\\
\\
$s_{\phi_{8,3''}}=-y_1/y_0\Phi_{6}(y_0/y_1)\cdot\allowbreak\Phi_{6}(x_0/x_1)\cdot
\allowbreak\Phi_{1}(x_0/x_1)\cdot\allowbreak\Phi_{1}(x_1/x_0)\cdot\allowbreak
\Phi_{3}(x_0/x_1)\cdot\allowbreak\Phi_{1}(x_0y_1^2/x_1y_0^2)\cdot
\allowbreak\Phi_{1}(x_0y_0^2/x_1y_1^2)$
\\
\\
$s_{\phi_{9,2}}=\Phi_{1}(y_0/y_1)\cdot\allowbreak\Phi_{1}(x_0/x_1)\cdot\allowbreak
\Phi_{1}(x_0y_1^2/x_1y_0^2)\cdot\allowbreak\Phi_{4}(x_0y_0/x_1y_1)\cdot
\allowbreak\Phi_{1}(x_1^2y_0/x_0^2y_1)\cdot\allowbreak\Phi_{2}(x_0y_0/x_1y_1)
\cdot\allowbreak\Phi_{2}(x_0y_0/x_1y_1)$
\\
\\
$s_{\phi_{12,4}}=6\Phi_{3}(y_0/y_1)\cdot\allowbreak\Phi_{3}(x_1/x_0)\cdot\allowbreak
\Phi_{2}(x_0y_1/x_1y_0)\cdot\allowbreak\Phi_{2}(x_0y_1/x_1y_0)\cdot
\allowbreak\Phi_{2}(x_0y_0/x_1y_1)\cdot\allowbreak\Phi_{2}(x_1y_1/x_0y_0)$
\\
\\
$s_{\phi_{16,5}}=2x_1y_1/x_0y_0\Phi_{6}(y_0/y_1)\cdot\allowbreak\Phi_{6}(x_1/x_0)\cdot
\allowbreak\Phi_{4}(x_0y_1/x_1y_0)\cdot\allowbreak\Phi_{4}(x_0y_0/x_1y_1)$}
\\
\\
Following theorem $\ref{Semisimplicity Malle}$, if we set
$$\begin{array}{cl}
    X_i^2:=(\zeta_2)^{-i}x_i & (i=0,1), \\
    Y_j^2:=(\zeta_2)^{-j}y_j & (j=0,1),
  \end{array}$$
then $\mathbb{Q}(X_0,X_1,Y_0,Y_1)$ is a splitting field for
$\mathcal{H}(G_{28})$. Hence the factorization of the Schur elements
over that field is as described by theorem $\ref{Schur element
generic}$.

\subsection*{The group $G_{32}$}

Let $\mathcal{H}(G_{32})$ be the generic Hecke algebra of the
complex reflection group $G_{32}$ over the ring
$\mathbb{Z}[x_0^\pm,x_1^\pm,x_2^\pm]$. We have
$$\begin{array}{rccl}
     \mathcal{H}(G_{32}) & = & <S_1,S_2,S_3,S_4 \,\,| &S_iS_{i+1}S_i=S_{i+1}S_iS_{i+1},  \\
     &  &  &  S_iS_j=S_jS_i \textrm{ when }|i-j|>1, \\
     &  &  & (S_i-x_0)(S_i-x_1)(S_i-x_2)=0>
  \end{array}$$
The Schur elements of all irreducible characters of
$\mathcal{H}(G_{32})$ have been calculated in \cite{Ma5} and they
are obtained by permutation of the parameters from the following
ones:
\\
\\
{\footnotesize
$s_{\phi_{1,0}}=\Phi_{1}(x_0/x_2)\cdot\allowbreak\Phi_{1}(x_0/x_2)\cdot
\allowbreak\Phi_{1}(x_0/x_1)\cdot\allowbreak\Phi_{1}(x_0/x_1)\cdot\allowbreak
\Phi_{1}(x_0^3/x_1x_2^2)\cdot\allowbreak\Phi_{1}(x_0^3/x_1^2x_2)
\cdot\allowbreak\Phi_{1}(x_0^5/x_1^3x_2^2)\cdot\allowbreak
\Phi_{1}(x_0^5/x_1^2x_2^3)
\cdot\allowbreak\Phi_{2}(x_0^4/x_1x_2^3)\cdot\allowbreak\Phi_{2}(x_0^4/x_1^3x_2)\cdot\allowbreak\Phi_{2}(x_0^2/x_1x_2)\cdot\allowbreak
\Phi_{2}(x_0^2/x_1x_2)\cdot\allowbreak\Phi_{6}(x_0/x_2)\cdot\allowbreak
\Phi_{6}(x_0/x_1)\cdot\allowbreak\Phi_{6}(x_0^3/x_1x_2^2)\cdot\allowbreak
\Phi_{6}(x_0^3/x_1^2x_2)\cdot\allowbreak\Phi_{6}(x_0^2/x_1x_2)\cdot\allowbreak
\Phi_{4}(x_0^2/x_1x_2)\cdot\allowbreak\Phi_{4}(x_0/x_2)\cdot\allowbreak
\Phi_{4}(x_0/x_1)\cdot\allowbreak\Phi_{3}(x_0^2/x_1x_2)\cdot\allowbreak
\Phi_{10}(x_0/x_2)\cdot\allowbreak\Phi_{10}(x_0/x_1)\cdot\allowbreak
\Phi_{5}(x_0^2/x_1x_2)$
\\
\\
$s_{\phi_{4,1}}=\Phi_{1}(x_0^4/x_1x_2^3)\cdot\allowbreak\Phi_{1}(x_0^3/x_1x_2^2)
\cdot\allowbreak\Phi_{1}(x_0^3/x_1^2x_2)\cdot\allowbreak\Phi_{1}(x_0^2x_1
/x_2^3)\cdot\allowbreak\Phi_{1}(x_1/x_0)\cdot\allowbreak\Phi_{1}(x_1/x_2
)\cdot\allowbreak\Phi_{1}(x_0/x_2)\cdot\allowbreak\Phi_{1}(x_0/x_2)\cdot
\allowbreak\Phi_{2}(x_0^5/x_1x_2^4)\cdot\allowbreak\Phi_{2}(x_0^3x_1/x_2^4)
\cdot\allowbreak\Phi_{2}(x_0^3/x_1x_2^2)\cdot\allowbreak\Phi_{2}(x_0^2
/x_1x_2)\cdot\allowbreak\Phi_{2}(x_0/x_2)\cdot\allowbreak\Phi_{2}(x_0x_1
/x_2^2)\cdot\allowbreak\Phi_{6}(x_0/x_2)\cdot\allowbreak\Phi_{6}(x_0/x_1)
\cdot\allowbreak\Phi_{4}(x_0/x_2)\cdot\allowbreak\Phi_{3}(x_0^2/x_1x_2
)\cdot\allowbreak\Phi_{10}(x_0/x_1)\cdot\allowbreak\Phi_{15}(x_0/x_2)$
\\
\\
$s_{\phi_{5,4}}=\Phi_{1}(x_0^3x_1^2/x_2^5)\cdot\allowbreak\Phi_{1}(x_0^2x_1/x_2^3
)\cdot\allowbreak\Phi_{1}(x_0/x_2)\cdot\allowbreak\Phi_{1}(x_0/x_2)\cdot
\allowbreak\Phi_{1}(x_1/x_0)\cdot\allowbreak\Phi_{1}(x_1/x_0)\cdot\allowbreak
\Phi_{1}(x_1/x_2)\cdot\allowbreak\Phi_{1}(x_1/x_2)\cdot\allowbreak
\Phi_{2}(x_0^3/x_1x_2^2)\cdot\allowbreak\Phi_{2}(x_0x_1^2/x_2^3)\cdot\allowbreak\Phi_{2}(x_1/x_2)\cdot\allowbreak\Phi_{2}(x_0^2/x_1x_2)\cdot
\allowbreak\Phi_{2}(x_0x_1/x_2^2)\cdot\allowbreak\Phi_{2}(x_0^4x_1/x_2^5)\cdot\allowbreak\Phi_{2}(x_0/x_2)\cdot\allowbreak\Phi_{2}(x_0/x_2)\cdot
\allowbreak\Phi_{6}(x_0/x_2)\cdot\allowbreak\Phi_{6}(x_0/x_2)\cdot\allowbreak\Phi_{6}(x_0/x_1)\cdot\allowbreak\Phi_{4}(x_0/x_1)\cdot\allowbreak\Phi_{3}(x_0x_1/x_2^2)\cdot\allowbreak\Phi_{12}(x_0/x_2)$
\\
\\
$s_{\phi_{6,8}}=x_1^2/x_0^2\Phi_{1}(x_0/x_1)\cdot\allowbreak\Phi_{1}(x_1/x_0)\cdot\allowbreak\Phi_{1}(x_0/x_2)\cdot\allowbreak
\Phi_{1}(x_1/x_2)\cdot\allowbreak\Phi_{1}(x_0/x_2)\cdot\allowbreak\Phi_{1}(x_1/x_2)\cdot\allowbreak\Phi_{1}(x_0x_1^2/x_2^3)\cdot\allowbreak
\Phi_{1}(x_0^2x_1/x_2^3)\cdot\allowbreak\Phi_{2}(x_0x_1^2/x_2^3)\cdot\allowbreak\Phi_{2}(x_0^2x_1/x_2^3)\cdot\allowbreak\Phi_{2}(x_0^2/x_1x_2)\cdot
\allowbreak\Phi_{2}(x_1^2/x_0x_2)\cdot\allowbreak\\
\Phi_{2}(x_0x_1/x_2^2)\cdot\allowbreak\Phi_{2}(x_0x_1/x_2^2)\cdot\allowbreak\Phi_{2}(x_1/x_2)\cdot
\allowbreak\Phi_{2}(x_0/x_2)\cdot\allowbreak\Phi_{6}(x_0x_1/x_2^2)\cdot
\allowbreak\Phi_{6}(x_0/x_2)\cdot\allowbreak\Phi_{6}(x_1/x_2)\cdot\allowbreak\Phi_{10}(x_0/x_1)\cdot\allowbreak\Phi_{5}(x_0x_1/x_2^2)$
\\
\\
$s_{\phi_{10,2}}=\Phi_{1}(x_0^2x_1/x_2^3)\cdot\allowbreak\Phi_{1}
(x_0^3/x_1x_2^2)\cdot\allowbreak\Phi_{1}(x_1/x_2)\cdot\allowbreak\Phi_{1}
(x_1/x_2)\cdot\allowbreak\Phi_{1}(x_1/x_0)\cdot\allowbreak\Phi_{1}(x_2/x_0)\cdot\allowbreak\Phi_{1}(x_0/x_2)\cdot\allowbreak\Phi_{1}(x_0/x_2)
\cdot\allowbreak\Phi_{2}(x_0^4/x_1^3x_2)\cdot\allowbreak\Phi_{2}(x_0^3/x_1x_2^2)\cdot\allowbreak\Phi_{2}(x_0x_1/x_2^2)\cdot\allowbreak\Phi_{2}
(x_0x_2/x_1^2)\cdot\allowbreak\Phi_{2}(x_0/x_2)\cdot\allowbreak\Phi_{2}(x_0/x_2)\cdot\allowbreak\Phi_{2}(x_1/x_2)\cdot\allowbreak\Phi_{2}(x_0^2/
x_1x_2)\cdot\allowbreak\Phi_{6}(x_0^2x_1/x_2^3)\cdot\allowbreak\Phi_{6}(x_0/x_2)\cdot\allowbreak\Phi_{6}(x_0/x_1)\cdot\allowbreak\Phi_{4}(x_0/x_2)
\cdot\allowbreak\Phi_{3}(x_0^2/x_1x_2)$
\\
\\
$s_{\phi_{15,6}}=\Phi_{1}(x_0^3/x_1x_2^2)\cdot\allowbreak\Phi_{1}(x_0^3/x_1^2x_2)\cdot\allowbreak\Phi_{1}(x_0/x_1)\cdot\allowbreak\Phi_{1}(x_0/x_1)
\cdot\allowbreak\Phi_{1}(x_2/x_1)\cdot\allowbreak\Phi_{1}(x_2/x_1)\cdot\allowbreak\Phi_{1}(x_2/x_0)\cdot\allowbreak\Phi_{1}(x_0/x_2)\cdot
\allowbreak\Phi_{2}(x_0^3x_2/x_1^4)\cdot\allowbreak\Phi_{2}(x_0^3x_1/x_2^4)\cdot\allowbreak\Phi_{2}(x_1^2/x_0x_2)\cdot\allowbreak\Phi_{2}(x_0x_1/x_2^2)
\cdot\allowbreak\Phi_{2}(x_0/x_2)\cdot\allowbreak\Phi_{2}(x_0/x_1)\cdot\allowbreak\Phi_{2}(x_0^2/x_1x_2)\cdot\allowbreak\Phi_{2}(x_0^2/x_1x_2)
\cdot\allowbreak\Phi_{6}(x_0^2/x_1x_2)\cdot\allowbreak\Phi_{6}(x_0/x_1)\cdot\allowbreak\Phi_{6}(x_0/x_2)\cdot\allowbreak\Phi_{4}(x_0^2/x_1x_2)$
\\
\\
$s_{\phi_{15,8}}=\Phi_{1}(x_1^2x_2/x_0^3)\cdot\allowbreak\Phi_{1}(x_0^2x_2/x_1^3)\cdot\allowbreak\Phi_{1}(x_0/x_2)\cdot\allowbreak\Phi_{1}(x_0/
x_2)\cdot\allowbreak\Phi_{1}(x_1/x_2)\cdot\allowbreak\Phi_{1}(x_1/x_2)\cdot\allowbreak\Phi_{1}(x_1/x_0)\cdot\allowbreak\Phi_{1}(x_0/x_1)\cdot
\allowbreak\Phi_{2}(x_1x_2/x_0^2)\cdot\allowbreak\Phi_{2}(x_0x_2/x_1^2)\cdot\allowbreak\Phi_{2}(x_0x_1/x_2^2)\cdot\allowbreak\Phi_{2}(x_0x_1/x_2^2)
\cdot\allowbreak\Phi_{2}(x_1/x_2)\cdot\allowbreak\Phi_{2}(x_1/x_2)\cdot\allowbreak\Phi_{2}(x_0/x_2)\cdot\allowbreak\Phi_{2}(x_0/x_2)\cdot
\allowbreak\Phi_{6}(x_0/x_2)\cdot\allowbreak\Phi_{6}(x_1/x_2)\cdot\allowbreak
\Phi_{4}(x_0x_1/x_2^2)$
\\
\\
$s_{\phi_{20,3}}=\Phi_{1}(x_0^2x_2/x_1^3)\cdot\allowbreak\Phi_{1}(x_0/x_1)\cdot\allowbreak\Phi_{1}(x_0/x_1)\cdot
\allowbreak\Phi_{1}(x_2/x_0)\cdot\allowbreak\Phi_{1}(x_0/x_2)\cdot\allowbreak\Phi_{1}(x_0^4/x_1x_2^3)\cdot\allowbreak\Phi_{1}(x_0^3/x_1x_2^2)\cdot
\allowbreak\Phi_{1}(x_1/x_2)\cdot\allowbreak\Phi_{2}(x_0x_1^2/x_2^3)\cdot\allowbreak\Phi_{2}(x_1^2/x_0x_2)\cdot\allowbreak\Phi_{2}(x_0^2/x_1x_2)
\cdot\allowbreak\Phi_{2}(x_0x_1/x_2^2)\cdot\allowbreak\Phi_{2}(x_0/x_2)\cdot\allowbreak\Phi_{2}(x_2/x_0)\cdot\allowbreak\Phi_{6}(x_0^3/x_1^2x_2)
\cdot\allowbreak\Phi_{6}(x_0/x_2)\cdot\allowbreak\Phi_{4}(x_0/x_2)\cdot\allowbreak\Phi_{3}(x_0x_1/x_2^2)$
\\
\\
$s_{\phi_{20,5}}=-\Phi_{1}(x_1^3/x_0^2x_2)\cdot\allowbreak\Phi_{1}(x_0^3/x_1^2x_2)\cdot\allowbreak\Phi_{1}(x_0x_1^2/x_2^3)\cdot
\allowbreak\Phi_{1}(x_0^2x_1/x_2^3)\cdot\allowbreak\Phi_{1}(x_2/x_1)\cdot\allowbreak\Phi_{1}(x_2/x_1)\cdot\allowbreak\Phi_{1}(x_0/x_2)\cdot\allowbreak
\Phi_{1}(x_0/x_2)\cdot\allowbreak\Phi_{2}(x_1^3/x_0x_2^2)\cdot\allowbreak\Phi_{2}(x_0^3/x_1x_2^2)\cdot\allowbreak\Phi_{2}(x_0x_1/x_2^2)\cdot
\allowbreak\Phi_{2}(x_0x_1/x_2^2)\cdot\allowbreak\Phi_{6}(x_0x_1/x_2^2)\cdot\allowbreak\Phi_{6}(x_1/x_0)\cdot\allowbreak\Phi_{6}(x_1/x_2)\cdot
\allowbreak\Phi_{6}(x_0/x_2)\cdot\allowbreak\Phi_{3}(x_0x_1/x_2^2)$
\\
\\
$s_{\phi_{20,7}}=\Phi_{1}(x_0^3x_1/x_2^4)\cdot\allowbreak\Phi_{1}(x_0x_1^2/x_2^3)\cdot\allowbreak\Phi_{1}(x_1/x_0)
\cdot\allowbreak\Phi_{1}(x_1/x_0)\cdot\allowbreak\Phi_{1}(x_2/x_0)\cdot\allowbreak\Phi_{1}(x_0/x_2)\cdot\allowbreak\Phi_{1}(x_1/x_2)\cdot
\allowbreak\Phi_{1}(x_1/x_2)\cdot\allowbreak\Phi_{2}(x_0/x_2)\cdot\allowbreak\Phi_{2}(x_2/x_0)\cdot\allowbreak\Phi_{2}(x_0^3x_2/x_1^4)\cdot
\allowbreak\Phi_{2}(x_0x_1^2/x_2^3)\cdot\allowbreak\Phi_{2}(x_0x_1/x_2^2)\cdot\allowbreak\Phi_{2}(x_1^2/x_0x_2)\cdot\allowbreak\Phi_{6}(x_0/x_2)\cdot
\allowbreak\Phi_{6}(x_0/x_2)\cdot\allowbreak\Phi_{6}(x_1/x_2)\cdot\allowbreak\Phi_{3}(x_0^2/x_1x_2)$
\\
\\
$s_{\phi_{20,12}}=2\Phi_{1}(x_2/x_1)\cdot\allowbreak\Phi_{1}(x_1/x_2)\cdot\allowbreak\Phi_{1}(x_2/x_0)\cdot\allowbreak\Phi_{1}(x_1/x_0)\cdot\allowbreak\Phi_{1}(x_2/x_0)\
\cdot\allowbreak\Phi_{1}(x_1/x_0)\cdot\allowbreak\Phi_{1}(x_0/x_2)\cdot\allowbreak\Phi_{1}(x_0/x_1)\cdot\allowbreak\Phi_{2}(x_0x_1^2/x_2^3)\cdot
\allowbreak\Phi_{2}(x_0x_2^2/x_1^3)\cdot\allowbreak\Phi_{2}(x_0^2/x_1x_2)\cdot\allowbreak\Phi_{2}(x_1x_2/x_0^2)\cdot\allowbreak\Phi_{2}(x_0/x_1)\
\cdot\allowbreak\Phi_{2}(x_0/x_2)\cdot\allowbreak\Phi_{2}(x_0/x_1)\cdot\allowbreak\Phi_{2}(x_0/x_2)\cdot\allowbreak\Phi_{6}(x_2/x_1)\cdot\allowbreak\Phi_{6}(x_1/x_2)\cdot\allowbreak\Phi_{3}(x_0^2/x_1x_2)$
\\
\\
$s_{\phi_{24,6}}=\Phi_{1}(x_1^3/x_0^2x_2)
\cdot\allowbreak\Phi_{1}(x_2^3/x_0^2x_1)\cdot\allowbreak\Phi_{1}(x_0/x_1)\cdot\allowbreak\Phi_{1}(x_0/x_2)\cdot\allowbreak\Phi_{1}(x_0/x_1)\cdot\
\allowbreak\Phi_{1}(x_0/x_2)\cdot\allowbreak\Phi_{1}(x_2/x_1)\cdot\allowbreak\Phi_{1}(x_1/x_2)\cdot\allowbreak\Phi_{2}(x_0/x_1)\cdot\allowbreak\
\Phi_{2}(x_0/x_2)\cdot\allowbreak\Phi_{2}(x_0x_2/x_1^2)\cdot\allowbreak\Phi_{2}(x_0x_1/x_2^2)\cdot\allowbreak\Phi_{6}(x_0/x_1)\cdot\allowbreak
\Phi_{6}(x_0/x_2)\cdot\allowbreak\Phi_{4}(x_0/x_1)\cdot\allowbreak\Phi_{4}(x_0/x_2)\cdot\allowbreak\Phi_{5}(x_0^2/x_1x_2)$
\\
\\
$s_{\phi_{30,4}}=\Phi_{1}(x_0^5/x_1^3x_2^2)\cdot\allowbreak\Phi_{1}(x_0/x_2)\cdot\allowbreak\Phi_{1}(x_1/x_0)\cdot\allowbreak\Phi_{1}(x_1/x_0)
\cdot\allowbreak\Phi_{1}(x_1/x_2)\cdot\allowbreak\Phi_{1}(x_1/x_2)\cdot\allowbreak\Phi_{1}(x_0/x_2)\cdot\allowbreak\Phi_{1}(x_0/x_2)\cdot
\allowbreak\Phi_{2}(x_1/x_0)\cdot\allowbreak\Phi_{2}(x_1/x_2)\cdot\allowbreak\Phi_{2}(x_2/x_0)\cdot\allowbreak\Phi_{2}(x_0^5/x_1x_2^4)\cdot\allowbreak
\Phi_{2}(x_0x_2^2/x_1^3)\cdot\allowbreak\Phi_{2}(x_0/x_2)\cdot\allowbreak\Phi_{2}(x_0x_2/x_1^2)\cdot\allowbreak\Phi_{2}(x_0^2/x_1x_2)\cdot
\allowbreak\Phi_{6}(x_0/x_1)\cdot\allowbreak\Phi_{6}(x_1/x_2)\cdot\allowbreak\Phi_{6}(x_0/x_2)\cdot\allowbreak\Phi_{4}(x_0/x_2)$
\\
\\
$s_{\phi_{30,12}'}=\Phi_{1}(x_1^5/x_0^3x_2^2)\cdot\allowbreak\Phi_{1}(x_1/x_2)\cdot
\allowbreak\Phi_{1}(x_0/x_1)\cdot\allowbreak\Phi_{1}(x_0/x_1)\cdot\allowbreak\Phi_{1}(x_0/x_2)\cdot\allowbreak\Phi_{1}(x_0/x_2)\cdot\allowbreak\Phi_{1}(x_1/x_2)\cdot\allowbreak\Phi_{1}(x_1/x_2)\cdot
\allowbreak\Phi_{2}(x_0/x_1)\cdot\allowbreak\Phi_{2}(x_0/x_2)\cdot\allowbreak\Phi_{2}(x_2/x_1)\cdot\allowbreak\Phi_{2}(x_1^5/x_0x_2^4)\cdot
\allowbreak\Phi_{2}(x_1x_2^2/x_0^3)\cdot\allowbreak\Phi_{2}(x_1/x_2)\cdot\allowbreak\Phi_{2}(x_1x_2/x_0^2)\cdot\allowbreak\Phi_{2}(x_1^2/x_0x_2)\cdot\allowbreak\Phi_{6}(x_1/x_0)\cdot\allowbreak\Phi_{6}(x_0/x_2)\cdot\allowbreak\Phi_{6}(x_1/x_2)\cdot\allowbreak\Phi_{4}(x_1/x_2)$
\\
\\
$s_{\phi_{36,5}}=\Phi_{1}(1/\zeta_3)\cdot\allowbreak\Phi_{1}(x_0/x_2)\cdot\allowbreak\Phi_{1}(x_1/x_0)\cdot\allowbreak\Phi_{1}(\zeta_3
x_0^2/x_1x_2)\cdot\allowbreak\Phi_{1}(\zeta_3^2x_0x_2/x_1^2)\cdot\allowbreak\Phi_{1}(x_2^2/\zeta_3^2x_0x_1)\cdot\allowbreak\Phi_{2}(x_1x_2/x_0^2)
\cdot\allowbreak\Phi_{2}(x_0^2x_2/\zeta_3x_1^3)\cdot\allowbreak\Phi_{2}(\zeta_3^2x_0^2x_1/x_2^3)\cdot\allowbreak\Phi_{6}(x_0^2/x_1x_2)\cdot
\allowbreak\Phi_{6}(x_0/x_1)\cdot\allowbreak\Phi_{6}(x_0/x_2)\cdot\allowbreak\Phi_{6}(x_1/x_2)\cdot\allowbreak\Phi_{5}(\zeta_3x_0/x_2)\cdot\allowbreak
\Phi_{5}(\zeta_3x_0/x_1)$
\\
\\
$s_{\phi_{40,8}}=\Phi_{1}(x_0^3x_1^2/x_2^5)\cdot\allowbreak\Phi_{1}(x_1/x_0)\cdot\allowbreak\Phi_{1}(x_0/x_1)\cdot\allowbreak\Phi_{1}(x_0/x_2)\cdot
\allowbreak\Phi_{1}(x_2/x_0)\cdot\allowbreak\Phi_{1}(x_0/x_2)\cdot\allowbreak\Phi_{1}(x_1^2x_2/x_0^3)\cdot\allowbreak\Phi_{1}(x_2/x_1)\cdot\allowbreak
\Phi_{2}(x_0/x_1)\cdot\allowbreak\Phi_{2}(x_0x_1/x_2^2)\cdot\allowbreak\Phi_{2}(x_0^2/x_1x_2)\cdot\allowbreak\Phi_{2}(x_0/x_2)\cdot\allowbreak\
\Phi_{6}(x_0/x_1)\cdot\allowbreak\Phi_{6}(x_1/x_2)\cdot\allowbreak\Phi_{4}(x_0/x_1)\cdot\allowbreak\Phi_{4}(x_0/x_2)\cdot\allowbreak\Phi_{3}(x_0x_2/x_1^2)$
\\
\\
$s_{\phi_{45,6}}=\Phi_{1}(\zeta_3)\cdot\allowbreak\Phi_{1}(\zeta_3^2x_0^2/x_1x_2)\cdot\allowbreak\Phi_{1}(\zeta_3x_0x_2/x_1^2)
\cdot\allowbreak\Phi_{1}(\zeta_3x_0x_1/x_2^2)\cdot\allowbreak\Phi_{1}(x_2/x_0)\cdot\allowbreak\Phi_{1}(x_1/x_2)\cdot\allowbreak\Phi_{2}(\zeta_3^2x_0^2x_2/x_1^3)\cdot\allowbreak\Phi_{2}(\zeta_3^2x_1^2x_2/x_0^3)
\cdot\allowbreak\Phi_{2}(x_0/\zeta_3^2x_2)\cdot\allowbreak\Phi_{2}(\zeta_3x_1/
x_2)\cdot\allowbreak\Phi_{2}(x_0x_1/x_2^2)\cdot\allowbreak\Phi_{6}(x_1/x_0)\cdot\allowbreak\Phi_{6}(x_1/x_2)\cdot\allowbreak\Phi_{6}(x_0/x_2)
\cdot\allowbreak\Phi_{6}(x_0x_1/x_2^2)\cdot\allowbreak\Phi_{4}(\zeta_3x_0/x_2)\cdot\allowbreak\Phi_{4}(\zeta_3x_1/x_2)$
\\
\\
$s_{\phi_{60,7}}=\Phi_{1}(x_0x_1^2/x_2^3)\cdot\allowbreak\Phi_{1}(x_0/x_1)\cdot\allowbreak\Phi_{1}(x_0/x_1)\cdot\allowbreak\Phi_{1}(x_1/x_2)\cdot\allowbreak\Phi_{1}(x_1/x_2)
\cdot\allowbreak\Phi_{1}(x_0/x_2)\cdot\allowbreak\Phi_{1}(x_0/x_2)\cdot\allowbreak\Phi_{1}(x_1x_2^2/x_0^3)\cdot\allowbreak\Phi_{2}(x_0/x_1)\cdot\allowbreak\Phi_{2}(x_0^4x_1/x_2^5)\cdot\allowbreak\Phi_{2}(x_1^2/x_0x_2)\cdot\allowbreak\Phi_{2}(x_1x_2/x_0^2)\cdot\allowbreak\Phi_{2}(x_1/x_2)\cdot\allowbreak
\Phi_{2}(x_1/x_2)\cdot\allowbreak\Phi_{6}(x_0/x_1)\cdot\allowbreak\Phi_{6}(x_2/x_1)\cdot\allowbreak\Phi_{4}(x_0/x_1)$
\\
\\
$s_{\phi_{60,11}''}=\Phi_{1}(x_1x_2^3/x_0^4)\cdot\allowbreak\Phi_{1}(x_0/x_1)\cdot\allowbreak\Phi_{1}(x_0/x_1)\cdot\allowbreak\Phi_{1}(x_1/x_2)\cdot\allowbreak\Phi_{1}(x_1/x_2)\cdot\allowbreak\Phi_{1}(x_2/x_0)\cdot\allowbreak\Phi_{1}(x_2/x_0)\cdot\allowbreak\Phi_{1}(x_0x_2^3/x_1^4)\cdot\allowbreak\Phi_{2}
(x_1^2/x_0x_2)\cdot\allowbreak\Phi_{2}(x_0x_1/x_2^2)\cdot\allowbreak\Phi_{2}(x_0x_1/x_2^2)\cdot\allowbreak\Phi_{2}(x_0^2/x_1x_2)\cdot\allowbreak\Phi_{2}(x_0/x_1)\cdot\allowbreak\Phi_{2}(x_0/x_1)\cdot\allowbreak\Phi_{6}(x_0x_1/x_2^2)\cdot\allowbreak\Phi_{6}(x_1/x_0)$
\\
\\
$s_{\phi_{60,12}}=2\Phi_{1}(x_0x_2^2/x_1^3)\cdot\allowbreak\Phi_{1}(x_0x_1^2/x_2^3)\cdot\allowbreak\Phi_{1}(x_1/x_0)\cdot\allowbreak\Phi_{1}(x_1/x_0)\cdot\allowbreak\Phi_{1}(x_1/x_0)\cdot\allowbreak\Phi_{1}(x_2/x_0)\cdot\allowbreak\Phi_{1}(x_2/x_0)\cdot\allowbreak\Phi_{1}(x_2/x_0)\cdot\allowbreak\Phi_{2}(x_2/x_1)\cdot\allowbreak\Phi_{2}
(x_1/x_2)\cdot\allowbreak\Phi_{2}(x_0^2/x_1x_2)\cdot\allowbreak\Phi_{2}
(x_0^2/x_1x_2)\cdot\allowbreak\Phi_{6}(x_0^2/x_1x_2)\cdot\allowbreak\Phi_{6}(x_1/x_2)\cdot\allowbreak\Phi_{6}(x_2/x_1)\cdot\allowbreak\Phi_{4}
(x_0/x_1)\cdot\allowbreak\Phi_{4}(x_0/x_2)$
\\
\\
$s_{\phi_{64,8}}=2\Phi_{1}(rx_1/x_2^2)\cdot
\allowbreak\Phi_{1}(x_2^2/rx_0)\cdot\allowbreak\Phi_{1}(x_0/x_2)\cdot\allowbreak\Phi_{1}(x_2/x_1)\cdot\allowbreak\Phi_{1}(x_0/x_1)\cdot\allowbreak
\Phi_{1}(x_1/x_0)\cdot\allowbreak\Phi_{1}(x_0^3/x_1x_2^2)\cdot\allowbreak\Phi_{1}(x_0x_2^2/x_1^3)\cdot\allowbreak\Phi_{2}(rx_0^2/x_1^2x_2)\cdot
\allowbreak\Phi_{2}(rx_1^2/x_0^2x_2)\cdot\allowbreak\Phi_{3}(x_0x_1/x_2^2)\cdot\allowbreak\Phi_{10}(r/x_2)\cdot\allowbreak\Phi_{15}(r/x_0)$\\
where $r=\root 2\of{x_0x_1}$
\\
\\
$s_{\phi_{80,9}}=2\Phi_{1}(x_0x_2^2/x_1^3)\cdot\allowbreak
\Phi_{1}(x_0x_1^2/x_2^3)\cdot\allowbreak\Phi_{1}(x_0/x_1)
\cdot\allowbreak\Phi_{1}(x_0/x_1)\cdot\allowbreak\Phi_{1}(x_2/x_0)\cdot\allowbreak\Phi_{1}(x_2/x_0)\cdot\allowbreak\Phi_{1}(x_2/x_1)\cdot
\allowbreak\Phi_{1}(x_2/x_1)\cdot\allowbreak\Phi_{2}(x_0/x_2)\cdot\allowbreak\Phi_{2}(x_0/x_1)\cdot\allowbreak\Phi_{4}(x_1x_2/x_0^2)\cdot\allowbreak
\Phi_{4}(x_0/x_1)\cdot\allowbreak\Phi_{4}(x_0/x_2)\cdot\allowbreak\Phi_{3}(x_0^2/x_1x_2)\cdot\allowbreak\Phi_{12}(x_1/x_2)$
\\
\\
$s_{\phi_{81,10}}=3\Phi_{2}(rx_2/x_0^2)\cdot\allowbreak\Phi_{2}(rx_2/x_1^2)\cdot\allowbreak\Phi_{2}(rx_0/x_2^2)\cdot\allowbreak\Phi_{2}
(rx_0/x_1^2)\cdot\allowbreak\Phi_{2}(rx_1/x_0^2)\cdot\allowbreak\Phi_{2}(rx_1/x_2^2)\cdot\allowbreak\Phi_{2}(x_0x_1/x_2^2)\cdot\allowbreak\Phi_{2}(x_0x_2/
x_1^2)\cdot\allowbreak\Phi_{2}(x_1x_2/x_0^2)\cdot\allowbreak\Phi_{2}(r/x_2)
\cdot\allowbreak\Phi_{2}(r/x_0)\cdot\allowbreak\Phi_{2}(r/x_1)\cdot
\allowbreak\Phi_{4}(r^2/x_0x_1)\cdot\allowbreak\Phi_{4}(r^2/x_0x_2)\cdot
\allowbreak\Phi_{4}(r^2/x_1x_2)\cdot\allowbreak\Phi_{5}(r/x_0)\cdot\allowbreak
\Phi_{5}(r/x_2)\cdot\allowbreak\Phi_{5}(r/x_1)$\\
where $r=\root 3\of{x_0x_1x_2}$}
\\
\\
Following theorem $\ref{Semisimplicity Malle}$ and \cite{Ma4}, Table
8.2, if we set
$$X_i^6:=(\zeta_3)^{-i}x_i  \,\,\,(i=0,1,2),$$
then $\mathbb{Q}(\zeta_3)(X_0,X_1,X_2)$ is a splitting field for
$\mathcal{H}(G_{32})$. Hence the factorization of the Schur elements
over that field is as described by theorem $\ref{Schur element
generic}$.

\subsection*{Some corrections on the article \cite{MaRo}}

\begin{description}
  \item [Example 3.17 ($G_{34}$)]:\\
        type $\mathcal{M}(Z_3)''$: the first character of the first block is
        $\phi_{105,8}''$.\\
        type $\alpha_{3 \times 5}''$: the second block is\\
        $(\phi_{84,41},\phi_{84,37} ,\phi_{336,34},\phi_{336,32},\phi_{420,31},\phi_{420,35},\phi_{504,33}).$
  \item [Example 3.19 ($G_{36}$)]: There are two families of type $\mathcal{M}(Z_2)$
  missing:\\
        $(\phi_{15,28},\phi_{105,26},\phi_{120,25})$ and
        $(\phi_{15,7},\phi_{105,5},\phi_{120,4})$ \\
        (the characters may not be in the correct order).
  \item [Example 4.1 ($G_{5}$)]: The decomposition matrix of the last family for
  $p=3$ is wrong; the characters of degree 1 aren't in the same
  3-block of the group algebra with the characters of degree 2.
  \item [Example 4.3 ($G_{13}$)]: The decomposition matrix of the last family for
  $p=2$ is wrong; the characters of degree 2 aren't in the same
  2-block of the group algebra with the characters of degree 4.
  \item [Example 6.4 ($G_{9}$)]: The blocks given here are completely wrong. The correct, probably, are:\\
  $( \phi_{2,7}', \phi_{2,11}' ), ( \phi_{2,7}'', \phi_{2,11}'' ),
  ( \phi_{1,12}', \phi_{1,12}'', \phi_{2,4}, \phi_{2,8}),
  ( \phi_{1,6}, \phi_{1,30}, \phi_{2,10}, \phi_{2,14} ),$\\
  $( \phi_{2,5}, \phi_{2,13}, \phi_{2,1}, \phi_{2,17}, \phi_{4,9}, \phi_{4,3},
      \phi_{4,7}, \phi_{4,5} ).$
  \item [Example 6.5 ($G_{10}$)]: The last three non-trivial blocks are:\\
  $( \phi_{1,16}, \phi_{2,10}, \phi_{3,12}', \phi_{3,8}', \phi_{3,4} ),
  ( \phi_{1,12}, \phi_{2,18}, \phi_{3,8}'', \phi_{3,12}'', \phi_{3,16}
  )$,\\
  $( \phi_{1,14}, \phi_{1,26}, \phi_{2,8}, \phi_{3,14}\, \phi_{3,2},
      \phi_{3,10}'', \phi_{3,10}', \phi_{3,6}'', \phi_{3,6}', \phi_{4,11},
      \phi_{4,5} ).$
\end{description}

\end{document}